%% file: RW_typeA7.tex
\documentclass[12pt]{amsart}
\usepackage[latin1]{inputenc}
\usepackage{xspace,amssymb,amsfonts}
\usepackage{amsthm,amsmath,amscd,stmaryrd,latexsym}
\usepackage[hmargin=3cm,vmargin=3cm]{geometry}

\usepackage{graphicx}
\usepackage[all]{xy}
\usepackage{tikz}
\usetikzlibrary{calc, matrix, arrows, decorations.pathmorphing}
\usepackage[dvips]{epsfig}
\usepackage{pinlabel}

\newcommand{\F}{\mathbb{F}}
\newcommand{\g}{\mathfrak{g}}
\newcommand{\Z}{\mathbb{Z}}
\newcommand{\Cat}{\mathcal{C}}
\newcommand{\SB}{\mathbb{S}\mathrm{Bim}}
\newcommand{\h}{\mathfrak{h}}
\newcommand{\OCat}{\mathcal{O}}
\newcommand{\C}{\mathbb{C}}
\newcommand{\Hom}{\operatorname{Hom}}
\newcommand{\slf}{\mathfrak{sl}}
\newcommand{\K}{\mathbb{K}}
\newcommand{\gr}{\operatorname{gr}}
\newcommand{\Ext}{\operatorname{Ext}}
\newcommand{\End}{\operatorname{End}}

\newcommand{\bs}{\underline{s}}
\newcommand{\Part}{\mathcal{P}}
\newcommand{\Ind}{\operatorname{Ind}}
\newcommand{\head}{\mathsf{head}}
\newcommand{\gl}{\mathfrak{gl}}
\newcommand{\GL}{\operatorname{GL}}

    \def\BM{{\mathbb{B}}}
    \def\CM{{\mathbb{C}}}
\def\DG{{\mathfrak D}}    
    
    \def\FM{{\mathbb{F}}}
  \def\gfrak{{\mathfrak g}}

  \def\hg{{\mathfrak h}}

    \def\LM{{\mathbb{L}}}

    \def\QM{{\mathbb{Q}}}
    \def\RM{{\mathbb{R}}}
    \def\SM{{\mathbb{S}}}

    \def\VM{{\mathbb{V}}}

    \def\ZM{{\mathbb{Z}}}

    \def\CC{{\mathcal{C}}}
    \def\DC{{\mathcal{D}}}
  \def\eb{{\mathbf e}}  \def\EC{{\mathcal{E}}}
  \def\fb{{\mathbf f}}  \def\FC{{\mathcal{F}}}
    \def\GC{{\mathcal{G}}}
\def\HB{{\mathbf H}}    \def\HC{{\mathcal{H}}}

    \def\MC{{\mathcal{M}}}
    
    \def\OC{{\mathcal{O}}}
    \def\PC{{\mathcal{P}}}
    \def\QC{{\mathcal{Q}}}
    
\def\SB{{\mathbf S}}    \def\SC{{\mathcal{S}}}
    
    \def\UC{{\mathcal{U}}}
    \def\VC{{\mathcal{V}}}

\def\a{\alpha}

\def\g{\gamma}

\def\D{\Delta}

\def\ep{\varepsilon}

\def\la{\lambda}
\def\La{\Lambda}

\let\phi=\varphi
\let\tilde=\widetilde

\usepackage{bbm}
\def\C{{\mathbbm C}}
\def\N{{\mathbbm N}}

\def\Z{{\mathbbm Z}}

\def\1{\mathbbm{1}}

\newcommand{\un}{\underline}

\newcommand{\ot}{\otimes}
\newcommand{\pa}{\partial}
\newcommand{\co}{\colon}

\renewcommand{\to}{\rightarrow}

\renewcommand{\sl}{\mathfrak{sl}}

\newcommand{\refequal}[1]{\xy {\ar@{=}^{#1}
(-1,0)*{};(1,0)*{}};
\endxy}

\newcommand{\Bim}{\textbf{Bim}}

\newcommand{\Res}{{\rm Res}}
\newcommand{\op}{{\rm op}}
\newcommand{\id}{{\rm id}}

\newcommand{\std}{{\rm std}}

\newcommand{\Rep}{\textbf{Rep}}
\newcommand{\Mod}{-\textbf{Mod}}

\newcommand{\Ring}{\mathcal{R}}
\newcommand{\Ringg}{\mathsf{R}}

\newcommand{\bari}{\overline{i}}

\newcommand{\Fock}{\FC}

\newcommand{\SBSBim}{\SC\BM\SM\textrm{Bim}}
\newcommand{\SSBim}{\SC\SM\textrm{Bim}}
\newcommand{\SBim}{\SM\textrm{Bim}}
\newcommand{\BSBim}{\BM\SM\textrm{Bim}}
\newcommand{\LLL}{\LM\LM}
\newcommand{\Flag}{\textbf{Flag}}

\newcommand{\tL}{\tilde{\La}}

\newcommand{\tgl}{\tilde{\gl}}

\newcommand{\fpoly}{a}
\newcommand{\epoly}{b}

\newcommand{\pre}{{\operatorname{pre}}}

\newcommand{\qchoose}[2]
{ \left [ \begin{array}{c}
    #1 \\ #2 \end{array} \right ]}

\newtheorem{Thm}{Theorem}[section]
\newtheorem{Prop}[Thm]{Proposition}
\newtheorem{Cor}[Thm]{Corollary}
\newtheorem{Lem}[Thm]{Lemma}
\theoremstyle{definition}

\newtheorem{Rem}[Thm]{Remark}

\newcommand{\ig}[2]{\vcenter{\xy (0,0)*{\includegraphics[scale=#1]{fig/#2}} \endxy}}

\RequirePackage{color}
\definecolor{myred}{rgb}{0.75,0,0}
\definecolor{mygreen}{rgb}{0,0.5,0}
\definecolor{myblue}{rgb}{0,0,0.65}

\RequirePackage{ifpdf}
\ifpdf
  \IfFileExists{pdfsync.sty}{\RequirePackage{pdfsync}}{}
  \RequirePackage[pdftex,
   colorlinks = true,
   urlcolor = myblue, 
   citecolor = mygreen, 
   linkcolor = myred, 
   pagebackref,
   bookmarksopen=true]{hyperref}
\else
  \RequirePackage[hypertex]{hyperref}
\fi

\newtheorem{thm}{Theorem}[section]
\newtheorem{lemma}[thm]{Lemma}

\newtheorem{prop}[thm]{Proposition}
\newtheorem{cor}[thm]{Corollary}
\newtheorem{claim}[thm]{Claim}

\newtheorem*{prop*}{Proposition}

\theoremstyle{definition}
\newtheorem{defn}[thm]{Definition}

\newtheorem{notation}[thm]{Notation}
\newtheorem{ex}[thm]{Example}
\newtheorem{example}[thm]{Example}

\theoremstyle{remark}
\newtheorem{remark}[thm]{Remark}

\numberwithin{equation}{section}

\title{Modular representation theory in type A via Soergel bimodules}
\author{Ben Elias and Ivan Losev}
\begin{document}
\begin{abstract}
In this paper we express the multiplicities of modular representation theoretic categories of type A
in terms of affine p-KL polynomials of Elias and Williamson. The representation theoretic categories we deal with
include the category of rational representations of $\operatorname{GL}_n$ and of the quantum group
$U_q(\mathfrak{gl}_n)$,  representations of (degenerate) cyclotomic Hecke and Schur algebras,
and the base field is an algebraically closed field of arbitrary prime characteristics.
In order to approach this problem we define Soergel theoretic version of parabolic  categories
$\mathcal{O}$ in characteristic $p$. We show that these categories have many common features with
the classical parabolic categories $\mathcal{O}$, for example, they are highest weight.
We produce a homomorphism from a (finite or affine) type A 2-Kac-Moody category
to the diagrammatic version of the category of singular Soergel bimodules (again, of finite or affine type A).
This leads to a categorical Kac-Moody action on the Soergel theoretic categories $\mathcal{O}$.
Then we relate the representation theoretic categories to Soergel theoretic ones by proving a uniqueness
result for highest weight categorical actions on Fock spaces.
\end{abstract}

\maketitle

\tableofcontents

\input Introduction5.tex

\input SoergelDiags2.tex

\input FiniteADiags2.tex

\input AffineADiags2.tex

\input Cellular2.tex

\input SoergelO.tex

\input FockUnique4.tex

\input decomp_pKL.tex

\end{document}

%% file: Introduction5.tex
\section{Introduction}

There are many categories of interest in classical representation theory whose decomposition numbers (certain numerical shadows of the structure of the category that allow to compute dimensions/characters of irreducible representations) in characteristic zero
are encoded by various kinds of Kazhdan-Lusztig polynomials. These same categories often have characteristic $p$ analogues (we will call them ``modular'' analogues) whose decomposition
numbers are essentially unknown, unless $p$ is very large. Kazhdan-Lusztig polynomials are known to encode decomposition numbers in category $\OC$, and one of the problems was
that, historically, there was no reasonable modular analogue of category $\OC$, whose decomposition numbers should play the role of Kazhdan-Lusztig polynomials in characteristic $p$.
Recent advances have constructed just such a modular analogue of $\OC$, using diagrammatic Soergel bimodules. Decomposition numbers in modular category $\OC$ go by the name of
$p$-Kazhdan-Lusztig polynomials.

The main goal of this paper is to show that certain kinds of $p$-Kazhdan-Lusztig polynomials do, in fact, encode decomposition numbers in classical modular representation categories.
We restrict our attention to type $A$ for a number of reasons.

In the remainder of the introduction we give a more detailed description of the results, and an overview of the proof. For precise statements of these theorems on decomposition numbers,
see Section \ref{sec-application}.

\subsection{Rep-theoretic and Soergel-theoretic categories}

\subsubsection{Representation theory} \label{sssec-repthry}

Given a Coxeter group $W$, one can define its Hecke algebra $\HB(W)$, a deformation over $\Z[v,v^{-1}]$ of the group algebra $\Z[W]$. It has a standard basis, analogous to the usual
basis of the group algebra using group elements. It also has a Kazhdan-Lusztig basis (shortened to \emph{KL basis}), as defined in the seminal paper \cite{KL79}. The change of basis
matrix between these two bases actually has entries which live in $\Z[v]$; these entries are called KL polynomials. They can be computed algorithmically inside $\HB(W)$, as was proven
by Kazhdan and Lusztig in \cite{KL79}. The Kazhdan-Lusztig conjectures, proven in \cite{BB, Br_Ka}, showed that when $W$ is a Weyl group, the KL polynomials evaluated at $v=1$
also encode decomposition numbers in category $\OC$ for the corresponding complex semisimple Lie algebra. This is wonderful, because at the time there was no known algorithmic way to
compute these decomposition numbers within category $\OC$.

There are many variants on KL polynomials, like singular and parabolic KL polynomials (with respect to parabolic subgroups of $W$), which can be algorithmically computed within certain
variants of $\HB(W)$. (For example, singular KL polynomials agree with certain ordinary KL polynomials.) When $W$ is a Weyl group, these encode decomposition numbers in variants on
category $\OC$, like singular and/or parabolic category $\OC$.

Now let $W$ be the affine Weyl group in type $A$. There are many interesting categories in classical representation theory whose decomposition numbers in characteristic zero are given
by values at $v = 1$ of parabolic KL polynomials (with respect to the finite Weyl group). This is the case for  parabolic categories $\OC$ over $\hat{\gl}_n$, by the work
of Kashiwara-Tanisaki \cite{KasTan}. This is also the case for the following categories.

\begin{enumerate} \item Modules over cyclotomic Hecke algebras (by work of Ariki \cite{Ariki}).
\item Finite dimensional representations of the Beilinson-Lusztig-Macpherson ``divided powers" form of the quantum group $U_\epsilon(\gl_n)$, where $\epsilon$ is a root of $1$ (by work of Kazhdan and Lusztig \cite{KazLusXX} and also Arkhipov-Bezrukavnikov-Ginzburg \cite{ABG}.
\item Categories $\OC$ over cyclotomic Rational Cherednik algebras (by work of Rouquier-Shan-Varagnolo-Vasserot \cite{RSVV},
Losev \cite{VV_proof}, and Webster \cite{Webster_Rouquier}). These include as a special case the categories of modules over the cyclotomic Schur algebras from \cite{DJM}. \end{enumerate}

The categories in this list, which we call \emph{rep-theoretic categories}, all have independent definitions and uses. What unites them, and explains why their decomposition numbers
are encoded by affine parabolic KL polynomials, is that they can all be realized as subquotients in a suitable  parabolic category $\OC$ of affine type $A$. For (2) this
follows from \cite{KazLusXX}, for (3) from \cite{RSVV,VV_proof}, and the categories in (1) are quotients of those from (3). This  parabolic category $\OC$ is something like a
highest weight category: it has simples, standards, costandards, and projectives indexed by a certain (infinite) poset, with projectives being filtered by standards, and standards by
simples, in a fashion respecting the partial order. (It is not truly a highest weight category, because the poset is not finite and the projectives are pro-objects rather than honest objects, but morally it behaves in the same way. For the rest of
this introduction we refer to such a category as a highest weight category; the details are in Section \ref{SS_HW}.) Consequently, a highest weight subquotient of $\OC$
(attached to an interval in the poset) will have the same decomposition numbers as the entire category (within that interval).

There are characteristic $p$ versions of each of the categories given in (1)-(3).
\begin{enumerate}
\item Modules over cyclotomic Hecke algebras in characteristic $p$, and their
degenerate analogs (for example, modular representations of the symmetric groups).
\item Modules over $U_\epsilon(\gl_n)$  and their degenerate analogs, the rational
representations of $\GL_n$.
\item Modules over (degenerate) cyclotomic Schur algebras, these categories are recalled in
Section \ref{SS_cycl_Schur}. \end{enumerate}
We may refer to any of these categories as a \emph{modular rep-theoretic category}.

Before now, there was no modular version of affine  category $\OC$ to relate these disparate constructions in representation theory. Recently, a reasonable candidate has appeared (modular diagrammatic Soergel bimodules). The goal of this paper is to prove that the reasonable candidate will fit the bill.

That modular diagrammatic Soergel bimodules should work in certain contexts was first conjectured by Williamson and Riche \cite{RW}. They have also proven their result in type $A$, using different techniques, and applying to a different set of modular rep-theoretic categories. We compare our work to theirs below in Section \ref{ssec-comparison}.

\subsubsection{Diagrammatic Soergel bimodules}

%

Let $\OC_0$ denote the principal block of category $\OC$, the block containing the trivial representation. In Soergel's alternate proof of the Kazhdan-Lusztig conjectures for finite Weyl
groups, he constructed an additive monoidal category $\SBim$ of bimodules over a polynomial ring $R$, now known as \emph{Soergel bimodules}, see \cite{Soergel, Soer07}. Taking the
quotient of these bimodules by the right action of positive degree polynomials yields a category of left $R$-modules called \emph{Soergel modules}. For (finite) Weyl groups in
characteristic zero, there are fully faithful functors from two other additive categories to Soergel modules: projective objects in category $\OC_0$, and semisimple perverse sheaves on
the flag variety. In these contexts, Soergel bimodules (acting on Soergel modules by tensor product) should be thought of as an algebraic encoding of projective functors (acting on
category $\OC_0$), or of perverse sheaves (acting by convolution on other perverse sheaves). Soergel bimodules categorify the Hecke algebra $\HB(W)$, and the indecomposable Soergel
bimodules descend to the KL basis (in characteristic zero). Consequently, KL polynomials encode multiplicity information in the Soergel category (in characteristic zero), and thus in
category $\OC_0$ as well.

Soergel generalized his category of bimodules to arbitrary Coxeter groups and arbitrary characteristic. In characteristic zero, they continue to categorify the Hecke algebra, and its
Kazhdan-Lusztig basis (as proven by Elias-Williamson \cite{EWHodge}). Thus, Soergel modules should be thought of as the correct analog of (projectives in) category $\OC_0$ for arbitrary
Coxeter groups, in characteristic zero. Even in positive characteristic $\SBim$ sometimes behaves well, categorifying the Hecke algebra, but the indecomposable bimodules may descend to
a different basis. However, when the reflection representation of the Coxeter group is not faithful, Soergel bimodules behave poorly: extra bimodule morphisms appear which did not
exist in characteristic zero, and the Grothendieck group will change. This is a major problem for affine Weyl groups (or any infinite Coxeter groups), as they admit no faithful
representations in finite characteristic.

To define a modular analogue of a category, it is not always the best approach to take an algebraic definition of the category and simply work in finite characteristic (this would be
the definition of modular Soergel bimodules). A better approach is to choose an integral form of the category, and specialize that to finite characteristic.

In \cite{EWGr4sb}, Elias and Williamson described the Soergel category by generators and relations, using the language of planar diagrammatics. More precisely, they defined a
$\Z$-linear category which we denote $\SB$, and call the category of \emph{diagrammatic Soergel bimodules}. There is a functor $\SB \to \SBim$, which is an equivalence (after base
change) when Soergel bimodules are well-behaved. Essentially, $\SB$ encodes those morphisms between Soergel bimodules which are generic, but does not contain the additional morphisms
which may appear in degenerate situations, such as when the reflection representation is not faithful. It is the ``correct'' integral form of $\SBim$; in particular, after base change
to any field $\SB$ will categorify the Hecke algebra. It can be used to construct the correct analog of (projectives in) $\OC_0$ for any Coxeter group in any characteristic.

We call the base change of $\SB$ to a field of finite characteristic $p$ the category of \emph{modular diagrammatic Soergel bimodules}, and continue to denote it $\SB$. The basis of
the Hecke algebra given by the indecomposable objects in $\SB$ is called the \emph{$p$-canonical} or \emph{$p$-KL basis}. Unlike the KL basis, the $p$-canonical basis can not be
computed algorithmically inside $\HB(W)$, but using the generators and relations description of $\SB$, it can be computed algorithmically (with much greater difficulty) within the
category $\SB$ itself. For an introduction to the $p$-canonical basis and further discussion of this algorithm, see \cite{ThoWil}.

\begin{remark} The geometric analog of modular diagrammatic Soergel bimodules are \emph{parity sheaves}, introduced by Juteau-Mautner-Williamson, \cite{JMW}. This is established in
\cite[Part 3]{RW}. Parity sheaves are the ``correct'' finite characteristic analog of perverse sheaves, for which a decomposition theorem still holds. \end{remark}

For reasons which will be clear to the reader soon, we will need to study not just Soergel bimodules but singular Soergel bimodules. The algebraic definition of singular Soergel
bimodules was introduced by Williamson in his thesis \cite{WillSingular}. Where Soergel bimodules are like projective functors acting on the principal block $\OC_0$, singular Soergel
bimodules are like projective functors acting on one block of each singular type.\footnote{The equivalence class of a block of category $\OC$ only depends on the stabilizer subgroup in $W$ of its highest weight, which we refer to as its \emph{singular type}.}

Williamson proved that singular Soergel bimodules (in characteristic zero) categorify the \emph{Hecke algebroid} $\HC(W)$, an idempotented algebra with one idempotent for each finite
parabolic subgroup of the Coxeter group. In precise analogy to Soergel bimodules above, Williamson's algebraic definition of singular Soergel bimodules will behave well in characteristic
zero and other nice situations, but in characteristic $p$ for affine Weyl groups will not actually categorify $\HC(W)$. Instead, we need the 2-category of \emph{diagrammatic singular
Soergel bimodules} $\DG$. This diagrammatic 2-category has yet to appear in the published literature, but is a long-standing work in preparation due to Elias-Williamson
\cite{EWSingular}. Though we hope a remedy will arrive soon, the reader may have to accept several statements about $\DG$ (analogous to the results of \cite{EWGr4sb} for $\SB$) on faith
for the moment. In Remark \ref{rmk:whatusedfromEWsing} we will explain which results from \cite{EWSingular} (apart from the existence of $\DG$) will be assumed.

\subsubsection{Producing singular parabolic categories from diagrammatic Soergel bimodules} \label{sssec-makingparabolic}
It will be important to keep in mind the difference between: the category $\OC$ itself, which splits into infinitely many blocks; the individual blocks of $\OC$, some equivalent to the
principal block $\OC_0$, others singular; and the (2-categories of) projective functors which act on these blocks. To reiterate, the monoidal category $\SB$ encodes projective functors
acting on $\OC_0$. It is not the same as $\OC_0$, but it can be used to produce a model for projectives in $\OC_0$ (c.f. Soergel modules), and thus it can also be used to reproduce the
abelian category $\OC_0$ itself. Analogously, the 2-category $\DG$ encodes projective functors acting on individual (singular) blocks of $\OC$, one of each singular type. It can be used
to reproduce any given block of $\OC$. Finally, by ``duplicating" the objects of $\DG$ to account for all the blocks of $\OC$ with the same singular type, one can obtain a variant of
$\DG$ which encodes projective functors acting on all blocks of $\OC$, and can extract from this the entire abelian category $\OC$ itself. So, when we refer to modular category $\OC$, or
a block thereof, we mean the abelian category extracted in the analogous fashion from the finite characteristic specialization of $\DG$.


Now let $W$ be an affine Weyl group. As noted in \S\ref{sssec-repthry}, it is not actually category $\OC$ itself one is interested in, but various (singular) blocks of a parabolic
category $\OC$. Moreover, we are interested in giving each block the structure of a highest weight category, in order that we can discuss highest weight subquotients. Highest weight
structures will be technically important for other reasons in this paper (categorical truncation, rigidity), which will be explained in due course.

In classical settings, one can actually recover all singular parabolic categories $\OC$ from the 2-category of projective functors (i.e. singular Soergel bimodules). For category $\OC$
itself, this was essentially done in the seminal work of Bernstein and Gelfand \cite{BG}. Constructing parabolic category $\OC$ from ordinary $\OC$ is straightforward: we just take the
Serre subcategory of all modules in $\OC$ that are integrable with respect to a given parabolic subgroup. One can expect similar constructions to work for (modular) diagrammatic
singular Soergel bimodules. Although many of the ideas are understood, there is unfortunately no formal discussion of how to explicitly reconstruct singular parabolic categories $\OC$
from projective functors in the literature, perhaps because there was no need. In Section \ref{sec-hwnonsense} we  give a formal construction which works for our setting of diagrammatic singular Soergel bimodules.

Fix a Coxeter group $W$, and assume one has a category $\DG$ which categorifies the Hecke algebroid. Section \ref{sec-hwnonsense} contains a formal construction which takes $\DG$ and
(under some technical assumptions) produces a family of categories which play the roles of singular blocks of parabolic category $\OC$ for various parabolic subgroups of $W$. First, one
produces all singular blocks of the ordinary category $\OC$ from $\DG$, as discussed above. Then one uses a cellular structure  on $\SB$ to equip the principal block with a highest weight
structure. One uses this highest weight structure and some formal properties of singular Soergel bimodules to produce highest weight structures on all the singular blocks. Then, for any
finite parabolic subgroup of $W$ one can produce a standardly stratified structure (in the sense of \cite{LW}) on the principal block of $\OC$. Using this standardly stratified
structure, one can produce parabolic subcategories $\OC$ in all singular blocks, and prove that they are highest weight. In Section \ref{sec-cell}, the groundwork is laid for the
technical assumptions needed in Section \ref{sec-hwnonsense}. In particular, the cellular structure on $\SB$ is described, using results of \cite{EWGr4sb}.

\begin{remark} There are many technical points about the base ring of the construction. Highest weight categories should live over fields or local complete rings, while $\DG$ lives
over various polynomial rings. We will ignore such questions in the introduction. \end{remark}

\begin{remark} The work in progress \cite{EWSingular} should eventually prove that $\DG$ is a fibered cellular category (in the sense of \cite{ELauda}). This would give a  simpler
proof and construction of the highest weight structures on singular categories $\OC$. We do not use this approach, to avoid reliance on the unavailable results of \cite{EWSingular} as
much as is possible. \end{remark}

In addition, Section \ref{sec-hwnonsense} discusses what happens when one takes the Ringel dual of these singular parabolic categories $\OC$. For affine Weyl groups, one obtains new
highest weight categories in this fashion. One should think that ordinary singular parabolic categories $\OC$ are always of \emph{positive level}, while their Ringel duals are of
\emph{negative level}. (It is not essential to understand what the level means in this context).

Thus, applying these constructions to modular diagrammatic Soergel bimodules yields a host of new categories at our disposal. We may refer to these as \emph{Soergel-theoretic} categories
$\OC$, adding the adjectives modular, parabolic, or singular as necessary. Multiplicities of simples in standards (or standards in projectives) in these categories are encoded by $p$-KL
polynomials.

\subsection{Integrable Kac-Moody representations}

Having introduced our major players, the modular rep-theoretic categories and the modular Soergel-theoretic catgeories $\OC$, we seek to explain why their decomposition numbers should
agree. The setting in which to organize this result is the categorical representation theory.


Here is a rough description of our results, so rough as to be quite false. It is one of the major observations of categorical representation theory that modular rep-theoretic categories
can be used to produce a categorification of a Fock space representation of the affine Lie algebra $\hat{\sl}_e$. Meanwhile, we prove that modular Soergel-theoretic categories $\OC$ can
also be used to produce a categorification of (a portion of) the same Fock space representation. This result is obtained by constructing a categorical action of $\hat{\sl}_e$ on
Soergel-theoretic categories $\OC$ for $\hat{\gl}_n$. 
Finally, we prove a rigidity result for categorifications of a Fock space.
This rigidity result is then used to prove the desired equivalence of rep-theoretic and Soergel-theoretic categories, which implies the equality of decomposition numbers.

The process of taking this rough overview and making it correct and precise is not only quite technical, but also extremely interesting! One significant issue is that Soergel-theoretic
categories $\OC$ do not actually categorify a Fock space representation -- they categorify the tensor product of suitable exterior powers of the tautological level $0$ representation of
$\hat{\mathfrak{sl}}_e$. So one must use some interesting tricks to extract something like a Fock space categorification therefrom. The main trick here, \emph{categorical truncation}
(that should be thought as some version of the semi-infinite wedge construction of the Fock spaces), requires keeping careful track of highest weight structures, which is one of the
reasons why they have been emphasized throughout. Even for the reader who does not intend to read the nasty details, it is worth reading through \S\ref{sssec-Fock} and
\S\ref{sssec-whatwedo} to see more of the interesting features of this story.

\subsubsection{Kac-Moody categorification, Fock space, and rigidity}

Let us briefly recall the well-known results of  \cite{Ariki,LLT,ChuRou06}.

We assume the reader is familiar with Young diagrams for partitions, and the notion of the content of a box in a Young diagram. Let $\PC_n$ denote the set of partitions of $n$, and $\PC
= \coprod_{n \in \N} \PC_n$ be the set of all partitions. For $i \in \Z$ one can define an operator $f_i$ on $\PC \coprod \{0\}$, where $0$ is a formal symbol, as follows: $f_i(\la) =
\mu$ when $\mu$ is obtained from $\la$ by adding a single box of content $i$, and $f_i(\la) = 0$ if no such $\mu$ exists. One can define the operator $e_i$ similarly, via the removal of
a box of content $i$. The \emph{(level 1, charge zero) Fock space} $\Fock$ is the vector space with a basis given by $\PC$. We will write $|\lambda\rangle$ for the basis vector labelled
by $\lambda$. Extending $e_i$ and $f_i$ to linear operators on $\Fock$, one obtains a representation of the affine Lie algebra $\sl_{\infty}$. This representation has level $1$, meaning
that the canonical central element in $\sl_\infty$ acts by $1$.

Fix an integer $e \ge 2$. For $i \in \Z$ we let $\bari$ denote its image in $\Z/e\Z$. Now define the operator $e_{\bari}$ on $\Fock$ to be the infinite sum $\sum_{j \equiv \bari} e_j$
(which acts by a finite sum on any element of $\Fock$), and define $f_{\bari}$ similarly. From this we obtain a representation of the affine Lie algebra $\hat{\sl}_e$ of affine type $A$.
This representation also has level $1$. Note also that we can twist this representation by a diagram automorphism of $\hat{\sl}_e$ that rotates the Dynkin diagram: for $d\in \Z$, we can
define operators $e_{\bari}$ and $f_{\bari}$ as $\sum_{j \equiv \bari+d} e_j, \sum_{j \equiv \bari+d} f_j$. We will denote the resulting representation by $\Fock_d$; collectively, such
representations are \emph{level 1 Fock spaces}, and the space $\mathcal{F}_d$ is said to have charge $d$.
Of course, $\Fock_d$ and $\Fock_{d+e}$ are the same representation.

The irreducible complex representations of the symmetric group $S_n$ are parametrized by $\PC_n$. Consequently, the Grothendieck group of the category of complex representations of
symmetric groups (of all sizes) can be identified with $\Fock$ as a vector space. Induction and restriction between $S_n$ and $S_{n+1}$, together with projection to the eigenvalues of a
Young-Jucys-Murphy operator, can be used to construct functors $E_i$, $F_i$ which categorify the operators $e_i$, $f_i$. In this way one obtains a (weak) categorification of $\Fock$ as
an $\sl_{\infty}$-module. If instead one considers representations of symmetric groups in characteristic $p$ (here, one should set $e = p$), the irreducibles are parametrized by a subset
of partitions, and the Grothendieck group is only a quotient of $\Fock$ (in fact, it is an irreducible $\hat{\sl}_p$-submodule generated by the vacuum vector $|\varnothing\rangle$).
However, if one considers representations of the Schur algebra in characteristic $p$ one recovers the entire Fock space as the Grothendieck group. Eigenvalues of the Young-Jucys-Murphy
operator live in $\Z/p\Z$ instead of $\Z$, so one can only construct functors $E_{\bari}$ and $F_{\bari}$ lifting the operators $e_{\bari}$ and $f_{\bari}$, and giving a (weak)
categorical representation of the affine lie algebra $\hat{\sl}_p$.

In the seminal paper of Chuang-Rouquier \cite{ChuRou06}, they studied the algebras of natural transformations between compositions of these functors, and defined what it means to be a
(strong) categorification of an $\sl_2$ representation. Later, independent work of Khovanov-Lauda \cite{KhoLau10} and Rouquier \cite{Rouq-2KM} defined what it means to be a
categorification of a $\mathfrak{g}$ representation for any Kac-Moody Lie algebra $\mathfrak{g}$, by defining a 2-category $\UC(\gfrak)$ which should act by natural transformations. \footnote{Their
definition was shown to match with algebras of natural transformations in representation theory by Brundan-Kleshchev \cite{BK_KLR}, and the different notions of categorification have
finally been shown to agree by Brundan \cite{Brundan}.} In particular, the categorification of Fock space above does admit an action of the 2-category $\UC(\sl_\infty)$ in characteristic 0, or $\UC(\hat{\sl}_p)$ in characteristic $p$.

In fact, all the rep-theoretic categories we are interested in take part in (strong) categorifications of Fock space, or variants thereof. The Lie algebra which acts is $\hat{\sl}_e$. Let $\SC_q(d,n)$ be the $q$-Schur algebra in characteristic $p$, with $q\in \F_p\setminus \{0,1\}$. That is, $\SC_q(d,n)$ is the image of $U_q(\gl_n(\F))$ in $\operatorname{End}_{\F}((\F^n)^{\otimes d})$,
where $\F$ is an algebraically closed field of characterstic $p$, and $\F^n$ is the tautological
representation of $U_q(\gl_n)$. So the category of $\SC_q(d,n)$-modules is nothing else but the category
of degree $d$ polynomial representations of $U_q(\gl_n)$, such categories for all $n\geqslant d$
are naturally equivalent. The category $\bigoplus_{d=0}^{\infty}\SC_q(d,d)\operatorname{-mod}$ is a categorification of the level one Fock space for $\hat{\sl}_e$, where $e$ is the order of $q$ in $\mathbb{F}^\times_p$.  Using classical Schur algebras (instead of $q$-Schur algebras) in characteristic $p$, one obtains a categorification of the level one Fock space for $\hat{\sl}_e$ with $e = p$. Studying instead the categories of modules over cyclotomic Schur algebras (for $e$ coprime to $p$) or their degenerate analogs (for $e=p$), one obtains
categorifications of higher level Fock spaces with various multicharges; these are, by definition, tensor products of level $1$ Fock spaces.

One useful and critical feature of categorical actions comes from rigidity statements, which can prove that two categories with categorical actions are actually equivalent. Such a
result was first proven for minimal categorifications of finite dimensional irreducible $\sl_2$-representations by Chuang and Rouquier. They also prove a result for isotypic
categorifications: two categorifications of an isotypic representation are equivalent if a small piece of each (the highest weight category) are equivalent. In \cite{Rouq-2KM} Rouquier
extended the first uniqueness result to isotypic categorifications of irreducible integrable highest weight representations for arbitrary Kac-Moody algebras.

In this paper (in Section \ref{S_Fock_unique}) we prove, roughly speaking, an analogous rigidity result for modular categorifications of Fock space, following earlier work of the second author, \cite{VV_proof}, in characteristic $0$.
This will be the tool we use to prove an equivalence between the rep-theoretic and Soergel-theoretic sides of the
story. Our techniques originate from Rouquier's paper, \cite{rouqqsch}, and employ highest weight
structures (and also deformations of categories of interest over formal power series) in an essential way.
We remark that there are several serious challenges one faces to prove rigidity results
in characteristic $p$ that are absent in characteristic $0$, we will briefly introduce the reader to
these nasty details in Section \ref{SS_Fock_rigidity}.

However, it is not at all obvious how one obtains a categorification of Fock space from the Soergel-theoretic categories $\OC$. In fact, we never categorify Fock space itself using
Soergel-theoretic categories! Instead, through a roundabout method, we are able to categorify portions of Fock space. Let us postpone a discussion of how one proves rigidity statements
until after we understand what statement should be proven.

\subsubsection{Fock space vs exterior powers} \label{sssec-Fock}

We begin with the decategorified picture.
We refer the reader to a survey paper of Leclerc \cite{LeclercFock} for background on Fock spaces
and references.

We have defined the Fock space $\Fock$ representation of $\hat{\sl}_e$ in the previous section, using partitions as a basis. A partition $\lambda$ can be thought of as a sequence of
non-increasing integers $\lambda = (\lambda_1 \ge \lambda_2 \ge \ldots \ge \lambda_m)$ where each $\la_i \ge 0$. Note that the empty partition is \emph{singular} in that it is killed by
all operators $e_{\bari}$.

Now let us recall the classical infinite wedge construction of the Fock spaces. For simplicity we will restrict ourselves to the level $1$ situation. The algebra $\hat{\sl}_e$ acts on
the space $\C^e[t,t^{-1}]$ via the epimorphism $\hat{\sl}_e\twoheadrightarrow \sl_e[t,t^{-1}]$. It is convenient to identify the space $\C^e[t,t^{-1}]$ with the space $\C^{\Z}$ with
basis $v_i, i\in \Z$: if $u_1,\ldots,u_e$ is the tautological basis of $\C^e$, then we send $u_i t^j\in \C^e[t,t^{-1}]$ to $v_{i+ej}$. So the operator $f_{\bari}$
sends $v_{j+1}$ to $v_{j}$ when $j \equiv \bari$ modulo $e$. 
The canonical central element of $\hat{\sl}_e$ tautologically acts on $\C^{\Z}$ by $0$, so this is a level $0$ representation.

Of course, we can form the $k$th exterior power $\Lambda^k \C^\Z$, this is also a level $0$ representation. But one also can form the semi-infinite wedge: $\Lambda^{+\infty/2}\C^{\Z}$,
which by definition has basis of semi-infinite wedges $v_{i_1}\wedge v_{i_2}\wedge\ldots\wedge v_{i_k}\wedge\ldots$, where $i_1<i_2<\ldots$, and eventually (for $k \gg 0$)
\begin{equation} \label{eq:eventuallyconstant} i_k=k+d \end{equation} for some integer $d$. This space has a natural action of $\hat{\sl}_e$, now of level $1$. For fixed $d$ all
semi-infinite wedges form a subrepresentation to be denoted by $(\Lambda^{+\infty/2}\C^{\Z})_d$. The representation $(\Lambda^{+\infty/2}\C^{\Z})_d$ is identified with the Fock space
$\Fock_d$, and $v_{i_1}\wedge v_{i_2}\wedge\ldots \mapsto |\lambda\rangle$, where the partition $\lambda$ is determined from the semi-infinite wedge via
\begin{equation}\label{eq:partition_from_wedge} \lambda_j=j+d-i_j. \end{equation} We note that $\lambda_k=0$ for $k \gg 0$ by \eqref{eq:eventuallyconstant}, so we indeed get a partition.

Unfortunately, we cannot categorify this construction, as we do not have a categorification of the semi-infinite
wedge. However, we can categorify all $\Lambda^m \C^{\Z}$ using Soergel-theoretic affine parabolic categories $\OC$. We obviously cannot realize the Fock space representation inside $\Lambda^m\C^{\Z}$ but we can realize a chunk of
it. Namely, let us write $\Fock_d(k)$ for the degree $k$ part of $\Fock_d$, i.e., the span of all
$|\lambda\rangle$ with $|\lambda|=k$. Further, we write $\Fock_d(\leqslant n)$ for $\bigoplus_{k\leqslant n}\Fock_d(k)$. Now fix $n<m$. We have an injective map of vector spaces $\Fock_m(\leqslant n) \to \Lambda^m \C^{\Z}$ by
sending $|\lambda\rangle$ to $v_{1-m-\lambda_1}\wedge v_{2-m-\lambda_2}\wedge\cdots\wedge v_0$. This embedding is
inspired by \eqref{eq:partition_from_wedge}; the final index $v_0$ arises because $\lambda_m=0$ whenever $\lambda$ is a partition of $k$ and $k < m$. We note that the operator $f_{\bari}:\mathcal{F}_m(k)\rightarrow
\mathcal{F}_m(k+1)$ for $k<n$ coincides with the restriction of $f_{\bari}$ from $\Lambda^m \C^{\Z}$. The same
is true for $e_{\bari}: \mathcal{F}_m(k)\rightarrow \mathcal{F}_m(k-1)$ for $k \leqslant n$
(we set $\mathcal{F}_m(-1)=0$), so long as $\bari\neq 0$. However, for $\bari\neq 0$, this is no longer
true, for in the action on $\bigwedge^{m}\C^{\Z}$, $e_{\bar{0}}(v_{i_1}\wedge\ldots \wedge v_{i_{m-1}}\wedge v_0)$
includes an additional monomial: $v_{i_1}\wedge v_{i_2}\wedge\ldots\wedge v_{i_{m-1}}\wedge v_{1}$.

Thus the embedding $\Fock_m(\leqslant n) \to \Lambda^m \C^{\Z}$ is not an $\hat{\sl}_e$ intertwiner, but it can be made into one by modifying or \emph{truncating} the action of
$e_{\bar{0}}$ to ignore this extra monomial. That is, one produces an alternate operator $\underline{e}_{\bar{0}}$ on the image of $\Fock_m(\leqslant n)$, and observes that now the embedding is an intertwiner of \emph{restricted} $\hat{\sl}_e$ representations. The term ``restricted'' refers to the fact that $\Fock_m(\leqslant n)$ is not actually an $\hat{\sl}_e$ representation, because it is not preserved by $f_{\bari}$, but nonetheless the embedding intertwines $f_{\bari}$ whenever that makes sense.

\begin{remark} Here is an analogous situation. Let $M$ be an $\sl_2$ representation, and $\C^2$ be the standard $\sl_2$ representation, with lowest weight $v_-$. Then there is a vector
space embedding $M \to M \ot \C^2$, $m \mapsto m \ot v_-$, which is not an intertwiner. By modifying the action of the raising operator $e$ on $M \ot v_-$, truncating it to an operator
$\underline{e}$ which ignores the final tensor factor, one can force this embedding to be an intertwiner. Note that $\underline{e}$ can be seen as the adjoint to $f$ on $M \ot v_-$ under the correct bilinear form. \end{remark}

The inclusion $\Fock_m(\leqslant n)\hookrightarrow \Lambda^m \C^{\Z}$ and the truncation procedure can be lifted to a categorical level following \cite{VV_proof}.
Namely, suppose that we have a {\it highest weight $\hat{\sl}_e$-categorification} $\VC$ of $\bigwedge^m \C^{\Z}$. What this basically means is that $\VC$
\begin{itemize}
\item[(i)] carries a categorical action of $\hat{\sl}_e$,
\item[(ii)] and a highest weight structure, where the standards are labelled by the increasing sequences of $m$ integers, with
a particular partial order;
\item[(iii)] the complexified Grothendieck group is identified (as a module over $\hat{\sl}_e$)
with $\bigwedge^m \C^{\Z}$ so that the standard object labelled by $i_1<i_2<\ldots<i_m$ corresponds to $v_{i_1}\wedge v_{i_2}\wedge\ldots\wedge v_{i_m}$,
\item[(iv)] and the categorification functors $E_{\bari},F_{\bari}$ map standard objects to standardly filtered
ones.
\end{itemize}
Over the base field $\C$, an example of $\VC$ is provided by the {\it Kazhdan-Lusztig category}
for $\hat{\gl}_m$. By definition, this is the parabolic subcategory consisting of all
$\operatorname{GL}[[t]]$-integrable modules in the (full) affine category $\OC$
with level $-m-e$.  We will produce $\VC$ in positive characteristic by taking a suitable
sum of various blocks of a Soergel-theoretic parabolic category $\OC$.

The choice of partial order guarantees that the labels corresponding to partitions form a poset ideal. So the image of
$\Fock_m(\leqslant n)$ in $\Lambda^m \C^{\Z}$ is categorified by a highest weight subcategory $\VC'=\bigoplus_{k=0}^n \VC'(n)$ inside $\VC$. We still have $F_{\bari}\VC'(k)\subset \VC'(k+1)$ for $k<n$ and $E_{\bari}\VC'(k)\subset \VC'(k-1)$ for $k \le n$ and $\bari \ne 0$, while, of course, $E_{\bar{0}}$ does not
preserve $\VC'$. However, it was shown in \cite{str} that the functor $F_{\bar{0}}:\VC'(k)\rightarrow \VC'(k+1)$ still has a biadjoint which we denote
$\underline{E}_{\bar{0}}$. Replacing $E_{\bar{0}}$ with $\underline{E}_{\bar{0}}$, we obtain a restricted categorical $\hat{\sl}_e$ action on $\VC'$.

We note that analogous constructions make sense for higher level Fock spaces and tensor products of wedges. Set $\underline{m}:=(m_1,\ldots,m_\ell)$ and
$\Fock_{\underline{m}}=\bigotimes_{j=1}^\ell \Fock_{m_\ell}$. We can embed $\Fock_{\underline{m}}(\leqslant n)$ into $V_{\underline{m}}:=\bigwedge^{m_1}\C^{\Z}\otimes \ldots\otimes
\bigwedge^{m_\ell}\C^{\Z}$ similarly to the above. The notion of a highest weight categorification, $\VC_{\underline{m}}$, of the tensor product of wedges still makes sense (over $\C$ an
example is provided by a parabolic affine category $\OC$). Inside, we can consider the subcategory $\VC'_{\underline{m}}(\leqslant n)$ corresponding to
$\Fock_{\underline{m}}(\leqslant n)$ (we emphasize that it will now depend on $\underline{m}$ itself not on its class modulo $e$) and equip it with a restricted categorical
$\hat{\sl}_e$-action.

\subsubsection{What we accomplish} \label{sssec-whatwedo}

Now we can finally explain what the main results of this paper are, and how they are proven.

As mentioned previously, Section \ref{sec-hwnonsense} constructs the Soergel-theoretic singular parabolic categories $\OC$, attached to for $\hat{\sl}_m$ for any $m$ (or more precisely,
attached to an algebra $\tgl_m$ that differs from $\hat{\sl}_m$ in the Cartan part). When $\un{m}$ satisfies $\sum m_i = m$, we can take what is essentially a large direct sum of such categories $\OC$ to produce a big category called
$\VC_{\un{m}}$. Because one understands the Grothendieck group of $\DG$, one can bootstrap this to show that the Grothendieck group of $\VC_{\un{m}}$ is $V_{\underline{m}}=\Lambda^{m_1}\C^{\Z}\otimes\ldots\otimes \Lambda^{m_\ell}\C^{\Z}$,
as a vector space.

\begin{remark} This category $\VC_{\un{m}}$ is the (modular) diagrammatic Soergel version of a parabolic category $\OC$ in level $-m-e$ for $\hat{\gl}_m$ (for the standard parabolic subalgebra of $\hat{\gl}_m$ corresponding to the composition
$\underline{m}$ of $m$).  However, do not be confused: what this categorifies is $V_{\un{m}}$, which is a level zero representation of $\hat{\sl}_e$. \end{remark}

Now we need to equip $\VC_{\un{m}}$ with the structure of a categorical action of $\hat{\sl}_e$  making it into a highest weight categorification of $V_{\underline{m}}$. In Section \ref{S_KM_affine_A}, a categorical action of $\hat{\sl}_e$ on a modification of $\DG$ is established, with the previous chapters laying the groundwork for this result. Here is where the diagrammatic, generators and relations description of $\DG$ really shines, making it possible to explicitly define the 2-functor
from Khovanov-Lauda's category $\UC$ (see Section \ref{SS_intro_2KM_action} for an overview).
This allows one to equip $\VC_{\un{m}}$ with a $\UC$-action. Using the properties
of the highest weight structure on the Soergel-theoretic category $\OC$, and its interplay
with a $\DG$-action, we see that $\VC_{\underline{m}}$ indeed becomes a highest weight categorification
of $V_{\underline{m}}$ (Section \ref{SSS_full_affine_cat}).
Then, using the categorical trunction construction described in \S\ref{sssec-Fock},
we cook up a restricted highest weight categorification of $\Fock_{\un{m}}(\leqslant n)$, which we denote $\VC'_{\un{m}}(\leqslant n)$.
This is done in Section \ref{SS_Fock_restr_categ}.

In the remainder of Section \ref{S_Fock_unique}, starting with Section \ref{SS_equiv_thms},
we prove our uniqueness result for highest weight (restricted)  categorifications of Fock
spaces. This allows us to identify  $\VC'_{\un{m}}(\leqslant n)$ with
the direct sum of categories of modules over cyclotomic $q$-Schur (or degenerate Schur)
algebras under the assumption that $m_1-m_2\gg m_2-m_3\gg\ldots m_{\ell-1}-m_\ell\gg n$.
This is the most technical part of the paper, and will be discussed further in Section \ref{SS_Fock_rigidity}.

Our equivalence theorem allows us to express the decomposition numbers in modular rep-theoretic categories
via $p$-KL polynomials, similarly to what is done in characteristic $0$. More precisely, we can express
the decomposition numbers for cyclotomic (degenerate) Hecke algebras via those for cyclotomic q- or
degenerate Schur algebras. Also the multiplicities for rational representations of $\operatorname{GL}_n$
or its q-analogs can be expressed via those for the usual Schur algebras or $q$-Schur algebras, this is basically due to the fact that any rational representation of $\GL_n$ becomes polynomial after twisting with a large enough power
of $\det$.
On the other hand,  $\VC'_{\un{m}}(\leqslant n)$ is a highest weight subcategory of $\VC_{\underline{m}}$
and the multiplicities in the latter are expressed via the p-KL polynomials in a standard way.
We will elaborate on this in Section \ref{sec-application}.

In the rest of this introduction, we give some additional background on the various pieces of this puzzle. For sake of clarity, we will not state precise theorems about which
decomposition numbers agree with which $p$-KL polynomials in this introduction, but defer these theorems to the very short Section \ref{sec-application}.


\subsection{2-Kac-Moody actions on Soergel theoretic categories}\label{SS_intro_2KM_action}

The idea that Kac-Moody quantum groups should act categorically (i.e. by functors) on sums of blocks of rep-theoretic parabolic categories $\OC$ is a ``classical'' one in type $A$, dating
back to work of Igor Frenkel and collaborators (e.g. \cite{BFK}) in the 90's. In these classical setups, the functors which lift the Kac-Moody generators are certain projective functors. This suggests that there should be a 2-functor from Khovanov-Lauda's category $\UC$ to the category of singular Soergel bimodules. We realize this functor in this paper.

\begin{remark} There is also a 2-Kac-Moody action on affine type A categories $\mathcal{O}$ (of negative level) \cite{RSVV,VV_proof}, but it has a different nature; the functors come
from Kazhdan-Lusztig fusion products. \end{remark}

More precisely, our immediate goal is to construct a 2-functor from a Kac-Moody 2-category attached to $\tgl_e$, to the 2-category $\DG$ of diagrammatic singular Soergel bimodules
attached to $\tgl_m$. Here we write $\tgl_?$ for the algebra $\hat{\sl}_?\oplus \C d\oplus \C I$, where $d$
is the grading element and $I$ is a central element (to be though as the unit matrix from $\gl_?$). A different but closely related $2$-functor is already in the literature when $e \ge m \ge 3$.

Because the history is somewhat complicated, let us explain several other 2-functors in the existing literature, beginning with finite type $A$. Khovanov and Lauda \cite{KhoLau10} define
an action of $\UC(\sl_e)$ (the Kac-Moody 2-category in finite type $A$) on what they call the \emph{equivariant flag $2$-category} for $\gl_m$. This flag $2$-category is essentially the
$2$-category of (algebraic) singular Soergel bimodules in disguise. Later, Mackaay-Stosic-Vaz \cite{MSV} constructed an action of $\UC(\gl_e)$ on categorified MOY diagrams. Their
categorified MOY diagrams are precisely the $2$-category of (algebraic) singular Soergel bimodules for $\gl_m$. In this paper, in Section \ref{S_KM_fin_A}, we will construct an action of
$\UC(\gl_e)$ on $\DG(\gl_m)$, on diagrammatic singular Soergel bimodules in finite type $A$. Philosophically, these actions are all one and the same. It will be important for us to
define the action diagrammatically, rather than algebraically, in order to guarantee that it works appropriately in generalized settings (like the modular setting). The relationship
between these 2-functors is described in the commutative diagram \eqref{eq:bigdiagram}, and these previous results are discussed in much more detail in Section \ref{subsec-previous}.

It is worth mentioning that the diagrammatic and algebraic approaches to singular Soergel bimodules are quite different. The algebraic approach is grounded in polynomials and their
manipulation. The diagrammatic approach is built from the \emph{Frobenius extension} structures between different rings of invariant polynomials; it rarely wants to examine polynomials
themselves, but only cares about the properties of various Frobenius trace and coproduct maps. Very complicated operations on polynomial rings can be encoded easily using Frobenius
structure maps, which in turn are encoded with rather simple diagrams. This simplicity is one of the major advantages of the diagrammatic approach. Our 2-Kac-Moody action on $\DG$ is an
excellent illustration of this idea, as the well-definedness of our 2-Kac-Moody action is vastly easier to check than the algebraic proof in Mackaay-Stosic-Vaz \cite[Section
4]{MSV}.\footnote{In fact, Mackaay-Stosic-Vaz never give a complete proof of this result, giving only example computations of several relations in \cite[Section 4.3.1]{MSV}, and leaving
the remaining computations as tedious exercises to the reader. It is less tedious with our diagrammatic technology, so we give a complete proof.}

Fix $e \ge 2$ and $m \ge 2$. The next goal, roughly, is to glue these actions of $\UC(\gl_e)$ on $\DG(\gl_m)$ into an action of $\UC(\tgl_e)$ on $\DG(\tgl_m)$. Two important subtleties
now appear.

The 2-category of (diagrammatic) singular Soergel bimodules has one object for each finite (standard) parabolic subgroup of the Coxeter group in question. The analogy is this: when
studying category $\OC$ in finite type, one splits it into blocks corresponding to orbits of the Weyl group, and an orbit is called \emph{singular} if it the action of the Weyl group has
a nontrivial stabilizer, which is a parabolic subgroup. Two blocks are equivalent when they have the same stabilizer, so one can choose to ignore most blocks and just consider one block
of each type. Then one has a category which splits into blocks, one for each parabolic subgroup. When projective functors act on these specific blocks, they can be organized into a
2-category with one object for each parabolic subgroup.

However, in order to create a Kac-Moody action in affine type, one should consider not just one block of each stabilizer type, but all the blocks at once. For any two blocks, the set of
projective functors from one to another yield a Hom category which is equivalent to the Hom categories between the corresponding stabilizer subgroups in $\DG(\tgl_m)$. Thus the
2-category we actually construct our Kac-Moody action on has many more objects than $\DG(\tgl_m)$ does, but all of its morphism categories are already encapsulated in the structure of
$\DG(\tgl_m)$. Let us not introduce extra notation in this introduction; we refer to this larger $2$-category abusively as $\DG$. \footnote{This 2-category $\DG$, or rather the quotient analogous to Soergel modules, is the diagrammatic analog of the full category $\OC$ (of level $-m-e$) for $\tgl_m$.} This was the first subtlety.

One fact about diagrammatic singular Soergel bimodules is that their relations only depend on finite type parabolic subgroups. So, for any finite type $A$ Dynkin diagram inside the
affine Dynkin diagram of $\tgl_m$, we can use our computations above to construct a $2$-functor from $\UC(\gl_e)$ to the corresponding part of $\DG$. The second subtlety is that these
actions do not ``glue'' together in the most naive way, because the imaginary root (which we shall denote as $y$) appears as ``monodromy'' as the parabolic subgroups travel in circles
around the affine Dynkin diagram. The way to fix this is to modify the 2-category $\UC(\tgl_e)$, and use a slightly different notion of a 2-Kac-Moody action, an idea that was first
introduced by Mackaay-Thiel \cite{MacThi}.

The affine analog of Mackaay-Stosic-Vaz is the work of Mackaay-Thiel \cite{MacThi}. They construct an extended version of (non-singular) Soergel bimodules in affine type $A$, analogous
to the construction of the extended affine Weyl group from the affine Weyl group, by formally adding a rotation operator. They then construct a $y$-deformation $\UC(\tgl_e)_{[y]}$ of the
2-Kac-Moody category in affine type $A$, and reproduce the results of Mackaay-Stosic-Vaz in this setting. They work algebraically with polynomials, and they do not express the singular
Soergel bimodule portion of the story in diagrammatic language, because the technology was not available at the time.

We treat the affine setting in a slightly different fashion than Mackaay-Thiel. They construct an action of $\UC(\tgl_e)_{[y]}$ on extended singular Soergel bimodules, which allows them
to give a simple formula for the action of the dots (certain generators of the Kac-Moody $2$-category), and let the rotation operator do a lot of work. In Section \ref{S_KM_affine_A} we
construct an action of $\UC(\tgl_e)_{[y]}$ on $\DG$ (which has no rotation operator, but has many additional objects). We must define the action of the dots in a careful way, depending
on the actual block and not just its stabilizer, to account for the monodromy.

These actions are different but philosophically the same. Our proof is new, because we use Frobenius-style arguments rather than computations with polynomials.


\begin{remark} We should mention that we treat all $m$ and $e$ in a uniform way. Because Mackaay-Thiel dealt with (non-singular) Soergel bimodules, they assume $e > m \ge 3$ in order to
avoid categories which were entirely singular. They dealt with the $m=e \ge 3$ case in a separate paper \cite{MacThi2}, though this has a different flavor. We have no issues working with
arbitrary $e$ and $m$, although the $e=2$ case does have some additional peculiarities. \end{remark}

Ironically, for the remaining results in this paper, i.e. the equivalence between Soergel-theoretic and rep-theoretic categories, we will use the degenerate $y=0$ specialization of
$\DG$. The second subtlety mentioned above - the monodromy which appears when glueing finite type $A$ actions into an affine one, and the corresponding deformation $\UC(\tgl_e)_{[y]}$ of
the Kac-Moody category - is entirely unnecessary in the specialization $y=0$!

\begin{remark} We do expect the generic, $y \ne 0$ case to have future applications, although we do not use it in this paper. Although it does create some additional complications
(largely confined to Section \ref{S_KM_affine_A}) to discuss the generic case, most of these complications are already necessitated by our other constructions anyway, such as the first
subtlety, the construction of a 2-category $\DG$ with a larger set of objects. So we felt it was worthwhile to prove the generic case here, to save the trouble of reproducing this
painstaking work elsewhere. \end{remark}

In summary, we have produced an action of (a modified version of) $\UC(\tgl_e)$ on $\DG$, a modified version of diagrammatic singular Soergel bimodules for $\tgl_m$. Fixing a parabolic
subgroup attached to a composition $m_1 + \ldots + m_\ell = m$, and applying the general construction of parabolic categories $\OC$ mentioned in  \S\ref{sssec-makingparabolic} to
each block of $\DG$, one obtains a category $\VC_{\un{m}}$ with a categorical action of $\UC(\tgl_e)$ which categorifies the level $0$ representation $V_{\un{m}}$.

\begin{remark} \label{rmk:whatusedfromEWsing} This remark is for the experts. Because the paper \cite{EWSingular} is not yet available, we wish to make clear which results from that
paper we use. We write down a list of generators for $\DG$, and a non-exhaustive list of relations. We need the 2-category $\DG$ to have these generators and relations (possibly with
more relations) so that our 2-functor is defined. We also need $\DG$ to categorify the Hecke algebroid and the Soergel-Williamson Hom formula, for our statements about the Grothendieck
group to hold, and we need the morphisms between Bott-Samelson objects in $\DG$ to agree with the morphisms in $\SB$, so that we can bootstrap known results for $\SB$ into results for
$\DG$. \end{remark}

\subsection{Rigidity for Fock spaces}\label{SS_Fock_rigidity}

Proving that two categories are equivalent is usually a difficult task, but a recent philosophy is that it becomes much easier in the presence of highest weight (or standardly stratified) structures assuming that the two categories have
a common quotient category that is large enough.

A general approach here was proposed by Rouquier, \cite{rouqqsch}. The setting is as follows. Let $\F$ be a base field, $\Ringg$ a regular complete local $\F$-algebra, and $\K$ its fraction field. Suppose that we have two highest weight categories $\OC^1_{\Ringg},\OC^2_{\Ringg}$ over $\Ringg$. For $i = 1, 2$ we let $\OC^i_\F$ and $\OC^i_{\K}$ denote corresponding categories after base change. Suppose that $\OC^1_R$ and $\OC^2_R$ come equipped with quotient functors to the same category, $\pi^i_R:\OC^i_R\to \mathcal{C}_R$. Rouquier proved that there is a highest weight equivalence $\OC^1_R \to \OC^2_R$ intertwining the functors $\pi^i_R$ provided the following two conditions hold:
\begin{itemize}
\item[(R1)] Everything is easy after base change to the fraction field $\K$. That is, $\OC^i_{\K}$
is split semisimple, the functor $\pi^i_{\K}$ is an equivalence, and
under the identification $\operatorname{Irr}(\OC^1_{\F})\cong \operatorname{Irr}(\mathcal{C}_{\K})
\cong \operatorname{Irr}(\OC^2_{\F})$, there is a common highest weight order for $\OC^i_{\F}$.
\item[(R2)] The specializations $\pi^i_{\F}$ are {\it 0-faithful}, i.e., fully faithful on standardly filtered objects.
\end{itemize}
In particular, this implies an equivalence $\OC^1_{\F} \to \OC^2_{\F}$.

The reason why the approach works is as follows. (R1) implies that the images of the standard objects under the
quotient functors coincide: $\pi^1_R(\Delta^1_R(\lambda))=\pi^2_R(\Delta^2_R(\lambda))$ for all labels $\lambda$
of irreducible objects. (R2) together with the condition that $\pi^i_\K$ are equivalences imply
that $\pi^i_{R}$ are {\it 1-faithful}, meaning that $\Ext^k_{\OC^i_R}(M,N)\xrightarrow{\sim}
\Ext^k_{\mathcal{C}_R}(\pi^i_R M,\pi^i_R N)$ for $k=0,1$ and any standardly filtered objects $M,N$.
The coincidence of images of the standard objects together with 1-faithfulness imply that
the images of the indecomposable projectives agree, which establishes an equivalence $\OCat^1_R\xrightarrow{\sim}
\OCat^2_R$. So the purpose of considering deformations over $R$ is two-fold: this helps to show that
the images of standards coincide and also improves faithfulness.

If $\OC^i_{\F}$ are the categories of interest, then one proceeds to proving  an equivalence by finding deformations $\OC^i_{\Ringg}$ (with compatible highest weight structures)
which are generically split semisimple (i.e. $\OC^i_{\K}$ is split semisimple). Extensions of Rouquier's approach have led to proving that the categories $\OC$ for cyclotomic
Rational Cherednik algebras are equivalent to categorical truncations of affine parabolic categories $\OC$, see \cite{RSVV,VV_proof}.

\begin{remark} As a warning, our situation will differ from Rouquier's setup in several aspects. For example, the deformations we use will not be generically semisimple. This adds a considerable amount of pain and requires several new techniques. \end{remark}

The functors $\pi_{\F}^i$ should be thought of as something like a ``Soergel functor.'' The original Soergel functor $\VM$ was a functor $\OC_0 \to C\operatorname{-mod}$, where $\OC_0$ is the principal block of the
original category $\OC$ associated to a Weyl group, and $C$ is the corresponding coinvariant ring. The functor itself is defined as $\VM = \Hom(P,-)$ for a particular projective object
$P$, and $C = \End(P)^{\op}$. Note that $P$ is not a projective generator ($P$ is the projective cover of the antidominant (simple) Verma module, so it has a
unique simple quotient), so that $\VM$ is very far from being faithful in general.
However, $\VM$ is faithful on standardly filtered objects and on costandardly filtered objects,
because the socle of any standard object (and the head of any costandard object) is the direct sum of several
copies of the antidominant simple. A corollary of this faithfulness is that any projective is included into the sum of several
copies of $P$ with a standardly filtered cokernel. This, in turn, implies that $\VM$ is fully faithful on projective
objects.


Note that $\VM$ is not fully faithful on all standardly filtered objects, as is already clear in the case of $\sl_2$. This can be remedied to an extent: we can enlarge the projective
$P$ to a projective $\bar{P}$ by adding all projective covers corresponding to \emph{sub-antidominant simples}, which are the simple objects $L$ such that
$\Ext^1(L,\Delta)\neq 0$ for some standard $\Delta$. The resulting functor $\bar{\VM}:=\operatorname{Hom}_{\OC}(\bar{P},-)$ is now fully faithful on the standardly filtered objects, and is analogous to Rouquier's functor $\pi$.


In Section \ref{S_Fock_unique}, we prove a rigidity result for restricted deformed categorifications of Fock space, that is, categorifications of $\Fock_{\un{m}}(\leqslant n)$. We define
a \emph{restricted $\Ringg$-deformed highest weight $\tgl_e$-categorification of $\Fock_{\un{m}}$} in Section \ref{SS_Fock_restr_categ}, to be a categorification over the base ring
$\Ringg:=\FM[[z_0,\ldots,z_{\ell-1}]]$, with highest weight structures on various subcategories compatible with the functors from the categorical action, and such that the action of the
``dots'' in $\UC(\tgl_e)$ is compatible with the ring $\Ringg$ in a particular way. We prove that our Soergel-theoretic construction can be equipped with such a structure. We think of it
as the prototypical restricted Fock space categorification (coming at least in a roundabout way from the explicit generators and relations description of $\DG$), to which other
categorifications are to be compared. We also prove that the main rep-theoretic setup, modules over cyclotomic $q$-Schur algebras, has this structure when the components of $\un{m}$ are
sufficiently far apart: $m_1-m_2\gg m_2-m_3\gg\ldots \gg\ m_{\ell-1}-m_{\ell}\gg n$ (this is sufficient for all the numerical data we seek, see Section \ref{sec-application}).


Our approach to prove the equivalence of categories is similar in spirit to that of \cite{VV_proof}, but is considerably more involved (the categories considered in \cite{VV_proof}
behave much nicer than ours). In our situation, we again have two quotient functors (to be denoted by $\pi^i$ and $\bar{\pi}^i$) constructed as $\Hom(P^i,-)$ and
$\Hom(\bar{P}^i,-)$ for some special projective objects $P^i$ and $\bar{P}^i$. The assumption in Rouquier's setup that $\pi^1$ and $\pi^2$ have the same target category corresponds to
the statement that $\End(\bar{P}^1)$ and $\End(\bar{P}^2)$ can be identified.

The construction of these projective objects is based on the categorical Kac-Moody action
on $\OCat^i_{\Ringg}$. Fix an integer $n$ and let $\OC^i_{\Ringg}(n)$ denote the part of $\OC$ which categorifies the degree $n$ part of the Fock space. We can consider the projective
object $P^i_{\Ringg}=F^n P^i_{\Ringg}(\varnothing)$, where $P^i_{\Ringg}(\varnothing)$ corresponds to $\Ringg\in \OCat^i_{\Ringg}(0)\cong \Ringg\operatorname{-mod}$. This is a direct
analog of Soergel's $P$, or of the projective $P_{KZ}$ in the category $\OC$ for Rational Cherednik algebras that was used in \cite{rouqqsch}. As with the functor $\VM$, the
functor $\Hom_{\OCat^i_{\Ringg}}(P^i_{\Ringg},-)$ is not 0-faithful for some choices of $e,m_1,\ldots,m_\ell$. So we enlarge $P^i_{\Ringg}$ to a projective $\bar{P}^i_{\Ringg}$ by adding
certain projective objects. We note that $\End(P^1_{\Ringg})$ is naturally identified with $\End(P^2_{\Ringg})$, which is basically Rouquier's uniqueness theorem for minimal
categorifications. The analogous isomorphism for $\End(\bar{P}^i_{\Ringg})$ is more complicated. Our treatment of this isomorphism can be found in Section \ref{SS_Ext_quot_equiv}. It is
similar to what was done in \cite{VV_proof}, but is more complicated roughly because  there are more ``induction functors'' available in \cite{VV_proof}. This is logically the first instance where the
assumption (for cyclotomic $q$-Schur algebras) that the $m_i$'s are far apart is used.

Now that we have the quotient functors $\bar{\pi}^i_R: \OCat^i_\Ringg(n)\twoheadrightarrow \mathcal{C}_\Ringg$
to the same category, we could try to apply Rouquier's approach to produce an equivalence between 
$\OCat^1_{\Ringg}(n)$ and $\OCat^2_{\Ringg}(n)$. Bad news is that
(R1) doesn't hold (the generic specializations of our categories are not semisimple), while establishing
(R2) seems to be out of reach. We proceed roughly as follows. We relax the conditions of the equivalence
theorem requiring that
\begin{itemize}
\item[(i)] $\pi^i_{\F}$, the functor defined by $P^i_{\F}$, is faithful on standardly filtered objects and on costandardly filtered objects,
\item[(ii)] $\bar{\pi}^i_{\K}$ is fully faithful on standardly filtered objects.
\item[(iii)] $\bar{\pi}^1_{\K}(\Delta^1_{\K}(\lambda))=\bar{\pi}^2_{\K}(\Delta^2_{\K}(\lambda))$ for all labels (=$\ell$-multipartitions of $n$) $\lambda$.
\end{itemize}
The result establishing an equivalence under (roughly) these assumptions is in Section \ref{SS_equiv_thms}.
The logic of the proof is basically as follows. (i) and (ii) imply that $\bar{\pi}^i_R$ is $0$-faithful.
This together with (iii) shows that  $\bar{\pi}^1_{\Ringg}(\Delta^1_{\Ringg}(\lambda))=
\bar{\pi}^2_{\Ringg}(\Delta^2_{\Ringg}(\lambda))$. Despite the fact that the functors
$\bar{\pi}^i_{\Ringg}$ are not known to be 1-faithful we still manage to show (Theorem \ref{Thm:equi_techn}) that the images of the
indecomposable projectives under our quotient functors agree.

Establishing (i),(ii),(iii) is very technical and is done in Section \ref{SS_-1_faith} (condition (i)),
Sections \ref{SS_companion} and \ref{SS_0_faith} (condition (ii), this is an especially technical part), and Section \ref{SS_equiv_ell1} (condition (iii)). Prominent ingredients of the proofs include the combinatorics of crystals, and the abstract nonsense of categorical Kac-Moody actions and highest weight categories.

\subsection{Related work} \label{ssec-comparison}

That a connection between the classical representation theoretic categories in characteristic $p$ and a suitable version of the category of Soergel bimodules exists seems to be long
expected by the experts in the subject (such as Rouquier and Soergel). An immediate inspiration for our work was a lecture by Geordie Williamson given at MSRI in November 2014. There
Williamson stated his conjecture with Riche that there is an equivalence between the category of tilting objects in the principal block of $\operatorname{Rep}(G)$, where $G$ is a
semisimple algebraic group over an algebraically closed field $\F$, and a parabolic category of diagrammatic Soergel modules for the corresponding affine Weyl group. Recently, Riche and
Williamson proved their conjecture for $G=\GL_n(\F)$, \cite[Part 2]{RW}. That work is independent from ours.

Let us compare our approach and results to those of Riche and Williamson. Perhaps, the only common idea used is a connection between 2-Kac-Moody category and the category of Soergel
bimodules. However, the approaches are in some sense opposite. We equip a Soergel-type category with a 2-Kac-Moody action, while Riche and Williamson equip the principal block of
$\operatorname{\Rep}(G)$ with an action of the diagrammatic category of Soergel bimodules.

Our results are stronger in two ways. First, in order for the principal block to exist, one needs to impose the inequality $p \ge m$, which is not needed in our work. Second, they
restrict their attention to $\operatorname{Rep}(\GL_n(\F))$ (the categories of modules over the quantum groups should also be accessible by their techniques), while we also deal with
``higher level categories'', that is, categories like modules over cyclotomic $q$-Schur algebras, which categorify higher level Fock space.  In particular, our theorems recover many more decomposition numbers than Riche and Williamson can.

On the other hand, our rigidity techniques do not apply to the categories of rational representations
(essentially because they categorify level $0$ representations and those are not highest weight).
 In particular, at this point we cannot prove that the principal block of $\operatorname{Rep}(\GL_n(\F))$ is equivalent to the principal block in a Soergel parabolic category $\mathcal{O}$.

{\bf Acknowledgements} The first author is supported by NSF CAREER grant DMS-1553032, and by the Sloan Foundation.
The second author is supported by NSF grant DMS-1501558. We would like to thank Jon Brundan, Sasha Kleshchev,
Raphael Rouquier, Ben Webster and Geordie Williamson for stimulating discussions.

%% file: SoergelDiags2.tex
\section{Singular Soergel bimodules and diagrammatics} \label{sec-ssbim}

\subsection{The upgraded Chevalley theorem}

Let $(W,S)$ be a Coxeter system, so that $W$ is a Coxeter group, and $S$ a set of simple reflections. We refer to a subset $I \subset S$ as a \emph{parabolic subset}. The reflections $s_i$ for $i \in I$ generate the \emph{parabolic subgroup} $W_I$. When $W_I$ is finite, the subset $I$ is called \emph{finitary}. In this case, we let $w_I$ denote the longest element of $W_I$, and $\ell(I) = \ell(w_I)$ its length.

\begin{example} When $W$ is finite, every parabolic subset is finitary. When $W$ is an affine Weyl group, every proper parabolic subset is finitary. \end{example}

A \emph{realization} of a Coxeter system $(W,S)$ over a commutative base ring $\Bbbk$ is, roughly speaking, a free $\Bbbk$-module $\hg$ equipped with an action of $W$, and a choice of
simple roots $\a_s \in \hg^*$ and simple coroots $\a_s^\vee \in \hg$ for each $s \in S$, which are compatible with the $W$ action. See \cite[\S 3.1]{EWGr4sb} for a precise
definition. A realization of $(W,S)$ gives rise to a realization of $(W_I, I)$ for any parabolic subset $I$, with the same space $\hg$ equipped with a subset of the simple roots and coroots. 

\begin{example} Let $W = S_n$ with the usual choice $S$ of simple reflections. Let $\Bbbk = \ZM$ and let $\hg$ be spanned by $\{\ep_1, \ldots, \ep_n\}$, with a dual basis $\{x_1, \ldots, x_n\}$ in $\hg^*$. Let $\a_i = x_i - x_{i+1}$ and $\a_i^\vee = \ep_i - \ep_{i+1}$, for $1 \le i \le n-1$. This is the \emph{$\gl_n$-realization}, and can be defined after base change for any commutative base ring $\Bbbk$. \end{example}

Given a realization, one can construct the polynomial ring $R = \textrm{Sym}(\hg^*)$. This is a $\Bbbk$-algebra with a natural $W$ action, and is graded such that $\deg \hg^* = 2$. For
each parabolic subset $I$ one can consider the invariant subring $$R^I = R^{W_I} = \{f \in R \;|\; s_i(f) = f \textrm{ for all } i \in I \}.$$ We often use shorthand like $R^s$ for
$R^{\{s\}}$, and $R^{st}$ for $R^{\{s,t\}}$. When $I$ is finitary, the Chevalley theorem states that $R^I$ is also a polynomial ring (generated in various degrees determined by the group $W_I$), and $R$ is free of finite rank over $R^I$.

\begin{example} For the $\gl_n$-realization, we get $R = \Bbbk[x_1, \ldots, x_n]$ with the usual action of $S_n$. Parabolic subgroups $W_I$ have the form $S_{n_1} \times S_{n_2} \times \cdots \times S_{n_d}$ where $\sum n_k = n$. The subrings $R^I$ are generated by elementary symmetric polynomials of the appropriate subsets of $\{x_i\}$. \end{example}

In fact, under mild assumptions on the realization, the ring extension $R^I \subset R$ can be equipped with the structure of a graded Frobenius extension. We refer to this fact as the
\emph{upgraded Chevalley theorem}.

\begin{defn} A \emph{graded Frobenius extension of degree} $d$ between two commutative graded rings $A$ and $B$ is a homogeneous ring extension $A \subset B$ together with the
\emph{Frobenius trace map}, an $A$-linear map $\pa^B_A \co B \to A$, homogeneous of degree $-2d$. This data satisfies: \begin{itemize} \item $B$ is free over $A$ of finite rank. \item
The map $\pa^B_A$ is a perfect pairing. That is, $B$ has dual bases $\{b_i\}$ and $\{b_i^*\}$ as an $A$-module such that $\pa^B_A(b_i b_j^*) = \delta_{ij}$. \end{itemize} \end{defn}

Let $(1)$ denote the grading shift operator, so that the degree $0$ part of $M(1)$ is the degree $1$ part of $M$, for a graded module $M$. One major consequence of the existence of a
graded Frobenius extension is the fact that the functors of induction $\Ind_A^B$ and shifted restriction $\Res^B_A(d)$ are biadjoint (adjoint on both the left and the right), where $d$
is the degree of the extension. Henceforth, when we mention the restriction functor in the context of a graded Frobenius extension, we always mean the shifted restriction $\Res^B_A(d)$.
The choice of inclusion $A \to B$ and Frobenius trace map $B \to A$ determine the units and counits of adjunction.

As functors between module categories, $\Ind_A^B$ and $\Res^B_A(d)$ can be thought of as bimodules, acting via tensor product. More precisely, $\Ind_A^B = {}_A A_B$ as an
$(A,B)$-bimodule, in the sense that tensoring on the left with ${}_A A_B$ will take a $B$-module $M$ to $A \ot_B M \cong \Ind_A^B M$. Similarly, $\Res^B_A(d) = {}_B A_A(d)$ as a
$(B,A)$-bimodule.

\begin{defn} We say that the upgraded Chevalley theorem holds for a given realization if, whenever $J, I$ are finitary parabolic subsets and $J \subset I$, then $R^I \subset R^J$ is a
graded Frobenius extension of degree $\ell(I) - \ell(J)$. More precisely, we require the trace maps $\pa^J_I \co R^J \to R^I$ to be given by Demazure operators, as discussed in
\S\ref{subsec-frobhyper} below. We write $\Ind_I^J$ instead of $\Ind_{R^I}^{R^J}$, and similarly for $\Res^J_I$. \end{defn}

The following is a folklore theorem. It appears to be well-known, but we do not know the proper attribution for it.

\begin{thm} The upgraded Chevalley theorem holds for all realizations of all Coxeter groups so long as $\Bbbk$ is a field of characteristic zero. \end{thm}

In general, however, the upgraded Chevalley theorem may fail when $\Bbbk$ does not have the appropriate units (e.g. a field of finite characteristic). This paper is concerned only
with (finite and affine) type $A$, and the following theorem will suffice for our needs.

\begin{thm} For the $\gl_n$-realization over $\ZM$, or after base change over any commutative ring $\Bbbk$, the upgraded Chevalley theorem holds. \end{thm}

To prove this, one can explicitly construct dual bases (the Schubert basis and dual Schubert basis) for the extension $R^{S_n} \subset R$. This is the topic of Schubert calculus. One
can adapt this construction to examine all the extensions $R^I \subset R^J$. A similar, direct construction of dual bases (over $\QM$) for other Coxeter groups should exist, but we have
not been able to find it in the literature.

\begin{remark} Because of this theorem, we have chosen to work with the $\gl_n$-realization (and the $\tgl_n$-realization for affine type $A$) in this paper, where the upgraded
Chevalley theorem holds in any characteristic. The reader interested in generalizations will note that our constructions in \S\ref{S_KM_fin_A} and \S\ref{sec-quantumaction-affine} will
work for any realization where the upgraded Chevalley theorem holds, and where the braid relations hold (see \S\ref{subsec-frobhyper}). \end{remark}

\subsection{Compatible Frobenius hypercubes and positive roots} \label{subsec-frobhyper}

For a finite set $S$, the collection of subsets form the vertices of a hypercube, which is a graded poset. When $J \subset I \subset S$ and $I \setminus J$ consists of a single element, this inclusion corresponds to an edge of the hypercube. For any Coxeter group $(W,S)$, the poset $\Gamma(W,S)$ of finitary subsets of $S$ need not be a hypercube, but it is closed below in the sense that whenever $I \in \Gamma(W,S)$ and $J \subset I$ then $J \in \Gamma(W,S)$. Consequently, the graded poset $\Gamma(W,S)$ is a union of hypercubes.

Given a realization of a Coxeter system $(W,S)$, one has a hypercube of graded ring extensions. Namely, for each parabolic subset $I \subset S$ one has a graded ring $R^I$, and whenever
$J \subset I$ one has a ring extension $R^I \subset R^J$. Tacit in this construction is the fact that this hypercube is \emph{compatible}, in the sense that whenever $K \subset J
\subset I$, the composition of the inclusion map $R^I \to R^J$ with the inclusion map $R^J \to R^K$ is the inclusion map $R^I \to R^K$. Consequently, if $K \subset I$, one can describe
the inclusion $R^I \to R^K$ as a composition of the inclusions along any sequence of edges giving a path from $I$ down to $K$, obtained by removing one simple reflection at a time, and
the result will not depend on the path chosen.

These inclusions give rise to adjoint functors $\Ind$ and $\Res$, and the choice of inclusion map $R^I \to R^J$ fixes the unit and counit of this (one-sided) adjunction. Compatibility
says that whenever $K \subset J \subset I$, the functors $\Ind^K_I$ and $\Ind^K_J \circ \Ind^J_I$ are naturally isomorphic (and similarly for $\Res$), and these isomorphisms are
compatible with the units and counits of adjunction. As a consequence, if one is interested in compositions of induction and restriction functors, one need only consider compositions of
induction and restriction along an edge.

Now consider the subposet $\Gamma(W,S)$. When the upgraded Chevalley theorem holds, this becomes a graded poset of Frobenius extensions. Note that being a Frobenius extension is a
structure, not a property. The specific units and counits of adjunction between $\Ind_I^J$ and $\Res^J_I(\ell(I)-\ell(J))$ depend on the choice of the trace map $\pa^J_I$. In fact, we
can choose our trace maps such that these Frobenius extensions are also \emph{compatible} in the sense that whenever $K \subset J \subset I$, one has the equality $\pa^J_K \circ \pa^I_J
= \pa^I_K$. Consequently, we need only remember the Frobenius structure for edges $J \subset I$, and this will be enough data to recover all induction and restriction functors, and
their biadjoint structures.

Hypercubes of Frobenius extensions satisfying the compatibility axioms above were called \emph{compatible Frobenius hypercubes} in \cite[Definition 1.2]{EWSFrob}, and were studied in
that paper. We will return to the results of \cite{EWSFrob} in \S\ref{subsec-hyperdiag}, but let us mention now that the results of that paper will easily apply to any graded poset of Frobenius extensions with the analogous compatibility axioms, when the poset is a union of hypercubes. We abusively refer to the compatible Frobenius poset attached to $\Gamma(W,S)$ as the \emph{Soergel hypercube}, even though it is not a hypercube in general.

Let us discuss how the trace maps $\pa^I_J$ are defined. The choice of simple roots and coroots made in the realization fixes the Frobenius extension structure on $R^s \subset R$. The
corresponding trace map $R \to R^s$ is given by the \emph{Demazure} or \emph{divided difference operator} $\pa_s$, for which \[ \pa_s(f) = \frac{f - sf}{\a_s}. \] For any Frobenius
extension there is an invariant $\mu$ of degree $2d$ called the \emph{product-coproduct element}; for $R^s \subset R$, this invariant $\mu_s$ is equal to $\a_s$.

\begin{defn} We say that a realization satisfies the \emph{braid relations} if the operators $\pa_s$ for $s \in S$ satisfy the braid relations. \end{defn}

\begin{prop} The $\gl_n$-realization satisfies the braid relations. \end{prop}

When a realization satisfies the braid relations, one obtains a \emph{Demazure operator} $\pa_w$ for each $w \in W$ by iterating $\pa_s$ over a reduced expression for $w$. Letting $w_I$
be the longest element of a finite parabolic subgroup $W_I$, one obtains a Frobenius structure on $R^I \subset R$ with trace map $\pa_I = \pa_{w_I}$. When the realization satisfies the
braid relations, there is also a consistent notion of \emph{positive roots} for any finitary parabolic subgroup. These are the elements of $\hg^*$ of the form $w(\a_s)$ where $ws > w$.
The product-coproduct element $\mu_I$ for the Frobenius extension $R^I \subset R$ is the product of the positive roots for $I$. More generally, for $J \subset I$ finitary, the relative
longest element $w_I^J = w_I w_J^{-1}$ leads to a Frobenius structure on $R^I \subset R^J$ with trace map $\pa_I^J = \pa_{w_I^J}$, and the product-coproduct element $\mu_I^J$ is the
product of the positive roots of $I$ which are not roots of $J$. Although they are not technically equivalent data, a choice of Frobenius extensions is almost the same information as a
choice of positive roots.

\begin{remark} We have used the words ``choice of positive roots'' because, for realizations where the Demazure operators do not satisfy the braid relations, there is a genuine choice
involved. Demazure operators $\pa_I$ and positive roots are only well-defined up to scalar, and one must choose scalars carefully so that the Frobenius hypercube is compatible. See
\cite{ECathedral} for more on when the braid relations hold. \end{remark}

\subsection{Algebraic Soergel bimodules} \label{subsec-algsbim}

Generalizing ideas of Soergel, Williamson defined certain $(R^I,R^J)$-bimodules called \emph{singular Soergel bimodules}. Together, singular Soergel bimodules form an additive
2-category $\SSBim$, where the objects are the finitary parabolic subsets $I \subset S$, the $1$-morphisms in $\Hom(J,I)$ are the $(R^I,R^J)$-bimodules which are singular Soergel
bimodules, and the $2$-morphisms are bimodule maps. Composition of $1$-morphisms is given by tensor product.

Whenever $J \subset I$, the bimodules $\Ind^J_I$ and $\Res^J_I(\ell(I)-\ell(J))$ are singular Soergel bimodules. By definition, singular Soergel bimodules are obtained from these
induction and restriction bimodules by applying tensor products and grading shifts, and taking direct sums and direct summands.

Let us describe this construction in a more combinatorial fashion. Let $\un{K}$ denote a path through the Soergel hypercube. That is, $\un{K} = (K_1, \ldots, K_d)$ is a sequence of
finitary parabolic subsets, such that $K_i$ and $K_{i+1}$ differ by a single index (and thus either $K_i \subset K_{i+1}$ or $K_{i+1} \subset K_i$) for each $i$. To such a path $\un{K}$
one can define a $(K_1, K_d)$-bimodule $BS(\un{K})$ as the tensor product of the corresponding induction and (shifted) restriction functors along this path. Such a bimodule is a
\emph{singular Bott-Samelson bimodule}, and they form a $2$-category $\SBSBim$. The Karoubi envelope of the (additive graded closure of) $\SBSBim$ is $\SSBim$.

\begin{example} Let $S = \{s,t,u\}$ be the simple reflections in type $A_3$, and $\un{K} = (\{s\},\{s,t\},\{t\},\{t,u\})$. Then $BS(\un{K}) = R^s \ot_{R^{s,t}} R^t(2)$ as an
$(R^s,R^{t,u})$ bimodule. The grading shift $(2)$ arises from the restriction $\Res^t_{s,t}(2)$, because $\ell(sts) - \ell(t) = 2$. \end{example}

\begin{defn} \label{defn:SWCT} Let $\CC$ be a 2-category with objects $I$ in bijection with finitary parabolic subsets of $S$. We say that the \emph{Soergel-Williamson categorification theorem holds} for $\CC$ if the Grothendieck group of $\CC$ is identified with the so called {\it Hecke algebroid},
so that the indecomposable objects in $\Hom(I,J)$ are parametrized by double cosets
$W_J \setminus W / W_I$. Moreover, morphism spaces in $\Hom(I,J)$ are free left $R^J$-modules and free right $R^I$-modules, whose ranks are governed by a pairing in the Hecke algebroid
(see \cite[Theorem 7.2.2]{WillSingular}). This last property is usually called the Soergel-Williamson Hom Formula. \end{defn}

\begin{thm}[The Soergel-Williamson categorification theorem] \label{thm:SWCT} Suppose that the realization of the Coxeter system $(W,S)$ is faithful and sufficiently
nice.\footnote{Soergel \cite{Soer07} and Williamson \cite{WillSingular} proved their results when the realization is \emph{reflection faithul}, a stronger notion than faithful, and when the
base ring $\Bbbk$ is an infinite field of characteristic not equal to $2$. Libedinsky \cite{LibRR} generalized Soergel's results to certain other realizations of affine Weyl groups in
characteristic zero. It is not known precisely what conditions on a realization guarantee that the algebraic bimodule category satisfies the Soergel-Williamson categorification theorem.
} Then the Soergel-Williamson categorification theorem holds for $\SSBim$. \end{thm}

It is not particularly relevant to this paper what exactly the Hecke algebroid is, and consequently, we have been vague in our statement of Definition \ref{defn:SWCT}. The point is
that, for many interesting realizations, the Soergel-Williamson categorification theorem will not hold for $\SSBim$, and so we are not particularly interested in $\SSBim$. We included
this theorem to provide context for the diagrammatic construction to come.

Let us quickly mention the original work of Soergel. If one considers only paths $\un{K}$ in $\Gamma(W,S)$ of the form $$\un{K} = (\emptyset, \{s_1\}, \emptyset, \{s_2\}, \emptyset,
\ldots, \{s_d\}, \emptyset),$$ one obtains a monoidal category of $(R,R)$-bimodules. These are the \emph{Bott-Samelson bimodules} $\BSBim$, and the Karoubi envelope of the (additive
graded closure of) $\BSBim$ is the category of \emph{Soergel bimodules} $\SBim$. Though it is not immediately obvious, the categories $\SBim$ and $\Hom_{\SSBim}(\emptyset,\emptyset)$
are the same (as subcategories of $(R,R)$-bimodules). According to Theorem \ref{thm:SWCT}, $\SBim$ categorifies the Hecke algebra. This is an earlier result of Soergel \cite{Soer07},
which Williamson built on to prove Theorem \ref{thm:SWCT}.

\subsection{Diagrammatics for Frobenius hypercubes} \label{subsec-hyperdiag}

Thus we have defined algebraic $2$-categories $\SBSBim$ and $\SSBim$ (and monoidal categories $\BSBim$ and $\SBim$) attached to a realization of a Coxeter system $(W,S)$. As discussed
in the introduction, there are too many morphisms between (singular) Soergel bimodules when the realization is not nice and the Soergel-Williamson categorification theorem does not
hold, such as for affine Weyl groups in characteristic $p$. Instead, we will need a diagrammatic version of (singular) Bott-Samelson bimodules, a presentation of the generic morphisms by
generators and relations, whose Karoubi envelope will satisfy the Soergel-Williamson categorification theorem in broader generality.

For the monoidal category $\BSBim$, this diagrammatic presentation is due in general to Elias-Williamson \cite{EWGr4sb}, building on earlier work of Elias-Khovanov \cite{EKho}, Elias
\cite{ECathedral}, and Libedinsky \cite{LibLL}. The presentation of $\BSBim$ is not immediately relevant, though we will need to recall it in \S\ref{sec-cell}. For the 2-category
$\SBSBim$, this presentation is due to Elias-Williamson in forthcoming work \cite{EWSingular}. The input to our construction is the Soergel hypercube of Frobenius extensions, with its
choice of Frobenius structure.

In the paper \cite{EWSFrob}, given any compatible hypercube of Frobenius extensions, a diagrammatic $2$-category was constructed by generators and relations, which was meant to encode
morphisms between iterated tensor products of induction and restriction bimodules. Note that, because of compatibility, one need only induct and restrict along edges of the hypercube.
The paper also provides relations which hold for any compatible hypercube. In any given compatible hypercube, there may be more generators and more relations required. It is easy to adapt these results to any poset built as a union of hypercubes.

For realizations in (finite or affine) type $A$ where $\SSBim$ satisfies the Soergel-Williamson categorification theorem, it turns out that the generating morphisms from \cite{EWSFrob}
are sufficient to generate all the morphisms in $\SBSBim$. Consequently, these should be the generators of a diagrammatic category which encodes the generic morphisms. Moreover, there
is only one more (relatively complicated) family of relations that needs to be enforced, beyond the relations from \cite{EWSFrob}. In this section we define a $2$-category $\DG^\pre =
\DG^\pre(W,S)$, which is just the diagrammatic $2$-category constructed in \cite{EWSFrob}. In the
Section \ref{subsec-finiteAsingdiags}
we define the quotient $\DG$ of $\DG^\pre$, which imposes this additional
family of relations.

We fix a Coxeter group $(W,S)$ and its Soergel hypercube.

\begin{defn} A \emph{singular $S$-diagram} is a collection of decorated, colored, oriented $1$-manifolds (with boundary) embedded in the planar strip $\RM \times [0,1]$ (with boundary
embedded in $\RM \times \{0,1\}$). There is one oriented $1$-manifold colored by each $s \in S$. For $s \ne t \in S$, the corresponding $1$-manifolds intersect transversely, and do not
intersect on the boundary. A \emph{region} is a connected component of the planar strip minus the union of the $1$-manifolds. Each region is labeled by a finitary subset $I \subset S$. These
region labels must be compatible with the orientation as follows: if $I$ and $J$ are on opposite sides of a $1$-manifold labeled $s$, then $I = J \cup \{s\}$ and $J = I \setminus \{s\}$,
and $I$ is on the right hand side of the $1$-manifold. Finally, one may place polynomials as floating symbols (0-manifolds) in each region: a polynomial $f \in R^I$ may live in a
region labeled $I$. Sometimes polynomials are placed in boxes for heightened visibility.

Associated to a singular $S$-diagram is its \emph{bottom} (resp. \emph{top}) sequence. This is a sequence $\un{K}$ of finitary subsets of $S$, obtained by intersecting the diagram with
$\RM \times \{0\}$ (resp. $\RM \times \{1\}$) and reading off the region labels from left to right. Note that the adjacent terms in $\un{K}$ are subsets which differ by a single
element; thus one can view $\un{K}$ as a sequence of elements of $S$, together with a finitary subset $K_1$ which labels the leftmost region.

The \emph{degree} of a singular $S$-diagram is computed by summing the degrees of its individual constituents: clockwise cups and caps, counterclockwise cups and caps, crossings of
various orientations, and polynomials. Recall that $\ell(I)$ is the length of the longest element of $W_I$. The degrees of these building blocks are encoded in the following chart.

\begin{equation*} {
\labellist
\small\hair 2pt
 \pinlabel {Generator} [ ] at 42 111
 \pinlabel {Degree} [ ] at 111 111
 \pinlabel {$\ell(J) - \ell(L)$} [ ] at 111 89
 \pinlabel {$\ell(L) - \ell(J)$} [ ] at 111 65
 \pinlabel {$0$} [ ] at 111 40
 \pinlabel {$\ell(I) + \ell(L) - \ell(J) - \ell(K)$} [ ] at 111 18
 \pinlabel {$J$} [ ] at 20 83
 \pinlabel {$L$} [ ] at 30 92
 \pinlabel {$L$} [ ] at 20 60
 \pinlabel {$J$} [ ] at 30 68
 \pinlabel {$L$} [ ] at 9 40
 \pinlabel {$J$} [ ] at 20 45
 \pinlabel {$K$} [ ] at 20 35
 \pinlabel {$I$} [ ] at 30 40
 \pinlabel {$K$} [ ] at 10 16
 \pinlabel {$L$} [ ] at 20 21
 \pinlabel {$I$} [ ] at 20 11
 \pinlabel {$J$} [ ] at 30 16
\endlabellist
\centering
\ig{2}{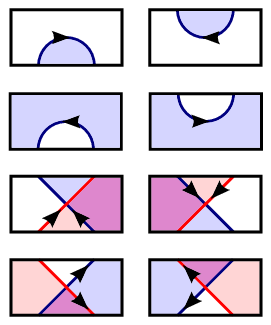}
} \end{equation*}

Above, white represents the label $L$, blue the label $J$, red the label $K$, and purple the label $I$; one has $L \subset J, K \subset I$. The more ``colors" a region has, the more
invariant a polynomial must be. Note that $\ell(I) + \ell(L) - \ell(J) - \ell(K)$ is non-negative, that clockwise cups and caps have positive degree, and that counterclockwise cups and caps have negative degree.

Two singular $S$-diagrams are \emph{isotopic} if there is a continuous family of singular $S$-diagrams which forms a homotopy from one to the other. Isotopic diagrams have the same
bottom and top. The degree of a singular $S$-diagram is an isotopy invariant. \end{defn}

We will often shorten the phrase ``singular $S$-diagram" to $S$-diagram or just diagram. Unless stated otherwise, we will abuse notation and say ``diagram" when we mean ``isotopy class
of diagrams."

\begin{notation} Let $L$ be a parabolic subset, and $\{s,t, \ldots, u\}$ a subset of $S$ which is disjoint from $L$. We write $Lst\ldots u$ for the disjoint union $L \cup \{s, t,
\ldots, u\}$. In particular, the notation $Lst$ will always imply that $\{s,t\}$ is disjoint from $L$. \end{notation}

\begin{defn} The \emph{connected components} of a parabolic subset are the connected components of the corresponding sub-Dynkin diagram. Let $K$ and $L$ be disjoint parabolic subsets. We say that $K$ is \emph{disconnected relative to $L$} if $K$ contains elements in at least two distinct connected components of $K \cup L$. For example, when $K = \{s,t\}$ and $L=\emptyset$, $K$ is disconnected relative to $L$ if and only if $s$ and $t$ are distant. When $L$ is understood, we simply say that $K$ is \emph{relatively disconnected}. \end{defn}

Consider the sideways crossing above, having degree $\ell(I) + \ell(L) - \ell(J) - \ell(K)$. Suppose that $I = Lst$ and that $\{s,t\}$ is disconnected relative to $L$. Then the sideways crossing has degree $0$. In fact, it will be an isomorphism in the category we are about to define. This will be a great simplifying factor.

Now we impose relations on $S$-diagrams.

\begin{defn} Let $\DG^\pre$ denote the following 2-category. The objects are finitary subsets $I \subset S$. The
$1$-morphisms are sequences $\un{K} = K_1 \cdots K_d$ of finitary subsets such that $K_i$ and $K_{i+1}$ differ by a single element of $S$. It is viewed as a $1$-morphism from $K_d$ to $K_1$.
The $2$-morphisms from $\un{K}$ to $\un{L}$ are the $\Bbbk$-span of (isotopy classes of) singular $S$-diagrams with bottom $\un{K}$ and top $\un{L}$, modulo the relations below.
$2$-morphisms are graded by the degree of the diagram (and all the relations below are homogeneous). \end{defn}

The first relations state that polynomials behave the way one would expect. That is, juxtaposition is multiplication, and polynomials slide freely across induction/restriction functors from
a smaller to a bigger ring. In these relations, the labels on the regions are arbitrary, so long as they are consistent.
\begin{subequations}
\begin{equation} \label{polymult} {
\labellist
\small\hair 2pt
 \pinlabel {$f$} [ ] at 14 27
 \pinlabel {$g$} [ ] at 34 27
 \pinlabel {$fg$} [ ] at 88 27
\endlabellist
\centering
\ig{1}{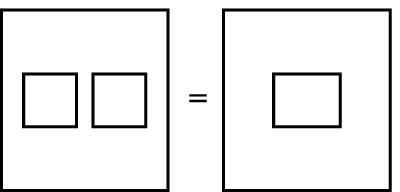}
} \end{equation}
\begin{equation} \label{slidepoly} {
\labellist
\small\hair 2pt
 \pinlabel {$f$} [ ] at 18 19
 \pinlabel {$f$} [ ] at 145 19
\endlabellist
\centering
\ig{1}{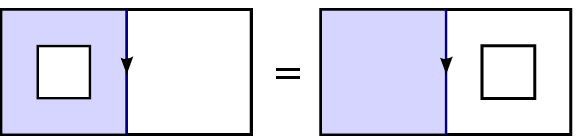}
} \end{equation}
\end{subequations}
Note that, in order for \eqref{slidepoly} to make sense, $f$ must be an element of the smaller, more invariant subring.

The next relations state that the cups and caps behave like the four canonical maps associated to a Frobenius extension. For more information on these canonical maps, see \cite{EWSFrob}.
White represents the label $J$ and blue the label $I$, for $J \subset I$.
\begin{subequations}
\begin{equation} \label{cwcirc} {
\labellist
\small\hair 2pt
 \pinlabel {$\mu_I^J$} [ ] at 70 20
\endlabellist
\centering
\ig{1}{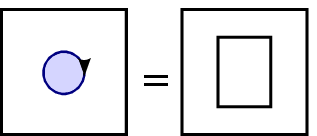}
} \end{equation}
\begin{equation} \label{ccwcirc} {
\labellist
\small\hair 2pt
 \pinlabel {$f$} [ ] at 21 21
 \pinlabel {$\pa_I^J(f)$} [ ] at 83 21
\endlabellist
\centering
\ig{1}{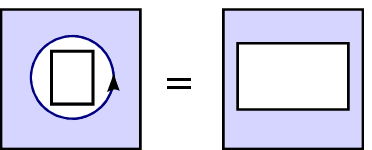}
} \end{equation}
\begin{equation} \label{coproduct} {
\labellist
\small\hair 2pt
 \pinlabel {$\D_I^J$} [ ] at 86 25
\endlabellist
\centering
\ig{1}{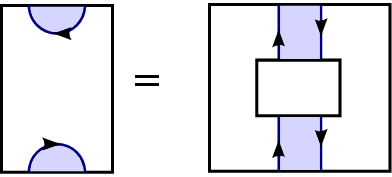}
} \end{equation}
Here, the tensor product $\D_I^J$ indicates a sum of diagrams, each of which has polynomials (corresponding to dual bases) in the two white regions.

To give an example of \eqref{coproduct}, consider the $\gl_n$-realization when $J = \emptyset$ and $I = \{s_i\}$. One can choose dual bases for $\pa_{s_i}$ to be $\{x_i, 1\}$ and $\{1, -x_{i+1}\}$. Thus the coproduct element is $x_i \ot 1 - 1 \ot x_{i+1} \in R \ot_{R^{s_i}} R$, and one would obtain the same coproduct element even had one chosen different dual bases. Then \eqref{coproduct} reads as follows.

\[ {
\labellist
\small\hair 2pt
 \pinlabel {$x_i$} [ ] at 69 27
 \pinlabel {$-$} [ ] at 123 27
 \pinlabel {$x_{i+1}$} [ ] at 180 27
\endlabellist
\centering
\ig{1}{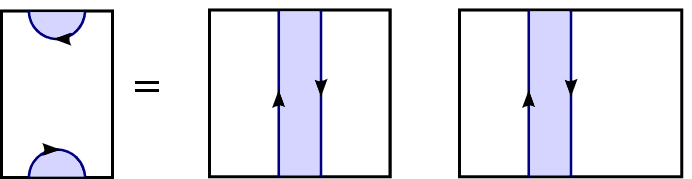}
} \]

One implication of \eqref{ccwcirc} is
\begin{equation} \label{ccwcircempty} \ig{1}{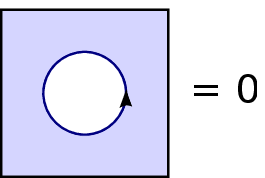}. \end{equation}
\end{subequations}

The next relations deal with a square of Frobenius extensions. Fix a parabolic subset $L$, colored white, and suppose that $I = Lst$ is purple, $J = Ls$ is blue,
and $K = Lt$ is red. A red (resp. blue) 1-manifold will separate two regions which differ in $s$ (resp. $t$).
\begin{subequations}
\begin{equation} \label{R2oriented} \ig{1}{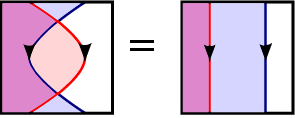} \end{equation}
\begin{equation} \label{R2nonoriented1} 	{
	\labellist
	\small\hair 2pt
	 \pinlabel {$\pa \D^L_I$} [ ] at 68 17
	\endlabellist
	\centering
	\ig{1.2}{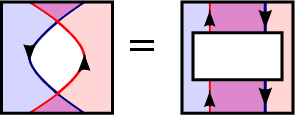}
	} \end{equation}
\begin{equation} \label{R2nonoriented2} 	{
	\labellist
	\small\hair 2pt
	 \pinlabel {$\mu^{J,K}_I$} [ ] at 72 17
	\endlabellist
	\centering
	\ig{1.2}{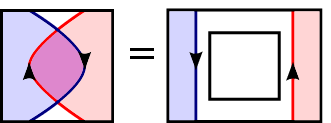}
	} \end{equation}
The symbol $\pa \D^L_I$ in \eqref{R2nonoriented1} represents the element \[\pa \D^L_I = \pa^L_J(\D^K_{I (1)}) \ot \D^K_{I (2)} = \D^J_{I (1)} \ot \pa^L_K (\D^J_{I (2)}) \in R^J \ot_{R^I}
R^K, \] using Sweedler notation. The (non-obvious) equality between these two expressions was proven in \cite{EWSFrob}. As in \eqref{coproduct} above, this use of $\pa \D^L_I$ is shorthand for a sum of diagrams where polynomials appear in the blue and red regions. The symbol $\mu^{J,K}_I$ in \eqref{R2nonoriented2} represents the product of the positive roots in $I$ which are in neither $J$ nor $K$. One has $\mu^{J,K}_I = \frac{\mu_I
\mu_L}{\mu_J \mu_K}$.

When $\{s,t\}$ is disconnected relative to $L$, then both sides of \eqref{R2nonoriented1} and \eqref{R2nonoriented2} must have degree $0$. In fact, $\pa \D^L_I = 1 \ot 1$ and
$\mu^{J,K}_I = 1$.

Some comments about isotopy are in order. In the definition above we have asserted that two singular $S$-diagrams that are isotopic represent equal $2$-morphisms in $\DG^\pre$. The two
sides of \eqref{R2oriented} are not isotopic as singular $S$-diagrams, because the isotopy would have to pass through a diagram where the red and blue $1$-manifolds do not intersect
transversely. Thus the relation \eqref{R2oriented} is not just an implication of isotopy invariance. Let us now discuss a more general notion of isotopy.

\begin{defn} Two singular $S$-diagrams are \emph{weakly isotopic} if the underlying oriented $1$-manifolds of each color are isotopic (independently). \end{defn}

It is a standard result in topology that two weakly isotopic $S$-diagrams are related by a sequence of simple transformations which we call \emph{Reidemeister moves} in analogy with the
corresponding isotopies of braids. The \emph{RII move} is illustrated by \eqref{R2oriented}: it allows two differently-colored $1$-manifolds to be ``pulled apart'' (if going from the
left hand side to the right hand side), or ``crossed over'' (if going in the other direction). The \emph{RIII move} is illustrated below by \eqref{R3oriented}. These are the only
transformations required.

Unfortunately, weakly isotopic diagrams need not be equal in $\DG^\pre$, as evidenced by \eqref{R2nonoriented1} and \eqref{R2nonoriented2}. In particular, the relation
\eqref{R2oriented} only allows one to apply the RII move when the strands are oriented in the same direction. When the strands are oppositely-oriented, \eqref{R2nonoriented1} and
\eqref{R2nonoriented2} allow one to pull strands apart at the cost of creating (linear combinations of diagrams with) polynomials in various regions. Neither relation makes it easy to
cross strands back over again, unless the requisite (linear combinations of) polynomials happen to be present already. However, as mentioned above, when the two strand colors $s$ and
$t$ are disconnected relative to the ambient parabolic subgroup $L$, the polynomials in question are trivial, and the strands may be pulled apart or crossed over without any trouble.

We call \eqref{R2oriented}, and applications of \eqref{R2nonoriented1} and \eqref{R2nonoriented2} in the relatively disconnected case, by the name of \emph{unobstructed RII moves}.
These are examples of weak isotopies which actually produce equal 2-morphisms in $\DG^\pre$. In contrast, we refer to other applications of \eqref{R2nonoriented1} and
\eqref{R2nonoriented2} as \emph{obstructed RII moves}, because the weak isotopy only leads to an equality of $2$-morphisms when particular polynomials are present.

We pause to note a consequence of the preceding relations, which is akin to \eqref{R2nonoriented1}.
\begin{equation} \label{R2nonoriented1var} {
\labellist
\small\hair 2pt
 \pinlabel {$f$} [ ] at 21 16
 \pinlabel {$\D^J_{I (1)}$} [ ] at 81 16
 \pinlabel {$\pa^L_K(f \D^J_{I (2)})$} [ ] at 149 16
\endlabellist
\centering
\ig{1.2}{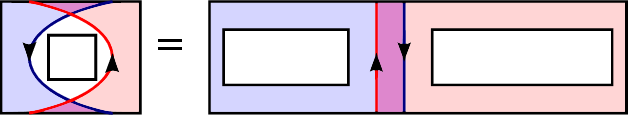}
}
\end{equation}
\end{subequations}

Because of this relation, we can always use (obstructed or unobstructed) RII moves to pull apart two strands, even in the presence of polynomials. However, it is difficult to cross over
two strands.

The final relations deal with a cube of Frobenius extensions. We have not bothered to color the regions.
\begin{subequations}
\begin{equation} \label{R3oriented} \ig{1}{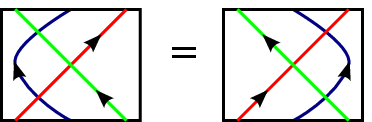} \end{equation}
\begin{equation} \label{R3unoriented} \ig{1}{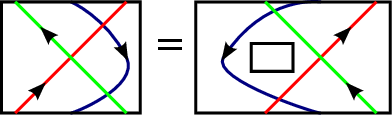} \end{equation}
\end{subequations}
The polynomial in the box in \eqref{R3unoriented} is the product of all the positive roots in the biggest parabolic subset $I = Lrbg$ (containing red, blue, and green, and an ambient parabolic
subset $L$) which are not positive roots in any of the other parabolic subgroups (namely $Lrb$, $Lrg$, and $Lbg$). 
Similar to $\mu^{J,K}_I$ above, it can be written as a ratio using the inclusion-exclusion principle.

In analogy with the discussion above, the RIII move (a weak isotopy) may hold in $\DG^\pre$ with or without obstruction. Relation \eqref{R3oriented} allows one to apply an RIII move
without any polynomials (that is, unobstructed), when the strands are all oriented the same direction. Relation \eqref{R3unoriented} involves polynomials, so it may be obstructed.
However, when the triple $\{s,t,u\}$ of strand colors is disconnected relative to $L$, then the polynomial in \eqref{R3unoriented} must be degree zero, and in fact it is $1$.
Thus the RIII move is unobstructed in the relatively disconnected case.

This ends the definition of $\DG^\pre$. Let us make one further comment before discussing the connection with bimodules.

\begin{defn} Given a singular $S$-diagram, let $L$ be the intersection of all the parabolic subsets labeling each region. We call $L$ the \emph{ambient parabolic subset}. Let $K$ denote
the collection of all the strand labels, so that $L \coprod K$ is the union of all the parabolic subsets labeling each region. Let $J \subset L$ be the set of all simple reflections
which do not share a connected component in $L \coprod K$ with any element of $K$. We call $J$ the \emph{irrelevant parabolic subset}. We will usually refer to $L \coprod K \setminus J$ as the \emph{interesting parabolic subset}. \end{defn}

\begin{lemma} \label{lem:removeirrelevant} A relation holds between singular $S$-diagrams if and only if the same relation would hold if the irrelevant parabolic subset were removed
from every region label. \end{lemma}

To give an example, suppose $s$ and $t$ are distant. The Frobenius structure map for $R^s \subset R$ is $\pa_s$, and $\mu_s = \a_s$. The Frobenius structure map for $R^{s,t} \subset
R^t$ is also $\pa_s$, restricted from $R$ to $R^t$, and $\mu^t_{st} = \a_s$ (now viewed as an element of $R^t$). The key point here is that adding or removing an irrelevant subset does
not ``change" any of the Frobenius structure maps, so it does not affect any of the relations above. This lemma will be proven in \cite{EWSingular}.

\subsection{The additional relations in finite type $A$} \label{subsec-finiteAsingdiags}

There is a $2$-functor from $\DG^\pre$ to $\SBSBim$, sending $1$-morphisms to the corresponding induction or restriction bimodules, sending cups and caps to the corresponding units or
counits of biadjunction, and sending crossings to the isomorphisms between iterated induction operators. That this $2$-functor is well-defined is the content of \cite{EWSFrob}.
Henceforth, whenever we mention induction or restriction, we refer to $1$-morphisms in $\DG^\pre$ or their image in $\DG$ (defined below), not the actual algebraic bimodules.

The functor from $\DG^\pre$ to $\SBSBim$ is full when $\SBSBim$ satisfies the Soergel-Williamson categorification theorem \ref{thm:SWCT}. It is not faithful, but the kernel is described
by one additional family of relations in type $A$. Imposing this additional family, we obtain our $2$-category $\DG$. These facts will be proven in \cite{EWSingular}.

\begin{defn} Let $\DG$, the \emph{diagrammatic 2-category of singular Bott-Samelson bimodules in type $A$}, be the quotient of $\DG^\pre$ by one more family of relations, partially
described below. The Karoubi envelope of $\DG$ is the \emph{diagrammatic 2-category of singular Soergel bimodules}. 
\end{defn}

The missing relation categorifies the following relation in the Hecke algebroid, which we draw here using MOY diagrams.

\begin{equation} \label{MOYrelation}
\begin{array}{c}
\begin{tikzpicture}[scale = 0.7,bull/.style={inner sep=0pt,minimum size=1mm,draw,circle,fill}]
\node [bull] (c) at (0,0) {};
\node [bull] (a) at (0,3) {};
\node [bull] (a+b) at (2,0) {};
\node [bull] (c+b) at (2,3) {};
\node [below] at (c) {$c$};
\node [above] at (a) {$a$};
\node [below] at (a+b) {$a+b$};
\node [above] at (c+b) {$c+b$};
\draw (c) -- (0,0.6) -- (1,1.1) -- (1,1.9) -- (0,2.4) -- (a);
\draw (a+b) -- (2,0.6) -- (3,1.1) to node[right] {$b$} (3,1.9) -- (2,2.4) -- (c+b);
\draw (1,1.1) -- (2,0.6) -- (3,1.1) -- (3,1.9) -- (2,2.4) -- (1,1.9);
\end{tikzpicture} \end{array}
=
\sum_{i = 0}^{\min(a,c)}
\qchoose{b}{i}
\begin{array}{c}
\begin{tikzpicture}[scale = 0.7,bull/.style={inner sep=0pt,minimum size=1mm,draw,circle,fill}]
\node [bull] (c) at (0,0) {};
\node [bull] (a) at (0,3) {};
\node [bull] (a+b) at (2,0) {};
\node [bull] (c+b) at (2,3) {};
\node [below] at (c) {$c$};
\node [above] at (a) {$a$};
\node [below] at (a+b) {$a+b$};
\node [above] at (c+b) {$c+b$};
\draw (c) -- (0,0.6) -- (1,1.1) -- (1,1.9) -- (0,2.4) -- (a);
\draw (a+b) -- (2,0.6) -- (1,1.1) -- (1,1.9) -- (2,2.4) -- (c+b);
\draw (1,1.1) -- (1,1.9) -- (0,2.4) -- (-1,1.9) to node[left] {$i$}  (-1,1.1) -- (0,0.6);
\end{tikzpicture}
\end{array}
\end{equation}
For the reader unfamiliar with MOY diagrams, the diagram above can be read as follows. A horizontal slice is a sequence of numbers adding to $n = a+b+c$, giving a parabolic subgroup of $S_n$. (However, this relation can be applied locally, so one need not assume $a+b+c=n$ in general.) Reading from bottom to top, a merge (resp. split) represents a $1$-morphism in the Hecke algebroid, categorified by the corresponding restriction (resp. induction) bimodule. Both sides of \eqref{MOYrelation} should be thought of as describing $1$-morphisms in $\DG$, so that \eqref{MOYrelation} is categorified if the $1$-morphism on the left hand side is a direct sum of the $1$-morphisms on the right hand side.

The general relation in $\DG$ looks something like this:
\begin{equation} \label{relationmysterious} {
\labellist
\tiny\hair 2pt
 \pinlabel {$\displaystyle \sum$} [ ] at 74 29
 \pinlabel {$\displaystyle + \sum$} [ ] at 161 30
 \pinlabel {$g_\beta$} [ ] at 115 11
 \pinlabel {$g_\beta^*$} [ ] at 115 42
 \pinlabel {$h_\gamma$} [ ] at 204 15
 \pinlabel {$h_\gamma^*$} [ ] at 204 37
\endlabellist
\centering
\ig{1.3}{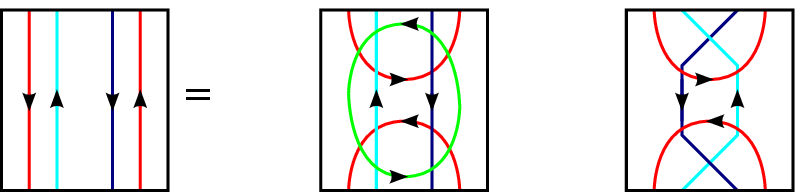}
} \end{equation}
In this relation, the identity morphism on the left hand side is rewritten as a sum of orthogonal idempotents. Depending on the various region labels, the terms on the right hand side will have differently colored ``green" circles, or perhaps no ``green" circle at all. For the terms of a given type, the polynomials $g$ or $h$ are dual bases for some particular combination of Frobenius operators, in order that the diagrams are orthogonal idempotents. One advantage of the $\gl_n$-realization is that these dual bases can be realized explicitly using Schur polynomials. However, the particular choice of dual bases is irrelevant. A degenerate version occurs when the blue and aqua colors are the same.
\begin{equation} \label{relationmysteriousdegen} {
\labellist
\small\hair 2pt
 \pinlabel {$\displaystyle \sum$} [ ] at 74 29
 \pinlabel {$\displaystyle + \sum$} [ ] at 161 30
 \pinlabel {\tiny $g_\beta$} [ ] at 115 11
 \pinlabel {\tiny $g_\beta^*$} [ ] at 115 42
 \pinlabel {\tiny $h_\gamma$} [ ] at 204 15
 \pinlabel {\tiny $h_\gamma^*$} [ ] at 204 37
\endlabellist
\centering
\ig{1.3}{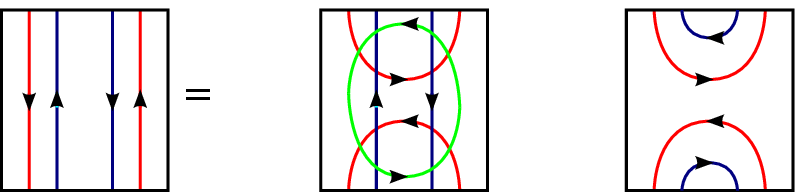}
} \end{equation}

To make this family of relations precise is not actually necessary for this paper. Instead, we point out the two special instances of this relation that we do use.

The first case we will use is the dihedral relation, which was also described in \cite[(6.8)]{ECathedral} (the $m=3$ case). \begin{equation} \label{dihedralrelation} \ig{1}{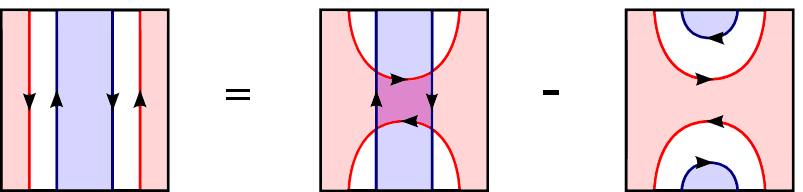}
\end{equation} Here, the white region is labeled by $\emptyset$, the red and blue regions are labeled by adjacent simple reflections $s$ and $t$, and the purple region is labeled by
$st$. We use a consequence of \eqref{dihedralrelation}, namely \cite[Claim 6.5]{ECathedral}, in the computation \eqref{NilHecke3} below, though we do not use \eqref{dihedralrelation} directly.

The second special case occurs when $a=b=1$ in the MOY diagrams above. We describe this relation after having removed the irrelevant parabolic subset, using Lemma \ref{lem:removeirrelevant}; this is like assuming that $a+b+c=n$. So, let $L$ denote a type $A$ configuration, such that $\{s,t,L,z\}$ is (in order) a larger $A_n$ configuration.
Then our relation takes the form below.
\begin{equation} \label{relation11n1} {
\labellist
\tiny\hair 2pt
 \pinlabel {$Lsz$} [ ] at 45 37
 \pinlabel {$Lz$} [ ] at 36 21
 \pinlabel {$Ltz$} [ ] at 23 37
 \pinlabel {$Lt$} [ ] at 11 22
 \pinlabel {$Lst$} [ ] at 4 37
 \pinlabel {$Lsz$} [ ] at 111 29
 \pinlabel {$Lstz$} [ ] at 92 27
 \pinlabel {$Lst$} [ ] at 74 28
 \pinlabel {$Ltz$} [ ] at 92 11
 \pinlabel {$Ltz$} [ ] at 92 45
 \pinlabel {$Lsz$} [ ] at 177 28
 \pinlabel {$Ls$} [ ] at 160 28
 \pinlabel {$Lst$} [ ] at 141 28
 \pinlabel {$Ltz$} [ ] at 161 4
 \pinlabel {$Ltz$} [ ] at 160 51
 \pinlabel {$L$} [ ] at 160 39
 \pinlabel {$L$} [ ] at 160 16
\endlabellist
\centering
\ig{2}{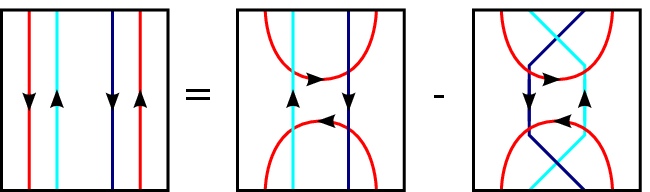}
} \end{equation}
We use this relation when we check the most complicated triple crossing relation in the KLR algebra.


The main theorem of \cite{EWSingular} will be:

\begin{thm} \label{thm:SWCTdiag} The $2$-category $\DG$ satisfies the Soergel-Williamson categorification theorem. That is, its Karoubi envelope categorifies the Hecke algebroid, and its
$2$-morphism spaces are free modules of the appropriate size. The $2$-functor $\DG^\pre \to \SBSBim$ descends to an essentially surjective, faithful (but not necessarily full)
$2$-functor $\DG \to \SBSBim$. This $2$-functor is an equivalence if and only if (the Karoubi envelope of) $\SBSBim$ satisfies the Soergel-Williamson categorification theorem. \end{thm}

\subsection{Diagrammatics in affine type $A$} \label{subsec-diagaffinetype}

Let $(W,S)$ denote the finite type $A_{n-1}$ Coxeter system. Let $(W_a,S_a)$ denote the affine type $A$ Coxeter system, with $S_a = S \cup \{s_0\}$.

\begin{defn} We now describe the \emph{$\tgl_n$-realization} of $(W_a,S_a)$, over the base ring $\Bbbk = \ZM$. The free $\Bbbk$-module $\hg^*$ has basis $\{x_i\}_{i = 1}^n \cup \{y\}$,
while $\hg$ has the dual basis $\{\ep_i\}_{i=1}^n \cup \{\phi\}$. The simple roots $\a_i$ for $i \in S$ are still given by $\a_i = x_i - x_{i+1}$, while the simple root $\a_0$ is given
by $x_n - x_1 + y$. The simple coroots for $i \in S$ are still given by $\a_i^\vee = \ep_i - \ep_{i+1}$, and the simple coroot $\a_0^\vee$ is given by $\ep_n - \ep_1$. With
these definitions, $R = \ZM[x_1, \ldots, x_n, y]$ where $y$ is fixed by $W_a$, $S_n$ acts as usual on the variables $x_i$, and $s_0(x_n) = x_1 - y$, $s_0(x_1) = x_n + y$. \end{defn}

The realization above looks asymmetric in that $\a_0 = x_n - x_1 + y$ while $\a_0^\vee = \ep_n - \ep_1$. It is more common to work in a larger realization with an extra generator $a \in
\hg^*$ with dual $b \in \hg$, such that $\a_0^\vee = \ep_n - \ep_1 + b$. Instead, we work in ``level 0'' by setting $b = 0$. This larger realization does not affect the ring $R$ or its
Demazure operators in any interesting way (c.f. \cite[Sectios 2.1]{MacThi}; they work with the larger realization, but have the same formulas for the $W_a$-action on $R$), so we
have chosen the smaller realization above.

\begin{remark} The Demazure operators for this realization satisfy the braid relations, so there is a well-defined notion of positive roots for each finite parabolic subgroup, and a
standard Frobenius trace map $\pa_I \co R \to R^I$ whenever $I$ is finitary. Moreover, one can explicitly construct dual bases for $\pa_I$ (again using Schubert calculus), showing that
the upgraded Chevalley theorem holds. \end{remark}

We now wish to use the results of \cite{EWSFrob} to construct a diagrammatic 2-category $\DG_a^\pre$ attached to the Soergel hypercube of the $\tgl_n$-realization. As previously
mentioned, the subposet $\Gamma(W_a,S_a)$ of finitary parabolic subgroups is a union of hypercubes. To obtain the corresponding $2$-category $\DG_a^\pre$, one just considers the
generating morphisms and relations to be those which lie in any hypercube inside the poset $\Gamma(W_a,S_a)$.

Analogously, in affine type $A$, \cite{EWSingular} will define the category of \emph{diagrammatic singular Bott-Samelson bimodules} $\DG_a = \DG(W_a,S_a)$ to be the quotient of
$\DG^\pre_a$ by the relations of \S\ref{subsec-finiteAsingdiags} for each finite type $A$ sub-hypercube. This will be equivalent to the algebraic category of
singular Bott-Samelson bimodules whenever the latter behaves nicely. 

%% file: FiniteADiags2.tex
\section{The Kac-Moody action in finite type $A$}\label{S_KM_fin_A}

Fix $n \ge 2$ and $e \ge 2$.\footnote{When comparing to \cite{MSV}, their $n$ is our $e$, and their $d$ is our $n$. When comparing to \cite{KhoLau10}, their $n$ is our $e$, and their $N$ is our $n$. We will use our conventions throughout, even when discussing their results.} We now move towards defining a $2$-functor $\FC \co \UC(\gl_e) \to \DG(\gl_n)$.

\subsection{Previous constructions} \label{subsec-previous}

The goal of this section is to explain the following commutative diagram of $2$-functors.

\begin{equation} \label{eq:bigdiagram} {
\labellist
\small\hair 2pt
 \pinlabel {$\UC_{\leftarrow}(\gl_e)$} [ ] at 60 130
 \pinlabel {$\UC_{\leftarrow}(\sl_e)$} [ ] at 120 130
 \pinlabel {$\SC(e,n)$} [ ] at 60 84
 \pinlabel {$\DG(\gl_n)$} [ ] at 60 44
 \pinlabel {$\SBSBim$} [ ] at 60 4
 \pinlabel {$\Flag_{n,e}$} [ ] at 120 4
 \pinlabel {$\SB$} [ ] at 0 44
 \pinlabel {$\BSBim$} [ ] at 0 4
 \pinlabel {$\psi_{e,n}$} [ ] at 93 103
 \pinlabel {$\Sigma_{e,n}$} [ ] at 27 72
 \pinlabel {$\FC$} [ ] at 54 61
 \pinlabel {$\FC_{MSV}$} [ ] at 103 45
 \pinlabel {$\FC_{KL}$} [ ] at 110 77
 \pinlabel {$\FC_{EK}$} [ ] at 12 26
 \pinlabel {$\FC_{EW}$} [ ] at 52 26
\endlabellist
\centering
\ig{1.5}{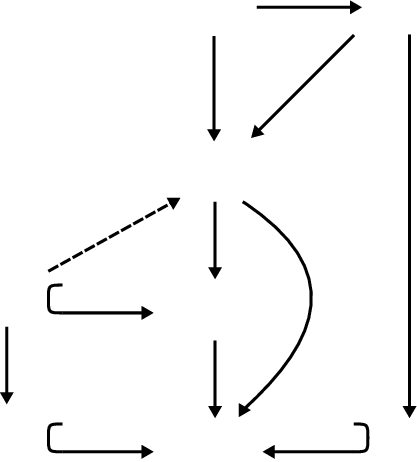}
} \end{equation}

Let $\Bim$ denote the 2-category whose objects are graded rings, and where $\Hom_{\Bim}(A,B)$ is the category of graded $(B,A)$-bimodules. The 2-categories in the bottom row of
\eqref{eq:bigdiagram} are all full sub-2-categories of $\Bim$, and thus are algebraic in nature. The 2-categories $\BSBim$ and $\SBSBim$ have already been defined in
\S\ref{subsec-algsbim}, and are attached to the $\gl_n$-realization. Accordingly, we let $R = \Bbbk[x_1, \ldots, x_n]$ be the polynomial ring of the $\gl_n$-realization, with base ring
$\Bbbk = \ZM$.\footnote{Note that $\BSBim$ is a monoidal category inside $R$-bimodules. In this diagram, monoidal categories are thought of as 2-categories with a single object. Thus
the map $\BSBim \to \SBSBim$ is an inclusion of 2-categories. Similarly, the map $\SB \to \DG$ is an inclusion of a monoidal category inside a 2-category.}

The 2-categories in the penultimate row of \eqref{eq:bigdiagram} are the diagrammatic versions of their algebraic analogs, as discussed in \S\ref{subsec-hyperdiag}. The functor
$\FC_{EK} \co \SB \to \BSBim$ is the monoidal functor from \cite{EKho}, and the 2-functor $\FC_{EW} \co \DG \to \SBSBim$ is the functor from \cite{EWSingular}, which sends a diagram to
the corresponding bimodule morphism. These 2-functors are equivalences when the algebraic categories behave nicely, such as after base change to $\Bbbk = \QM$.

Now let us recall what Khovanov-Lauda accomplish in \cite{KhoLau10}, which is the definition of the 2-functor $\FC_{KL}$ on the right side of \eqref{eq:bigdiagram}. The target of their
2-functor is the \emph{equivariant flag 2-category} $\Flag_{n,e}$. The objects of $\Flag_{n,e}$ are in bijection with sequences $$\un{k} = (0 = k_0 \le k_1 \le \ldots \le k_e = n)$$ of
non-decreasing integers, having length $e$. To this sequence one associates a (partial) flag variety, with $k_i$ representing the dimension of the $i$-th flag inside $\CM^n$. If $P =
P_{\un{k}}$ is the parabolic subgroup in $G=GL_n$ stabilizing the standard partial flag, then the $P$-equivariant cohomology of a point is none other than $R^{I}$, where $I=I_{\un{k}}$
parametrizes the simple reflections inside $P_{\un{k}}$. The $G$-equivariant cohomology of this flag variety is also $R^I$, now viewed as a module for the $G$-equivariant cohomology of
a point, which is $R^W$.

Let $\un{k}$ be a sequence as above, and let $\un{k}_{+i}$ denote the identical sequence except with $k_i$ replaced by $k_i+1$ (if possible). Consider the flag variety for the length
$e+1$ sequence $$(0 = k_0 \le k_1 \le \ldots \le k_i \le k_i + 1 \le k_{i+1} \le \ldots \le n).$$ This $(e+1)$-step flag variety admits forgetful maps to the $e$-step flag varieties for
$\un{k}$ and $\un{k}_{+i}$, and thus its (equivariant) cohomology ring becomes a bimodule over the cohomology rings of these two $e$-step flag varieties. Viewed as a functor from
$R^{I_{\un{k}}}$-modules to $R^{I_{\un{k}_{+i}}}$-modules, this is the action of $E_i$; viewed as a functor in the other way, this is the action of $F_i$ (we have ignored the grading
shifts). These bimodules and their tensor products are the 1-morphisms in $\Flag_{n,e}$. The 2-morphisms are bimodule maps. The Khovanov-Lauda 2-functor $\FC_{KL} \co
\UC_{\leftarrow}(\sl_e) \to \Flag_{n,e}$ is constructed explicitly using operations on polynomials.

We have just discussed where $\FC_{KL}$ sends the $1$-morphisms in $\UC(\sl_e)$, but neglected to discuss where it sends objects. There is a bijection between valid sequences $\un{k}$
as above, and a certain subset $\La_{n,e}$ of the $\gl_n$-weights or $\sl_n$-weights. See Section \ref{subsec-prelims-finite} for more details. A weight in this subset $\La_{n,e}$ is sent to
the corresponding sequence $\un{k}$, and a weight not in this subset is sent to zero. \footnote{Here, zero is a formal symbol. In an additive $2$-category, one may as well include an
object named zero, for which all morphism categories to and from zero are the zero category. To say that an object is sent to zero is the same as to say that all $1$-morphisms to or
from that object are sent to the zero $1$-morphism. This latter condition makes sense even when the zero object is not included.}

Let us quickly remark on what the subscript arrow means in $\UC_{\leftarrow}(\sl_e)$. In \cite[page 15]{KhoLau11} they define a diagrammatic category attached to an oriented
simply-laced Dynkin diagram. The arrow in $\UC_{\leftarrow}(\sl_e)$ indicates that the orientation on the Dynkin diagram of $\sl_e$ points from vertex $j+1$ to $j$, for $1 \le j \le
e-2$. The effect of this orientation is to determine a sign in the ``double-crossing'' and ``triple-crossing'' relations (c.f. \cite[equations (28), (33), and (34)]{KhoLau11}). Let us
ignore the arrows for the moment, to avoid cluttering the page. \footnote{The other reason we ignore the arrows temporarily is that there is an error in the literature, see \S\ref{subsec-signerror}.}

Now we segue to the work of Mackaay-Stosic-Vaz \cite{MSV}, who constructed a 2-functor $\FC_{MSV}$ from a quotient $\SC(e,n)$ of $\UC(\gl_e)$ to $\SBSBim$. The difference between
$\UC(\sl_e)$ and $\UC(\gl_e)$ is just the parametrization of the objects by either $\sl_e$-weights or $\gl_e$-weights. Otherwise the categories are identical, modulo the following
remark.

\begin{remark} \label{rmk:signerror} The original version of Khovanov and Lauda's paper \cite{KhoLau10} had a sign error in the definition of the category $\UC(\sl_e)$ relating to
certain cups and caps, which was noticed by Mackaay-Stosic-Vaz when they constructed their version in \cite{MSV}. See the erratum \cite{KhoLau10erratum}. There are two ways to fix this
sign error. In \cite{MSV}, one extends the action to $\UC(\gl_e)$, the main effect being that objects are labeled by $\gl_e$-weights rather than $\sl_e$-weights; this extra data allows
one to choose a sign for the cups and caps consistently. With the sign convention from \cite{MSV}, the category $\UC(\gl_e)$ is cyclic. Cyclicity in this context refers to the fact that
isotopy classes of diagrams all represent the same $2$-morphisms. Meanwhile, in \cite{CauLau}, Cautis and Lauda provide a different renormalization of the category $\UC(\sl_e)$ or
$\UC(\gl_e)$ which is not cyclic, but which acts correctly on the equivariant flag 2-category. Both approaches have their merits\footnote{There are certain degree zero diagrams, called
degree zero bubbles, which are philosophically thought of as the zeroth Chern class of certain vector bundles. One of the merits of Cautis and Lauda's normalization is that these degree
zero bubbles are set equal to the identity. In Mackaay-Stosic-Vaz's normalization, degree zero bubbles are set equal to a sign.}.

In the rest of the paper we will stick with these conventions of \cite{MSV}, working entirely with $\UC(\gl_e)$ and not $\UC(\sl_e)$. The main reason for this choice is because $\DG$ is
a cyclic 2-category, and thus it is easier to define a $2$-functor from another cyclic 2-category. We never bother to choose a sign convention for $\UC(\sl_e)$, and only discussed the
2-functor $\FC_{KL}$ for context. \end{remark}

Mackaay-Stosic-Vaz's description of the target of their $2$-functor goes as follows. The objects are still sequences $\un{k}$ as above, identified with the same rings $R^I$ as above.
Given a sequence $\un{k}$, the differences $k_i - k_{i-1}$ are the labels on the strands in a MOY diagram. However, their 2-category has more 1-morphisms, as their generating
1-morphisms come from splitting and merging maps in a MOY diagram: these are interpreted as bimodules coming from induction and restriction from $R^I$ to $R^J$ where $I$ and $J$ differ
by a single element. Thus the $1$-morphisms and $2$-morphisms in the categorified MOY setup are the same as for singular Bott-Samelson bimodules. In truth, Mackaay-Stosic-Vaz construct
a 2-functor $\FC_{MSV} \co \UC(\gl_e) \to \SBSBim(\gl_n)$, which they denote $\FC_{Bim}$ in their paper. It is constructed explicitly using operations on polynomials, analogously to the
construction of $\FC_{KL}$.

They define the \emph{categorified $q$-Schur algebra} $\SC(e,n)$, which is the quotient of $\UC(\gl_e)$ by those 2-morphisms factoring through objects which are not in the subset
$\La_{n,e}$, and thus go to zero under $\FC_{MSV}$. Clearly the functor $\FC_{MSV}$ descends to $\SC(e,n)$. They also discuss the connections between $\UC(\gl_e)$ and $\UC(\sl_e)$.
There is no map from $\UC(\sl_e)$ to $\UC(\gl_e)$, but nonetheless the quotient functor $\UC(\gl_e) \to \SC(e,n)$ descends to a map $\UC(\sl_e) \to \SC(e,n)$, see \cite[page 8]{MSV}.
These connections we include only for sake of background, as we will stick to $\UC(\gl_e)$ in this paper. The functor $\FC_{KL}$ also descends to $\SC(e,n)$. Clearly $\Flag_{n,e}$ is a
sub-2-category of $\SBSBim$, and $\FC_{MSV}$ is just $\FC_{KL}$ composed with the inclusion of $\Flag_{n,e}$ into $\SBSBim$.

What we do is to lift $\FC_{MSV}$ to a purely diagrammatic setting, constructing a functor $\FC \co \UC(\gl_e) \to \DG(\gl_n)$. Complicated operations involving polynomials in $\SBSBim$
are encoded using rather simple $S$-diagrams, making the construction of this functor significantly easier. The verification that the action is well-defined is still a lengthy
computation, but only because there are many relations in $\UC(\gl_e)$ to check, and many special cases; the individual computations for each relation are much easier. In fact, most of
the computations are ``topological," involving the ``Reidemeister-like" relations above of $\DG$. Thus our proof has a very different flavor than theirs, and our result is 
 stronger than theirs. \footnote{To be fair to Mackaay, Stosic, and Vaz, the technology for diagrammatic singular Soergel bimodules was not available when they wrote their paper. And it is still not available, technically!} Composing $\FC \co \UC(\gl_e)
\to \DG(\gl_n)$ with the $2$-functor $\FC_{EW} \co \DG(\gl_n) \to \SBSBim(\gl_n)$, one obtains $\FC_{MSV}$.

\begin{remark} It should be noted that Khovanov-Lauda and Mackaay-Stosic-Vaz both work over the base ring $\QM$. However, no coefficients appear in their calculations except signs, and
their functors are defined over $\ZM$, so they remain actions after base change to other fields. One minor advantage of defining the functor to diagrammatic 2-categories, rather than
algebraic ones, is that one has guaranteed connection to parity sheaves and geometry, even in positive characteristic. \end{remark}

We have now covered the entire diagram \eqref{eq:bigdiagram} except for the dashed arrow $\Sigma_{e,n}$. This $2$-functor appears in \cite{MSV} only when $e>n$, as a partial inverse to
$\FC_{MSV}$. Namely, consider the $\gl_e$-weight $1^n = (1,\ldots,1,0,\ldots,0)$, with $n$ ones and $e-n$ zeroes. This object corresponds to the ring $R$ itself.\footnote{The
corresponding sequence $\un{k}$ is $(0 \le 1 \le 2 \le \ldots \le n \le n \le n \le \ldots \le n)$. The corresponding partial flag variety is actually the full flag variety, with some
duplication.} Thus the image of $\FC_{MSV}$ applied to the monoidal category $\Hom_{\UC(\gl_e)}(1^n,1^n)$ will be a subcategory of $R$-bimodules. In fact, it is the category  $\BSBim$
of Bott-Samelson bimodules. They construct a functor which they call $\Sigma_{e,n}$ from $\SB$, the diagrammatic category of Bott-Samelson bimodules, to $\Hom_{\SC(e,n)}(1^n,1^n)$, and
they prove that this is an equivalence (when $e>n$) in \cite[Proposition 6.9]{MSV}. Note that $\FC_{EK}$, the functor defined in \cite{EKho} from $\SB$ to $\BSBim$, is an equivalence
when the Soergel categorification theorem holds (such as over $\QM$, where they work), which shows that these three categories are all equivalent over $\QM$. In \cite[Lemma 6.8]{MSV}
they show that $\FC_{MSV} \circ \Sigma_{e,n} \cong \FC_{EK}$.

\begin{remark} Although the computations of \cite{MSV} are specific to the weight $1^n$, the existence of an inverse functor like $\Sigma_{e,n}$ result should apply to any $\gl_e$
weight whose corresponding ring is $R$. \end{remark}

\subsection{Problems with signs} \label{subsec-signerror}

The only unfortunate thing about \eqref{eq:bigdiagram} is that it does not match what was stated in the literature!

The commutative diagram \eqref{eq:bigdiagram} is correct as we have described it in the previous section. However, it does not match what is stated in \cite{MSV} and \cite{KhoLau10}
because of a sign error we discovered in \cite{MSV} (unrelated to the sign error discussed in Remark \ref{rmk:signerror}). This sign error arose from differing
conventions, and it is difficult to notice, so we will try to be clear and precise. The same error propogated to the work of Mackaay-Thiel \cite{MacThi} in the affine case.

Recall that the definition of $\UC = \UC_{\leftarrow}(\gl_e)$ depended on a choice of orientation on the Dynkin diagram. In reality, both Khovanov-Lauda in \cite{KhoLau10} and
Mackaay-Stosic-Vaz in \cite{MSV} work with the other orientation, using $\UC_{\rightarrow}(\gl_e)$ instead of $\UC_{\leftarrow}(\gl_e)$. The difference is some signs in the quiver Hecke
relations. However, the diagram \eqref{eq:bigdiagram} does not commute when using the functors they originally used.

Let us discuss the problem with $\FC_{MSV}$ as originally defined, using $\UC_{\rightarrow}(\gl_e)$. Under the functor $\FC_{EK}$, a certain $i$-colored diagram ($1 \le i \le e-1$)
called a ``barbell'' or ``double dot'' is sent to multiplication by by the simple root $\a_i = x_i - x_{i+1}$ (c.f. \cite[relation (6.19)]{MSV}). Under the functor $\Sigma_{e,n}$, this
barbell is sent to a counterclockwise $i$-colored bubble inside the region label $1^n$. However, by \cite[equation (4.13)]{MSV}, the counterclockwise bubble is sent by $\FC_{MSV}$ to
$-\a_i = x_{i+1} - x_i$! \footnote{Looking at \cite[equation (4.13)]{MSV} in this special case, one sets $a = b = 1$, $r=0$. The $t$ variable and the $u$ variable correspond to
$x_{i+1}$ and $x_i$ respectively, since by the conventions on \cite[page 20]{MSV}, one reads the variable indices from right to left.} Thus, \cite[Lemma 6.8]{MSV} does not actually hold
as stated.

To reiterate, the functors $\FC_{MSV}$ and $\FC_{KL}$ are perfectly well-defined, but $\FC_{MSV} \circ \Sigma_{e,n}$ is not actually the same as $\FC_{EK}$, differing by a sign on
certain morphisms.

There are multiple ways one could fix this sign error. Signs can be sprinkled in various places by: changing the orientation on the Dynkin diagram and thus the conventions for
$\UC(\gl_e)$ itself; changing the Frobenius structure and thereby the conventions for $\SB$ and $\DG$; applying one of the symmetries in \cite[Section 3.3.2]{KhoLau10}; modifying the
functor $\FC_{MSV}$ by introducing signs into the image of a dot or a crossing; etcetera. In this paper, we have chosen to keep all the standard conventions for Soergel bimodules and
their Frobenius structures intact, and we have chosen to keep the standard conventions for the decategorified action (e.g. so we can not swap $E_i$ with $F_i$). We have also chosen to
keep signs away from polynomials as much as possible (e.g. the image of a dot is multiplication by $x_i$, not multiplication by $-x_i$). Our solution is to use $\UC_{\leftarrow}(\gl_e)$
to begin with, and to modify the functor $\FC_{MSV}$ by adding signs. One can repeat the proof of Mackaay-Stosic-Vaz with our conventions instead, though this is unnecessary. Instead,
we will confirm that our functor $\FC \co \UC_{\leftarrow}(\gl_e) \to \DG(\gl_n)$ is well-defined, and simply define $\FC_{MSV}$ to be the composition $\FC_{EW} \circ \FC$.

In the rest of this chapter, we will work with the 2-category $\UC = \UC_{\leftarrow}(\gl_e)$. This is the 2-category defined precisely as in \cite[Definition 3.1]{MSV}, except that in
equations (3.17) and (3.20), the sign $(i-j)$ is replaced with $(j-i)$. The reader can find the relations written out explicitly below: \eqref{newXX1}, \eqref{newXX2},
\eqref{newtriple1}, and \eqref{newtriple2}; just ignore any terms with $\delta_{i0}$, these occur only in the affine case. \footnote{When comparing with \cite{MSV}, note that their relations involve upward-oriented strands, and ours are downwards oriented, so one must rotate by 180 degrees.}

\subsection{Prologue} \label{subsec-prologue}

The rest of this chapter is an in-depth construction of the 2-functor $\FC$ from $\UC$ to $\DG$. Because it is no more difficult, we also extend this $2$-functor to the partial
idempotent completion $\dot{\UC}$ which includes as 1-morphisms the divided powers $E_i^{(k)}$ and $F_i^{(k)}$. This partial idempotent completion was given a diagrammatic presentation
by Khovanov-Lauda-Mackaay-Stosic in \cite{KLMS}.

Let us also briefly mention the sub-2-categories $\UC^-$ (resp. $\dot{\UC}^-$), which contain 1-morphisms generated by the lowering operators $F_i$ (and the divided powers $F_i^{(k)}$),
without the raising operators. We will define our functor from $\UC^-$ first, and check that it is well-defined there. Then we will check that the functor extends from $\UC^-$ to
$\dot{\UC}^-$, and from $\UC^-$ to $\UC$. Together, this will imply that the functor extends to $\dot{\UC}$ as well, as a minimal set of relations for $\dot{\UC}$ live in either $\UC$
or $\dot{\UC}^-$.

\subsection{Preliminaries: Littelmann operators} \label{subsec-prelims-finite}

\begin{defn} \label{defn:Lane} Let $\La_{n,e}$ denote the set $\{ \la = (\la_1, \ldots, \la_n) \in \ZM^n \; | \; 1 \le \la_1 \le \ldots \le \la_n \le e\}$. \end{defn}
	
This set is in bijection with several other sets which have already played a role.

\begin{defn} To each $\la \in \La_{n,e}$ we associate a $\gl_e$ weight, that is, an element $(r_1, r_2, \ldots, r_e)$ of $\ZM^e$. One lets $r_j$ be the number of $1 \le i \le n$ such
that $\la_i = j$. Note that $\La_{n,e}$ is in bijection with the $\gl_e$ weights for which $\sum r_i = n$. We also let $\La_{n,e}$ denote this subset of $\gl_e$ weights. To each $\gl_e$
weight $(r_1, \ldots, r_e)$ in $\La_{n,e}$, we associate a sequence $\un{k} = (0 = k_0 \le k_1 \le \ldots \le k_e = n)$, where $k_j = \sum_{i \le j} r_i$. \end{defn}

The set of such sequences $\un{k}$ is also in bijection with $\La_{n,e}$. These sequences indexed the objects in the Khovanov-Lauda flag 2-category $\Flag_{n,e}$.

In the rest of this paper we will always refer to elements of $\La_{n,e}$ via the formulation in Definition \ref{defn:Lane}, with the tacit awareness that its elements are also $\gl_e$
weights.

\begin{defn} Let $W = S_n$ act on $\ZM^n$ in the usual way. For $\la \in \La_{n,e}$ let $I(\la)$ be its \emph{stabilizer subset}, so that $W_I$ is the stabilizer of $\la$. \end{defn}

That is, for each inequality in $\la_1 \le \ldots \le \la_n$ which is actually an equality, one has a simple reflection in $I(\la)$.

\begin{defn} For $1 \le j \le e-1$ and $k \ge 1$, we define an operator $F_j^{(k)}$ on $\La_{n,e} \cup \{0\}$, the extension of $\La_{n,e}$ by a formal symbol $0$. It will take $k$
copies of $j$ inside $\la$ and replace them with $j+1$. When there are not enough instances of $j$ in $\la$, the operator yields $0$. Also, $F_j^{(k)}$ sends the formal symbol $0$ to
itself. Let $F_j = F_j^{(1)}$. \end{defn}

\begin{remark} The operator $F_j^{(k)}$ is just the $k$-th iterate of the operator $F_j$. The reason for the extra parentheses in the notation will become clear in the next section.\end{remark}
	
\begin{remark} Note that, when interpreted as $\gl_e$ weights, $F_j$ will add $-\a_i = (0,\ldots,0,-1,1,0 \ldots 0)$ to the weight, as expected. \end{remark}

\begin{example} $F_2$ applied to $(12233)$ is $(12333)$, $F_2^{(2)}$ applied to $(12233)$ is $(13333)$, and $F_2^{(3)}$ applied to $(12233)$ is $0$. \end{example}

To distinguish between numbers $1 \le i \le n$ which index entries of $\la \in \La_{n,e}$, and numbers $1 \le j \le e-1$ which index operators $F_j$, we refer to the former as \emph{indices} and the latter as \emph{colors}.

\begin{defn} A sequence of weights in $\La_{n,e}$ is an \emph{$F_j$-string} if it is a maximal chain of (nonzero) elements where each is obtained from the previous by applying $F_j$. \end{defn}

For example, \begin{equation} \label{F2string} (12222) \to (12223) \to (12233) \to (12333) \to (13333) \end{equation} is an $F_2$-string. We say this particular $F_2$-string has length
$4$; with \emph{start} $(12222)$ and \emph{finish} $(13333)$. The start and the finish are the two \emph{ends} of the $F_j$-string, and the remaining weights are on the \emph{interior}.
Note that the start of an $F_j$-string is a weight killed by $F_{j+1}$, and can not be the start of an $F_{j-1}$-string (not counting strings of length zero).

There are four possibilities for how the operator $F_j^{(k)}$ will affect the stabilizer subset $I(\la)$ (assuming that $F_j^{(k)} \la \ne 0$), depending on the position of $\la$ in an $F_j$-string. These possibilities are easily comprehended from the following example.

\begin{example} Consider the stabilizers in \eqref{F2string}: $$\{s_2, s_3, s_4\} \to \{s_2, s_3\} \to \{s_2, s_4\} \to \{s_3, s_4\} \to \{s_2, s_3, s_4\}.$$ \end{example}

In general, the start and the finish of an $F_j$-string will have stabilizer $W_L \times W_I$ where $L$ and $I$ are disjoint, and $W_I \cong S_m$ for some $m$. In the above example,
$m=4$ and $L = \emptyset$. The $k$-th term in the string will have stabilizer $W_L \times S_k \times S_{m-k}$. Thus $F_j^{(k)}$ will affect the stabilizer as follows: \begin{itemize}
\item Replace a simple reflection with the missing simple reflection in $I$, if going between two interior parts of the $F_j$-string. \item Remove a simple reflection, if going from the
start to the interior. \item Add a simple reflection, if going from the interior to the finish. \item Do nothing, if going from the start to the finish. \end{itemize} One should think
of the first case as being the most general, and the other cases as being ``degenerations" where one forgets to add and/or remove a simple reflection.

\begin{defn} Suppose that $\la$ and $\mu$ lie on the same $F_j$-string. Let $I(\la,\mu) = I(\la) \cap I(\mu)$.  Whenever the notation $I(\la, \mu)$ is used, the color $j$ is implicit, and it is understood that they lie on an $F_j$-string. \end{defn}

Note that $I(\la, \mu)$ differs from $I(\la)$ in at most one simple reflection.

\begin{defn} For $1 \le j \le e-1$ and $k \ge 1$, we define the operator $E_j^{(k)}$ similarly to $F_j^{(k)}$. It replaces $k$ copies of $j+1$ with $j$, if possible. \end{defn}
	
\begin{remark} One can think that these operators $E_j$ and $F_j$ give $\La_{n,e}$ the structure of a graph, with edges of different colors $1 \le j \le e-1$, analogous to the crystal
of a representation of $\gl_e$. In the literature these operators are often called \emph{Littelmann operators}, and the graph the \emph{Littelmann graph}, see \cite[Appendix]{GreenBook}
for more details. \end{remark}



\subsection{Definition of the categorical action: objects and 1-morphisms}

We recall the $\gl_n$-realization. Let $(W,S)$ denote the Coxeter system of type $A_{n-1}$, so that $W = S_n$, and $S = \{1, \ldots, n-1\}$ is identified with the set of adjacent
transpositions $\{s_i\}$. The realization $\hg$, defined over the base ring $\Bbbk = \ZM$, will be the (rank $n$) permutation representation of $S_n$. Let $\hg$ have basis $\{\ep_i\}_{1
\le i \le n}$, and $\hg^*$ have dual basis $\{x_i\}_{1 \le i \le n}$. As usual, the simple roots $\a_i$ for $i \in S$ are given by $\a_i = x_i - x_{i+1}$, and the simple coroots are
$\a^\vee_i = \ep_i - \ep_{i+1}$. Thus, $R = \ZM[x_1, \ldots, x_n]$ is the polynomial ring in $n$ generators, with the usual action of $S_n$. Let $\DG$ denote the diagrammatic 2-category
of singular Soergel bimodules attached to the realization above.

We now begin to define a 2-functor $\FC \co \dot{\UC} \to \DG$, by defining it on objects and 1-morphisms. We abusively let $\FC$ denote various related subfunctors (e.g. the functor
from $\UC^-$ to $\DG$). The objects of $\dot{\UC}$ are the $\gl_e$ weights, some of which have a corresponding element in $\La_{n,e}$. In this case, we send the $\gl_e$ weight $\la$ to the
stabilizer subset $I(\la)$. If a $\gl_e$ weight does not have a corresponding element in $\La_{n,e}$ then all 1-morphisms attached to this object will be sent to the zero, so it is
irrelevant where the object is sent. We say it is sent to the formal symbol $0$.

Let $1 \le j \le e-1$, and $k \ge 1$. To the divided power $1$-morphism $F_j^{(k)}$ in $\dot{\UC}$ we associate a $1$-morphism in $\DG$: letting $\mu = F_j^{(k)} \la$, we induce from
$I(\la)$ to $I(\la, \mu)$, and then restrict to $I(\mu)$. As always, by induction and restriction we refer to the corresponding $1$-morphisms in $\DG$, not the algebraic functors. Note
that $F_j^{(k)}$ is a very different $1$-morphism from $F_j \circ \cdots \circ F_j$ taken $k$ times. Meanwhile, to $E_j^{(k)}$ we attach the biadjoint 1-morphism to $F_j^{(k)}$, which
is again induction followed by restriction, though in the opposite direction.

In degenerate cases, when $\la$ is the start of an $F_j$-string (resp. when $\mu$ is the finish), restriction (resp. induction) will be simply the identity $1$-morphism.

\begin{ex} The $1$-morphism attached to $F_2$ from $(12234)$ to $(12334)$ is induction from $R^{s_2}$ to $R$, followed by (shifted) restriction to $R^{s_3}$. By induction and restriction, we mean the corresponding $1$-morphisms in the category $\DG$. After applying the functor to bimodules, $F_2$ is sent to $R(1)$ viewed as an $(R^{s_3}, R^{s_2})$-bimodule. \end{ex}

\begin{ex} The 1-morphism attached to $F_2^{(4)}$ from $(12222)$ to $(13333)$ is $R^{\{s_2, s_3, s_4\}}$, viewed as a bimodule over itself. \end{ex}

In the examples below, the picture on the left hand side (LHS) is a diagram in $\DG$, while the picture on the right hand side (RHS) is a diagram in $\UC$ or $\dot{\UC}$. The
``equality'' is shorthand for saying that $\FC$ sends the RHS to the LHS. We also identify the simple reflections $\{s_1, s_2, s_3, \ldots\}$ with the alphabet $\{s,t,u,\ldots\}$.

\begin{ex} \label{F1example} $F_1$ from $\la = (11123) \sim st$ to $\mu = (11223) \sim su$.
	\[{
\labellist
\small\hair 2pt
 \pinlabel {$su$} [ ] at 9 16
 \pinlabel {$s$} [ ] at 25 16
 \pinlabel {$st$} [ ] at 41 16
 \pinlabel {$\mu$} [ ] at 83 16
 \pinlabel {$\la$} [ ] at 108 16
\endlabellist
\centering
\ig{1}{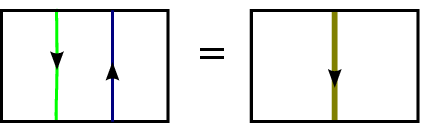}
}\] Technically, this is a drawing of the identity 2-morphism of $F_1$. \end{ex}

\begin{ex} \label{F1example2} $F_1^{(2)}$ from $\la = (11123) \sim st$ to $\rho = (12223) \sim tu$.
	\[{
\labellist
\small\hair 2pt
 \pinlabel {$tu$} [ ] at 9 16
 \pinlabel {$t$} [ ] at 25 16
 \pinlabel {$st$} [ ] at 41 16
 \pinlabel {$\rho$} [ ] at 83 16
 \pinlabel {\tiny{$2$}} [ ] at 92 8
 \pinlabel {$\la$} [ ] at 108 16
\endlabellist
\centering
\ig{1}{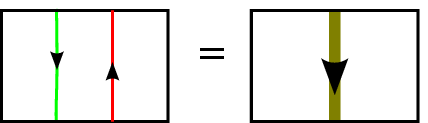}
}\] \end{ex}

(We will always denote simple reflections $s,t,u,v,w$ by the colors red, blue, green, pink, aqua respectively. When writing $\UC$ diagrams, olive represents $i$, teal represents $j=i+1$, and purple represents a distant color $k$.)

A general principle in this paper is that weights on the interior of an $F_i$-string are more complicated than weights on the end. We will generally consider weights on the interior first.
Pictures for ends of a string are typically obtained from interior pictures by deleting some strands.

\begin{ex} \label{F1examplestart} $F_1$ from $\nu = (11113) \sim stu$ to $\la = (11123) \sim st$.
	\[{
\labellist
\small\hair 2pt
 \pinlabel {$st$} [ ] at 9 16
 \pinlabel {$stu$} [ ] at 41 16
 \pinlabel {$\la$} [ ] at 83 16
 \pinlabel {$\nu$} [ ] at 108 16
\endlabellist
\centering
\ig{1}{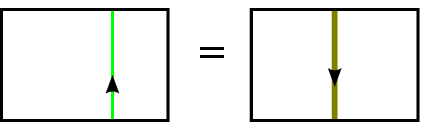}
}\] \end{ex}

\begin{ex} \label{F1examplefinish} $F_1$ from $\rho = (12223) \sim tu$ to $\xi = (22223) \sim stu$.
	\[{
\labellist
\small\hair 2pt
 \pinlabel {$stu$} [ ] at 9 16
 \pinlabel {$tu$} [ ] at 41 16
 \pinlabel {$\xi$} [ ] at 83 16
 \pinlabel {$\rho$} [ ] at 108 16
\endlabellist
\centering
\ig{1}{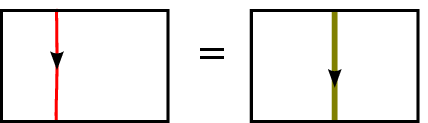}
}\] \end{ex}

\begin{ex} \label{F1example4} $F_1^{(4)}$ from $\nu = (11113) \sim stu$ to $\xi = (22223) \sim stu$.
	\[{
\labellist
\small\hair 2pt
 \pinlabel {$stu$} [ ] at 25 16
 \pinlabel {$\xi$} [ ] at 83 16
 \pinlabel {$\nu$} [ ] at 108 16
 \pinlabel {\tiny{$4$}} [ ] at 92 8
\endlabellist
\centering
\ig{1}{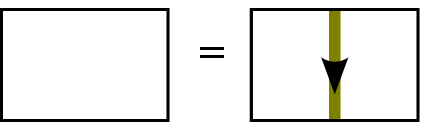}
}\] \end{ex}

We will abusively omit $\FC$ from the notation. So, we may write $\1_{\la}$ for the identity 1-morphism of an object in $\UC$, or abusively for the identity morphism $\FC(\1_{\la})$ of
the corresponding object in $\DG$. We write $\1_I$ instead when we want to emphasize its corresponding parabolic subset $I$. Thus, for instance, $\1_J F_i \1_I$ represents a 1-morphism
from $I$ to $J$ in $\DG$ coming from the action of $F_i$, although it is ambiguous which weight $\la$ corresponds to $I$.

Our next goal is to define the action of $\UC^-$ and $\dot{\UC}^-$, and check the relations.

\subsection{Definition of the categorical action: crossings}

The crossing in the quiver Hecke algebra is a map $ F_j F_i \1_\la \to  F_i F_j \1_\la$. There are three interesting separate cases: when $i=j$, when $i$ and $j$ are adjacent, and when $i$ and $j$ are distant.

\subsubsection{Same-color crossings}
First consider the endomorphism ring of $\1_K F_j \1_J F_j \1_I$. Suppose that all weights are on the interior of the $F_j$-string. Let $L = I \cap J \cap K$ be the ambient parabolic subset.
Note that the interesting parabolic subset $(I \cup J \cup K) \setminus L$ forms an $A_3$ configuration inside $S$; label its elements by $s, t, u$ in order. We may rewrite $I, J, K$ as
$Lst$, $Lsu$, and $Ltu$ respectively (which is shorthand for the disjoint union). Then we can define the following morphism inside $\End_{\DG}(\1_K F_i \1_J F_i \1_I)$, which we denote with
a crossing.

\[ {
\labellist
\small\hair 2pt
 \pinlabel {$Lst$} [ ] at 48 26
 \pinlabel {$Ls$} [ ] at 40 9
 \pinlabel {$Lsu$} [ ] at 30 4
 \pinlabel {$Lu$} [ ] at 18 11
 \pinlabel {$Ltu$} [ ] at 7 26
 \pinlabel {$Ls$} [ ] at 40 42
 \pinlabel {$Lsu$} [ ] at 31 46
 \pinlabel {$Lu$} [ ] at 18 41
 \pinlabel {$L$} [ ] at 30 16
 \pinlabel {$L$} [ ] at 30 32
 \pinlabel {$Lt$} [ ] at 30 25
\endlabellist
\centering
- \quad \ig{2}{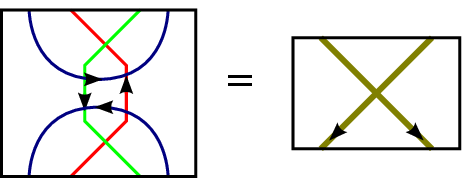}
} \]
Note the sign on the left hand side.

Note that $\{s,u\}$ is disconnected relative to $L$, so the sideways red-green crossing has degree $0$. However, $\{s,t,u\}$ is not disconnected relative to $L$, so the Reidemeister
III move can not be applied to either triangle in the diagram (there is a polynomial obstruction, see \eqref{R3unoriented} and the discussion afterwards). As the blue cup and cap have
degree $-1$ (since $t$ is not adjacent to any element of $L$), this degree of this crossing is $-2$.

When either $I$ or $K$ is on the end of the $F_j$-string, the crossing morphism is even easier. Here is what happens when the rightmost region is the start of the string.
\[ {
\labellist
\small\hair 2pt
 \pinlabel {$Lstu$} [ ] at 48 26
 \pinlabel {$Lst$} [ ] at 40 9
 \pinlabel {$Ls$} [ ] at 18 11
 \pinlabel {$Lsu$} [ ] at 7 26
 \pinlabel {$Lst$} [ ] at 40 42
 \pinlabel {$Ls$} [ ] at 18 41
\endlabellist
\centering
- \quad \ig{2}{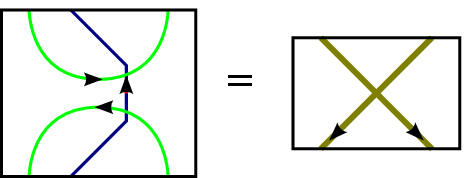}
} \]
Here is what happens for a string of length $2$.
\[ {
\labellist
\small\hair 2pt
 \pinlabel {$Lu$} [ ] at 48 26
 \pinlabel {$L$} [ ] at 40 9
 \pinlabel {$L$} [ ] at 40 42
\endlabellist
\centering
- \quad \ig{2}{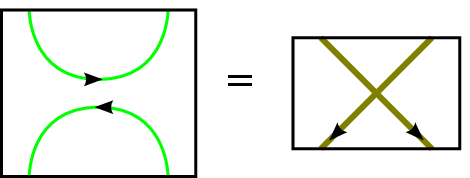}
} \]
The final edge case occurs when applying $F_j$ twice will leave the $F_j$-string. In this case the crossing is assigned to the zero map, as both the source and the target are the zero 1-morphism.

\subsubsection{Merge and split}

In the thick calculus, the crossing is a composition of two maps of degree $-1$: the merge (of two strands labeled $1$ into a strand labeled $2$), together with its flip, the split. These are easy to write down as well. Here is the merge when all weights are on the interior of the $F_j$-string.
\[ {
\labellist
\small\hair 2pt
 \pinlabel {\tiny{$2$}} [ ] at 104 30
\endlabellist
\centering
- \quad \ig{2}{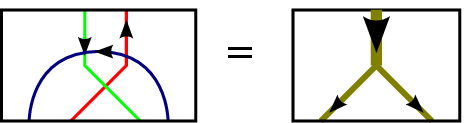}
} \]
In fact, the merge of two strands labeled $a$ and $b$ into a single strand labeled $k$ has the same picture, with different region labeling, when all weights are interior. The split has the upside-down picture, without the sign.\footnote{Whether one chooses to put the sign on the merge or the split is ultimately irrelevant. We  prefer this convention.}  The picture for weights on the ends of the $F_j$-string are obtained by removing red and/or green strands; we leave it to the reader.

\subsubsection{Adjacent-color crossings}

Now we consider a $2$-morphism from $F_j F_i \1_I$ to $F_i F_j \1_I$, for $j = i \pm 1$, and suppose that every weight is on the interior of its $F_i$- and $F_j$-strings. A picture of this crossing in $\DG$ is as follows.
\[
{
\labellist
\small\hair 2pt
 \pinlabel {$Lsu$} [ ] at 49 18
 \pinlabel {$Lu$} [ ] at 36 10
 \pinlabel {$Ltu$} [ ] at 29 4
 \pinlabel {$Lt$} [ ] at 18 9
 \pinlabel {$Ltv$} [ ] at 8 19
 \pinlabel {$Ls$} [ ] at 37 27
 \pinlabel {$Lsv$} [ ] at 28 30
 \pinlabel {$Lv$} [ ] at 18 28
 \pinlabel {$L$} [ ] at 28 17
 \pinlabel {$\mu$} [ ] at 93 16
 \pinlabel {$\la$} [ ] at 118 16
\endlabellist
\centering
\ig{2.5}{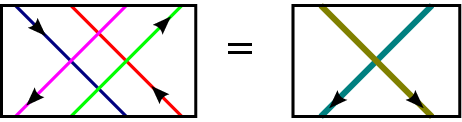}
} \]
Here, we let olive denote $F_i$ and teal denote $F_j$, and we have drawn the example where $j=i+1$. Once again, $L$ is the ambient parabolic subset. It contains an irrelevant part, together with an interesting part which interacts with the four strand colors $\{s,t,u,v\}$ as follows: $\{s,t, \textrm{stuff in }L, u, v\}$ is an $A_m$ configuration, in this order.  For example: $i=2$ and $j=3$, $\la = (2233333334)$ and $\mu = (2333333344)$, where $s$ and $t$ are the first two simple reflections, $u$ and $v$ are the last two, and $L$ contains everything in between.
Another example: $i=2$ and $j=3$, $\la = (2233477)$ and $\mu = (2334477)$, where $\{s,t,u,v\}$ are the first four simple reflections, and $L$ contains only the last reflection, which is irrelevant.

One can compute that the crossing has degree $+1$, coming from the blue-green crossing above.

Now we discuss what happens on the fringes of the $F_i$- and $F_j$-strings, and we continue with similar examples, where $i=2$ and $j=3$. We warn the reader that while the diagram above
looks symmetric, the region labels break the symmetry. For instance, the blue-green crossing has degree $+1$, while the red-pink crossing has degree $0$. The red and/or the pink strand may
cease to exist in certain edge cases, but the blue and green strands may not. For example, when $\la = (233344)$, it is as though the simple reflection $s$ did not exist, and the red strand
is removed (we are near the end of the $F_i$-string; there are too few 2's). When $\la = (22333)$ it is as though $v$ did not exist, and the pink strand is removed (we are near the start of
the $F_{i+1}$-string; there are too few 4's). When $\la = (2333)$ both red and pink strands are removed.

However, when we are near the start of the $F_i$-string and simultaneously near the end of the $F_{i+1}$-string, the situation does change! Consider when $\la = (2223444)$, so that $F_3$ can
be applied at most once. Now the simple reflections $t$ and $u$ effectively overlap, and the crossing appears different.
\begin{equation} \label{FiFjXspecial} \ig{2}{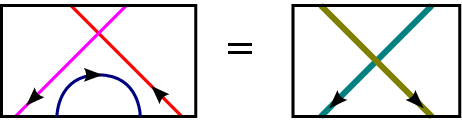} \end{equation}
The degree is still $+1$, coming from the blue cap. Once again, the red and the pink strands may be removed if there are too few $2$'s or too few $4$'s.

Finally, consider when $\la = (2224444)$, so that $F_3$ can not be applied at all. Then the crossing is assigned to the zero map, as its target is the zero 1-morphism (although the source is
not).

We leave the reader to explore the analogous possibilities when $j = i-1$. The diagrams will look similar, although the degree $+1$ crossing or cup will appear on top, rather than on
bottom.

\subsubsection{Distant-color crossings}

Now, consider a $2$-morphism from $F_k F_i \1_I$ to $F_i F_k \1_I$ for $k$ and $i$ distant. The crossing in $\DG$ is as follows.
\[
{
\labellist
\small\hair 2pt
 \pinlabel {$Lsv$} [ ] at 49 18
 \pinlabel {$Lv$} [ ] at 36 10
 \pinlabel {$Ltv$} [ ] at 29 4
 \pinlabel {$Lt$} [ ] at 18 9
 \pinlabel {$Ltw$} [ ] at 8 19
 \pinlabel {$Ls$} [ ] at 37 27
 \pinlabel {$Lsw$} [ ] at 28 30
 \pinlabel {$Lw$} [ ] at 18 28
 \pinlabel {$L$} [ ] at 28 17
\endlabellist
\centering
\ig{2.5}{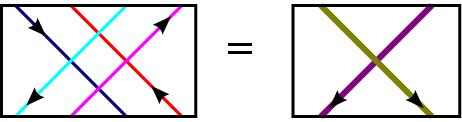}
} \]
Once again, $L$ is the ambient parabolic subgroup. The pairs $\{s,t\}$ and $\{v,w\}$ are both $A_2$ configurations, which are in different connected
components after adding $L$. This morphism has degree $0$.

Edge cases can remove any of the four strands. The target will be the zero morphism if and only if the source is as well.

\subsection{Definition of the categorical action: the dots}\label{SS_categ_act_dots}

\subsubsection{Dots are sent to polynomials}

We begin by determining the image under $\FC$ of the ``dot" from $\UC^-$. For each $\la$ with $F_j \la = \mu$, the $1$-morphism $\1_\mu F_j \1_\la$ in $\UC^-$ has a degree $2$ endomorphism known as a \emph{dot}. Let $I = I(\la)$ and $J = I(\mu)$. The image of the dot under $\FC$ will be a degree two endomorphism of (the diagrammatic version of) the $(R^J,R^I)$-bimodule $R^{I \cap J}$ (ignoring the grading shift). In both algebraic and diagrammatic contexts, this endomorphism space is spanned by multiplication by a polynomial, i.e. placing a polynomial in the middle region of the diagram corresponding to $\FC(F_j)$.

\begin{ex} \label{F1examplewdot} A dot on Example \ref{F1example}.
	\[ {
\labellist
\small\hair 2pt
 \pinlabel {$\fpoly_{\la \to \mu}$} [ ] at 24 17
\endlabellist
\centering
\ig{1.2}{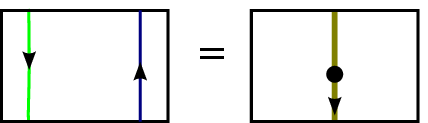}
} \] \end{ex}

We wish to choose a linear polynomial $\fpoly_{\la \to \mu} \in R^{I \cap J}$, for which multiplication by this element is $\FC$ applied to the dot. Note that this polynomial may depend
on $\la$, not just on the parabolic subgroups $I$ and $J$ in question.

The picture for a weight on the end of an $F_j$ string is obtained by removing strands.

\begin{ex} \label{F1examplestartwdot} A dot on Example \ref{F1examplestart}.
	\[{
\labellist
\small\hair 2pt
 \pinlabel {$\fpoly_{\nu \to \la}$} [ ] at 20 17
 \pinlabel {$\la$} [ ] at 83 16
 \pinlabel {$\nu$} [ ] at 108 16
\endlabellist
\centering
\ig{1}{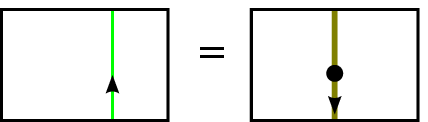}
}\] \end{ex}

When extending the functor $\FC$ to $\dot{\UC}^-$, a dot on a thick strand $F_j^{(k)}$ is defined in same way by a polynomial $\fpoly_{\la \to \mu}$. In the thick calculus, this dot
represents multiplication by the first symmetric polynomial in dots on the smaller strands, see \cite[\S 2.3]{KLMS}.

We call this family $\fpoly_{\la \to \mu}$ the collection of \emph{dot polynomials}. Let us define the standard choice for the dot polynomials.

\begin{defn} Let $\epoly_{\la \to \mu}$ be given by the sum $\sum x_i$, over indices $i$ which were changed from $\la$ to $\mu$. For example, if $\la = (1222234)$ and $\mu = (1223334)$,
so that positions 4 and 5 were changed by $F_2^{(2)}$, then $\epoly_{\la \to \mu} = x_4 + x_5$. In Example \ref{F1example}, $\epoly_{\la \to \mu} = x_3$. \end{defn}

In Khovanov-Lauda and Mackaay-Stosic-Vaz's construction, they set $\fpoly_{\la \to \mu} = \epoly_{\la \to \mu}$ by default. In particular, $\epoly_{\la \to \mu}$ only depends on the
parabolic subgroups $I = I(\la)$ and $J = I(\mu)$.

One of the goals of this chapter is to determine precisely which choices of dot polynomials will actually give rise to a functor from $\UC$, in an arbitrary realization of $(W,S)$. That
is, we extract the essential properties of the polynomials $\fpoly_{\la \to \mu}$ which will be needed for the Khovanov-Lauda relations to hold. This knowledge will be very useful when
we pass to the more complicated affine case in \S\ref{sec-quantumaction-affine}. It is also useful for other realizations, such as the realization spanned by the simple roots; the
polynomials $x_i$ and the standard choice of polynomials $\epoly_{\la \to \mu}$ are not well-defined in this context, because they do not lie within this span. 

\begin{prop} A collection of dot polynomials $\fpoly_{\la \to \mu}$ will give rise to a well-defined functor $\FC \co \UC_{\leftarrow}(\gl_e) \to \DG(\gl_n)$ if and only if they satisfy
the properties \eqref{f-properties} listed below. \end{prop}

As we confirm that $\FC$ is well-defined, we will use precisely the properties \eqref{f-properties}, which will prove the above proposition. Before listing the properties, let us describe the consequences.

As we will prove in Proposition \ref{onlychoices} below, the properties \eqref{f-properties} imply that, for the $\gl_n$-realization, one must choose $\fpoly_{\la \to \mu} = \epoly_{\la \to
\mu} + y$ for some fixed linear symmetric polynomial $y \in R^W$, independent of the choice of $\la$ and $\mu$. In particular, the polynomial $\fpoly_{\la \to \mu}$ only depends on the
parabolic subsets $I(\la)$ and $I(\mu)$. \footnote{When we pass to the affine case, the polynomials $\fpoly_{\la \to \mu}$ will depend on $\la$ and $\mu$, not just on their parabolic
subsets, which is why we keep this unneccessarily complicated notation.} The reader is welcome to accept $\fpoly_{\la \to \mu} = \epoly_{\la \to \mu} + y$ as a definition of the functor
$\FC$ on dots, and skip to Section \ref{subsec-conf-QHA-finite}.

Our construction and description of $\DG$ does not depend in any way on the choice of the basis $\{x_i\}$ of $\hg^*$. It depends only on the choice of roots and coroots. Because $\a_i =
x_i - x_{i+1} = (x_i + y) - (x_{i+1} + y)$, renaming every variable $x_i$ to $x_i + y$ will affect nothing. Thus by Proposition \ref{onlychoices} we may assume for the rest of the chapter
that $\fpoly_{\la \to \mu} = \epoly_{\la \to \mu}$ (we do this only when examining fake bubbles). However, the flexibility of adding $y$ will be become important when we treat the
affine case.

\subsubsection{Properties of the polynomials} \label{sssec-polyprops}

Now, let us state the properties $\fpoly_{\la \to \mu}$ must satisfy, and the relations of $\UC$ which use each property.

\begin{subequations} \label{f-properties}

{\bf Invariance:} When $F_j^{(k)}$ sends $\la$ to $\mu$, we have
\begin{equation} \label{fis-invariant} \fpoly_{\la \to \mu} \in R^{I(\la, \mu)}. \end{equation}
	
Invariance is required for the definition of the image of a dot on $F_j^{(k)}$ to even make sense, as the polynomial $\fpoly_{\la \to \mu}$ is placed in a region with this label.	

{\bf Additivity:} Suppose that $F_j^{(a)}\la = \nu$, and $F_j^{(b)}\nu = \mu$. Then
\begin{equation} \label{fis-additive} \fpoly_{\la \to \mu} = \fpoly_{\la \to \nu} + \fpoly_{\nu \to \mu}. \end{equation}

Additivity is equivalent to the compatibility between thick dots and ordinary dots in the thick calculus, see \S\ref{sssec-mergedot}. Because of additivity, we need only
worry about defining ordinary dot polynomials $\fpoly_{\la \to \mu}$ when $\mu = F_j \la$ for some $j$.

{\bf Dual basis condition:} Assume that $F_j \la = \nu$ and $F_j \nu = \mu$.  For example, $j=1$ and $\la = (11122)$, $\nu = (11222)$, $\mu = (12222)$. In particular, $\la$ and $\mu$ differ in precisely two adjacent indices, and there is a single simple reflection $t$ ($s_2$ in the previous example) which is in both $I(\la)$ and $I(\mu)$ but not $I(\nu)$. Our next conditions are:
\begin{equation} \label{fis-swapped} t(\fpoly_{\la \to \nu}) = \fpoly_{\nu \to \mu}, \end{equation}
\begin{equation} \label{fis-tunit} \pa_t(\fpoly_{\nu \to \mu}) = - \pa_t(\fpoly_{\la \to \nu}) = 1. \end{equation}

Alternatively, these two conditions can be replaced with:
\begin{equation} \label{fis-dualbases} \{\fpoly_{\nu \to \mu}, 1\} \textrm{ and } \{1, -\fpoly_{\la \to \nu}\} \textrm{ form dual bases for } R \textrm{ over } R^t \textrm{ w.r.t. } \pa_t. \end{equation}
Condition \eqref{fis-dualbases} is equivalent to \eqref{fis-swapped} and \eqref{fis-tunit} whenever it is known that $\fpoly_{\la \to \nu} + \fpoly_{\nu \to \mu}$ is $t$-invariant, which is implied by \eqref{fis-additive} and \eqref{fis-invariant}. An implication of \eqref{fis-dualbases} is that $\fpoly_{\nu \to \mu} - \fpoly_{\la \to \nu}$ is the simple root $\a_t$.

The dual basis condition is equivalent to the same color dot sliding relations, see \S\ref{sssec-samecolordotslide}.

{\bf (One-step) Locality:} Suppose that $i \ne j$ are two unequal colors, and $F_i \la = \mu$, $F_j \la = \nu$, $F_j \mu = F_i \nu = \rho$, and all these weights are not the formal object zero. For example, $i=1$, $j=2$, $\la = (11122)$, $\mu = (11222)$, $\nu = (11123)$, and $\rho = (11223)$. Then
\begin{equation} \label{fis-localpiece} \fpoly_{\la \to \mu} = \fpoly_{\nu \to \rho}. \end{equation}
By symmetry, $\fpoly_{\la \to \nu} = \fpoly_{\mu \to \rho}$ as well. One might say that parallel edges in the weight crystal yield equal polynomials.

One-step locality is equivalent to the different color dot sliding relations, see \S\ref{sssec-diffcolordotslide}.

\begin{claim} Condition \eqref{fis-localpiece} above is equivalent to the following condition, which we call \emph{locality}. Fix an index $1 \le k \le n$ and a color $1 \le j \le e-1$.
Suppose that $\la_k = j$ and $\la_{k+1} > j$, so that $F_j$ affects the $k$-th index of $\la$. Suppose that $\nu$ also satisfies $\nu_k = j$ and $\nu_{k+1} > j$. Then \begin{equation}
\label{fis-local} \fpoly_{\la \to F_j \la} = \fpoly_{\nu \to F_j \nu}. \end{equation} \end{claim}

\begin{proof} Clearly \eqref{fis-localpiece} is a special example of \eqref{fis-local}, so we need only show that \eqref{fis-localpiece} implies \eqref{fis-local}. By inspection of the graph $\La_{n,e}$, $\la$ can be related to $\nu$ by a sequence of various $E_k$ and $F_k$ for $k \ne j$, so that the edges $\la \to F_j \la$ and $\nu \to F_j \nu$ are actually parallel. Consequently, one can repeatedly apply \eqref{fis-localpiece} to obtain \eqref{fis-local}. \end{proof}
	
{\bf Repeated index condition:} Fix an index $1 \le k \le n$ and a color $1 \le i \le e-2$, and let $j = i+1$. Suppose that $\la_k = i$ and $\la_{k+1} > i+1$, so that $F_{j} \la = 0$. For example,
one can let $\la = (111333)$, with $i=1$ and $k=3$. Let $\mu = F_i \la$ and $\nu = F_{j} \mu$. Then \begin{equation} \label{fis-equalbyzero} \fpoly_{\la \to \mu} - \fpoly_{\mu \to \nu} = 0. \end{equation}

The repeated index condition is a consequence of the adjacent colored double crossing relation in this special case where $F_j \la = 0$, see \S\ref{sssec-adjacentXX}.
	
\begin{remark} \label{remark-locality} Thus locality is the statement that whenever $F_j$ acts on the $i$-th index of a weight, a dot on $F_j$ yields the same polynomial. The repeated
index condition below will compare when two different operators $F_j$ and $F_k$ act on the same index $i$, and deduce that they also give rise to the same polynomial. \end{remark}

There are several other conditions, but let us pause here. The above conditions are all relatively straightforward, while the remaining conditions use more complicated properties of
Frobenius extensions. Moreover, only the above conditions are required to pin down the dot polynomials uniquely.

\begin{prop} \label{onlychoices} Let $\{\fpoly_{\la \to \mu}\}$ be a family of linear polynomials satisfying invariance, additivity, the dual basis condition, locality, and the repeated
index condition. Then the family $\{\fpoly_{\la \to \mu} \}$ is determined uniquely by $\fpoly_{\nu \to F_j \nu}$ for any given $\nu \in \La_{n,e}$ and $1 \le j \le e-1$ with $F_j \nu
\ne 0$. Moreover, in the $\gl_n$ realization or any extension thereof, one has $\fpoly_{\la \to \mu} = \epoly_{\la \to \mu} + y$ for some linear $y \in R^W$, independent of the weights
$\la$ and $\mu$. \end{prop}

We give the proof in detail because we will need to repeat it in the affine case with more complexities.
	
\begin{proof} Fix some $\nu \in \La_{n,e}$ and $1 \le j \le e-1$ with $F_j \nu \ne 0$, and suppose that $f = \fpoly_{\nu \to F_j \nu}$ is given. Suppose that the $k$-th index of $\nu$
is incremented by $F_j$. Then by locality, $\fpoly_{\mu \to F_j \mu} = f$ whenever the $k$-th index of $\mu$ is incremented by $F_j$. In particular, let $\xi_m = (jj\ldots j(j+1) \ldots (j+1))$ where the final $j$ is in the $m$-th spot, so that $F_j \xi_m = \xi_{m-1}$. Then $\fpoly_{\xi_k \to \xi_{k-1}} = f$. However, $\xi_k$ is part of a large $F_j$-string involving the various $\{\xi_m\}$, and the other dot polynomials in this string are determined by the dual basis condition. So, for example, $\fpoly_{\xi_{k-1} \to \xi_{k-2}} = s_{k-1}(f)$ by \eqref{fis-swapped}, and by locality this determines the polynomial $\fpoly_{\la \to F_j \la}$ whenever $F_j$ increments the $(k-1)$-th index of $\la$. Repeating this argument, we see that additivity, invariance, the dual basis condition, and locality combine to determine $\fpoly_{\la \to F_j \la}$ for all $\la$, knowing just a single such polynomial.

In particular, in the $\gl_n$ realization, if $f = x_k + y$ for $y \in R^W$ then $\fpoly_{\la \to F_j \la} = x_\ell + y$ whenever the $\ell$-th index of $\la$ is incremented by $F_j$.

Now the repeated index condition can be used to determine the value of some $\fpoly_{\mu \to F_{j \pm 1} \mu}$ given the value of some specific $\fpoly_{\la \to F_j \la}$. Consequently, every dot polynomial attached to $F_{j+1}$ or $F_{j-1}$ is determined by $f$ as well. Repeating this, we see that every dot polynomial is determined uniquely by $f$.

Now suppose we work in the $\gl_n$ realization, and let $\la$ be such that $F_j$ increments the $n$-th index of $\la$, and let $g = \fpoly_{\la \to F_j \la}$. By \eqref{fis-dualbases},
$g$ must satisfy $\pa_{s_{n-1}}(g) = -1$, so $\pa_{s_{n-1}}(g - x_n) = 0$ and $g = x_n + y$ for some $y \in R^{s_{n-1}}$. But both $g$ and $x_n$ are invariant in $s_i$ for $1 \le i \le
n-2$, and thus so is $y$. Hence $y \in R^W$. Thus $\fpoly_{\la \to F_j \la} = \epoly_{\la \to F_j \la} + y$ for this particular $\la$ and $j$. But by uniqueness, we have $\{\fpoly\} =
\{\epoly\}$. (Note that $\{\epoly\}$ does satisfy all the conditions.) \end{proof}

Now let us conclude with two more conditions on the dot polynomials, which are equivalent to the adjacent colored double crossing relations in the remaining cases, see
\ref{sssec-adjacentXX}.

{\bf Generic Double crossing relations:} Here is the most complex requirement, which describes certain Frobenius invariants in terms of the polynomials $\fpoly$. Suppose that $1 \le i
\le e-2$ and let $j = i+1$. Suppose $F_i \la = \mu$, $F_{j} \la = \nu$, $F_{j} \mu = F_i \nu = \rho$, and all these weights are not the formal object zero. In particular, there are indices $a < b$ such
that $F_i$ increments the $a$-th index and $F_{j}$ increments the $b$-th index. We assume further that $a < b - 1$ (now we have avoided certain degenerate cases). For example, one could let $i=1$ and $\la = (111222333)$, with $a=3$ and $b=6$. Let $L$ be the
ambient parabolic subgroup $L = I(\la) \cap I(\rho) \cap I(\mu) \cap I(\nu)$. Let $t=s_a$ and $u = s_{b-1}$ (two distinct simple reflections by our assumptions) and note that
$s_c \in L$ whenever $a < c < b-1$. Then one has \begin{equation} \label{fis-padelta} \pa \D^L_{Ltu} = 1 \ot \fpoly_{\la \to \mu} - \fpoly_{\mu \to \rho} \ot 1 \in R^{Lt} \ot_{R^L} R^{Lu},
\end{equation} \begin{equation} \label{fis-mu} \mu^{Lt,Lu}_{Ltu} = \fpoly_{\nu \to \rho} - \fpoly_{\la \to \nu} \in R^L. \end{equation}

To put this condition in context, observe that if $a=1$ and $b=n$ then $\pa \D^L_{Ltu} = 1 \ot x_1 - x_n \ot 1$ and $\mu^{Lt,Lu}_{Ltu} = x_1 - x_n$ for the usual
Frobenius structure. In particular, $\epoly_{\la \to \mu}$ satisfies \eqref{fis-padelta} and \eqref{fis-mu}.

{\bf Redundant dual basis condition:} Continue the setup above, except instead of assuming that $a < b-1$, let $a = b-1$. For example, one can let $\la = (1112333)$. Now the simple reflections $t$ and $u$ coincide, so the equations \eqref{fis-padelta} and \eqref{fis-mu} no longer make sense. Instead, one has
\begin{equation} \label{fis-padeltaRedundant} \D^L_{Lt} = 1 \ot \fpoly_{\la \to \mu} - \fpoly_{\mu \to \rho} \ot 1 \in R^{Lt} \ot_{R^L} R^{Lt},
\end{equation} \begin{equation} \label{fis-muRedundant} \mu^{L}_{Lt} = \fpoly_{\nu \to \rho} - \fpoly_{\la \to \nu} \in R^L. \end{equation}
Actually, since the simple reflections adjacent to $t$ are not in $L$, $\mu^L_{Lt} = \a_t$, and these conditions are equivalent to:
\begin{equation} \label{fis-dualbasesRedundant} \{\fpoly_{\mu \to \rho}, 1\} \textrm{ and } \{1, -\fpoly_{\la \to \mu}\} \textrm{ form dual bases for } R \textrm{ over } R^t \textrm{ w.r.t. } \pa_t. \end{equation}
This looks just like \eqref{fis-dualbases}, although they apply to different situations. In fact, \eqref{fis-dualbasesRedundant} is redundant, following from \eqref{fis-dualbases} and locality.

This ends the list of conditions.
\end{subequations}


\begin{prop} \label{prop:thesearethepolyrelations} After applying the purported 2-functor $\FC$, the relations of $\UC^-$ and $\dot{\UC}^-$ which involve dots are equivalent to the properties of the polynomials $\fpoly_{\la \to \mu}$ given by \eqref{f-properties}. \end{prop}
	
The proof of this proposition will be found throughout the rest of this chapter, as the relations of $\UC^-$ and $\dot{\UC}^-$ are checked one by one.

\subsubsection{Summary}

We have now assigned maps in $\DG$ to the generators of $\UC^-$, namely dots and crossings.  We have also assigned maps in $\DG$ to the extra generators of the ``thick calculus'' from \cite{KLMS}, the split and merge maps.

\begin{thm} The assignments above induce a $2$-functor from $\dot{\UC}^-$ to $\DG$. \end{thm}

This theorem will be proven by the computations in the next section.

\subsection{Confirmation of the quiver Hecke action}
\label{subsec-conf-QHA-finite}

We now proceed to check the quiver Hecke relations, the relations of $\UC^-$, one by one. We continue to use the same conventions in our diagrams from before: there is some parabolic
subset $L$ which does not contain simple reflections named $s,t,u,v,w$, which are denoted by the colors red, blue, green, pink, aqua respectively. Every region in a diagram is labeled
by $L$ together with a subset of these simple reflections. When writing $\UC$ diagrams, olive represents $i$, teal represents $j=i+1$, and purple represents a distant color $k$.

\subsubsection{Same color dot sliding} \label{sssec-samecolordotslide}

There are two same-color dot sliding relations, which are the two equalities in \cite[(3.15)]{MSV}. We will check one, as the other follows by a mirrored argument. Note that we prefer to examine $F_i$ rather than $E_i$, so we will be checking the 180 degree rotation of \cite[(3.15)]{MSV}.

Assume that $\{s,t,u\}$ is an $A_3$ configuration. Let $f = \fpoly_{\la \to \nu}$ and $g = \fpoly_{\nu \to \mu}$.
\begin{equation} {
\labellist
\small\hair 2pt
 \pinlabel {$f$} [ ] at 194 170
 \pinlabel {$g$} [ ] at 294 138
 \pinlabel {\tiny $f$} [ ] at 44 98
 \pinlabel {\tiny $g$} [ ] at 121 81
 \pinlabel {$\la$} [ ] at 77 154
 \pinlabel {$\nu$} [ ] at 62 144
 \pinlabel {$\nu$} [ ] at 62 164
 \pinlabel {$\mu$} [ ] at 46 154
 \pinlabel {$\la$} [ ] at 144 154
 \pinlabel {$\nu$} [ ] at 129 144
 \pinlabel {$\nu$} [ ] at 129 164
 \pinlabel {$\mu$} [ ] at 115 154
\endlabellist
\centering
\ig{1.3}{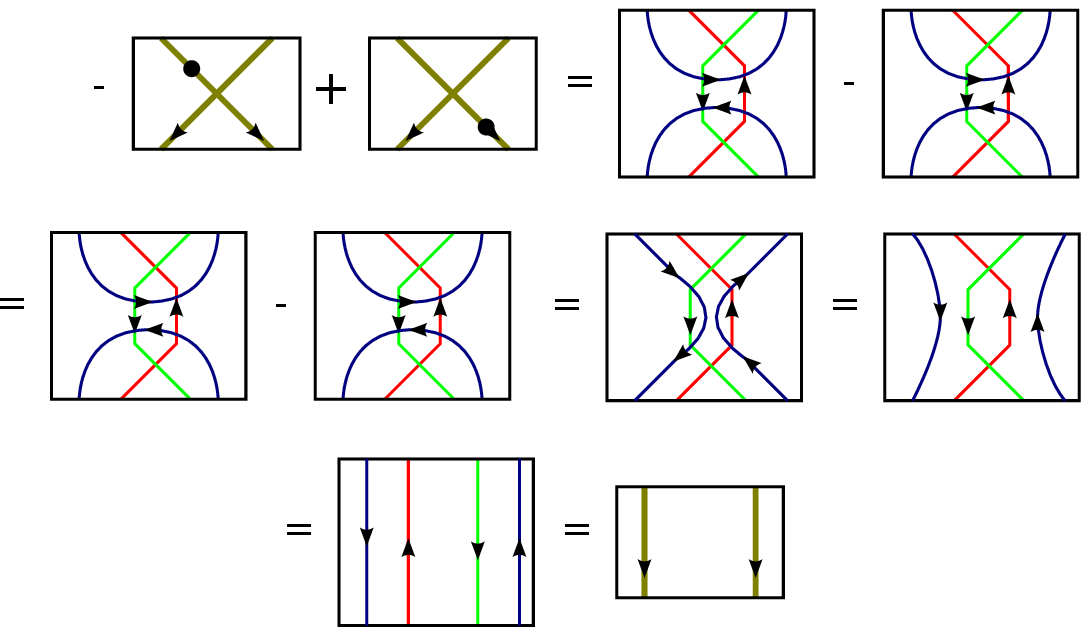}
} \end{equation}

The first equality follows from \eqref{slidepoly}, the second from \eqref{coproduct} and \eqref{fis-dualbases}, and the rest by unobstructed RII moves. The red-green sideways crossings are degree $0$, because $\{s,u\}$ is disconnected relative to $L$.

This was the computation on the interior of the $F_i$-string. Near the ends, one must remove either the red strand, the green strand, or both. However, it is easy to observe that the
computation above goes through in the same way, in the absence of any of these strands.

\subsubsection{Same color double crossing}

The double crossing relation \cite[(3.14)]{MSV} is easy.
\begin{equation} \label{NilHecke2} \ig{1.1}{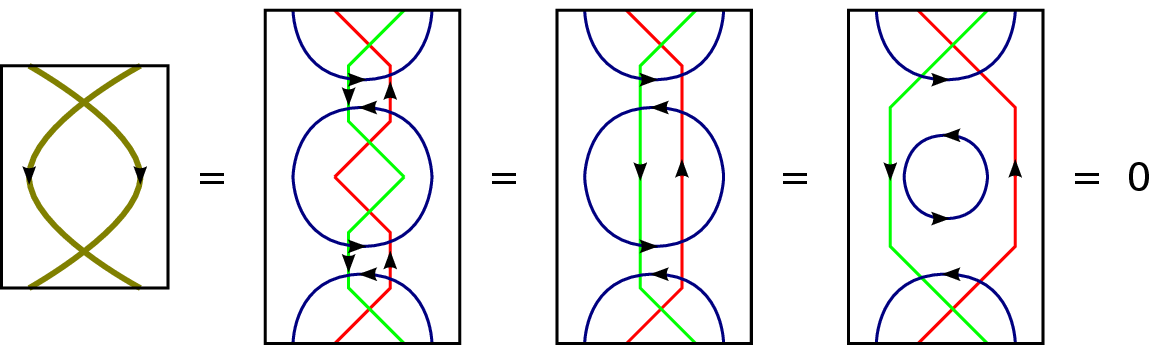} \end{equation}
We use unobstructed RII moves, and then \eqref{ccwcircempty}. The computation is the same in the absence of the green and/or red strands.

\subsubsection{Same color triple crossing}

The triple crossing relation \cite[(3.14)]{MSV} takes some work. Assume that $\{s,t,u,v\}$ is an $A_4$ configuration. \begin{equation} \label{NilHecke3} \ig{1.1}{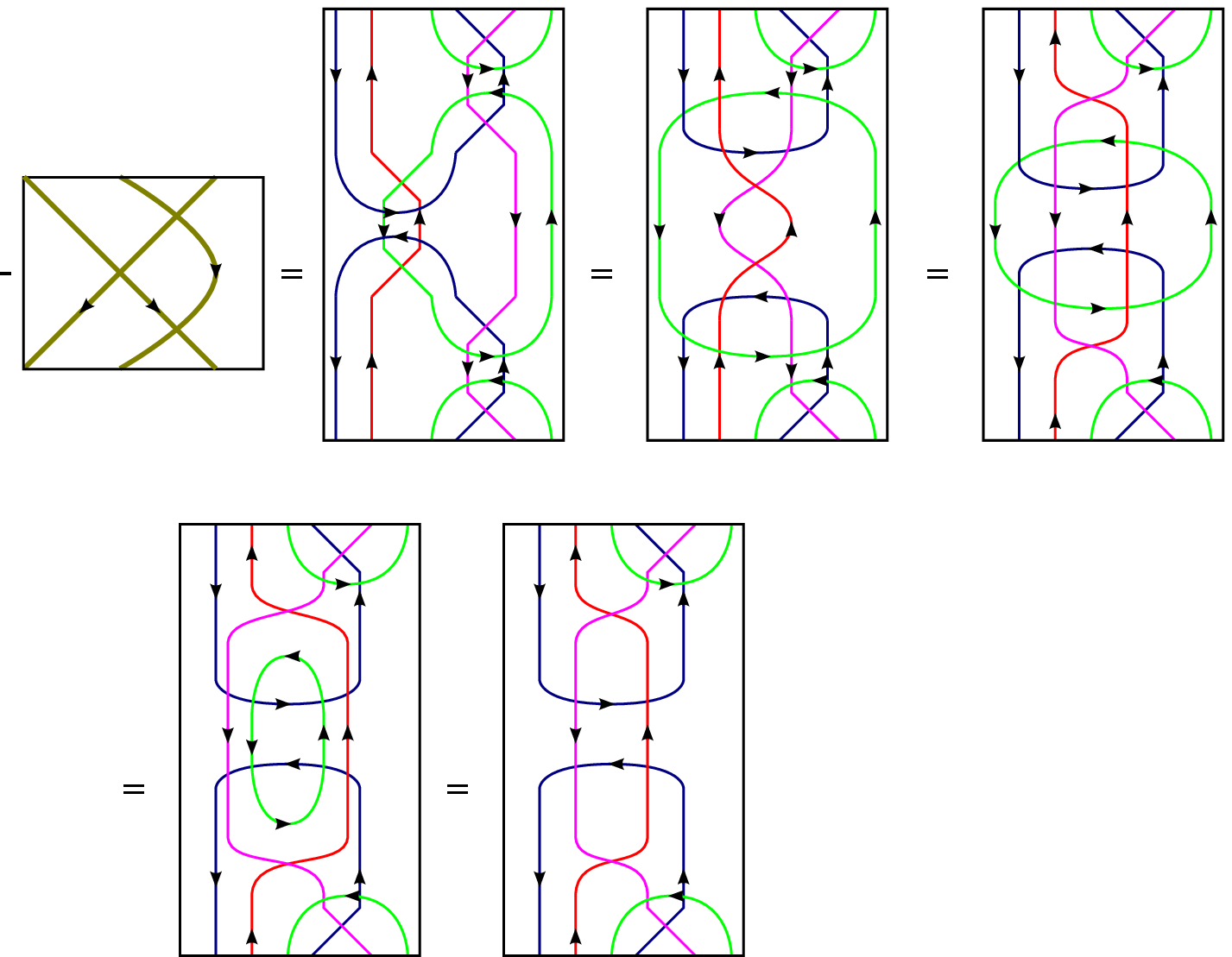} \end{equation} We use an
unobstructed red-pink RII move to add crossings, and then unobstructed RII and RIII moves to reach the second row. The final equality (removing the green circle) is a computation
involving obstructed RII moves for the $A_2$ configuration $\{t,u\}$ (the region does not contain either $s$ or $v$, so there is no interference from $L$). This equality was shown in
\cite[Claim 6.5]{ECathedral}.

Note that the final diagram in \eqref{NilHecke3} splits naturally into a top and bottom half, each of degree $-3$. Continuing the computation on the bottom half, we have
\begin{equation} \label{NilHecke3bottom} \ig{1.1}{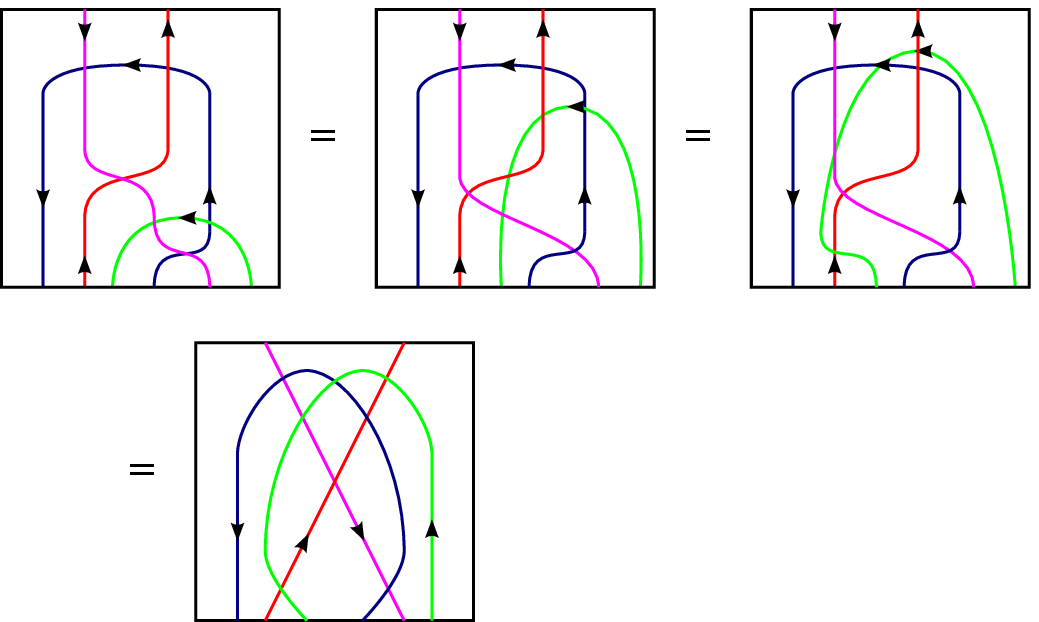} \end{equation}
We began with an unobstructed red-green RII move, followed by some oriented RIII moves. The final equality is just redrawing the diagram so it appears more symmetric.

A similar computation to \eqref{NilHecke3} shows that
\begin{equation} \label{NilHecke3-2} - \quad \ig{1.1}{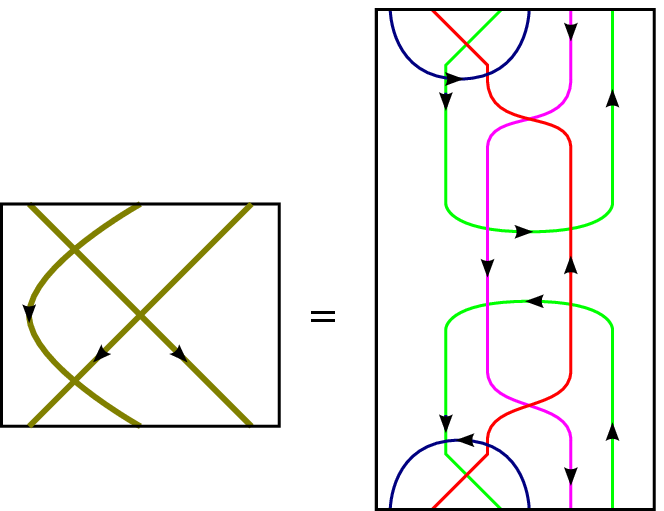} \end{equation}
and a similar computation to \eqref{NilHecke3bottom} shows that
\begin{equation} \label{NilHecke3bottom-2}  \ig{1.1}{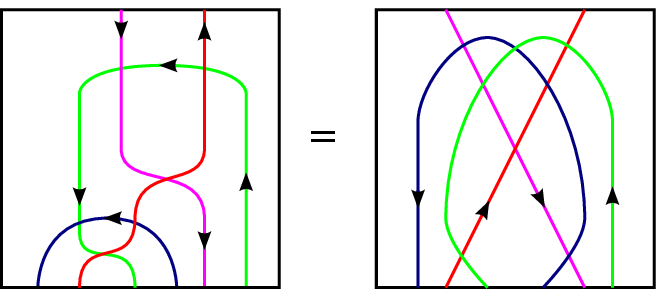}. \end{equation}
Thus \eqref{NilHecke3} and \eqref{NilHecke3-2} are equal, proving the triple crossing equality.

If the weights involved include the start and/or finish of the $F_i$-string, then the pink and/or red strands will need to be removed. However, this does not affect the computation at
all, as the reader can confirm.

\subsubsection{Different color dot sliding} \label{sssec-diffcolordotslide}

We now check \cite[(3.18)]{MSV}. However, using \eqref{slidepoly}, this relation is equivalent to the condition \eqref{fis-localpiece}, which is a very easy calculation.

\subsubsection{Distant color double crossing}

We now check the first row of \cite[(3.17)]{MSV}. Let $i$ and $k$ (olive and purple) be distant colors. Assume that $\{s,t,v,w\}$ is an $A_2 \times A_2$ configuration relative to $L$.
\begin{equation}\label{KLR3} \ig{1}{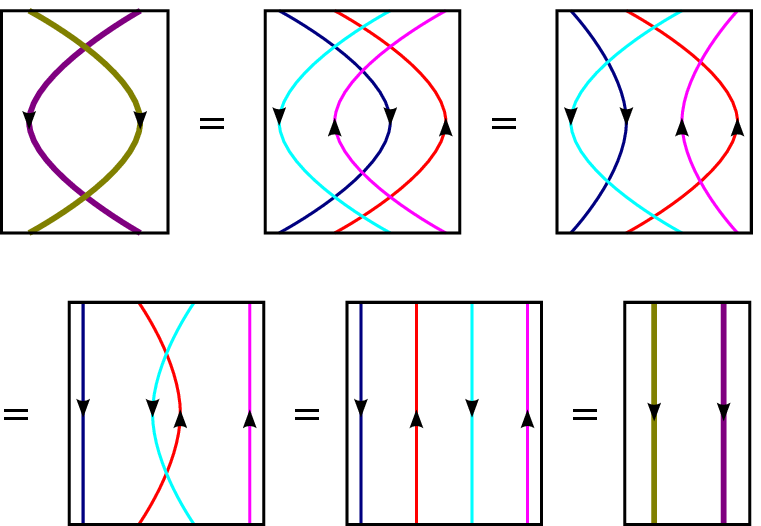} \end{equation}
We use only unobstructed RII moves. Edge cases can be dealt with by removing any given strand.

\subsubsection{Adjacent color double crossing} \label{sssec-adjacentXX}

We now check the second row of \cite[(3.17)]{MSV}. Recall that we use $\UC_{\leftarrow}(\gl_e)$ instead of $\UC_{\rightarrow}(\gl_e)$, so that the sign $(i-j)$ appearing in \cite{MSV} should be replaced by $(j-i)$. Let $i$ and $j$ (olive and teal) be adjacent colors, with $j = i+1$.

We first check the most generic case, on the interior of all the $F_i$ and $F_j$ strings. Recall the following setup which preceded \eqref{fis-padelta}. Suppose that $1
\le i \le e-2$ and $F_i \la = \mu$, $F_{j} \la = \nu$, $F_{j} \mu = F_i \nu = \rho$, and all these weights are not the formal object zero. In particular, there are indices $a < b$
such that $F_i$ increments the $a$-th index and $F_{j}$ increments the $b$-th index. We assume further that $a < b - 1$. Let $L$ be the ambient parabolic subgroup $L = I(\la) \cap
I(\rho) \cap I(\mu) \cap I(\nu)$. Let $t=s_a$ and $u = s_{b-1}$ (two distinct simple reflections by our assumptions) and note that $s_c \in L$ whenever $a < c < b-1$.

We further suppose that the $(a-1)$-th index is equal to the $a$-th index in $\la$, and the $(b+1)$-th index is equal to the $b$-th index in $\rho$. Thus $s = s_{a-1} \in I(\la)$ and $v
= s_b \in I(\rho)$. So $I(\la) = Lsu$, $I(\mu) = Ltu$, $I(\nu) = Lsv$, and $I(\rho)=Ltv$.

For example, the reader can ponder $\la = (111222333)$ and $i=1$.

There are actually two separate cases that need to be checked: when $i$ is on the left, and when $j$ is on the left. First consider $\1_{\rho} F_{j} \1_\mu F_i \1_{\la}$.
\begin{equation} \label{doublecrossing1} {
\labellist
\tiny\hair 2pt
 \pinlabel {$\fpoly_{\la \to \mu}$} [ ] at 80 33
 \pinlabel {$\fpoly_{\mu \to \rho}$} [ ] at 128 34
\endlabellist
\centering
\ig{1.4}{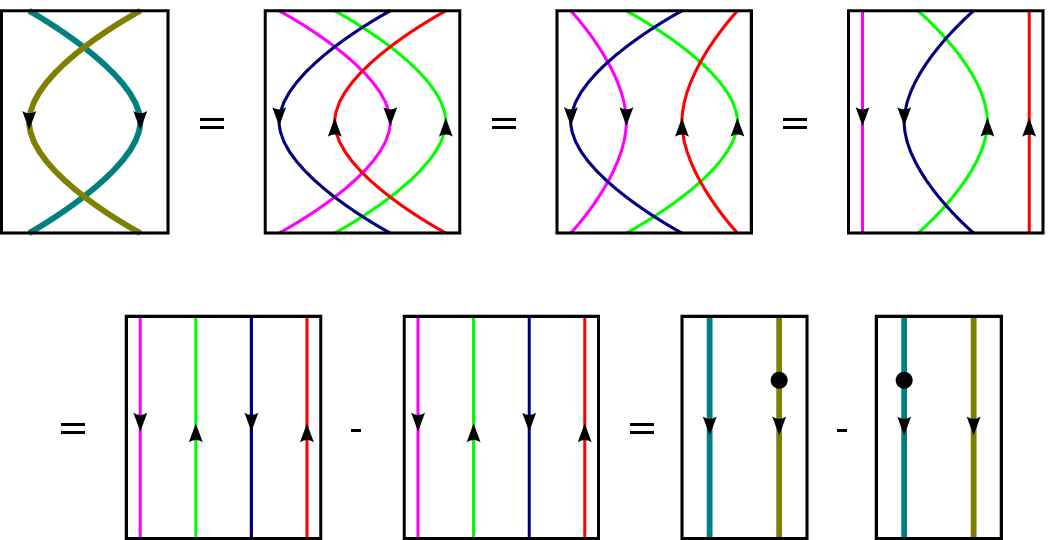}
} \end{equation}
We began with unobstructed RII moves. Then we applied \eqref{R2nonoriented1} together with \eqref{fis-padelta}.

Now consider $\1_{\rho} F_i \1_\nu F_{j} \1_{\la}$.
\begin{equation} \label{doublecrossing2} {
\labellist
\tiny\hair 2pt
 \pinlabel {$\fpoly_{\nu \to \rho}$} [ ] at 32 34
 \pinlabel {$\fpoly_{\la \to \nu}$} [ ] at 137 34
 \pinlabel {{\tiny $\fpoly_{\nu \to \rho} - \fpoly_{\la \to \nu}$}} [ ] at 228 123
\endlabellist
\centering
\ig{1.4}{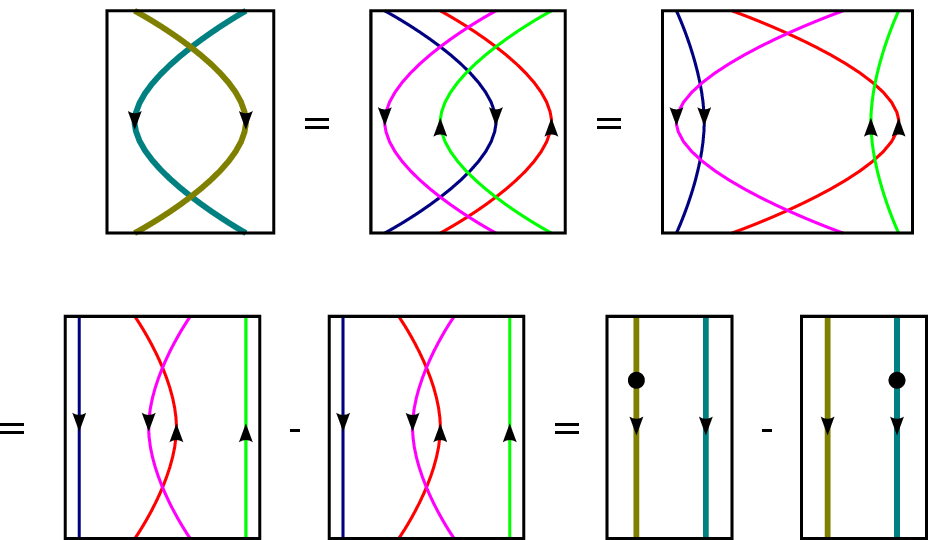}
} \end{equation}
We began with \eqref{R2nonoriented2} and \eqref{fis-mu}. The rest was polynomial sliding \eqref{slidepoly} and unobstructed RII moves.

We now remove some of our genericity assumptions. Leaving out the assumption that $s \in I(\la)$ (resp. that $v \in I(\rho)$) will correspond to removing the red (resp. pink) strand, which does not affect the above computations. For example, the reader can ponder $\la = (1222333)$ or $(111222)$. This handles the easiest edge cases. Now we must consider more degenerate edge cases.

The first difficult edge case occurs when $F_{j}$ can be applied at most once. The example to consider is $\la = (1112333)$, where the simple reflections $t$ and $u$ ``overlap." The derivation now goes as follows. Consider $\1_{\rho} F_{j} \1_\mu F_i \1_{\la}$.
\begin{equation} \label{doublecrossingweird1} {
\labellist
\tiny\hair 2pt
 \pinlabel {$\fpoly_{\la \to \mu}$} [ ] at 80 33
 \pinlabel {$\fpoly_{\mu \to \rho}$} [ ] at 128 34
\endlabellist
\centering
\ig{1.35}{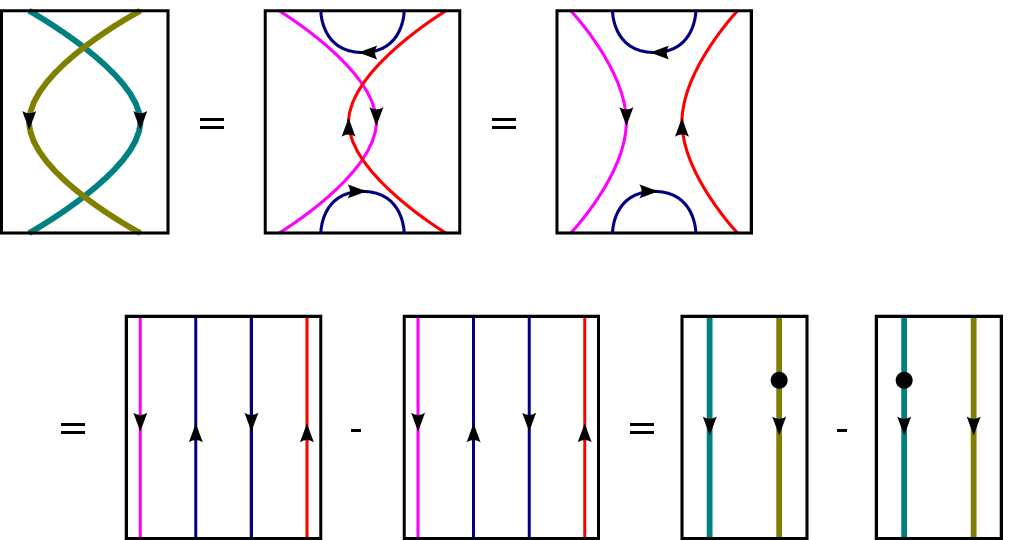}
} \end{equation}
This time we use \eqref{coproduct} and \eqref{fis-dualbases}. Now consider $\1_{\rho} F_i \1_\nu F_{j} \1_{\la}$.
\begin{equation} \label{doublecrossingweird2} {
\labellist
\tiny\hair 2pt
 \pinlabel {$\fpoly_{\nu \to \rho}$} [ ] at 32 34
 \pinlabel {$\fpoly_{\la \to \nu}$} [ ] at 137 34
 \pinlabel {{\tiny $\fpoly_{\nu \to \rho} - \fpoly_{\la \to \nu}$}} [ ] at 219 123
\endlabellist
\centering
\ig{1.3}{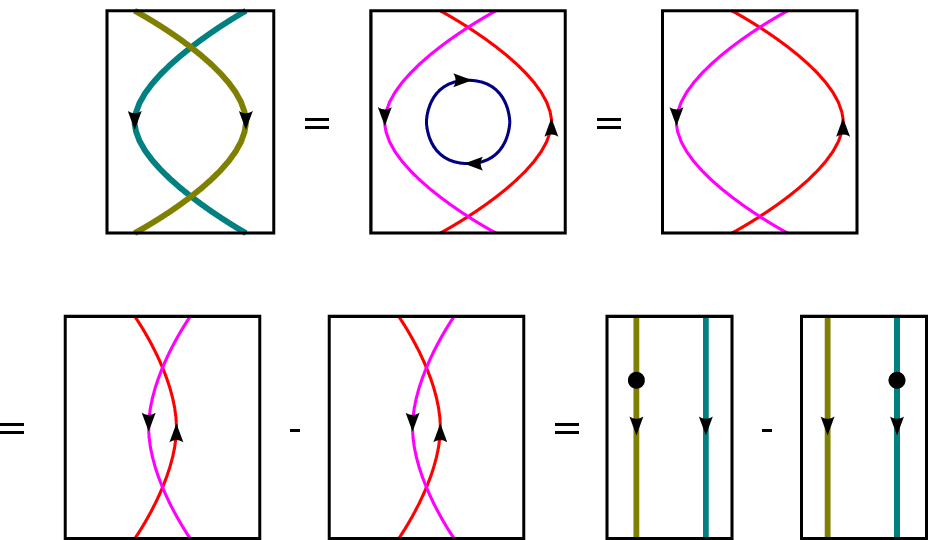}
} \end{equation}
This time we use \eqref{cwcirc} and \eqref{fis-dualbases}.  Both computations hold in the absence of the red and/or pink strand (e.g. too few 1s or too few 3s).

The final difficult edge case occurs when $F_{j}$ can not be applied at all, so that the crossing of $F_i$ and $F_{j}$ is the zero map, while $F_{j} F_i \1_{\la}$ is not the zero
1-morphism. For example, consider $\la = (111333)$, with $i=1$. Thus we are obligated to prove the following.
\begin{equation} \ig{1}{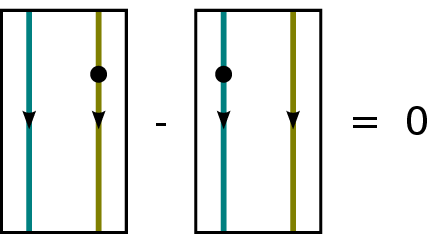} \end{equation}
In the corresponding $S$-diagrams, the polynomials appear in the same region, so the statement is that $\fpoly_{\la \to \mu} - \fpoly_{\mu \to \nu} = 0.$
However, this is exactly the content of \eqref{fis-equalbyzero}.

\begin{remark} \label{howpolysused} The reader can observe that we have now required all the properties \eqref{f-properties}, with the exception of \eqref{fis-additive}. In other words,
we have used the properties of $\fpoly_{\la \to \mu}$ when $\mu = F_i \la$, but not when $\mu = F_i^{(k)} \la$. The remaining property \eqref{fis-additive} will be required (in a
straightforward way) when we check the relations for the thick calculus. \end{remark}

\subsubsection{Typical triple crossing}

We now check the relation for the triple crossings $F_{i_1} F_{i_2} F_{i_3} \to F_{i_3} F_{i_2} F_{i_1}$. When $i_1 \ne i_3$, or $i_1 = i_3$ and $i_2$ is
distant, the desired relation is an equality between two KLR diagrams, which is \cite[(3.19)]{MSV}. It is possible that the functor sends one of the region labels (weights) in one of the triple crossing diagrams to
zero, because it is not in $\La_{n,e}$, but if this happens, one can check that the other triple crossing diagram also has a weight which is sent to zero. When both sides are non-zero,
the corresponding equality in $\DG$ will follow easily from Reidemeister II and Reidemeister III moves, all of which are unobstructed. We leave this verification as an exercise to the
reader.

\subsubsection{Interesting triple crossing}

We now check \cite[(3.20)]{MSV}, the triple crossing relation when $i_1 = i_3 = i$ (olive) and $i_2 = j = i+1$ (teal). The case when $i_1 = i_3 = j$ and $i_2 = i$ follows by an analogous computation. Recall that we use $\UC_{\leftarrow}(\gl_e)$ instead of $\UC_{\rightarrow}(\gl_e)$, so that the sign $(i-j)$ appearing in \cite{MSV} should be replaced by $(j-i)$. Also note that we prefer to examine $F_i$ rather than $E_i$, so we will be checking the 180 degree rotation of \cite[(3.20)]{MSV}.

Assume that $\{s,t,u\}$ is an $A_3$ configuration, $\{x,y\}$ is an $A_2$ configuration, and $L$ connects them (i.e. $L$ contains all simple reflections between $u$ and $x$). We now let pink denote $x$ and aqua denote $y$. We also omit the ambient region $L$ from the region labeling. One has
\begin{equation}\label{hardreln1} {
\labellist
\small\hair 2pt
 \pinlabel {$stx$} [ ] at 180 218
 \pinlabel {$sty$} [ ] at 151 235
 \pinlabel {$tuy$} [ ] at 117 221
 \pinlabel {$sx$} [ ] at 173 186
 \pinlabel {$su$} [ ] at 151 182
 \pinlabel {$uy$} [ ] at 122 193
 \pinlabel {$sy$} [ ] at 147 210
 \pinlabel {$sx$} [ ] at 172 302
 \pinlabel {$su$} [ ] at 149 310
 \pinlabel {$uy$} [ ] at 123 299
 \pinlabel {$sy$} [ ] at 147 280
 \pinlabel {$s$} [ ] at 158 201
 \pinlabel {$s$} [ ] at 160 286
 \pinlabel {$suy$} [ ] at 138 186
 \pinlabel {$suy$} [ ] at 137 304
 \pinlabel {$st$} [ ] at 168 238
 \pinlabel {$ty$} [ ] at 134 245
 \pinlabel {$stx$} [ ] at 241 45
 \pinlabel {$st$} [ ] at 224 60
 \pinlabel {$t$} [ ] at 208 70
 \pinlabel {$tu$} [ ] at 190 59
 \pinlabel {$tuy$} [ ] at 176 45
 \pinlabel {$sx$} [ ] at 233 16
 \pinlabel {$s$} [ ] at 218 36
 \pinlabel {$su$} [ ] at 207 16
 \pinlabel {$u$} [ ] at 195 36
 \pinlabel {$uy$} [ ] at 181 16
 \pinlabel {$sx$} [ ] at 233 131
 \pinlabel {$s$} [ ] at 218 111
 \pinlabel {$su$} [ ] at 207 131
 \pinlabel {$u$} [ ] at 195 111
 \pinlabel {$uy$} [ ] at 181 131
 \pinlabel {$-$} [ ] at -10 231
\endlabellist
\centering
\ig{1.1}{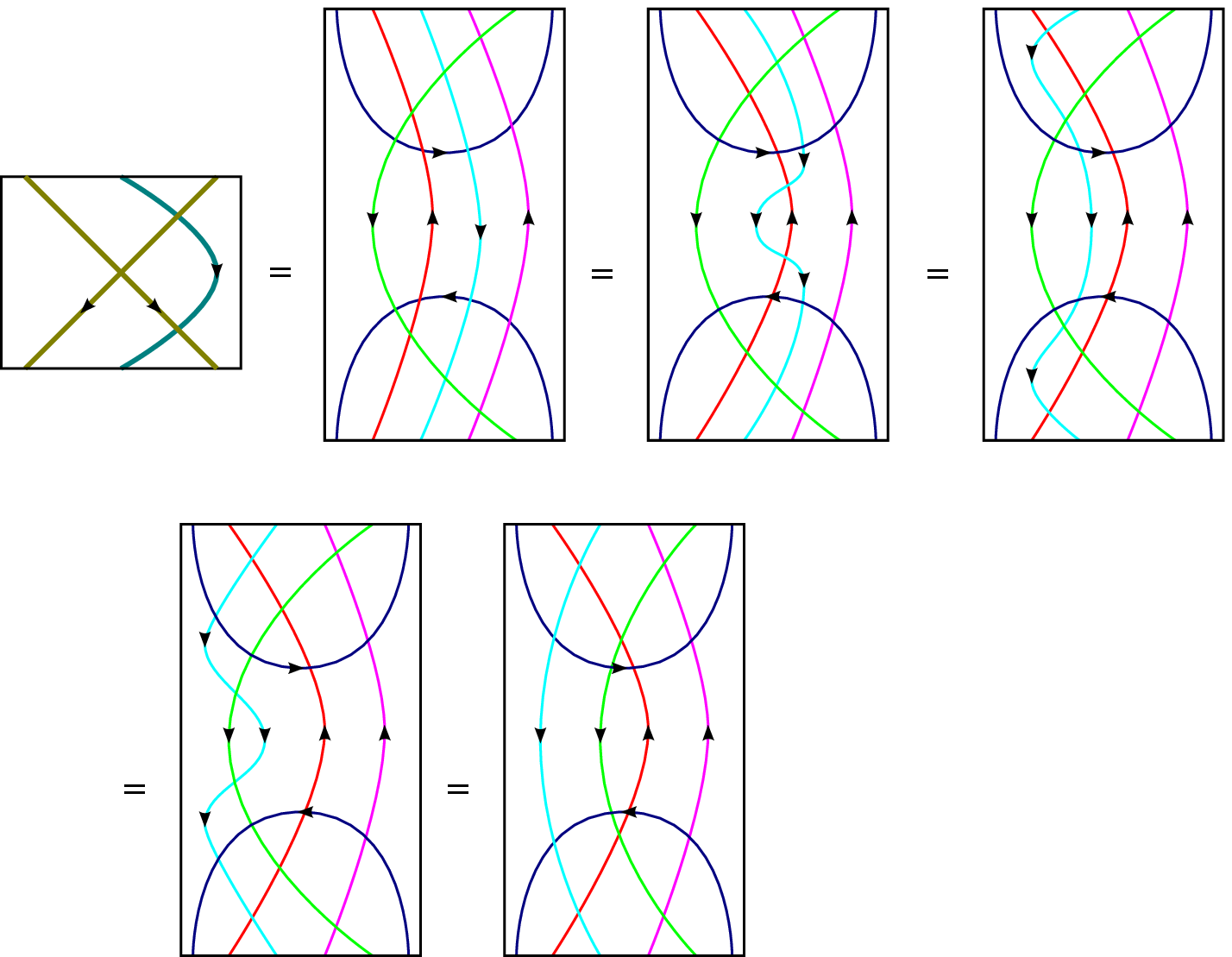}
} \end{equation}
This is just a series of unobstructed RII and RIII moves, used to pull the aqua strand to the left. They are unobstructed because $\{s,y\}$ is disconnected in the absence of either $t$, $u$, or $x$.

Meanwhile
\begin{equation} \label{hardreln2} {
\labellist
\small\hair 2pt
 \pinlabel {$stx$} [ ] at 178 224
 \pinlabel {$tx$} [ ] at 161 242
 \pinlabel {$tux$} [ ] at 144 242
 \pinlabel {$tu$} [ ] at 127 242
 \pinlabel {$tuy$} [ ] at 116 219
 \pinlabel {$sx$} [ ] at 171 200
 \pinlabel {$sux$} [ ] at 158 192
 \pinlabel {$ux$} [ ] at 150 213
 \pinlabel {$su$} [ ] at 145 185
 \pinlabel {$u$} [ ] at 136 204
 \pinlabel {$uy$} [ ] at 124 192
 \pinlabel {$sx$} [ ] at 171 292
 \pinlabel {$sux$} [ ] at 158 303
 \pinlabel {$ux$} [ ] at 150 284
 \pinlabel {$su$} [ ] at 145 313
 \pinlabel {$u$} [ ] at 136 292
 \pinlabel {$uy$} [ ] at 123 306
 \pinlabel {$stx$} [ ] at 400 224
 \pinlabel {$tx$} [ ] at 385 237
 \pinlabel {$t$} [ ] at 368 245
 \pinlabel {$tu$} [ ] at 348 242
 \pinlabel {$tuy$} [ ] at 336 219
 \pinlabel {$stx$} [ ] at 289 228
 \pinlabel {$x$} [ ] at 268 239
 \pinlabel {$ux$} [ ] at 254 232
 \pinlabel {$u$} [ ] at 241 221
 \pinlabel {$tuy$} [ ] at 224 219
 \pinlabel {$-$} [ ] at -10 231
\endlabellist
\centering
\ig{1.1}{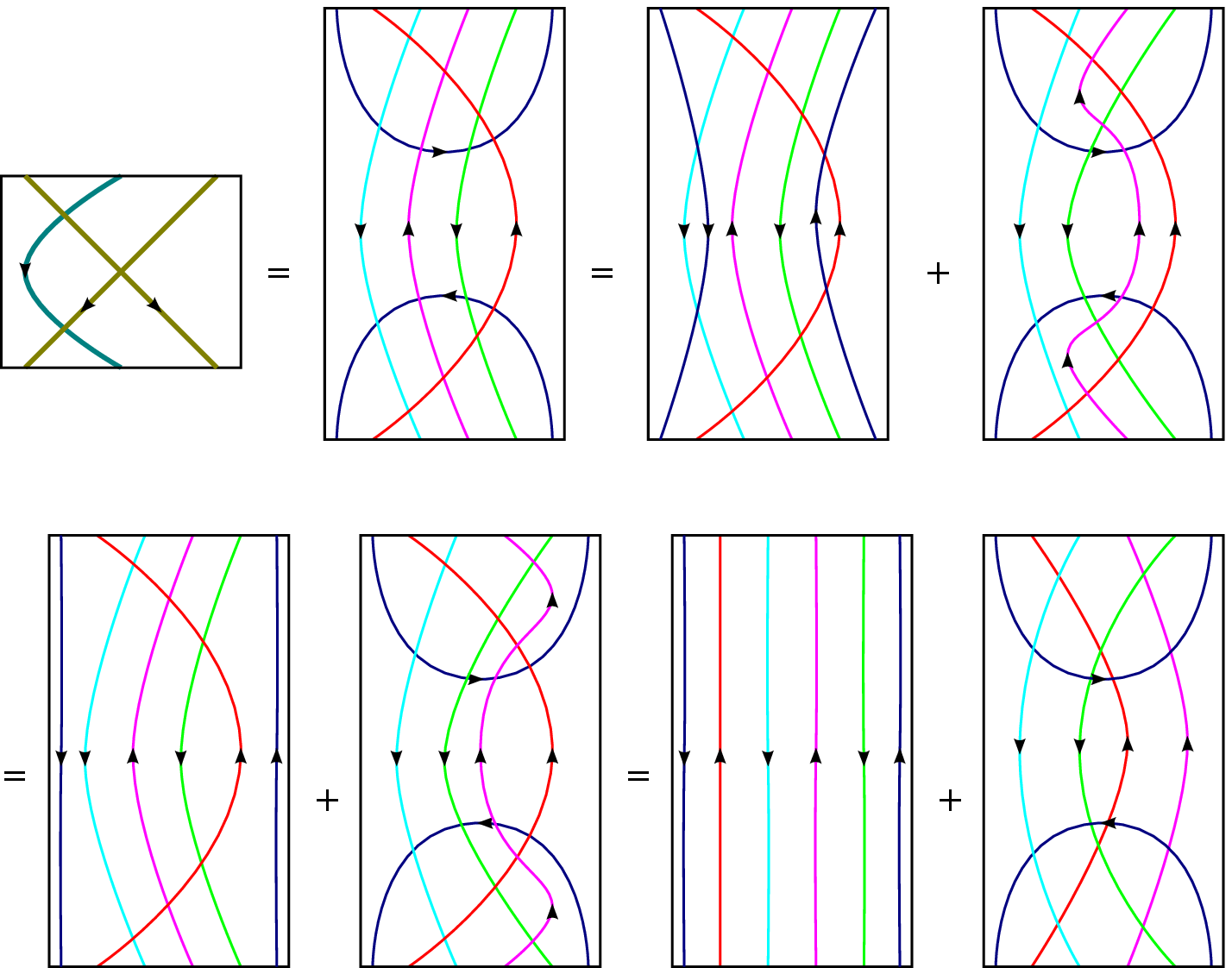}
} \end{equation}
To begin, we apply the interesting relation \eqref{relation11n1} inside the red-aqua bigon, where $s$ and $y$ are absent. The remaining simplifications use only unobstructed RII and RIII moves.

Thus we obtain
\begin{equation*} \ig{1}{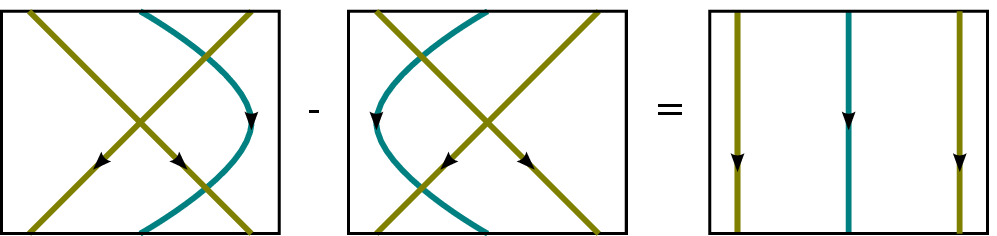} \end{equation*}
as desired. Note that the computation still holds in the absence of the red or aqua strands, which handles the easiest edge cases.

The remaining edge cases can be dealt with in similar fashion. We leave the confirmation to the reader.

\subsubsection{Conclusion}

We have checked all the relations of $\UC^-$, so the functor $\FC \co \UC^- \to \DG$ is well-defined.

\subsection{Confirmation of the thick calculus}
\label{subsec-conf-thick}

When extending a presentation for a category by generators and relations to a partial idempotent completion, one need only identify the idempotents whose images are being added as
objects. This is discussed in \cite[Section 2.4]{EThick}. Thus, generators and relations for $\dot{\UC}^-$, the partial idempotent completion of $\UC^-$ which adds the divided powers
$F_i^{(k)}$ as new objects, are not hard to come by. The ``thick calculus'' developed by Khovanov-Lauda-Mackaay-Stosic in \cite{KLMS} does this and more, giving a host of relations
which hold for morphisms between compositions of divided powers $E_i^{(k)}$ and $F_i^{(k)}$ for a single index $i$. Most of their relations are redundant, however, and we only check a
minimal set of relations here.

The thick calculus for $\dot{\UC}^-$, as constructed by \cite{KLMS}, includes the merge and split, and the dot on a thick strand, as new generators. We have already defined how $\FC$
acts on these new generators. Unfortunately, \cite{KLMS} is not written to make it clear which relations are foundational, and which are implied by others. We check some relations here,
and leave it as an exercise to prove that they are sufficient.

\subsubsection{Crossing is merge-split}

The first relation states that the crossing (in the usual quiver Hecke algebra) is a composition of a merge and a split. Checking this relation after applying the functor is trivial.
\begin{equation}\label{crossingismergesplit} {
\labellist
\small\hair 2pt
 \pinlabel {\tiny{$2$}} [ ] at 84 31
\endlabellist
\centering
\ig{1}{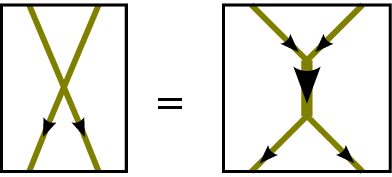}
} \end{equation}

\subsubsection{Split-merge}

The next relations also deal with a strand of thickness 2 splitting into strands of thickness 1.
\begin{equation} \label{splitmerge} {
\labellist
\small\hair 2pt
 \pinlabel {\tiny{$2$}} [ ] at 19 72
 \pinlabel {\tiny{$2$}} [ ] at 91 72
 \pinlabel {\tiny{$2$}} [ ] at 19 11
\endlabellist
\centering
\ig{1}{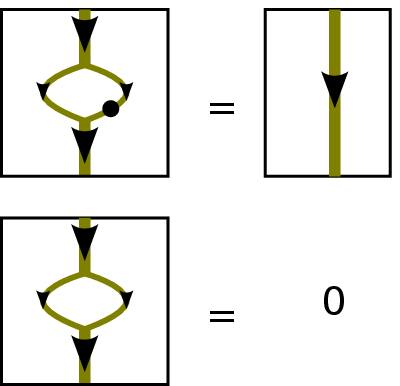}
} \end{equation}
The second relation checks out analogously to the computation in \eqref{NilHecke2}.
\begin{equation} \label{splitmergecheck} \ig{1}{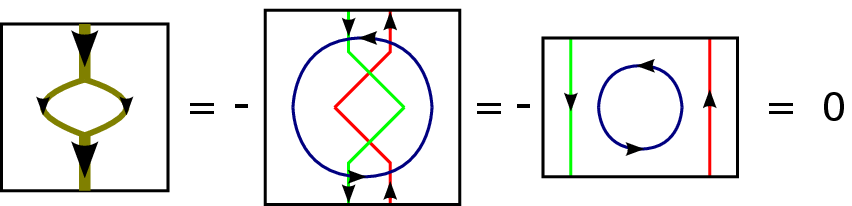} \end{equation}
Meanwhile, if a dot is placed on one of the strands, the corresponding polynomial ends up inside the blue circle. Thus, the first relation follows from \eqref{fis-tunit}.

Together, \eqref{crossingismergesplit} and \eqref{splitmerge} imply that the merge (resp. the split) of two thin strands into a thick strand is equal to the projection from $F_i F_i$ to
$F_i^{(2)}(-1)$ (resp. the inclusion of $F_i^{(2)}(+1)$ into $F_i F_i$) in the familiar decomposition $F_i F_i \cong F_i^{(2)}(-1) \oplus F_i^{(2)}(+1)$.

\subsubsection{Associativity}

The next ``associativity" relation shows that thicker strands are defined correctly. The line thicknesses in the source 1-morphism are arbitrary.
\begin{equation} \label{mergeassoc} \ig{1}{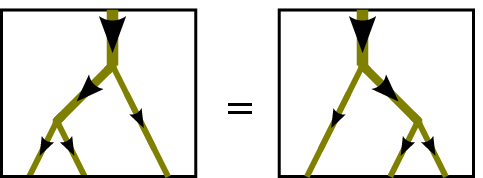} \end{equation}
This equality is checked exactly as \eqref{NilHecke3bottom} and \eqref{NilHecke3bottom-2}. Having checked that splits and merges for thickness 2 behave correctly, \eqref{mergeassoc} implies that all the thick merges and splits behave correctly.

\subsubsection{Thick dot} \label{sssec-mergedot}

Finally, we check that the dot on a thick strand is defined correctly. Again, the line thicknesses in the source 1-morphism are arbitrary.
\begin{equation} \label{mergedot} \ig{1}{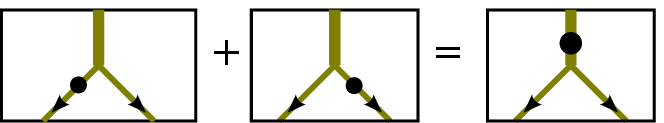} \end{equation}
After applying the functor and multiplying by a sign, we have
\begin{equation} \label{mergedotcheck} {
\labellist
\tiny\hair 2pt
 \pinlabel {$\fpoly_{\nu \to \mu}$} [ ] at 113 10
 \pinlabel {$\fpoly_{\la \to \nu}$} [ ] at 18 10
 \pinlabel {$\fpoly_{\la \to \mu}$} [ ] at 182 28
\endlabellist
\centering
\ig{2}{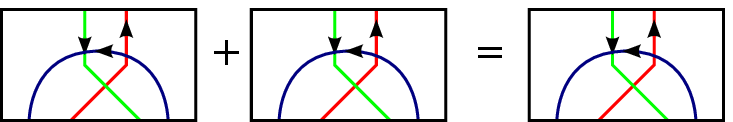}
} \end{equation}
All three polynomials can be slid into the central triangle, where the result is equivalent to \eqref{fis-additive}.

The modifications for the ends of an $F_i$-string are obvious, and we leave them to the reader.

\subsection{Definition of the categorical action: cups, caps, and bubbles}

Now we extend the functor $\FC$ from $\UC^-$ to $\UC$. The images of the cups and caps are obvious; here is an example.
\[ \ig{1}{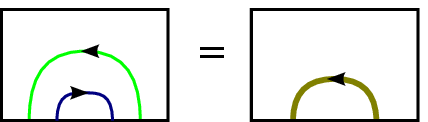} \]
Near the ends of the $F_i$-string, one simply removes one of the strands, as usual.

Because bubbles, both fake and real, play an essential role in the category $\UC$, let us discuss the images of the bubbles. We begin with a clockwise bubble labeled $i$, with $m$ dots, in a region $\la$.
\begin{equation} \label{CCWbub} {
\labellist
\small\hair 2pt
 \pinlabel {$\la$} [ ] at 39 11
 \pinlabel {$\mu$} [ ] at 24 34
 \pinlabel {\tiny{$m$}} [ ] at 41 32
 \pinlabel {$Ls$} [ ] at 117 7
 \pinlabel {$L$} [ ] at 96 20
 \pinlabel {$Lt$} [ ] at 92 34
 \pinlabel {$\fpoly^m$} [ ] at 110 35
 \pinlabel {$Ls$} [ ] at 189 7
 \pinlabel {$\pi_m$} [ ] at 171 38
\endlabellist
\centering
\ig{1}{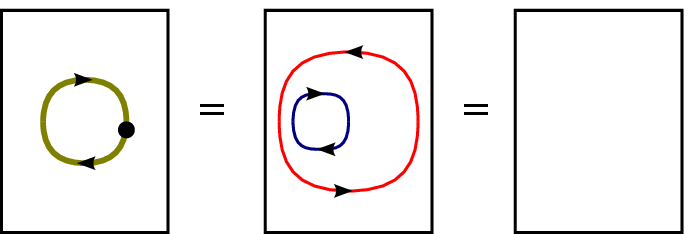}
} \end{equation}
Here, $\fpoly = \fpoly_{\la \to \mu}$. One deduces that $\pi_m$ is the polynomial
\begin{equation} \label{defnpi} \pi_m = \pa^L_{Ls} (\fpoly_{\la \to \mu}^m \mu^L_{Lt}). \end{equation}

Let us compute this in the familiar language of symmetric polynomials, in order to define $\pi_m$ for $m < 0$. We assume that $\fpoly_{\la \to \mu} = \epoly_{\la \to \mu}$ (otherwise, replace each
instance of $x_i$ with $x_i + y$). Suppose that $\la$ has $a+1$ copies of $i$ and $b$ copies of $i+1$. For example, when $a = 3$, $b=6$ and $i=3$ we have \[\la = (\ldots
333{\color{red}3}444444\ldots), \qquad \mu = (\ldots 333{\color{red}4}444444\ldots). \] Let $y_1, \ldots, y_a$ denote the variables $x_k$ for those $k$ with $\mu_k = i$; let $z_1, \ldots, z_b$ denote the variables $x_k$ for those $k$ with $\la_k = i+1$, and let $f$ denote $x_k$ for the red index above. Then one has
\[ \mu^L_{Lt} = \prod_{j=1}^{b} (f - z_j) = \sum_{p + q = b} (-1)^q h_p(f) e_q(\un{z}), \]
where $h_p$ and $e_q$ denote complete and elementary symmetric polynomials, respectively. Thus
\[ f^m \mu^L_{Lt} = \sum_{p + q = b} (-1)^q h_{p+m}(f) e_q(\un{z}) = \sum_{p + q = m+b} (-1)^q h_p(f) e_q(\un{z}). \]
Because
\[ \pa_k (h_p(x_{k+1}, x_{k+2}, \ldots, x_{k+d})) = - h_{p-1} (x_k, x_{k+1}, \ldots, x_{k+d}) \]
we have that
\begin{equation} \label{calcpi} \pi_m = \pa^L_{Ls} (f^m \mu^L_{Lt}) = (-1)^a \sum_{p + q = m + b} (-1)^q h_{p-a}(\un{y}, f) e_q(\un{z}) = (-1)^{m+b-2a} \sum_{p+q = m+b-a} (-1)^p h_p(\un{y}, f) e_q(\un{z}). \end{equation}
In this latter formula, there is clearly no problem defining $\pi_m$ for $m < 0$, and the result will be nonzero so long as $m+b-a \ge 0$. Note that when $m+b-a = 0$, one obtains $\pi_m = (-1)^a$.

Similarly, let us consider a counterclockwise bubble labeled $i$, with $m$ dots, in region $\la$. The parabolic subset $L$ from the previous paragraph is $Iu$ where $u$ is the simple reflection after $t$. Let $K = Is$.
\begin{equation} \label{CWbub} {
\labellist
\small\hair 2pt
 \pinlabel {$\la$} [ ] at 39 11
 \pinlabel {$\nu$} [ ] at 24 34
 \pinlabel {\tiny{$m$}} [ ] at 41 32
 \pinlabel {$Ku$} [ ] at 117 7
 \pinlabel {$K$} [ ] at 96 20
 \pinlabel {$Kt$} [ ] at 92 34
 \pinlabel {$g^m$} [ ] at 110 35
 \pinlabel {$Ku$} [ ] at 189 7
 \pinlabel {$\xi_m$} [ ] at 171 38
\endlabellist
\centering
\ig{1}{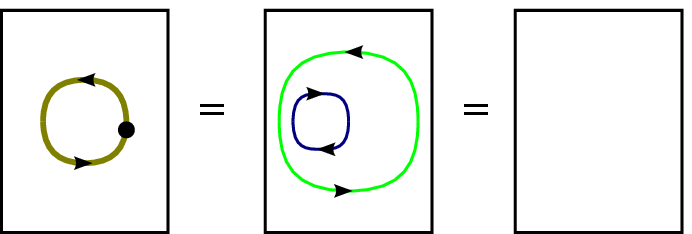}
} \end{equation}
where $g = \fpoly_{\la \to \nu}$ is equal to $t(f)$ or to $z_1$ from the previous paragraph. Repeating the argument, we get
\begin{equation} \label{defnxi} \xi_m = \pa^K_{Ku} (g^m \mu^K_{Kt}) \end{equation}
and
\begin{equation} \label{calcxi} \xi_m = (-1)^{m-b} \sum_{p+q = m+a-b} h_p(\un{z}) e_q(\un{y}, f). \end{equation}
Once again, when $m+a-b = 0$, the result is just the sign $(-1)^a$.

From \eqref{calcpi} and \eqref{calcxi}, the infinite Grassmannian relation is clear, which implies that the fake bubbles are well-defined. We leave the reader to confirm
that these formulas still hold at the ends of the string.

\begin{remark} \label{toobad} One of our goals was to describe the action of $\UC$ on $\DG$ using Frobenius-theoretic language, rather than the language of symmetric groups. While this
was successful for defining and checking the action of $\UC^-$, we have been unable to find a good formula for fake bubbles using only the Frobenius structures and the chosen
polynomials $\fpoly_{\la \to \mu}$. \end{remark}

\subsection{Confirmation of the full categorical action}

The first relation we treat from $\UC$ states that $E_j$ and $F_i$ can be ``pulled apart'' so long as $i \ne j$, which is \cite[(3.16)]{MSV}. This is easy to confirm, being very similar to the computation in \eqref{KLR3}.
\begin{equation} \label{EiFjnormal} \ig{1}{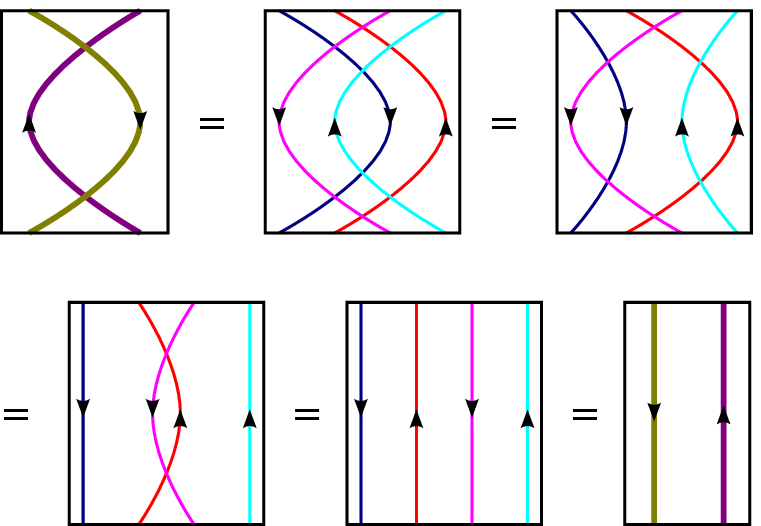} \end{equation}
The edge cases when $j = i \pm 1$ can be slightly different, but equally easy. Consider, for example, the case when $\la = (2244)$, $j=2$, $i=3$, and the simple reflections are $s,t,u$.
\begin{equation} \label{EiFjedge} {
\labellist
\small\hair 2pt
 \pinlabel {$\la$} [ ] at 472 241
 \pinlabel {$su$} [ ] at 553 243
 \pinlabel {$s$} [ ] at 539 233
 \pinlabel {$u$} [ ] at 539 259
 \pinlabel {$\emptyset$} [ ] at 525 245
 \pinlabel {$t$} [ ] at 511 246
\endlabellist
\centering
\ig{1}{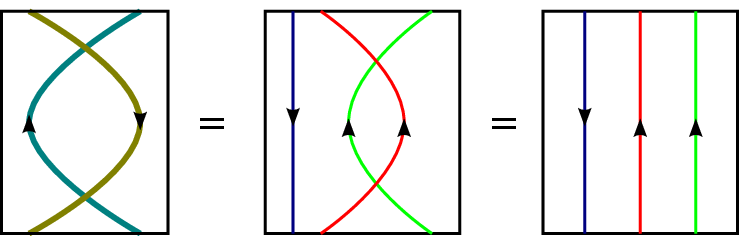}
} \end{equation}
In this example, the vertical blue line is actually made from the cup and cap of \eqref{FiFjXspecial}, turned sideways.

The remaining relations deal with $\sl_2$: straightening a curl, and pulling apart $E_i$ and $F_i$. Both relations involve bubbles, and as discussed in Remark \ref{toobad}, the
Frobenius perspective does not add much. Ultimately, these are checked by manipulating polynomials, precisely as in the original computations of Lauda in \cite{LauSL2}, and we do not
see a need to repeat the argument with our minor sign variations.

%% file: AffineADiags2.tex
\section{The Kac-Moody action in affine type $A$}\label{S_KM_affine_A}
\label{sec-quantumaction-affine}

\subsection{Preliminaries}\label{SS_affine_prelim}

Fix $n \ge 2$ and $e \ge 2$.\footnote{When comparing with \cite{MacThi}, the $n$ in this paper is the $r$ in \cite{MacThi}, and the $e$ in this paper is the $n$ in \cite{MacThi}.} Let $(W,S)$ continue to denote the finite type $A$ Coxeter system. Let $(W_a,S_a)$ denote the affine type $A$ Coxeter system, with $S_a = S \cup \{s_0\}$.

\begin{defn} Define an action of $W_a$ on $\La = \ZM^n$ as follows. The parabolic subgroup $S_n$ acts in the usual way. The affine reflection $s_0$ acts by $s_0(\la_1, \la_2, \ldots,
\la_{n-1}, \la_n) = (\la_n - e, \la_2, \ldots, \la_{n-1}, \la_1 + e)$. Let $\tL_{n,e}$ denote the set of $\la \in \ZM^n$ which satisfy the following condition: \begin{equation}
\label{orbitrepcondition} \la_1 \le \la_2 \le \ldots \le \la_n \le \la_1 + e. \end{equation} \end{defn}

The orbits of $W_a$ on $\La$ are in bijection with $\tL_{n,e}$, as each orbit has one distinguished representative in this set. Each inequality in \eqref{orbitrepcondition} which is
an equality corresponds to a simple reflection in $S_a$ which stabilizes $\la$. Given $\la \in \tL_{n,e}$, let $I(\la)$ denote the set of simple reflections in the stabilizer. Note that
$I(\la)$ is always a proper subset of $S_a$, and $I(\la)$ has at least $n-e$ elements.

One has $\La_{n,e} \subset \tL_{n,e}$. The same is true for $\La_{n,e} + c \cdot (1,1,\ldots,1)$ for any $c \in \ZM$. Together, the union $\cup_{c \in \ZM} \Big( \La_{n,e} + c \cdot (1,\ldots,1) \Big)$ is
precisely the subset of $\tL_{n,e}$ for which the stabilizer does not contain the affine reflection.

\begin{defn} For $0 \le j \le e-1$, let $F_j$ denote the operator on $\tL_{n,e} \cup \{0\}$ defined as follows. If $\la$ has no indices $\la_i$ equivalent to $j$ modulo $e$, then
$F_j(\la)$ is the formal symbol $0$. Otherwise, there is a unique $i$ such that $\la_i$ is equivalent to $j$ modulo $e$, and for which replacing $\la_i$ with $\la_i + 1$ still yields an
element of $\tL_{n,e}$, which we call $F_j(\la)$. As before, let $F_j^{(k)}$ denote the $k$-th iterate of $F_j$, and let $I(\la, \mu)$ denote the intersection $I(\la) \cap
I(\mu)$, when both $\la$ and $\mu$ are nonzero in the same $F_j$-string. \end{defn}

We also occassionally let $F_e$ denote $F_0$. In general, colors are considered as elements of $\ZM / e \ZM$.

\begin{example} \label{F0chain} Here is an $F_0$-chain when $e = 3$: $(000233) \to (001233) \to (011233) \to (111233) \to (111234) \to (111244)$. \end{example}

The way in which the stabilizing parabolic subgroups change within an $F_j$-chain is the same as in the finite case. There is an interesting copy of $S_k$ inside the stabilizer of the
start and finish, together with some distant simple reflections. As the chain progresses, the interesting part of the stabilizer varies: $S_k \to S_{k-1} \times S_1 \to S_{k-2} \times
S_2 \to \ldots \to S_1 \times S_{k-1} \to S_k$. The interesting copy of $S_k$ is just some type $A$ subdiagram of the affine Dynkin diagram, and may involve the affine reflection. The
identification with $S_k$ uses the cyclic ordering on the affine Dynkin diagram.

\begin{example} The stabilizers in Example \ref{F0chain} are $\{s_5, s_0, s_1, s_2\} \to \{s_5, s_0, s_1\} \to \{s_5, s_0, s_2\} \to \{s_5, s_1, s_2\} \to \{s_0, s_1, s_2\} \to \{s_5,
s_0, s_1, s_2\}$. \end{example}

\begin{defn}
 To $\la \in \tL_{n,e}$ we associate a $\tgl_e$-weight $\la$, as follows. This will be a level 0 weight
 meaning that the pairing of $\la$ with the affine root $\delta$ is zero. The value of $\la$ on $E_{ii}\in \gl_e\subset
 \tilde{\gl}_e$ equals to the number of the indices $j$ such that $\lambda_j-i$ is divisible by $e$. 
 The set of these weights we also denote $\tL_{n,e}$.  \end{defn}

\subsection{A modification of the diagrammatic singular Soergel category} \label{subsec-weirdSoergel}


Let us recall the $\tgl_n$-realization from \S\ref{subsec-diagaffinetype}. The base ring $\Bbbk$ is $\ZM$. The free $\Bbbk$-module $\hg^*$ has basis $\{x_i\}_{i = 1}^n \cup \{y\}$,
while $\hg$ has the dual basis $\{\ep_i\}_{i=1}^n \cup \{\phi\}$. The simple roots $\a_i$ for $i \in S$ are still given by $\a_i = x_i - x_{i+1}$, while the simple root $\a_0$ is given
by $x_n - x_1 + y$. The simple coroots for $i \in S$ are still given by $\a_i^\vee = \ep_i - \ep_{i+1}$, and the simple coroot $\a_0^\vee$ is given by $\ep_n - \ep_1$. With
these definitions, $R = \ZM[x_1, \ldots, x_n, y]$ where $y$ is fixed by $W_a$, $S_n$ acts as usual on the variables $x_i$, and $s_0(x_n) = x_1 - y$, $s_0(x_1) = x_n + y$.

As constructed in Section \ref{subsec-diagaffinetype}, the 2-category $\DG_a$ of diagrammatic singular Soergel bimodules in affine type $A$ has objects given by the poset
$\Gamma(W_a,S_a)$ of finitary subsets of $S_a$. Let us now define a very similar $2$-category.

\begin{defn} Let $\DG(\tL_{n,e})$ denote the following 2-category. The objects are given by the set $\tL_{n,e}$. For $\la, \mu \in \tL_{n,e}$, the morphism category $\Hom(\la, \mu)$ is
precisely the morphism category $\Hom(I(\la), I(\mu))$ in $\DG_a$, and the horizontal composition of $1$-morphisms is defined analogously.

There is an obvious 2-functor $\QC \co \DG(\tL_{n,e}) \to \DG_a$, which sends each object $\la \in \tL_{n,e}$ to the parabolic subset $I(\la)$. This 2-functor induces equivalences on
morphism categories. \end{defn}

Loosely speaking, $\DG(\tL_{n,e})$ is the same 2-category as $\DG_a$, only with the objects renamed and duplicated. When $e \ge n$, the 2-functor $\QC$ induces a Morita equivalence of
2-categories. However, when $e<n$ every element of $\tL_{n,e}$ has a nontrivial stabilizer, so some objects of $\DG_a$ (such as the empty parabolic subset) will never appear in the image
of $\QC$, and $\QC$ will not be a Morita equivalence.

\begin{lemma} Let $\CC$ be a 2-category with objects parametrized by a set $X$. Then the following data are equivalent. \begin{itemize} \item A 2-functor $\tilde{\GC} \co \CC \to
\DG(\tL_{n,e})$. \item A 2-functor $\GC \co \CC \to \DG_a$, and a map $\tilde{g}: X \to \tL_{n,e}$ such that $q \circ \tilde{g} = g$. Here, $q \co \tL_{n,e} \to \Gamma(W_a,S_a)$ is the
map induced by $\QC$, sending $\la \in \tL_{n,e}$ to the parabolic subset $I(\la)$, and $g \co X \to \Gamma(W_a,S_a)$ is the map on objects induced by $\GC$. \end{itemize} Under this
equivalence of data, one has $\FC = \QC \circ \tilde{\FC}$. \end{lemma}

\begin{proof} This is obvious. \end{proof}

Our next goal (roughly) is to define a 2-functor $\FC \co \UC(\tgl_e) \to \DG_a$, which will send a $\tgl_e$-weight $\la \in \tL_{n,e}$ to the parabolic subset $I(\la)$. By the above
lemma, it is immediate that there will be a 2-functor $\tilde{\FC} \co \UC(\tgl_e) \to \DG(\tL_{n,e})$ lifting $\FC$, sending the object $\la$ to $\la$. It is easier to define the
functor $\FC$ (see the remark below), which we do in the rest of this chapter. However, it is the 2-functor $\tilde{\FC}$ which is natural to study from the viewpoint of representation
theory, and this is the functor we will abusively refer to as $\FC$ in subsequent sections.

\begin{remark} There is one slight subtlety/annoyance in the definition of $\DG(\tL_{n,e})$. Consider for example the 1-morphism in $\DG(\tL_{n,e})$ from $\la = (11222)$ to $\mu =
(12222)$, which induces from $I(\la)$ to $I(\la, \mu)$, and then restricts to $I(\mu)$; this would be the image of $F_1 \1_{\la}$ under the functor from $\UC(\tgl_e)$. It is clearly a
1-morphism in $\DG(\tL_{n,e})$ because it arises from a 1-morphism in $\DG_a$, but the object $I(\la, \mu)$ which it factors through does not naturally come from any weight in
$\tL_{n,e}$ (the non-integral weight $(1,1.5,2,2,2)$ has the correct stabilizer). It factors as a composition of two 1-morphisms in $\DG_a$ (induction and then restriction), but it does
not factor in $\DG(\tL_{n,e})$.

It would be nice if one could picture morphisms in $\DG(\tL_{n,e})$ as $S_a$-diagrams with regions labeled by $\tL_{n,e}$ instead of $\Gamma(W_a, S_a)$. This is not strictly valid
because not every region can be consistently labeled by elements of $\tL_{n,e}$. Nonetheless, it is still effective to use $S_a$-diagrams to compute with 2-morphisms in $\DG(\tL_{n,e})$,
because one is simply computing in $\DG_a$ after applying the fully faithful 2-functor $\QC$.

It is not correct to think of $\tL_{n,e}$ as indexing a poset of compatible Frobenius extensions! The poset $\tL_{n,e}$ is too small to be built from hypercubes. It is perhaps possible
to come up with a natural combinatorial poset $\Gamma$ extending $\tL_{n,e}$, for which one could define $\DG(\Gamma)$ as the diagrammatic category associated to a poset of Frobenius
extensions. We have not bothered to attempt this. Instead, we perform all our computations in $\DG_a$. \end{remark}

\subsection{A modification of $\UC$} \label{subsec-objective}

As we have noted in the introduction, there is no 2-functor $\UC(\tgl_e) \to \DG_a$. Instead, we must define a modification of $\UC(\tgl_e)$ which will admit such a 2-functor. The problem is that the relations in $\UC(\tgl_e)$ impose conditions on the polynomials $\fpoly_{\la \to \mu}$ assigned to dots, and these conditions are contradictory. We discuss this in detail in \S\ref{sssec-failure} below.

\subsubsection{Properties of polynomials} \label{sssec-morepolys}
Let us begin by exploring some of the conditions on $\fpoly_{\la \to \mu}$ implied by the one-color relations and the different colored dot sliding relation. These are relations which should hold in any variant on $\UC$.

\begin{prop} \label{prop:junique} Let $\{\fpoly_{\la \to F_j^{(k)} \la}\}$ be a family of polynomials, indexed by $\la \in \tL_{n,e}$, $k \ge 1$, and $1 \le j \le e$ such that $F_j^{(k)} \la \ne 0$.
\begin{enumerate}
\item Suppose that this family satisfies invariance, additivity, the dual basis condition, and locality from \S\ref{sssec-polyprops}. Then, for any fixed $j$, if one knows the value of any given $\fpoly_{\nu \to F_j \nu}$, then one knows the values of all $\fpoly_{\la \to F_j^{(k)} \la}$.
\item	
Moreover, suppose that $F_j$ increments the $k$-th index of $\nu$ to $j$ to $j+1$, for $1 \le j \le e$, and consider the special case when $\fpoly_{\nu \to F_j \nu} = x_k + z$ for $z \in R^{W_a}$. Then, if $F_j$ increments the $\ell$-th index of $\la$ from $j + re$ to $j + 1 + re$ for $r \in \ZM$, then $\fpoly_{\la \to F_j \la} = x_\ell + z + ry$.
\end{enumerate}
\end{prop}

\begin{proof} We follow the outlines of the proof of Proposition \ref{onlychoices}. Fix $1 \le j \le e$ and $\nu \in \tL_{n,e}$ such that $F_j$ increments the $k$-th index of $\nu$ from
$d$ to $d+1$. Note that $d = j + re$ for some $r \in \ZM$. Let $f = \fpoly_{\nu \to F_j \nu}$. By locality, it is easy to deduce that $\fpoly_{\la \to F_j \la} = f$ whenever $F_j$ also
increments the $k$-th index of $\la$ from $d$ to $d+1$, because $\la \to F_j \la$ is parallel to $\nu \to F_j \nu$. Just as in the proof of Proposition \ref{onlychoices}, one can choose
(part of) an $F_j$-string which increments from $d$ to $d+1$ in various indices, and the dot polynomials along this string will be determined from $f$ by applying a sequence of simple
reflections, thanks to \eqref{fis-swapped}. For example, if $f = x_k + z$, then a dot polynomial $\fpoly_{\la \to F_j \la}$ where $F_j$ increments the $\ell$-th index from $d$ to $d+1$ will be $x_\ell + z$.

Unlike the finite case in Proposition \ref{onlychoices}, there are $F_j$ strings which increment some indices from $d$ to $d+1$, and increment other indices from $d+e$ to $d+e+1$. For
example, consider an element $\la \in \tL_{n,e}$ with first index $d$, second index $>d$, and last index $(d+e)$. If $g = \fpoly_{\la \to F_j \la}$, then $h = \fpoly_{F_j \la \to F_j F_j
\la} = s_0(g)$, and this dot polynomial increments the last index from $d+e$ to $d+e+1$. For example, if $g = x_1 + z$, then $h = s_0(g)= x_n + z + y$. This proves (1).

Repeating this argument, we see that any application of $F_j$ is determined from our original choice $f$. (2) of the proposition also follows. \end{proof}

\begin{lemma} \label{lem:jchoice} Again, suppose $\{\fpoly\}$ is a family of polynomials satisfying invariance, additivity, the dual basis condition, and locality from
\S\ref{sssec-polyprops}. In the $\tgl_n$ realization or any extension thereof, if $F_j$ increments the $k$-th index of $\la$, then $\fpoly_{\la \to F_j \la} = x_k + z_j$ for some $z_j
\in R^{W_a}$ only depending on $j$. \end{lemma}

\begin{proof} Let $\nu = (jj\ldots j)$. Just as in the proof of Proposition \ref{onlychoices}, one can show that $f = \fpoly_{\nu \to F_j \nu}$ is equal to $x_n + z$ for some $z \in
R^{W}$, that is, $z$ is symmetric under the finite Weyl group. By Proposition \ref{prop:junique}, one discovers that $\fpoly_{\la \to F_j \la}$ is determined for all $\la \in \tL_{n,e}$
with $F_j \la \ne 0$, and it has the form $x_k + z + ry$. In particular, $z$ was independent of the choice of $\la$, and $z \in R^{W_a}$ by invariance. \end{proof}

\subsubsection{Naivet\'e and failure}\label{sssec-failure}

Khovanov and Lauda have defined a $2$-category $\UC(\tilde{\sl}_e)$, as they have for any oriented Dynkin diagram. To define a $2$-category $\UC(\tgl_e)$, one should have objects given by
$\tgl_e$ weights instead of $\tilde{\sl}_e$ weights, and should fix the signs on cups and caps as per Mackaay-Stosic-Vaz \cite{MSV}. The $1$-morphisms and $2$-morphisms are defined so
that there is an inclusion $\UC(\gl_e) \to \UC(\tgl_e)$ for any finite Dynkin diagram inside the affine Dynkin diagram (appropriately oriented). We think of this 2-category as being ``vanilla,'' the simplest flavor of a construction of $\UC(\tgl_e)$.

A first hope would be to define a $2$-functor from $\UC(\tgl_e)$ on $\DG(\tgl_n)$. On objects, this $2$-functor should send any $\tgl_e$ weight not in the set $\tL_{n,e}$ to the formal
zero object, and should send any weight $\la \in \tL_{n,e}$ to its stabilizing parabolic $I(\la)$. Unfortunately, there is no ``vanilla" 2-functor $\FC \co \UC(\tgl_e) \to \DG(\tgl_n)$
extending the finite case of the previous chapter, unless one specializes to $y=0$.

The problem has to do with ``monodromy" in the definition of the polynomials $\{\fpoly_{\la \to \mu}\}$. In order to define the vanilla 2-functor $\FC$, one must choose dot polynomials,
and they must satisfy the properties \eqref{f-properties} from \S\ref{sssec-polyprops}.

\begin{example} \label{ex:failure} Let $\nu_m = (mm \ldots m)$ for various $m \in \ZM$. By Proposition \ref{prop:junique} and Lemma \ref{lem:jchoice}, if $1 \le j \le e$ is equal to $m$ modulo $e$, then $\fpoly_{\nu_m \to F_j \nu_m} = x_n + z_m$ for some $z_m \in R^{W_a}$, and moreover, $z_{m+e} = z_m + y$. However, if the repeated index condition is supposed to hold, then it is easy to deduce that $z_{m+1} = z_m$ for all $m$, and this is clearly a contradiction. \end{example}

The properties of dot polynomials used in Proposition \ref{prop:junique} rely on relations of $\UC$ that are quite important and ``immutable'': the nilHecke relations from the $\sl_2$
case, and different colored dot sliding. The repeated index condition, on the other hand, relies on the (adjacent-colored) double crossing relations in $\UC$, which are more
``malleable'': Rouquier \cite{Rouq-2KM} defines an entire family of possible double crossing relations one could use to define $\UC$. Thus, one expects to fix this contradiction by
changing the double crossing relations in $\UC$, and then defining a 2-functor. There is a compatibility condition between the (adjacent-colored) double crossing and the (interesting)
triple crossing relations, so that one may need to adjust the triple crossing relations as well.

\begin{remark} For our applications to the modular representation theory we only use the specialization $y=0$, so that a vanilla functor will suffice. However, the general case is still worth having in the literature, and may have  its practical uses: modulo some technical results including an extension of the Brundan-Kleshchev isomorphism from
\cite{BK_KLR}, this will allow to significantly simplify arguments in Section \ref{S_Fock_unique}.
\end{remark}


\subsubsection{Tweaking $\UC$}

In order to obtain a 2-category which actually acts on $\DG(\tgl_n)$, one must tweak the double crossing relation in $\UC(\tgl_e)$ in order to tweak the requirements \eqref{f-properties},
so that the functions $\fpoly_{\la \to \mu}$ can pick up some monodromy as they voyage through $\tL_{n,e}$. This tweaking was accomplished by Mackaay and Thiel.

When $e>2$, Mackaay and Thiel in \cite[Definition 3.18]{MacThi} construct a category $\UC(\tgl_e)_{[y]}$, linear over an extra central parameter $y$ of degree $2$. It agrees with the
category $\UC(\tgl_e)$ discussed above except in the adjacent-color double crossing relation, which is (3.48) in \cite{MacThi}. They also construct a variant on singular Soergel
bimodules in the algebraic setting (categorifying the extended affine Hecke algebra). Then, under the assumption $e > n > 2$, they construct in \cite[Section 5.2]{MacThi} a 2-functor
$\FC_{MT}$ from $\UC(\tgl_e)_{[y]}$ to their extended singular Soergel bimodules.

The other results from \cite{MSV} are lifted to the affine case as well. They define the affine Schur category $\hat{\SC}(e,n)_{[y]}$ as a quotient of $\UC(\tgl_e)_{[y]}$ by those
objects which are sent to zero by $\FC_{MT}$. They construct a functor $\hat{\Sigma}_{n,e}$ from diagrammatic extended (non-singular) Bott-Samelson bimodules back to $\hat{\SC}(e,n)_{[y]}$.

We construct a 2-functor $\FC \co \UC(\tgl_e)_{[y]} \to \DG(\tgl_n)$, which lifts their functor $\FC_{MT}$ in some sense. After evaluating by applying the functor $\FC_{EW}$ from
diagrammatics to bimodules, one should have $\FC_{EW} \circ \FC \cong \FC_{MT}$, although this is not really true. Let us record the various functors in a large commutative diagram. We put $\EC$ before a category to indicate it is Mackaay-Thiel's extended version of that category.

\vskip .5cm

\begin{equation} \label{eq:bigdiagaffine} {
\labellist
\small\hair 2pt
 \pinlabel {$\UC_{\leftarrow}(\tgl_e)_{[y]}$} [ ] at 122 175
 \pinlabel {$\hat{\SC}(e,n)_{[y]}$} [ ] at 122 127
 \pinlabel {$\DG(\tL_{n,e})$} [ ] at 124 88
 \pinlabel {$\DG(\tgl_n)$} [ ] at 122 47
 \pinlabel {$\SBSBim$} [ ] at 122 3
 \pinlabel {$\EC\DG(\tgl_n)$} [ ] at 65 90
 \pinlabel {$\EC\SBSBim$} [ ] at 65 3
 \pinlabel {$\EC\SB$} [ ] at -2 90
 \pinlabel {$\EC\BSBim$} [ ] at -2 3
 \pinlabel {$\hat{\Sigma}_{n,e}$} [ ] at 50 118
 \pinlabel {$\tilde{\FC}$} [ ] at 128 108
 \pinlabel {$\QC$} [ ] at 112 69
 \pinlabel {$\FC$} [ ] at 147 76
 \pinlabel {``$\FC_{MT}$''} [ ] at 161 33
 \pinlabel {$\FC_{EW}$} [ ] at 111 25
 \pinlabel {$\FC_{MT}$} [ ] at 82 59
\endlabellist
\centering
\ig{1.5}{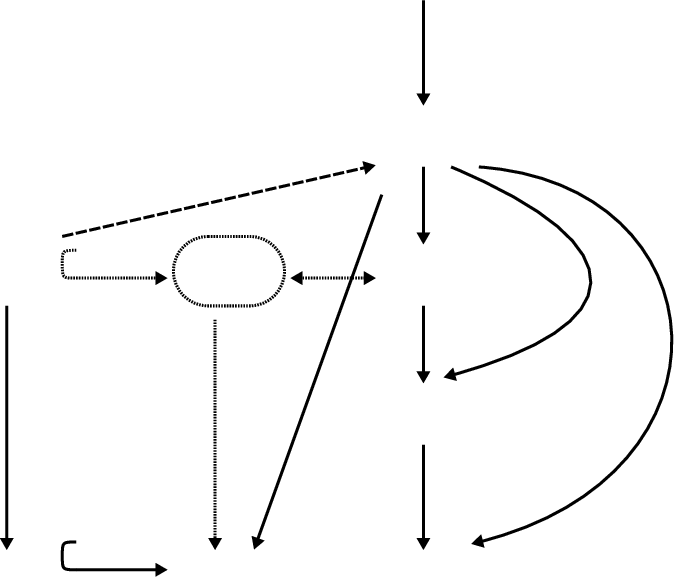}
} \end{equation}

\begin{remark} The same sign error from \S\ref{subsec-signerror} repeats itself here. We use the opposite orientation on the affine Dynkin diagram than was used by Mackaay-Thiel, and
consequently what we mean by $\FC_{MT}$ in this commutative diagram is a sign-modified version of their functor. \end{remark}

The main complication is that Mackaay-Thiel's extended (singular) Soergel bimodules are not quite the same 2-category as the singular Soergel bimodules for the affine Weyl group $\SSBim$
that were defined in \S\ref{subsec-algsbim}. Their 2-category is larger, including an extra twisting bimodule which acts as a rotation operator on the Dynkin diagram. Their 2-functor
$\FC_{MT}$ sends $F_0$ to a bimodule which has this twisting bimodule as a tensor factor. However, our category $\DG(\tgl_e)$ does not have this twisting operator, and thus the image of
$F_0$ under $\FC$ will not actually agree with the image under $\FC_{MT}$! A consequence of our different definition for $1$-morphisms is that our choice of polynomials $\fpoly_{\la \to
\mu}$ will depend on the weights $\la$ and $\mu$ and not just on the ambient parabolic subsets, while their polynomials do not have this dependency, making their 2-functor somewhat
simpler. We desired a categorical quantum group action of on a categorification of the affine Hecke algebra, not the extended affine Hecke algebra, so the extra complication was
necessary.

Here, quotes around $\FC_{MT}$ is defined to be the composition $\FC_{EW} \circ \FC$, but should be thought of being analogous to $\FC_{MT}$ in the non-extended world.

One guesses vaguely that there is a functor between $\DG(\tL_{n,e})$, the modified Soergel category from \S\ref{subsec-weirdSoergel}, and a (not yet constructed) diagrammatic version
of Mackaay-Thiel's extended singular Soergel bimodules, yielding the dashed part of the commutative diagram above. We have not attempted to make this precise.

\subsubsection{Improvements over Mackaay-Thiel}

Mackaay and Thiel assume throughout that $e>n>2$, and they do not even construct $\UC(\tgl_e)_{[y]}$ when $e=2$. We provide a construction for $e=2$ below in Section \ref{subsec-e2}, modifying
the double crossing relation via \eqref{KLRmodified} and \eqref{KLRmodified2} and the triple crossing relation via \eqref{funkyR3} and \eqref{funkyR3alt}. One feature of the
Frobenius-theoretic approach is that the correct version of the double crossing relation is easy to find, coming directly from computations analogous to \eqref{doublecrossing1} and
\eqref{doublecrossing2}.

The case $n=2$, omitted by Mackaay and Thiel, does not produce any new or interesting wrinkles in our perspective. Weights are almost always near the ends of various $F_j$-strings, which
makes computations easier, not harder.

Thus we can generalize Mackaay-Thiel's result to all $e,n \ge 2$, without the requirement $e>n$, providing a 2-functor from $\UC = \UC(\tgl_e)_{[y]}$ to $\DG_a = \DG(\tgl_n)$. Our
computations are quite different than theirs.

\subsection{Definition and confirmation of the categorical action}\label{subsec:e_bigger_2}

\subsubsection{Our version of $\UC$}

To save the reader some mental headache, let us record here the interesting relations in the category $\UC = \UC_{\leftarrow}(\tgl_e)_{[y]}$. These relations account for the sign changes
coming from the swap of Dynkin diagram orientation, and are written using the $1$-morphisms $F_i$ rather than $E_i$. The reader can take these relations as replacements for
\cite[(3.17)]{MSV} or \cite[(3.48)]{MacThi}, and \cite[(3.20)]{MSV} or \cite[(3.51)]{MacThi}. In other regards, the reader should follow the definition of $\UC(\tgl_e)_{[y]}$ in \cite[Section 3.5.1]{MacThi}.

In these diagrams, the symbol $y$ represents this new central parameter $y$, a degree $+2$ endomorphism of $\1_\la$ in $\UC$, for any weight $\la$. It is central in the sense that $F_i \ot y
= y \ot F_i$ and $E_i \ot y = y \ot E_i$. Diagrammatically, this means that $y$ can slide freely across all strands, so its placement is irrelevant.

When $e>2$:
\begin{equation} \label{newXX1} {
\labellist
\small\hair 2pt
 \pinlabel {$\delta_{i0}$} [ ] at 203 34
 \pinlabel {$y$} [ ] at 234 34
\endlabellist
\centering
\ig{1}{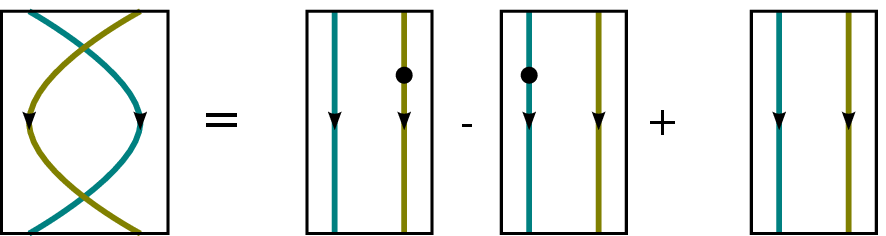}
} \end{equation}

\begin{equation} \label{newXX2} {
\labellist
\small\hair 2pt
 \pinlabel {$\delta_{i0}$} [ ] at 203 34
 \pinlabel {$y$} [ ] at 234 34
\endlabellist
\centering
\ig{1}{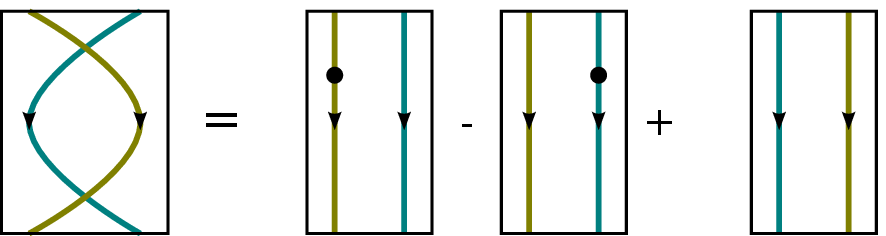}
} \end{equation}

\begin{equation} \label{newtriple1} \ig{1}{hardreln3} \end{equation}
	
\begin{equation} \label{newtriple2} \ig{1}{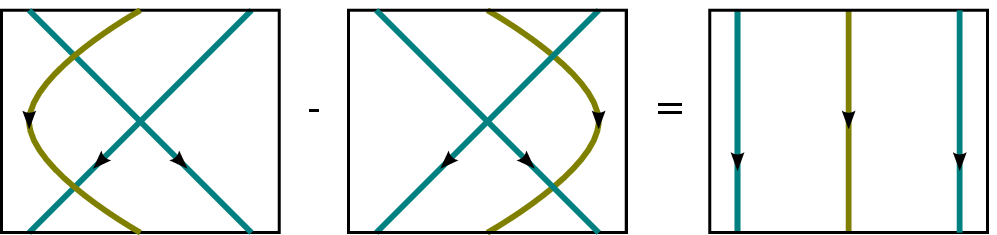} \end{equation}
In the relations above, olive is $i$ and teal is $j = i+1$. The term $\delta_{i0}$ indicates that this term only appears when $i$ is equal to $0$ modulo $e$. Note that \eqref{newtriple1} and \eqref{newtriple2} are unchanged by the $y$-deformation.

\begin{remark} Note that we have a plus sign in front of the $\delta_{i0} y$ term in \eqref{newXX1} and \eqref{newXX2}, while \cite[(3.48)]{MacThi} has a minus sign. We are unsure whether this sign arises from the same sign error that led to the decision to change the Dynkin diagram orientation, or whether it arises from a conventional change along the lines of Remark \ref{rmk:conventioncare} below. \end{remark}

When $e=2$:
\begin{equation} \label{KLRmodified} {
\labellist
\small\hair 2pt
 \pinlabel {$1$} [ ] at 40 13
 \pinlabel {$0$} [ ] at 10 13
 \pinlabel {$y$} [ ] at 350 38
\endlabellist
\centering
\ig{1}{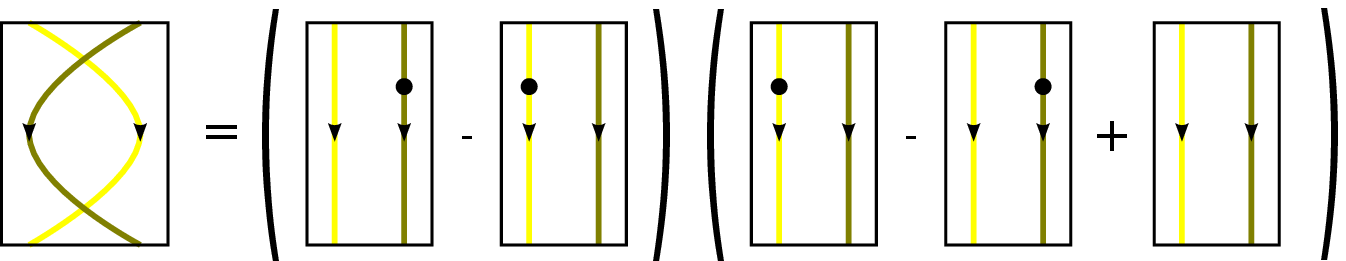}
} \end{equation}
\begin{equation} \label{KLRmodified2} {
\labellist
\small\hair 2pt
 \pinlabel {$0$} [ ] at 40 13
 \pinlabel {$1$} [ ] at 10 13
 \pinlabel {$y$} [ ] at 350 38
\endlabellist
\centering
\ig{1}{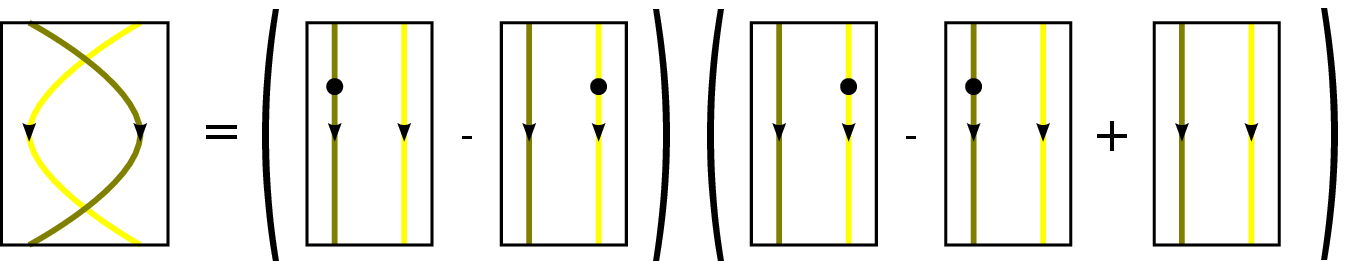}
} \end{equation}

\begin{equation} \label{funkyR3} {
\labellist
\small\hair 2pt
 \pinlabel {$y$} [ ] at 417 34
\endlabellist
\centering
\ig{1}{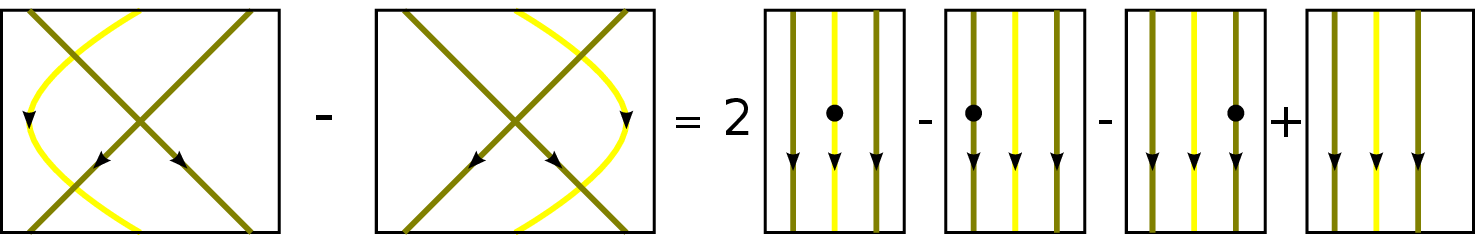}
} \end{equation}

\begin{equation} \label{funkyR3alt} {
\labellist
\small\hair 2pt
 \pinlabel {$y$} [ ] at 417 34
\endlabellist
\centering
\ig{1}{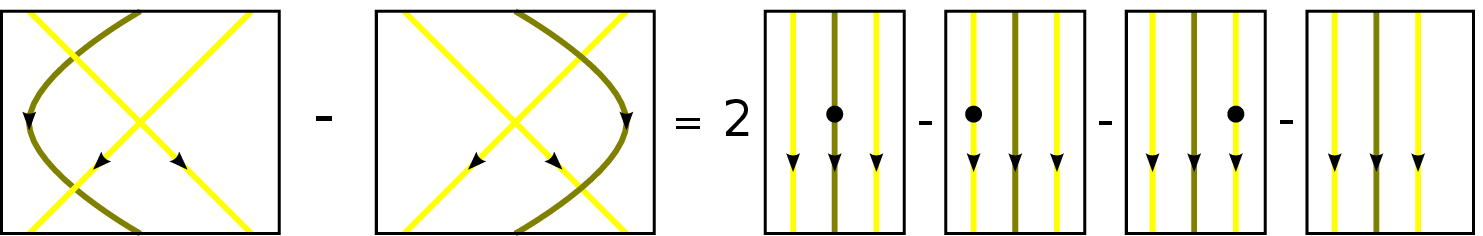}
} \end{equation}
In the relations above, olive is $1$ and yellow is $0$.

\begin{remark} The polynomial which appears in the double crossing relation and the polynomial appearing in the triple crossing relation must satisfy a consistency
requirement.\footnote{This consistency is part of Rouquier's definition in \cite{Rouq-2KM}. It can be seen in diagrammatic language in \cite[(2.6e) and (2.6g)]{WebsKIHRT}.} We describe how this consistency condition applies to
\eqref{KLRmodified}, \eqref{KLRmodified2}, \eqref{funkyR3}, and \eqref{funkyR3alt}. Namely, let $Q(\fpoly_0, \fpoly_1)$ is the polynomial in the dots appearing on the right hand side of
\eqref{KLRmodified} or \eqref{KLRmodified2}, where $\fpoly_0$ is a dot on the $0$-strand and $\fpoly_1$ is a dot on the $1$ strand. Then the polynomial on the right hand side of \eqref{funkyR3}
should be $$\frac{Q(\fpoly_0,\fpoly_{1,R}) - Q(\fpoly_0,\fpoly_{1,L})}{\fpoly_{1,R} - \fpoly_{1,L}},$$ where $\fpoly_{1,R}$ and $\fpoly_{1,L}$ represent dots on the $1$ strand on the right and left respectively. This is, in
fact, the case.
\end{remark}

Note that one also has related 2-categories $\UC^-$, $\dot{\UC}^-$, and $\dot{\UC}$, just as in \S\ref{subsec-prologue}. 

\subsubsection{Definition of the 2-functor}

We now define the $2$-functor $\FC$ from  $\UC$ to $\DG: = \DG(\tgl_n)$. On objects, a $\tgl_e$-weight $\la \in \tL_{n,e}$ will be sent to the parabolic subset $I(\la)$,
an object in $\DG$. If the weight $\la$ is not in $\tL_{n,e}$, it is sent to the formal zero object.

To the divided power $1$-morphism $F_j^{(k)}$ in $\UC$ we associate a $1$-morphism in $\DG$ defined exactly as in the finite case, as an induction from $I(\la)$ to $I(\la, \mu)$
followed by a restriction to $I(\mu)$, where $\mu = F_j^{(k)} \la$. To its biadjoint $E_j^{(k)}$ we assign the biadjoint $1$-morphism in $\DG$, as in the finite case.

To define the action of the 2-morphisms in $\UC$, we can follow the rubrick developed in the finite case. We need only deal with two subtleties: how to define the polynomial $\fpoly_{\la
\to \mu}$ which corresponds to the dot, and how to define the crossing of $F_0$ and $F_1$ in the special case when $e=2$ (this being the only generating morphism in $\UC$ for any $e$
which does not come from a finite type Dynkin subdiagram).

In \S\ref{sssec-polyprops} we listed precisely the properties of the polynomials $\fpoly_{\la \to \mu}$ which were needed for the functor to be well-defined. Now, the relations in $\UC$
have been modified, meaning that these properties will have to be modified accordingly. This is done below in Section \ref{subsec-polysredux}, giving the list of properties
\eqref{f-properties2}.

Most of the work needed to confirm that the functor is well-defined was accomplished in the previous chapter. So long as the relation in question depends only on a finite type
sub-Dynkin-diagram of $\tgl_e$ acting on a finite type sub-Dynkin-diagram of $\tgl_n$, and so long as our chosen polynomials $\fpoly_{\la \to \mu}$ satisfy the properties \eqref{f-properties2}, then relation is already confirmed. The only relations in $\UC$ for any $e$ which do not come from finite type subconfigurations are: several relations for the
special case $e=2$, and the triple crossing relation for the special case $e=3$ when using all three colors.

Let us immediately note that the triple crossing relation, when $e=3$ and the colors are all distinct, is an equality between two 2-morphisms $F_{i_1} F_{i_2} F_{i_3} \1_{\la} \to
F_{i_3} F_{i_2} F_{i_1} \1_{\la}$. Like before, it is an immediate consequence of the unobstructed RIII moves.

In the Section \ref{subsec-polysredux}, we deal with the choice of the polynomials $\fpoly_{\la \to \mu}$ when $e \ge 3$. In the subsequent section, we treat the special case when $e=2$.

\subsection{Polynomials and dots} \label{subsec-polysredux}


\begin{prop}\label{prop:thepolyrelns2} After applying the purported 2-functor $\FC$, the relations of $\UC^-$ and $\dot{\UC}^-$ which involve dots are equivalent to the following properties of the polynomials
$\fpoly_{\la \to \mu}$: invariance \eqref{fis-invariant}, additivity \eqref{fis-additive}, the dual basis condition \eqref{fis-dualbases} and its consequences \eqref{fis-swapped} and
\eqref{fis-tunit}, and one-step locality \eqref{fis-localpiece}, together with the modified versions of the repeated index condition, generic double crossing relation, and redundant dual
basis condition, found below in \eqref{f-properties2}. \end{prop}

This proposition can be compared to Proposition \ref{prop:thesearethepolyrelations}, which was proven by checking the relations of $\UC$ one by one. Note that only those properties of
the polynomials $\fpoly_{\la \to \mu}$ which depended on the adjacent color double crossing relation or the triple crossing relation when $e=2$ will need adjustment. Therefore, the
proofs of several properties can be quoted verbatim, namely: invariance \eqref{fis-invariant}, additivity \eqref{fis-additive}, the dual basis condition \eqref{fis-dualbases} and its
consequences \eqref{fis-swapped} and \eqref{fis-tunit}, and one-step locality \eqref{fis-localpiece}. For each of the remaining conditions, which are stated below, we prove the
correspondence with one of the relations of $\UC$ immediately afterwards.

We begin by stating the modified version of the repeated index condition.

\begin{subequations} \label{f-properties2} Fix an index $1 \le k \le n$ and a color $1 \le i \le e-2$. Suppose that $\la_k = i$ and $\la_{k+1} > i+1$, so that $F_{i+1} \la = 0$. Then, if $e > 2$, one has
\begin{equation} \label{fis-equalbyzero2} \fpoly_{\la \to E_i \la} - \fpoly_{E_i \la \to E_{i+1} E_i \la} + \delta_{i0} y = 0. \end{equation}
If $e = 2$, one has
\begin{equation} \label{fis-equalbyzero3} (\fpoly_{E_i \la \to E_{i+1} E_i \la} - \fpoly_{\la \to E_i \la} + \delta_{i1} y)(\fpoly_{\la \to E_i \la} - \fpoly_{E_i \la \to E_{i+1} E_i \la} + \delta_{i0} y) = 0. \end{equation}

\begin{proof}[Proof of part of Proposition \ref{prop:thepolyrelns2}] When $e> 2$, let $\la$ label the rightmost region of \eqref{newXX1}. Then after applying $\FC$, the left hand side of \eqref{newXX1} is zero. Consequently so is
the right hand side. In the corresponding $S$-diagram the dot polynomials are all in the same region, and consequently one must have \eqref{fis-equalbyzero2}.

When $e = 2$, let $\la$ label the rightmost region of \eqref{KLRmodified} if $i=0$, or \eqref{KLRmodified2} if $i=1$. The same argument yields \eqref{fis-equalbyzero3}.
\end{proof}

Before giving the remaining modified conditions, we state the implications for the $\tgl_n$-realization.

\begin{prop} \label{prop:affineuniquepoly} In the $\tgl_n$-realization, one can define the standard choice $\epoly_{\la \to \mu}$ of polynomials as follows. If $\mu = F_j \la$ for $1 \le
j \le e$, and $F_j$ increments the $k$-th index from $j + re$ to $j + 1 + re$ for some $r \in \ZM$, then we let $\epoly_{\la \to \mu} = x_k + ry$. Then these polynomials satisfy
\eqref{f-properties2}. Moreover, in any extension of the $\tgl_n$ realization, when $e>3$ the only possible choices of $\fpoly_{\la \to \mu}$ have the form $\epoly_{\la \to \mu} + z$ for
some linear polynomial $z \in R^{W_a}$. When $e = 2$ the only possible choices are $\fpoly_{\la \to F_j \la} = \epoly_{\la \to F_j \la} + z$, or $\fpoly_{\la \to F_j \la} = \epoly_{\la \to F_j \la} + z + \delta_{j0} y$. \end{prop}

\begin{proof} Define $\nu_m$ and $z_m$ as in Example \ref{ex:failure}. By Proposition \ref{prop:junique} and Lemma \ref{lem:jchoice}, these polynomials $z_m$ determine all the polynomials $\fpoly_{\la \to \mu}$. We know that $z_{m+e} = z_m + y$.
	
When $e > 2$, the new repeated index condition states that $z_{m+1} = z_m + \delta_{m0} y$, where $\delta_{m0}$ is one if and only if $m$ is zero modulo $e$. Letting $z = z_1$, we see that $\fpoly_{\la \to \mu} = \epoly_{\la \to \mu} + z$.

When $e = 2$, the new repeated index condition is the vanishing of a quadratic which factors into two linear terms, so one of these two linear terms is zero. The consequence is that
$z_{m+1}$ is either equal to $z_m$, or is equal to $z_m + y$. However, since $z_{m+2} = z_m + y$, the values of $z_m$ are determined by $z_1$ and $z_2$. So long as either $z_2 = z_1$ or
$z_2 = z_1 + y$, this will imply the same fact about $z_{m+1}$ compared to $z_m$, for all $m$. Consequently, these are the two possibilities. Letting $z = z_1$ we obtain the desired
result. \end{proof}

\begin{remark} \label{rmk:conventioncare} To reiterate, in the standard choice $\epoly$, incrementing from $e$ to $(e+1)$ gives a polynomial $x_k$, while incrementing from $(e+1)$ to $(e+2)$ gives $x_k + y$, and incrementing from $0$ to $1$ gives $x_k - y$. It is very important to respect this convention! The ``other'' convention, where incrementing from $0$ to $1$ gives $x_k$, and incrementing from $e$ to $e+1$ gives $x_k + y$, is incorrect. One would have to modify the definition of $\UC$ by sprinking signs in the $y$ terms of \eqref{newXX1}, \eqref{funkyR3}, etcetera, in order for $\epoly$ to work with the other convention. \footnote{ To be honest, we were unsure which convention was used by Mackaay-Thiel in \cite[Section 5]{MacThi}, and so we just chose one convention, and apologize if we did not guess right.}
\end{remark}

We now state the modified versions of the remaining properties, when $e>2$. The corresponding conditions for $e=2$ could also be written down (and there is also a condition for the triple crossing relation), but we have not bothered. Instead, when $e=2$, we just check directly in the next section that the polynomials $\fpoly_{\la \to \mu}$, already determined by Proposition \ref{prop:affineuniquepoly}, do satisfy the desired relations.
	
{\bf Generic Double crossing relations:} Fix a color $i$ with $j = i+1$, and a weight $\la$ and set $F_i \la = \mu$, $F_{j} \la = \nu$, $F_{j} \mu = F_i \nu = \rho$. Suppose that none of
these weights are zero. Let $L$ be the ambient parabolic subgroup $L = I(\la) \cap I(\rho) \cap I(\mu) \cap I(\nu)$. Suppose that there are distinct simple reflections $t$ and $u$ such
that $t \in I(\rho) \cap I(\mu) \setminus I(\la)$ and $u \in I(\la) \cap I(\mu) \setminus I(\rho)$. (This is our nondegeneracy condition.)

Then instead of \eqref{fis-padelta} and \eqref{fis-mu} we have:
\begin{equation} \label{fis-padelta2} \pa \D^L_{Ltu} = 1 \ot \fpoly_{\la \to \mu} - \fpoly_{\mu \to \rho} \ot 1 + \delta_{i0} y(1 \ot 1) \in R^{Lt} \ot_{R^L} R^{Lu}, \end{equation} \begin{equation} \label{fis-mu2} \mu^{Lt,Lu}_{Ltu} = \fpoly_{\nu \to \rho} - \fpoly_{\la \to \nu} + \delta_{i0} y \in R^L. \end{equation}
That is, when $i=0$, there is an extra $y$ term appearing.

\begin{proof}[Proof of part of Proposition \ref{prop:thepolyrelns2}] Property \eqref{fis-padelta2} follows from a computation entirely analogous to \eqref{doublecrossing1}. Namely, when $i=0$ and the blue-green strands are pulled apart using \eqref{R2nonoriented1}, \eqref{fis-padelta2} says that an extra term $y(1 \ot 1)$ appears, while \eqref{newXX1} says that this extra $y$ term is required.

Similarly, property \eqref{fis-mu2} follows from a computation analogous to \eqref{doublecrossing2}. When $i=0$ and the blue-green strands are pulled apart using \eqref{R2nonoriented2}, \eqref{fis-mu2} says that there is an extra term $y$, while \eqref{newXX2} says this extra $y$ term is required. \end{proof}

{\bf Redundant dual basis condition:} Continue the setup above, except instead of assuming the existence of $t$ and $u$ as above, assume the existence of $t \in I(\mu)$ which is in neither $I(\rho)$ nor $I(\la)$. Instead of \eqref{fis-dualbasesRedundant} we have:
\begin{equation} \label{fis-padeltaRedundant2} \D^L_{Lt} = 1 \ot \fpoly_{\la \to \mu} - \fpoly_{\mu \to \rho} \ot 1 + \delta_{i0} y(1 \ot 1) \in R^{Lt} \ot_{R^L} R^{Lt},
\end{equation} \begin{equation} \label{fis-muRedundant2} \mu^{L}_{Lt} = \fpoly_{\nu \to \rho} - \fpoly_{\la \to \nu} + \delta_{i0} y \in R^L. \end{equation}

\begin{proof}[Proof of part of Proposition \ref{prop:thepolyrelns2}] As above, this is analogous to the computations \eqref{doublecrossingweird1} and \eqref{doublecrossingweird2}. \end{proof}

\end{subequations}

This concludes the proof of Propostion \ref{prop:thepolyrelns2}, when $e>2$. As the remaining relations of $\UC$ are checked precisely as in the previous chapter, we have proven the following.

\begin{thm} When $e>2$, the 2-functor $\FC$ is well-defined. \end{thm}

\subsection{The case $e=2$} \label{subsec-e2}

We now try to define the 2-functor $\FC$ and prove it is well-defined in the special case $e=2$. This case has a number of idiosyncrasies, because the Dynkin diagram of $\tgl_e$ is no longer simply-laced. We can no longer refer to the colors $0$ and $1$ as being adjacent or distant, so we refer to them as being \emph{funky}.

\subsubsection{The funky crossing}

We begin by defining the image under $\FC$ of the crossing of a $0$-colored strand with a $1$-colored strand. Let us first consider a non-degenerate case. Let $\la$ be a weight in the
middle of its $F_0$- and $F_1$-strings, and let $\mu = F_1 \la$, $\nu = f_0 \la$, and $\rho = F_1 F_0 \la$ as usual. For example, $\la = (111222233)$, $\mu = (112222233)$, $\nu =
(111222333)$, and $\rho = (112222333)$. We let $L$ denote the ambient parabolic, which in this example includes the affine reflection $s_0$, but in general need not. Label the other four
relevant simple reflections by the letters $s,t,u,v$, so that $I(\la) = Lsu$, $I(\mu) = Ltu$, $I(\nu) = Lsv$, and $I(\rho) = Ltv$. If one labels the simple reflections by $s_i$ for $i
\in \ZM/e \ZM$, then there are some $k$ and $l$ for which $s = s_k$, $t = s_{k+1}$, $u = s_l$, and $v = s_{l+1}$ are all distinct, and $L$ contains all the remaining simple reflections.
That is, reading in a circle around the affine Dynkin diagram, our simple reflections are $\{L,L,s,t,L,L,u,v,L\}$.

The $2$-morphism in $\DG$ corresponding to the crossing from $\1_\rho F_0 F_1 \1_{\la} \to \1_\rho F_1 F_0 \1_{\la}$ is below. We color $F_0$ yellow to emphasize the special nature of $e=2$ colors.
\[
{
\labellist
\small\hair 2pt
 \pinlabel {$Lsu$} [ ] at 49 18
 \pinlabel {$Lu$} [ ] at 36 10
 \pinlabel {$Ltu$} [ ] at 29 4
 \pinlabel {$Lt$} [ ] at 18 9
 \pinlabel {$Ltv$} [ ] at 8 19
 \pinlabel {$Ls$} [ ] at 37 27
 \pinlabel {$Lsv$} [ ] at 28 30
 \pinlabel {$Lv$} [ ] at 18 28
 \pinlabel {$L$} [ ] at 28 17
 \pinlabel {$\rho$} [ ] at 93 16
 \pinlabel {$\la$} [ ] at 118 16
 \pinlabel {$\mu$} [ ] at 105 6
 \pinlabel {$\nu$} [ ] at 105 28
\endlabellist
\centering
\ig{2.5}{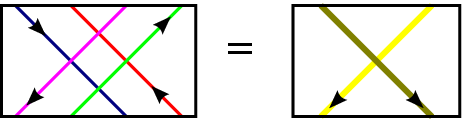}
} \]
While this morphism has a familiar appearance, there is a significant difference between this 2-morphism and the crossings which appeared in the previous chapter. Unlike before, both the blue-green crossing and the red-pink crossing have degree $+1$, because $\{s,v\}$ is no longer disconnected relative to $L$. Thus the overall degree is $+2$. Consequently, there are no ``easy" fringe cases for the ends of the $F_0$- or $F_1$-string, as no strands can be harmlessly removed.

Let $a$ denote the number of even entries of $\la$, and $b$ the number of odd entries. Note that $a+b = n$. Within the interior of the strings, $a>1$ and $b>1$. If $b=1$ but $a>1$ then the picture resembles the difficult edge case from before, see \eqref{FiFjXspecial}.
\begin{equation} \label{F1F0Xspecial} \ig{2}{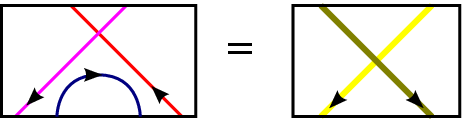} \end{equation}
Unlike \eqref{FiFjXspecial}, both the red-pink crossing and the blue cap in \eqref{F1F0Xspecial} have degree $+1$.  If $a=1$ but $b>1$ then the diagram has a cup on top and a crossing below instead. If either $a=0$ or $b=0$ then the crossing of $F_1$ and $F_0$ is the zero 2-morphism.

The most unique case is when $a=1$ and $b=1$, and $n=2$. Thus for example $\la = (01)$ and $\mu = (12)$, or some $\ZM$-translation of this. There are only two simple reflections in $S_a = \{s,t\}$, and the crossing appears as follows.
\begin{equation} \label{superspecial} {
\labellist
\small\hair 2pt
 \pinlabel {$\emptyset$} [ ] at 47 17
 \pinlabel {$s$} [ ] at 27 29
 \pinlabel {$t$} [ ] at 27 7
\endlabellist
\centering
\ig{1}{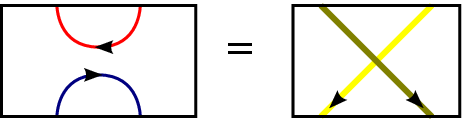}
} \end{equation}
Both the cup and the cap have degree $+1$.

We now check the relations of $\UC$. The one-color relations hold as before.

\subsubsection{Easy relations}

The different-colored dot sliding relation, \cite[(3.18)]{MSV}, follows by \eqref{fis-localpiece} exactly as before.

The typical triple crossing relations, \cite[(3.19)]{MSV}, also follow as before. That is, sometimes both sides of the equality are zero (as must always happen for the special case of
$e=n=2$), and the rest of the time equality follows from a sequence of unobstructed RII and RIII moves.

The commutativity of $F_i$ and $E_j$, \cite[(3.16)]{MSV}, also follows immediately by unobstructed RII moves and isotopy just as in \eqref{EiFjnormal}, and the special case $e=n=2$ is
actually trivial (since rotating the left hand side of \eqref{superspecial} by 90 degrees yields the identity map in $\DG$).

\subsubsection{Funky double crossing relation}

We now compute the double crossing relation. We have already seen how the most degenerate version, where the left hand side of \eqref{KLRmodified} is sent to zero, is equivalent to
\eqref{fis-equalbyzero3}. We now begin with the most generic situation.

\begin{equation} \label{doublecrossinge2}
{
\labellist
\tiny\hair 2pt
 \pinlabel {$\la$} [ ] at 40 102
 \pinlabel {$\mu$} [ ] at 23 94
 \pinlabel {$\nu$} [ ] at 23 121
 \pinlabel {$\rho$} [ ] at 5 101
 \pinlabel {$Lsu$} [ ] at 125 101
 \pinlabel {$Lu$} [ ] at 113 94
 \pinlabel {$Ltu$} [ ] at 104 86
 \pinlabel {$Lt$} [ ] at 93 94
 \pinlabel {$Ltv$} [ ] at 82 102
 \pinlabel {$Ls$} [ ] at 115 110
 \pinlabel {$Lsv$} [ ] at 104 119
 \pinlabel {$Lv$} [ ] at 92 112
 \pinlabel {$L$} [ ] at 103 101
 \pinlabel {$L$} [ ] at 187 103
 \pinlabel {$\mu^{Ls,Lv}_{Lsv}$} [ ] at 189 124
 \pinlabel {$L$} [ ] at 271 106
 \pinlabel {$\mu^{Ls,Lv}_{Lsv}$} [ ] at 271 124
 \pinlabel {$Lu$} [ ] at 283 97
 \pinlabel {$Lt$} [ ] at 257 97
 \pinlabel {$\pa \Delta \mu_{(1)}$} [ ] at 61 33
 \pinlabel {$\pa \Delta \mu_{(2)}$} [ ] at 109 33
 \pinlabel {$1$} [ ] at 40 147
 \pinlabel {$0$} [ ] at 8 147
\endlabellist
\centering
\ig{1.5}{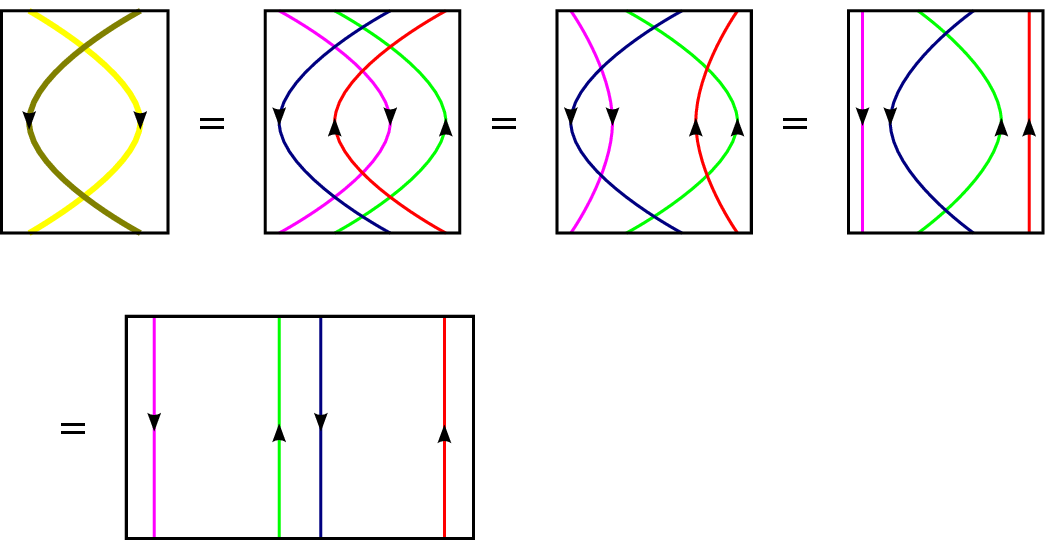}
} \end{equation}
First we apply \eqref{R2nonoriented2}, then an unobstructed RII, and then \eqref{R2nonoriented1var}, letting $\pa \Delta \mu$ be shorthand for the result which is
\begin{equation} \label{eq:nastyresult} \D^{Lt}_{Ltu (1)} \ot \pa^L_{Lu}(\mu^{Ls,Lv}_{Lsv} \D^{Lt}_{Ltu (2)}) \in R^{Lt} \ot_{R^{Ltu}} R^{Lu}. \end{equation}

This complicated-looking expression is actually not so bad. Consider the example above where $\la = (111222233)$. Then $\mu^{Ls, Lv}_{Lsv}$ is the longest root in the connected component
of $Lsv$ containing $s$ and $v$. This is $x_7 - x_3 + y$ in this example. Consequently, $\mu^{Ls, Lv}_{Lsv}$ is a sum of two terms, one in $R^t$ and one in $R^u$, and
each term can be slid out of the middle of the blue-green bigon. Thus, instead of applying \eqref{R2nonoriented1var} as we did in the last step of \eqref{doublecrossinge2} above, we can
empty the bigon and apply \eqref{R2nonoriented1} instead. The empty bigon is equal to $\pa \D^L_{Ltu}$, which is related to the longest root in the connected component of $Ltu$
containing $t$ and $u$. In this example, $\pa \D^L_{Ltu} = 1 \ot x_3 - x_7 \ot 1$. Consequently, in this example \eqref{eq:nastyresult} is equal to
\begin{equation} \label{thispoly} \pa \Delta \mu = (x_7 \ot 1 - 1 \ot x_3 + y(1 \ot 1)) (-x_7 \ot 1 + 1 \ot x_3). \end{equation}

Let $\epoly_{\la \ot \mu}$ be defined just as in Proposition \ref{prop:affineuniquepoly}. Then, in the example $\la = (111222233)$, $\epoly_{\la \to \mu} = x_3$ and $\epoly_{\mu \to
\rho} = x_7$. Consequently, one has \begin{equation} \label{thisexpression} \pa \Delta \mu = (\epoly_{\mu \to \rho} \ot 1 - 1 \ot \epoly_{\la \to \mu} + y(1 \ot 1)) (-\epoly_{\mu \to
\rho} \ot 1 + 1 \ot \epoly_{\la \to \mu}). \end{equation}

If all the indices are shifted up by two, for example $\la = (333444455)$, then the computation \eqref{thispoly} is unchanged. The values of $\epoly$ are changed: $\epoly_{\la \to \mu} = x_3 + y$ and $\epoly_{\mu \to \rho} = x_7 + y$. However, these extra $y$ terms cancel out, and \eqref{thisexpression} still holds. In particular, in this generic situation, \eqref{KLRmodified} holds when the dot polynomials $\fpoly$ are chosen to agree with the standard choices $\epoly$, or with $\epoly + z$ for any $z \in R^{W_a}$ (the extra $z$ terms cancel out as well).

Meanwhile, consider the situation when all indices are shifted up by $1$, so that $\la = (222333344)$, $\mu = (223333344)$, etcetera. The parabolic subgroups involved are identical, but the roles of $F_1$ and $F_0$ are swapped. The resulting polynomial is still \eqref{thispoly} by the identical computation. However, $\epoly_{\la \to \mu} = x_3$ and $\epoly_{\mu \to \rho} = x_7 + y$. The result is now expressed as
\begin{equation} (\epoly_{\mu \to \rho} \ot 1 - 1 \ot \epoly_{\la \to \mu}) (- \epoly_{\mu \to \rho} \ot 1 + 1 \ot \epoly_{\la \to \mu} + y(1 \ot 1)). \end{equation}
Now we see that \eqref{KLRmodified2} holds for the standard choice of dot polynomials, or their shift by $z$.

Now we must check these modified double crossing relations in the edge cases. Consider, for example, the case $\la = (111233)$; this computation is analogous to \eqref{doublecrossingweird1}. Now our
circle around the affine Dynkin diagram is $\{L,L,s,t,v,L\}$, and $I(\la) = Ls$, $I(\mu) = Lt$, $I(\nu) = Lsv$, and $I(\rho) = Lv$.
\begin{equation} {
\labellist
\small\hair 2pt
 \pinlabel {$Ls$} [ ] at 121 34
 \pinlabel {$Lv$} [ ] at 85 34
 \pinlabel {$Lt$} [ ] at 103 7
 \pinlabel {$\mu^{Ls, Lv}_{Lsv}$} [ ] at 188 33
\endlabellist
\centering
\ig{1}{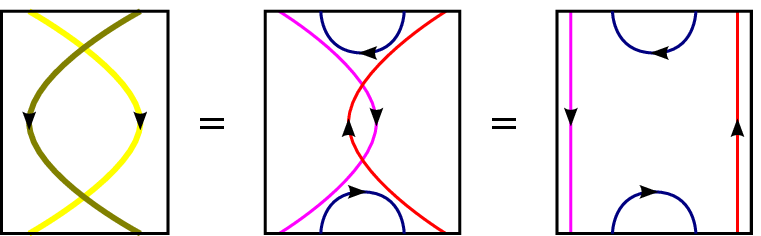}
} \end{equation}
We applied \eqref{R2nonoriented2} then an unobstructed RII move. Our next step (not pictured) is to apply \eqref{coproduct} to the two blue caps. The overall result is
\begin{equation} \label{nasty2} \mu^{Ls,Lv}_{Lsv} \Delta^L_{Lt}. \end{equation}
Again, this nasty-looking formula is not bad in any example. When $\la = (111233)$ we have $\mu^{Ls,Lv}_{Lsv} = x_4 - x_3 + y$, while $\Delta^L_{Lt} = -x_4 \ot 1 + 1 \ot x_3$. Once again, one can confirm that \eqref{KLRmodified} holds, in any example.

The $e=2$ computation which is analogous to \eqref{doublecrossingweird2} follows by exactly the same arguments.

The $e=n=2$ special case looks as follows. Let $\la = (12)$.
\begin{equation} {
\labellist
\small\hair 2pt
 \pinlabel {$\emptyset$} [ ] at 124 20
 \pinlabel {$t$} [ ] at 104 8
 \pinlabel {$s$} [ ] at 103 33
 \pinlabel {$\mu_s$} [ ] at 187 35
\endlabellist
\centering
\ig{1}{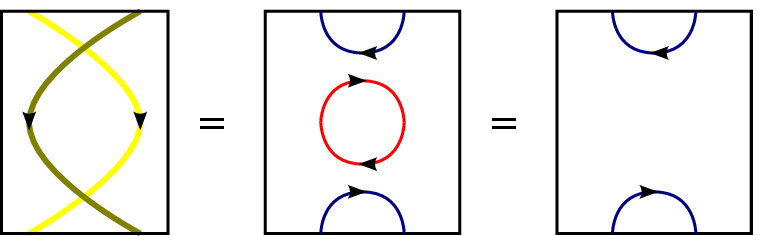}
} \end{equation}
Again, the final result is $\mu_s \Delta_t$, which agrees with \eqref{KLRmodified} when the dot polynomials are standard, by analogous arguments.

\subsubsection{Interesting triple crossing relation}

Now we must check the interesting triple crossing relation. Reading around our affine Dynkin diagram, label the simple reflections as $\{L, L, s, t, u, L, L, x, y, L \}$.  The computation is analogous to \eqref{hardreln2}.
\begin{equation} \label{hardrelne2} {
\labellist
\small\hair 2pt
 \pinlabel {$stx$} [ ] at 178 224
 \pinlabel {$tx$} [ ] at 161 242
 \pinlabel {$tux$} [ ] at 144 242
 \pinlabel {$tu$} [ ] at 127 242
 \pinlabel {$tuy$} [ ] at 116 219
 \pinlabel {$sx$} [ ] at 171 200
 \pinlabel {$sux$} [ ] at 158 192
 \pinlabel {$ux$} [ ] at 150 213
 \pinlabel {$su$} [ ] at 145 185
 \pinlabel {$u$} [ ] at 136 204
 \pinlabel {$uy$} [ ] at 124 192
 \pinlabel {$sx$} [ ] at 171 292
 \pinlabel {$sux$} [ ] at 158 303
 \pinlabel {$ux$} [ ] at 150 284
 \pinlabel {$su$} [ ] at 145 313
 \pinlabel {$u$} [ ] at 136 292
 \pinlabel {$uy$} [ ] at 123 306
 \pinlabel {$stx$} [ ] at 400 224
 \pinlabel {$tx$} [ ] at 385 237
 \pinlabel {$t$} [ ] at 368 245
 \pinlabel {$tu$} [ ] at 348 242
 \pinlabel {$tuy$} [ ] at 336 219
 \pinlabel {$stx$} [ ] at 289 228
 \pinlabel {$x$} [ ] at 268 239
 \pinlabel {$ux$} [ ] at 254 232
 \pinlabel {$u$} [ ] at 241 221
 \pinlabel {$tuy$} [ ] at 224 219
 \pinlabel {$uy$} [ ] at 236 27
 \pinlabel {$su$} [ ] at 261 27
 \pinlabel {$u$} [ ] at 248 72
 \pinlabel {\tiny $suy$} [ ] at 249 6
\endlabellist
\centering
\ig{1.1}{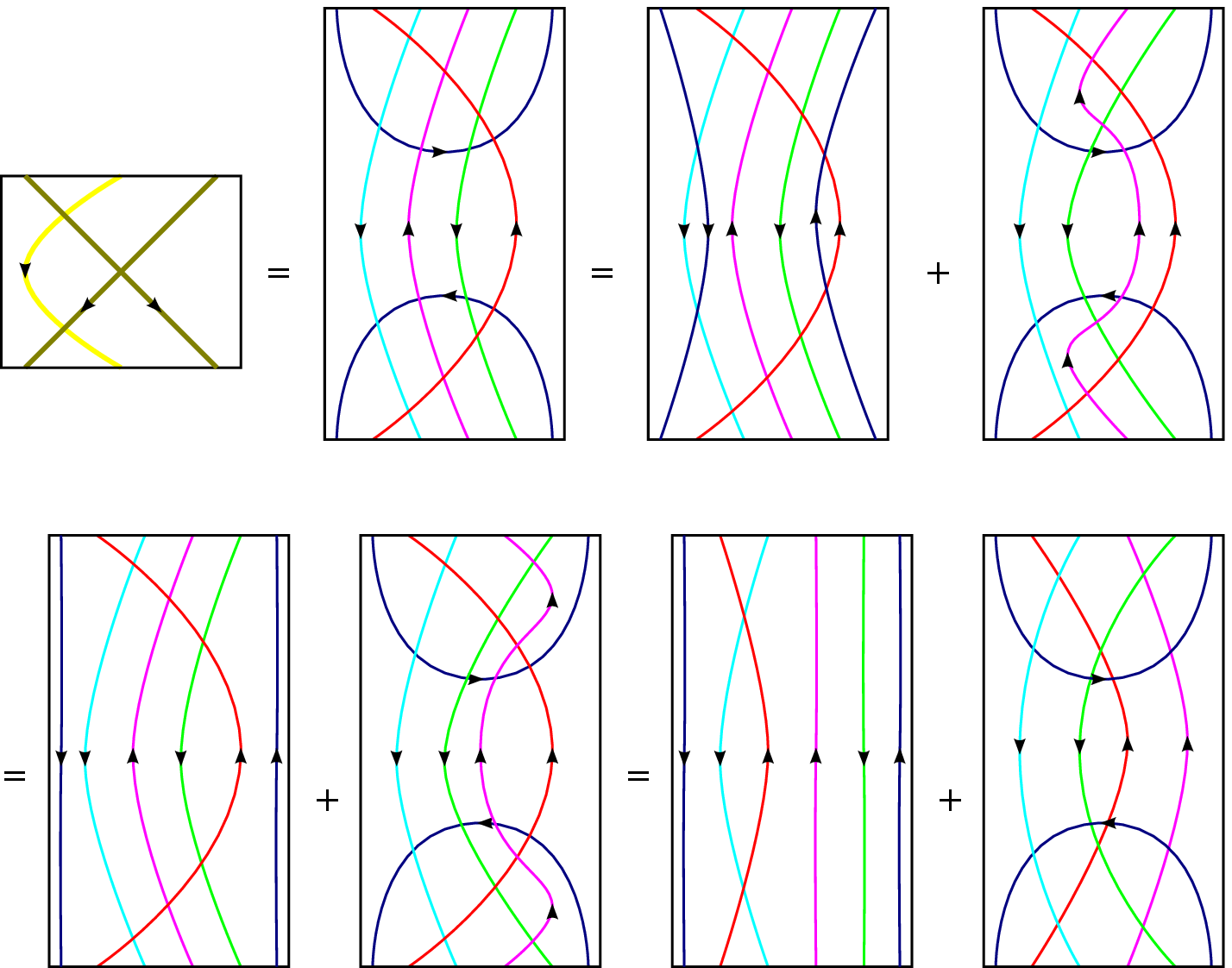}
} \end{equation}
To begin, we apply the interesting relation \eqref{relation11n1} inside the large red-aqua bigon, where $s$ and $y$ are absent. The remaining simplifications use only unobstructed RII and RIII moves.

Unlike the analogous computation in \eqref{hardreln2}, the penultimate diagram has a red-aqua RII move which is obstructed. For this obstructed RII move, one should think that $u$ is in the ambient parabolic, and the interesting stuff is happening on the finite type $A$ subdiagram of the form $\{y, L, L, L, s\}$. Applying \eqref{R2nonoriented1} one obtains $\pa \Delta_{Ksy}^{K}$, where $K = Lu$, which by \eqref{fis-padelta} or \eqref{fis-padelta2} can be expressed in terms of the dots and the polynomial $y$.

For example, suppose the rightmost region $\la$ is $(2223333444)$. Let $\fpoly_{1,R}$ denote a dot on the right strand labeled $1$, $\fpoly_{1,L}$ be a dot on the left strand labeled $1$, and $\fpoly_0$ be a dot on the strand labeled $0$. Then $\fpoly_{1,R} = x_7 + y$, $\fpoly_0 = x_3$, and $\fpoly_{1,L} = x_6 + y$. Moreover, $s = s_5$, $t = s_6$, $u = s_7$, $x = s_2$, and $y = s_3$. One has $\pa \Delta_{Ksy}^{K} = -x_6 \ot 1 + 1 \ot x_3 = -\fpoly_{1,L} + \fpoly_0 + y$.

By an exactly analogous computation one has
\begin{equation} \label{hardrelne2part2} {
\labellist
\small\hair 2pt
 \pinlabel {$sx$} [ ] at 192 26
 \pinlabel {$su$} [ ] at 165 27
 \pinlabel {$s$} [ ] at 180 70
 \pinlabel {\tiny $sux$} [ ] at 180 6
\endlabellist
\centering
\ig{1}{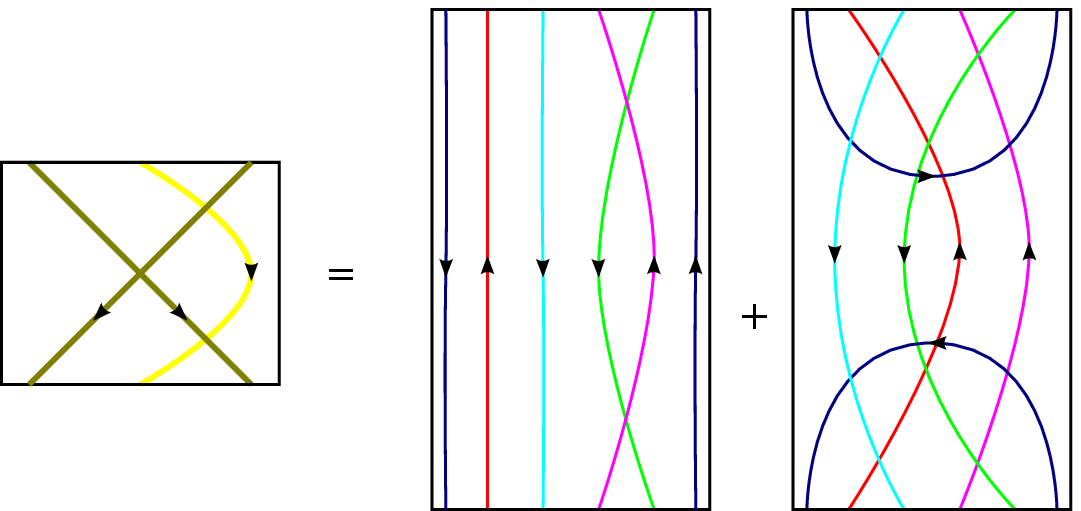}
} \end{equation}
Again, we apply an obstructed RII move to obtain $\pa \Delta_{K'ux}^{K'}$, where $K' = Ls$, which by \eqref{fis-padelta} or \eqref{fis-padelta2} can be expressed in terms of the dots and the polynomial $y$.

In the example of $\la = (2223333444)$, one has $\pa \Delta_{K'ux}^{K'} = - x_3 \ot 1 + 1 \ot x_7 + y(1 \ot 1) = - \fpoly_0 + \fpoly_{1,R}$.

In the end, one obtains
\begin{equation*} {
\labellist
\small\hair 2pt
 \pinlabel {$y$} [ ] at 417 34
\endlabellist
\centering
\ig{1}{funkyR3}
} \end{equation*}
The second terms of \eqref{hardrelne2} and \eqref{hardrelne2part2} cancel out. The first term of \eqref{hardrelne2} accounts for middle dot minus left dot, possibly with a $y$ term (as in the example above). The first term of \eqref{hardrelne2part2} accounts for middle dot minus right dot, possibly with the $y$ term (not in the example above), and has a $y$ term if and only if \eqref{hardrelne2} does not.

For example, when $\la = (1112222333)$, then $s = s_1$, $t = s_2$, $u = s_3$, $x = s_6$ and $y = s_7$. We have $\fpoly_{1,R} = x_3$, $\fpoly_0 = x_7$, and $\fpoly_{1,L} = x_2$. We have $\pa \Delta_{Ksy}^K = -x_2 \ot 1 + 1 \ot x_7 + y(1 \ot 1) = -\fpoly_{1,L} + \fpoly_0 + y$, and $\pa \Delta_{K'ux}^{K'} = -x_7 \ot 1 + 1 \ot x_3 = -\fpoly_0 + \fpoly_{1,R}$. Again, the difference between these two things is $2\fpoly_0 - \fpoly_{1,R} - \fpoly_{1,L} + y$.

A similar computation with colors reversed gives
\begin{equation*} {
\labellist
\small\hair 2pt
 \pinlabel {$y$} [ ] at 417 34
\endlabellist
\centering
\ig{1}{funkyR3alt}
} \end{equation*}

For example, when $\la = (2223333444)$, one has $s = s_1$, $t = s_2$, $u = s_3$, $x = s_6$, and $y = s_7$. One has $\fpoly_{0,R} = x_3$, $\fpoly_1 = x_7 + y$, and $\fpoly_{0,L} = x_2$. One has $\pa \Delta_{Ksy}^K = -x_2 \ot 1 + 1 \ot x_7 + y(1 \ot 1) = -\fpoly_{0,L} + \fpoly_1$, and $\pa \Delta_{K'ux}^{K'} = -x_7 \ot 1 + 1 \ot x_3 = -\fpoly_1 + \fpoly_{0,R} + y$. The difference between these two things is $2\fpoly_1 - \fpoly_{0,R} - \fpoly_{0,L} - y$.

Having checked all the relations, we have proven the following. 

\begin{thm} When $e=2$, the 2-functor $\FC$ is well-defined. \end{thm} 

%% file: Cellular2.tex
\section{Cellular Structures and Soergel bimodules} \label{sec-cell}
%
%

Recall from the introduction that our eventual goal is to prove a rigidity result for categorifications of Fock space representations. This rigidity result will rely on techniques
involving highest weight categories. In this chapter, we discuss the cellular structure on $\SB$, and the connections between cellular structures and highest weight categories.
We also prove some technical lemmas about the cellular structure on $\SB$.

%
%

%

\subsection{Cellular algebras and highest weight categories}

It is one of the key results of Graham and Lehrer \cite{GraLeh} that for cellular algebras with certain properties, their module categories are highest weight categories. Our first goal
is to justify the analogous results for object-adapted 
cellular categories like $\SB$. These results are mostly in the literature, although not stated in the precise form we use.

Let us first recall the results of Graham and Lehrer.

A \emph{cellular algebra} $C$ over a commutative base ring $A$ is a filtered algebra with a special basis and an antiinvolution. Attached to this algebra is a poset of \emph{cells}
$\PC$. For each $\la \in \PC$ there is a finite set $M(\la)$, and the algebra has an $A$-basis $\{c^\la_{S,T}\}$ attached to a cell $\la$ and a pair $S,T \in M(\la)$. It is typically
assumed that $\PC$ is finite, in which case the cellular algebra is also finite rank over $A$. This data satisfies several conditions: \begin{itemize} \item The $A$-linear map $\iota$
which sends $c^\la_{S,T}$ to $c^\la_{T,S}$ for each $\la \in \PC$ and $S, T \in M(\la)$ is an anti-involution of $C$. \item The span of $c^\mu_{S,T}$ for $\mu \le \la$ forms an ideal
$J_{\le \la} \subset C$. Similarly the span for $\mu < \la$ is an ideal $J_{< \la}$, which one kills when working \emph{modulo lower cells} with respect to $\la$. \item Finally,
multiplication by $c^\la_{S,T}$ satisfies a particular ``cellular formula" in the quotient by lower cells. More precisely, if $x \in C$ then \begin{equation} \label{eq:cellformula} x
c^\la_{S,T} = \sum_{S' \in M(\la)} \ell(x,S,S') c^\la_{S',T} \textrm{ modulo lower cells.} \end{equation} Here, $\ell$ is some function which does not depend on the choice of $T$.
Applying $\iota$, a similar formula holds for right multiplication. \end{itemize}
	
The function $\ell$ from \eqref{eq:cellformula} allows one to define the \emph{cell module} $V(\la)$ for $\la \in \PC$ as follows. It has a basis $v_S$ for $S \in M(\la)$, and the action
of $A$ is given by \begin{equation} x v_S = \sum_{S' \in M(\la)} \ell(x,S,S') v_{S'}. \end{equation} In particular, it is obvious that the $C$-module $J_{\le \la} / J_{< \la}$, a
quotient of ideals in $C$, is isomorphic to a direct sum of copies of $V(\la)$ for each $T \in M(\la)$, via the map sending $v_S \mapsto c^\la_{S,T}$.

For each cell $\la \in \PC$ one has a \emph{cellular pairing}, a map $\phi^\la \co M(\la) \times M(\la) \to A$ determined from \begin{equation} \label{eq:cellmult} c^\la_{S,T} c^\la_{U,V} =
\phi^\la(T,U) c^\la_{S,V} \textrm{ modulo } J_{< \la}.\end{equation} This extends by $A$-linearity to a form on $V(\la)$. Many of the properties of a cellular algebra can be understood using
the cellular pairings, or equivalently (and often in the literature) by studying the squared ideal $J_{\le \la}^2$. Here are three properties a cellular pairing might have.
\begin{enumerate} \item The pairing on $V(\la)$ is non-degenerate. For example, when $\phi^\la(T,U) = \delta_{T,U}$ is the Kronecker delta (or, more generally, is non-degenerate), then
$J_{\le \la} / J_{< \la}$ is a matrix algebra over $A$. By non-degenerate, we mean that every $v \in V(\la)$ pairs against some $v' \in V(\la)$ to be a unit in $A$. \item The pairing is
non-zero. Equivalently by \eqref{eq:cellmult}, $J_{\le \la}^2$ is not contained in $J_{< \la}$. \item In between these, the pairing is \emph{unit-worthy}, that is, some $v \in V(\la)$ pairs against some $v' \in V(\la)$ to be a unit in $A$. Equivalently by \eqref{eq:cellmult}, $J_{\le \la}^2$ is equal to $J_{\le \la}$ modulo $J_{< \la}$. In this case, note that if
$J_{< \la}^2 = J_{< \la}$ then also $J_{\le \la}^2 = J_{\le \la}$. \end{enumerate} When $A$ is a field, being non-zero is the same as being unit worthy.

Suppose that $A = \Bbbk$ is a field. Graham and Lehrer prove that when the cellular pairing $\phi^\la$ is nonzero, then the radical of the pairing on $V(\la)$ is a maximal proper
submodule, and the quotient by this radical is an absolutely irreducible module $L(\la)$. They prove that these simple modules $L(\la)$, ranging over the cells for which the pairing is
nonzero, enumerate all the simple modules of $C$. Furthermore, they prove the following proposition.

\begin{prop} \label{prop:whencqh} (See \cite[Theorem 3.8 and Remark 3.10]{GraLeh}) Suppose that $A = \Bbbk$ is a field and $\PC$ is finite. Then the cellular algebra $C$ is
quasi-hereditary if and only if each cellular pairing $\phi^\la$ is nonzero (equivalently, unit-worthy). It is semisimple if and only if each cellular pairing is non-degenerate.
\end{prop}


A finite-dimensional algebra is quasi-hereditary if it admits a filtration by ideals where each descends to a heredity ideal in the quotient by lower terms. The details of heredity
ideals are not relevant to us. It is important that an algebra is quasi-hereditary if and only if its module category is a highest weight category. The paper \cite{KonXiWhenis} discusses precisely when a cellular algebra is quasi-hereditary, and \cite[Lemma 2.1(3)]{KonXiWhenis}
contains the proposition above, stated in terms of $J_{\le \la}^2$.

The axioms of a highest weight category will be reviewed in Section \ref{SS_HW}, but one major theme is the definition of standard modules $\Delta$. In the context of cellular algebras
over a field, the standard module $\Delta(\la)$ is the same as the cell module $V(\la)$.

\subsection{Generalizations}\label{SS_cell_hw}

We now proceed to discussing some generalizations of cellular algebras and Proposition \ref{prop:whencqh}.

First, suppose that $C$ is a cellular algebra over a base ring $A$, but that $A$ is not a field. For example, suppose that $A$ is an infinite-dimensional commutative $\Bbbk$-algebra for
some field $\Bbbk$. The definition of quasi-hereditary algebras implies that they are finite dimensional, but $C$ is not finite dimensional. Thus $C$ is not quasi-hereditary, and its
module category is not highest weight, using the original definitions of these terms. One could attempt to study what happens after base change from $A$ to various fields, although this
is an awkward approach (improved by Kleshchev, see \S\ref{SSS_cell_affine}). Thankfully, the literature contains generalizations of the concepts of quasi-hereditary rings and highest weight
categories which apply when $A$ is not a field.

Secondly, we are interested in categories, not algebras. There is a notion of a cellular category, due to Westbury \cite{WestCC}. We will need a generalization of Proposition
\ref{prop:whencqh} to this setting. It becomes easier to prove when one restricts to a special kind of cellular category, known as an object-adapted cellular category, defined in
\cite{ELauda}.

Ultimately, we are interested in studying $\SB$ and related categories; these are (graded) object-adapted cellular categories over the base ring $R$, which is a polynomial ring over a field. Moreover, we are interested in graded modules over this category. This explains why we need these generalizations.

\subsubsection{Affine cellular algebras}\label{SSS_cell_affine}


We are interested in cellular algebras over a base ring $A$ which is not a field. The results we plan to quote from the literature, however, are stated in a much broader level of generality: that of affine cellular algebras. This is why we introduce affine cellular algebras, which are not otherwise necessary in this paper. \footnote{Rouquier \cite{rouqqsch} has definitions and results which apply to the case when $A$ is a complete ring. For example, he has a definition of a highest weight category over a complete ring.}

Let $\Bbbk$ be a field. An \emph{affine cellular algebra}, as introduced by Koenig and Xi \cite{KonXi12}, is analogous to a cellular algebra defined over a different $\Bbbk$-algebra
$A_{\la}$ for each cell $\la$. We do not give a definition here, see \cite[Definition 2.1]{KonXi12}. The word \emph{affine} indicates that these rings $A_{\la}$ are assumed to be
quotients of finitely generated polynomial rings by a homogeneous ideal. Although affine cellular algebras are not usually cellular algebras, they possess many of the same structural
properties. The analog of Proposition \ref{prop:whencqh} above is \cite[Theorem 4.1(1)]{KonXi12}, which declares that the condition $J_{\le \la}^2 = J_{\le \la}$ is equivalent to a
statement about heredity ideals in quotients by the maximal ideals of $A_{\la}$.

\begin{remark} A cellular algebra over a base ring $A$ can be viewed as an affine cellular algebra where each of the rings $A_\la$ is equal to $A$. The ability to change the ring for each cell is not relevant for us, though see Remark \ref{rmk:fibredDG}. \end{remark}

A streamlined approach to studying module categories over affine cellular algebras was provided by Kleshchev \cite{Klesh}, who defined the analogous notion of an \emph{affine
quasi-hereditary algebra} and an \emph{affine highest weight category}. He proves in \cite[Theorem C]{Klesh} that a module category is affine highest weight if and only if the algebra is
affine quasi-hereditary. Kleshchev is also careful to handle graded rings, defining the analog of an affine cellular algebra for graded $\Bbbk$-algebras where the grading is bounded
below (see \cite[Definition 9.3]{Klesh}). The analog of Proposition \ref{prop:whencqh} is the following.

\begin{prop} \label{prop:whenacqh} A (graded) affine cellular algebra $C$ with finite\footnote{This finiteness is assumed in Kleshchev's definition of an affine cellular algebra.} cell
poset $\PC$ is affine quasi-hereditary if each cellular pairing $\phi^\la$ is unit-worthy. \end{prop}

One expects this proposition to be an if and only if statement. Unfortunately, neither direction is stated directly in \cite{Klesh}, although Chapter 9 discusses ideas closely related to
the converse. The proof below merely hacks together Kleshchev's results.

\begin{proof} We seek to show that $J_{\le \la}$ descends to an affine heredity ideal modulo $J_{< \la}$, as defined in \cite[Definition 6.2]{Klesh}. This requires two conditions, called
(SI1) and (SI2) in \cite[Definition 6.1]{Klesh}, where (SI2) is modified in \cite[Definition 6.2]{Klesh}. At this point in Kleshchev's paper, one is working with arbitrary algebras, not
with affine cellular algebras, and these conditions are simplified somewhat in an affine cellular algebra. For example, the modified condition (SI2) holds automatically in an affine
cellular algebra, following from condition (Stand) in \cite[Definition 9.3]{Klesh}. Meanwhile, by \cite[Lemma 6.5]{Klesh}, (SI1) is equivalent to $J_{\le \la}^2 = J_{\le \la}$, given
another condition (that $J_{\le \la}/J_{< \la}$ is projective over $A_\la$) which holds automatically in an affine cellular category. By induction, one has $J_{\le \la}^2 = J_{\le \la}$
for all $\la$ if and only if $\phi^\la$ is unit-worthy for all $\la$. \end{proof}

\begin{remark} Note also that an affine highest weight category over a base ring $A$, after base change from $A$ to a field, becomes a highest weight category in the usual sense. After base change to a complete ring, it becomes a highest weight category in the sense of Rouquier \cite{rouqqsch}. \end{remark}

\subsubsection{Cellular categories}

Now we segue to discussing categories instead of algebras. It is a well-known philosophy that (certain) categories can be thought of as locally-unital algebras, or \emph{algebroids}.
Given a small category $\CC$ whose morphism spaces are abelian groups, with composition being bilinear (e.g. a full subcategory of an additive category), one can consider $C = C_{\CC} =
\oplus_{X,Y \in \CC} \Hom(X,Y)$, where we have abusively written $X \in \CC$ to say that $X$ is an object of $\CC$. The abelian group $C$ inherits a multiplication structure, with $f
\cdot g$ being equal to $f \circ g$ if the maps are composable, and $0$ otherwise. If there are finitely many objects, then $1 = \sum_X \id_X$ is a unit in $C$, but otherwise there is no
unit. However, the infinite collection $\{\id_X\}_{X \in \CC}$ of orthogonal idempotents will serve as a replacement for a unit, since the infinite sum $\sum_X \id_X$ acts locally
finitely on any element of $C$, and acts by the identity.

By a module over a category, we mean a locally-unital module over the locally-unital algebra $C_{\CC}$. This is the same data as the following definition.

\begin{defn} Let $\CC$ be a $\Bbbk$-linear category. Then $\CC\Mod$ is the category of $\Bbbk$-linear functors from $\CC$ to $\Bbbk\Mod$. This is an abelian category. \end{defn}
	
In particular, all representable functors $\Hom(M,-)$ are projective as $C_{\CC}$-modules, being the image of the idempotent $\id_M$.

A \emph{cellular category} over a commutative base ring $A$ is the analogue of a cellular algebra for categories or locally-unital algebras, introduced by Westbury \cite{WestCC}. The
data consists of \begin{itemize} \item an $A$-linear category $\CC$, \item a poset $\PC$, \item a finite set $M(X,\la)$ for each object $X$ and each cell $\la$, \item and an $A$-basis
for the space $\Hom(X,Y)$ given by $\{c^\la_{S,T}\}$ for elements $S \in M(Y,\la)$ and $T \in M(X,\la)$. \end{itemize} This data satisfies the following conditions. \begin{itemize} \item There is a contravariant isomorphism of categories $\iota \co \CC
\to \CC$ sending $c_{S,T}^\la$ to $c_{T,S}^\la$. \item The span of various $c^\mu_{S,T}$ for all $\mu \le \la$ is an ideal $J_{\le \la}$ in the category. \item A cellular formula holds. Namely, if $x \in \Hom(Y,Z)$ and $S \in M(Y,\la)$, $T \in M(X,\la)$, then \begin{equation} \label{eq:cellformula2} x
c^\la_{S,T} = \sum_{S' \in M(Z,\la)} \ell(x,S,S') c^\la_{S',T} \textrm{ modulo lower cells.} \end{equation} Here, $\ell$ is some function which does not depend on the choice of $T$. \end{itemize} 

There is still a \emph{cellular pairing} $\phi^\la_X \co M(X,\la) \times M(X,\la) \to A$, for which \[c^\la_{S,T} \circ c^\la_{U,V} = \phi^\la_Y(T,U) c^\la_{S,V} \textrm{ modulo } J_{<
\la} \] holds in $\Hom(X,Z)$ for any $V \in M(X,\la)$, $S \in M(Z,\la)$, and $T,U \in M(Y,\la)$. If $\CC$ has finitely many objects, the algebra $C_\CC$ is a cellular algebra, with the
same poset $\PC$, and with $M(\la) = \coprod_X M(X,\la)$. Moreover, the cellular pairing in $C_\CC$ is automatically zero between any $S \in M(X,\la)$ and any $T \in M(Y,\la)$ for $X \ne
Y$, so that it is the product of the cellular pairings in $\CC$ on each $M(X,\la)$.

Thus, if there are finitely many objects, one can deduce from Proposition \ref{prop:whencqh} that a cellular category over a field is semisimple if for each $\la$ the pairing on each set
$M(X,\la)$ is non-degenerate, but it is quasi-hereditary so long as for each $\la$ the pairing on some $M(X,\la)$ is nonzero. The same conclusion holds even when the cellular category
has infinitely many objects, as we state in the proposition below.

\begin{prop} \label{prop:whenccqh} The category of modules for a cellular category over a field (with finite poset $\PC$) is highest weight, or equivalently, the cellular category is a
quasi-hereditary algebroid, if and only if for each $\la$ there is an object $X \in \CC$ such that the cellular pairing $\phi^\la_X$ is unit-worthy. \end{prop}

This proposition is the direct analog of Proposition \ref{prop:whencqh} for algebroids. Unfortunately, the representation theory of locally-unital rings is not nearly as well represented
in the literature; almost every result that holds for unital algebras has an analog for algebroids, and because the proof is almost identical, no one bothers to write it down. We will
not bother to prove this proposition either (or to define quasi-hereditary algebroids), but the proof should be an imitation of the proof of Proposition \ref{prop:whencqh}. See also
Remark \ref{rem:itsokay} below. It is worth mentioning that while there may be infinitely many objects (idempotents) in the algebroid, there are finitely many cells.

\subsubsection{Object-adapted cellular categories}

The notion of a \emph{(strictly) object-adapted cellular category} was introduced in \cite[Definition 2.4]{ELauda}, as a cellular category (over a base ring $A$) where the cells are associated with objects,
and the cellular basis element $c^\la_{S,T} \co X \to Y$ is a composition of maps $c_T \co X \to \la$ and $c_S \co \la \to X$ which go to and from the object $\la$ respectively. This
factorization property implies the mysterious ``cellular formula," and underlies most (but not all) of the cellular structures in the literature. These maps $c_S$ and $c_T$ can be
thought of as two ``halves'' of the cellular basis, although they are also members of the cellular basis of $\Hom(X,\la)$ and $\Hom(\la,Y)$.

To help distinguish between maps to $\la$ and maps from $\la$, we use slightly different notation. Instead of a set $M(X,\la)$, one has two sets $M(\la,X)$ and $E(X,\la)$, in bijection
via a map $\iota$. Then given $T \in E(X,\la)$ one has a map $c_T \in \Hom(X,\la)$, and given $S \in M(\la,X)$ one has a map $c_S \in \Hom(\la,X)$. These satisfy $c_{\iota(S)} =
\iota(c_S)$.

One of the axioms of an object-adapted cellular category is that $M(\la,\la) = E(\la,\la) = \{*\}$, and $c_{*}$ is the identity of $\la$. A consequence is that the identity of $\la$ spans $\End(\la)$ modulo lower terms. Another consequence is that the ideal $J_{< \la}$ consists precisely of those morphisms which (are linear combinations of morphisms which) factor through the objects $\mu$ for $\mu < \la$.

The cellular pairing can be calculated as follows: compose $c_T \circ c_S$ to obtain an endomorphism of the object $\la$, and then compute the coefficient of the identity (modulo lower
terms); this is $\phi(T,S)$. Because the identity of $\la$ always pairs with itself to be $1$, object-adapted cellular categories are always unit-worthy! Thus, their module categories
are always affine highest weight.

\begin{remark} The condition that $J_{\le \la}^2 = J_{\le \la}$ is often equivalent to $J_{\le \la}$ being generated as a two-sided ideal by some idempotent $e$. In this case, the
idempotent is $\sum \id_{\mu}$ where the sum is over all the objects $\mu$ with $\mu \le \la$. \end{remark}

\begin{defn} Given a cell $\la \in \PC$, the \emph{standard module} or \emph{cell module} $\Delta(\la)$ is the functor $\CC \to \Bbbk\Mod$ which sends an object $X$ to $\Hom(\la,X) / J_{< \la}$. \end{defn}

The standard module is a quotient of a representable functor by the image of a direct sum of representable functors, because the ideal $J_{< \la}$ is the sum, over all $\mu < \la$ and
all maps $\la \to \mu$, of the corresponding image of $\Hom(\mu,X)$ inside $\Hom(\la,X)$.

\begin{remark} \label{rem:itsokay} An object-adapted cellular category is Morita equivalent to the subcategory whose objects are the cells; that is, these two algebroids have equivalent
module categories. This easy fact is proven in Corollary \ref{lem:morita} below. Thus, if there are finitely many cells, one may as well assume there are finitely many objects for the
purposes of studying the module category. This gives a proof of Proposition \ref{prop:whenccqh} above for object-adapted cellular categories with finitely many cells, without needing any
technical results about algebroids. \end{remark}


\begin{remark} \label{rmk:fibredDG} Also introduced in \cite{ELauda} is the notion of a \emph{fibered cellular category}, which combines the idea of an affine cellular category with an
object-adapted cellular category. It is like an object-adapted cellular category but with a different base ring $A_{\la}$ for each cell $\la$. We mention this because $\DG$ is
conjectured to be a fibered cellular category (and the ability to choose different base rings for each object is essential), a result we expect will be proven in \cite{EWSingular}.
Because this result has not yet appeared, we have chosen to stick with $\SB$ and its known cellular structure. \end{remark}


%


\subsection{Cellular subquotients} \label{subsec-cellsubquot}

One limitation seen in all the propositions above is that the cellular poset is assumed to be finite. For the diagrammatic category of Soergel bimodules $\SB$, the poset is the Coxeter
group $W$ with its Bruhat order, which is not finite when $W$ is not finite. Section \ref{SS_HW} will deal with many technicalities involving infinite posets. In this section, we discuss a few useful operations involving object-adapted cellular categories, namely quotients for coideals and subcategories for ideals, which can help to reduce from the infinite to the finite setting.

Let $\CC$ be an object-adapted 
cellular category, with poset $\PC$. Choose $\la \in \PC$. Let $\PC_{\le \la} = \{ \mu \in \PC \; | \; \mu \le \la\}$ be the ideal
generated by $\la$ in the poset $\PC$, and $\PC_{\ge \la}$ be the similarly-defined coideal. Similarly, one has an ideal $\PC_{\ngeq \la}$, etcetera.

Since $J_{\ngeq \la}$ is an ideal in $\CC$, one can consider the quotient category $\CC^{\ge \la} = \CC / J_{\ngeq \la}$. This will also be an object-adapted  cellular category, with
poset $\PC_{\ge \la}$. Note that any object $\mu \in \CC$ corresponding to $\mu \in \PC_{\ngeq \la}$ becomes isomorphic to the zero object in this quotient.

Now we recall the Morita equivalence mentioned in Remark \ref{rem:itsokay}, which is mostly tautological.

\begin{lemma}\label{lem:morita} Let $\Bbbk$ be a field. Let $\CC$ be a $\Bbbk$-linear category, and $\PC$ be a set of objects for which the ideal generated by $\{\id_P\}$, $P \in \PC$, is the entire category
$\CC$. That is, up to linear combinations, every morphism factors through an object in $\PC$. Then $\CC$ is Morita-equivalent to the full subcategory with objects $\PC$. That is, $\CC$
and this full subcategory have equivalent module categories. \end{lemma}

\begin{proof} Any $\Bbbk$-linear category includes fully faithfully into its Karoubi envelope, and is Morita-equivalent to its Karoubi envelope, so we may as well assume that $\CC$ is
Karoubian. Yoneda-think states that any Karoubian $\Bbbk$-linear category is equivalent to the category of projective modules over itself, by identifying an object $M \in \CC$ with the
representable functor $\Hom(M,-)$. If every morphism factors through an object in $\PC$ (up to linear combinations) then the representable functor $\Hom(\oplus_{P \in \PC} P,-)$ is a
projective generator, and its endomorphism ring is precisely the full subcategory with objects in $\PC$. \end{proof}

\begin{cor} \label{cor:morita} An object-adapted 
cellular category is Morita equivalent to the full subcategory with objects $\PC$. \end{cor}

Now we claim that the full subcategory $\CC_{\le \la}$ with objects in $\PC_{\le \la}$ is an object-adapted 
cellular category with poset $\PC_{\le \la}$, which is easy to confirm.
It was shown in \cite{ELauda} that this full subcategory is entirely contained in the ideal $J_{\le \la}$. Although the ideal $J_{\le \la}$ contains other morphisms as
well, these factor through the objects in $\PC_{\le \la}$ (up to linear combinations). Philosophically, the module category over $\CC_{\le \la}$ is the Serre subcategory of the module category over $\CC$ corresponding to the ideal $\PC_{\le \la}$.

Finally, one can combine these operations to consider the subquotient category $\CC_{[\mu,\la]} = \CC_{\le \la} / J_{\ngeq \mu}$, which is an object-adapted  cellular
category having poset $\PC_{[\mu,\la]} = \{ \nu \; | \; \mu \le \nu \le \la\}$. When these intervals are finite, one can filter the module category over $\CC$ by affine highest weight
subquotients.

\subsection{The cellular structure on Soergel bimodules}\label{subsec-lightleaves}

Having concluded our discussion of the abstract nonsense of cellular structures and highest weight categories, we now focus on the cellular structure on diagrammatic Soergel bimodules.

In \cite{EWGr4sb}, the category $\SB = \SB(W,S)$ of \emph{diagrammatic Bott-Samelson bimodules} was defined. The category depends not only on a Coxeter system $(W,S)$, but also on a choice of realization, giving rise to a polynomial ring $R$ over a commutative base ring $\Bbbk$. An object is a sequence $\un{w}$ of simple reflections, also known as an expression. The Karoubi envelope of $\SB$ is the category of \emph{diagrammatic Soergel bimodules}.

In the remainder of this chapter we assume the reader is familiar with the diagrammatic calculus for $\SB$, see \cite{EWGr4sb} for more details.

In \cite{EWGr4sb}, it was proven that $\SB$, for any Coxeter group $W$, is an object-adapted cellular category where the base ring is the polynomial ring $R$.
The cells are parametrized by the elements of the Coxeter group $W$, with partial order
given by the Bruhat order. To an expression $\un{w}$ and an element $x \in W$, the
set $E(\un{w},x)$ is given by the set of all subexpressions $\eb \subset \un{w}$ which multiply to give the element $x$. To such a subexpression $\eb$, following an algorithm of
Libedinsky \cite{LibLL}, one associates a diagram called a \emph{light leaf}, living in the morphism space $\Hom(\un{w}, \un{x})$ for some reduced expression of $x$. This is the morphism $c_{\eb}$ in the object-adapted cellular structure, here denoted $LL_{\eb}$, which gives half of a cellular basis element. Drawing this morphism upside-down gives the other half $\overline{LL}_{\eb}$ of a cellular basis element associated to $\iota(\eb) \in M(x,\un{w})$. For two sequences $\eb \in E(\un{w},x)$ and $\fb \in M(x,\un{y})$, the composition $\overline{LL}_{\fb} \circ LL_{\eb} = \LLL_{\fb,\eb}$ is a cellular basis element, called a \emph{double leaf}.

Because $\SB$ is an object-adapted cellular category, and intervals in the poset $W$ are finite sets, one deduces from Section \ref{subsec-cellsubquot} that $\SB\Mod$ is filtered by affine highest weight categories over $R$.

Now we discuss the details of the light leaf construction. They will be needed to prove two technical lemmas about the cellular structure on $\SB$. Light leaves are constructed
inductively. If $\un{w}$ is the empty expression, then a subexpression $\eb$ must also be empty, and $LL_{\eb}$ is the identity endomorphism of the monoidal identity. This handles the base case.

Let $\eb \subset \un{w}$ be a subexpression, thought of as an element of $E(\un{w},x)$ for some $x$. Consider $\un{w}^-$, which omits the last simple reflection in the expression, and
$\eb^-$, the corresponding subexpression, and let $x^-$ be the element expressed by $\eb^-$. Then, inductively, we have a morphism $LL_{\eb^-} \co \un{w} \to \un{x^-}$ for some reduced
expression of $x^-$. Let $s$ be the final simple reflection of $\un{w}$, so that $\un{w} = \un{w}^- \ot s$ as objects in $\SB$. Let us decorate this final index with either U1, U0, D1,
or D0, depending on whether: \begin{itemize} \item The simple reflection is contained (1) or not contained (0) in the subexpression $\eb$. That is, we have (1) if $x = x^- s$ and (0) if $x = x^-$. \item The simple reflection goes up (U) or down
(D) in the Bruhat order, relative to $x^-$. That is, we have (U) if $x^- s > x^-$, and (D) if $x^- s < x^-$. \end{itemize} For each of these four possibilities, there is a corresponding morphism $\phi_{\rm{last}} \in \Hom(\un{x^-} \ot s, \un{x})$, pictured below. Now we define $LL_{\eb} = \phi_{\rm{last}} \circ (LL_{\eb^-} \ot \id_s)$.

\begin{equation} \label{eq:LLconst} {
\labellist
\small\hair 2pt
 \pinlabel {$\un{x^-}$} [ ] at 13 30
 \pinlabel {$\un{x^-}$} [ ] at 89 30
 \pinlabel {$\un{x^-}$} [ ] at 168 30
 \pinlabel {$\un{x^-}$} [ ] at 225 30
 \pinlabel {$\un{x}$} [ ] at 13 80
 \pinlabel {$\un{x}$} [ ] at 89 80
 \pinlabel {$\un{x}$} [ ] at 168 95
 \pinlabel {$\un{x}$} [ ] at 225 95
 \pinlabel {$s$} [ ] at 36 30
 \pinlabel {$s$} [ ] at 108 30
 \pinlabel {$s$} [ ] at 188 30
 \pinlabel {$s$} [ ] at 248 30
\endlabellist
\centering
\ig{1}{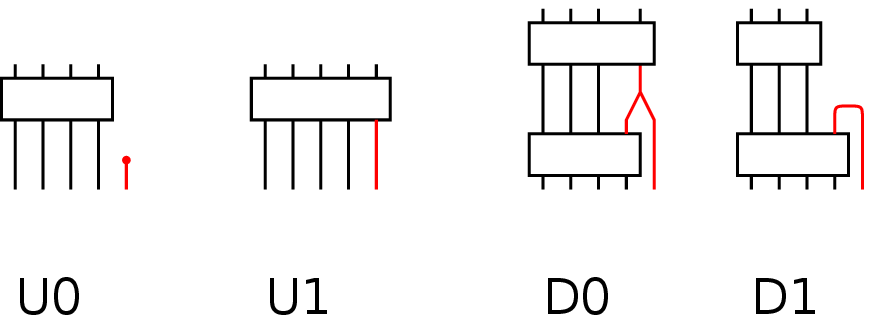}
} \end{equation}
Red represents the simple reflection $s$. The shaded boxes are \emph{rex moves}, morphisms built entirely of $2m$-valent vertices (no dots or trivalent vertices). Which rex move is used
to define the map is not actually relevant, so long as it has the correct source and target. (The neophyte reader can pretend that reduced expressions are chosen serendipitously, and can
pretend that each box is the identity map.)

Note that, inductively, every index in $\un{w}$ has been labeled with U1, U0, D1, D0 based on the subexpression $\eb$. Although $\eb$ is just the data of the zeroes and ones, it helps to
think of the Us and Ds as also being specified by $\eb$. For a longer exposition on this construction, see \cite[Section 6.1]{EWGr4sb}.

Here is an important consequence of the light leaf construction.

\begin{lemma} \label{lem:dotkernel} Let $\un{w}$ be a reduced expression, and $\un{z}$ be an arbitrary expression. Consider the quotient map $q \co \Hom(\un{w}, \un{z}) \to \Hom(\un{w},
\un{z})/J_{<w}.$ Then the kernel of $q$ is the span of all diagrams where a dot meets one of the strands on the bottom of the diagram. \end{lemma}

\begin{proof} Using the light leaves basis, a morphism in $J_{<w}$ factors through a light leaf corresponding to $\eb \in E(\un{w},v)$ for some $v < w$. Now $\un{w}$ is a reduced
expression, and $\eb$ is a subexpression containing some zero, and consequently the first zero in $\eb$ must be a $U0$. The corresponding light leaf therefore has a dot on this strand.
\end{proof}

\subsection{Standard modules and tensoring}\label{standardpreserved}

Note that $\SB$ is not just an object-adapted cellular category, it is also a monoidal category.

\begin{defn} If $\CC$ is a $\Bbbk$-linear monoidal category, then there is a monoidal action of $\CC$ on $\CC\Mod$ as follows. Given $\MC \co \CC \to \Bbbk\Mod$ and $B \in \CC$, let $B
\ot \MC$ be the functor which sends $X \in \CC$ to $\MC(B \ot X)$. It is clear how to define where $B \ot \MC$ sends a morphism $X \to X'$, and also how a morphism $B \to B'$ yields a
morphism $B \ot \MC \to B' \ot \MC$. \end{defn}

The goal of this section to prove the following proposition. Although the generating objects of $\SB$ are the simple reflections $s \in S$, let us denote these objects $B_s$ so that it is easier to distinguish between the simple reflection in $W$ and the corresponding object in $\SB$.

\begin{prop}\label{Prop:stand_action_preserv} Let $\SB\Mod_\D$ denote the full subcategory of $\SB\Mod$ consisting of functors with filtrations by standard modules $\Delta(x)$. Then $\SB\Mod_\D$ is preserved by the
action of $B_s$ under tensor product on the right (or left), for any simple reflection $s$. \end{prop}

For an element $w \in W$, let $\Delta(w)$ denote the corresponding standard module in $\SB\Mod$. By construction, this functor sends $X \in \SB\Mod$ to $\Hom(\un{w},X)/J_{< w}$, where
$\un{w}$ is some reduced expression of $w$. Let $\un{z}$ be an arbitrary expression. Note that, while $\Hom(\un{w},\un{z})$ has a basis of double leaves (our cellular basis), the
quotient $\Hom(\un{w},\un{z})/J_{<w}$ has a basis given by (upside-down) light leaves $\overline{LL}_{\eb}$ for $\eb \in M(\un{w},\un{z})$.

To prove the proposition it is enough to prove the following two lemmas.

\begin{lemma} \label{lem:stdpreserve1} Suppose that $w < ws$. Then $\Delta(w) \ot B_s$ has a submodule $\Delta(w)$ with quotient $\Delta(ws)$. \end{lemma}
	
\begin{lemma} \label{lem:stdpreserve2} Suppose that $ws < w$. Then $\Delta(w) \ot B_s$ has a submodule $\Delta(ws)$ with quotient $\Delta(w)$. \end{lemma}

In particular, the smaller element in the Bruhat order appears as the submodule, and the larger element as the quotient. These filtrations on $\Delta(w) \ot B_s$ arise in essentially the
same way that the light leaves basis is constructed.

\begin{proof}[Proofs of Lemmas \ref{lem:stdpreserve1},\ref{lem:stdpreserve2}] Let $\un{z}$ be any object of $\SB$. Suppose that one is to construct a light leaf map $\overline{LL}_{\eb}$ corresponding to $\eb \in M(\un{w},\un{z}s)$. That is, $\eb$ is
a subexpression of $\un{z}s$ which multiplies to $w$. Let $\fb$ be the subexpression of $\un{z}$ inside $\eb$, multiplying to an element $v$. Then either $w = vs$ or $w = v$, depending
on whether the final index in $\eb$ is labeled with 1 or 0. Whether this final index is labeled with U or D depends on whether $w < ws$ or $w > ws$.

Suppose that $w < ws$. Let $\un{z}$ be any object of $\SB$. We define a map $U0_{\un{z}}$ from $\Hom(\un{w},\un{z})$ to $\Hom(\un{w}, \un{z}s)$ as follows. The trapezoid is an arbitrary map from $\un{w}$ to $\un{z}$.
\begin{equation} \label{eq:natlU0} {
\labellist
\small\hair 2pt
 \pinlabel {$\un{w}$} [ ] at 32 6
 \pinlabel {$\un{w}$} [ ] at 124 6
 \pinlabel {$\un{z}$} [ ] at 32 45
 \pinlabel {$\un{z}$} [ ] at 124 45
 \pinlabel {$s$} [ ] at 159 45
 \pinlabel {$\mapsto$} [ ] at 77 25
\endlabellist
\centering
\ig{1}{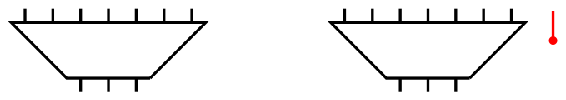}
} \end{equation}
This will send a light leaf in $\Hom(\un{w},\un{z})$ to the corresponding light leaf in $\Hom(\un{w},\un{z}s)$ with the final index $U0$. Clearly the family of maps $\{U0_{\un{z}}\}$ gives rise to a natural transformation $U0$ from $\Hom(\un{w},-)$ to $\Hom(\un{w},- \ot s)$.
Analogously, we define a natural transformation $D1$ from $\Hom(\un{w}s,-)$ to $\Hom(\un{w},- \ot s)$ as follows, which extends a light leaf by adding a final index $D1$.
\begin{equation} \label{eq:natlD1} {
\labellist
\small\hair 2pt
 \pinlabel {$\un{w}$} [ ] at 33 6
 \pinlabel {$s$} [ ] at 48 6
 \pinlabel {$\un{w}$} [ ] at 120 6
 \pinlabel {$\un{z}$} [ ] at 35 60
 \pinlabel {$\un{z}$} [ ] at 123 60
 \pinlabel {$s$} [ ] at 160 60
 \pinlabel {$\mapsto$} [ ] at 77 29
\endlabellist
\centering
\ig{1}{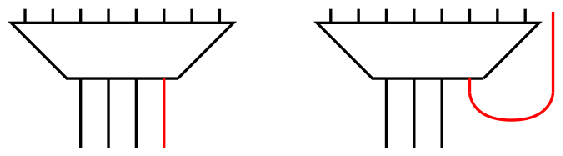}
} \end{equation}
Thus, applying $U0$ to the subspace spanned by light leaves, and applying $D1$ to the subspace spanned by light leaves, gives the subspace of $\Hom(\un{w},\un{z}s)$ spanned by light leaves.

Consider the quotient map $q \co \Hom(\un{w}, \un{z}s) \to \Hom(\un{w},\un{z}s)/J_{<w}$. The kernel of $q$ is spanned by diagrams with dots on bottom, by Lemma \ref{lem:dotkernel}. Now
consider the restriction of $q$ to the image of $U0$. If a morphism in $\Hom(\un{w},\un{z})$ had a dot on bottom (i.e. there is a dot on the bottom of the trapezoid) then after applying
U0 there is still a dot on bottom. Vice versa, if a diagram has a dot on bottom after applying U0, then there was a dot on the bottom of the trapezoid. Consequently, the preimage under
$U0$ of the ideal $J_{<w}$ (the kernel of $q$) is precisely the ideal $J_{<w}$. Thus $U0$ descends to an injective map $\psi \co \Hom(\un{w},M)/J_{< w} \to \Hom(\un{w},M \ot
B_s)/J_{<w}$. Another way to deduce that $\psi$ is injective is that it sends the light leaves basis of the source to a set of light leaves forming a subbasis of the target.

Now we claim that $D1$ induces an isomorphism $\phi$ from $\Hom(\un{w}s,\un{z})/J_{< ws}$ to the quotient of $\Hom(\un{w},\un{z}s)/J_{<w}$ by the image of $U0$. If we can show that
$\phi$ is well-defined, it will clearly be an isomorphism. This is because, once again, light leaves span the source, and the images of these light leaves are the complimentary subbasis
to the subbasis hit by $U0$.

Note that the map $D1$ from $\Hom(\un{w}s,\un{z})$ to $\Hom(\un{w},\un{z}s)$ does not have $J_{<ws}$ in the kernel, but $J_{<ws}$ will be in the kernel after taking the quotient both by
$J_{<w}$ and by the image of $U0$. The reason is that $\un{w}s$ is a reduced expression, so any lower term (i.e. map in $J_{<ws}$) will have a dot on the bottom (of the trapezoid) by
Lemma \ref{lem:dotkernel}. If this dot is on one of the strands in $\un{w}$, then it dies in $J_{<w}$ (it is on the bottom of the diagram after applying $D1$). If the dot meets the final $s$-colored strand, the result will be in the image of
$U0$.
\begin{equation} \ig{1}{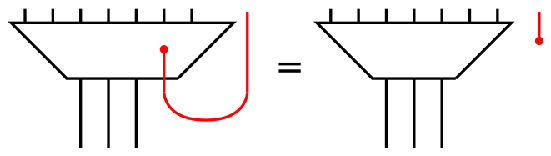} \end{equation} (See \cite[(3.9)]{EKho}.) Thus $\phi$ is well-defined. This concludes the proof of Lemma \ref{lem:stdpreserve1}.

Now suppose that $ws < w$. It does no harm to assume that our chosen reduced expression $\un{w}$ ends in the reflection $s$. Let $\un{y} = \un{w}^-l$ denote a reduced expression for $ws$. We now define natural transformations $U1 \co \Hom(\un{y},-) \to \Hom(\un{w},- \ot s)$ and $D0 \co \Hom(\un{w},-) \to \Hom(\un{w},- \ot s)$ as follows.
\begin{equation} \label{eq:natlU1} {
\labellist
\small\hair 2pt
 \pinlabel {$\un{y}$} [ ] at 32 6
 \pinlabel {$\un{y}$} [ ] at 124 6
 \pinlabel {$s$} [ ] at 161 7
 \pinlabel {$\un{z}$} [ ] at 31 45
 \pinlabel {$\un{z}$} [ ] at 124 45
 \pinlabel {$s$} [ ] at 161 44
 \pinlabel {$\mapsto$} [ ] at 76 23
\endlabellist
\centering
\ig{1}{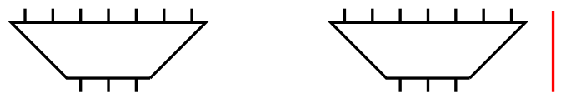}
} \end{equation}
\begin{equation} \label{eq:natlD0} {
\labellist
\small\hair 2pt
 \pinlabel {$\un{y}$} [ ] at 32 6
 \pinlabel {$s$} [ ] at 48 6
 \pinlabel {$\un{y}$} [ ] at 119 8
 \pinlabel {$s$} [ ] at 147 7
 \pinlabel {$\un{z}$} [ ] at 35 60
 \pinlabel {$\un{z}$} [ ] at 122 60
 \pinlabel {$s$} [ ] at 159 59
 \pinlabel {$\mapsto$} [ ] at 80 30
\endlabellist
\centering
\ig{1}{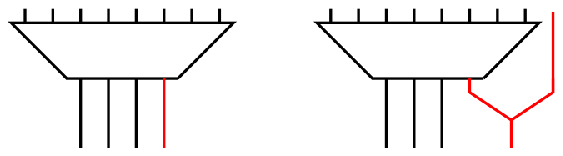}
} \end{equation}
Once again, by construction, these maps send light leaves to light leaves. Again, $U1$ descends to an injective map $\Hom(\un{y},-)/J_{< ws} \to \Hom(\un{w},- \ot s)/J_{<w}$. After all, if a diagram in the image of $U1$ has a dot on the bottom, that dot can not be on the final $s$ strand, so it must have come from a dot on $\un{y}$. Again, $D0$ descends to an isomorphism from $\Hom(\un{w},-)/J_{<w}$ to the quotient of $\Hom(\un{w},- \ot s)$ by both $J_{<w}$ and the image of $U1$. This is because a dot on the bottom of the trapezoid either corresponds to a dot on the bottom of the target, or to a dot entering the trivalent vertex, yielding a morphism in the image of $U1$.
\begin{equation} \ig{1}{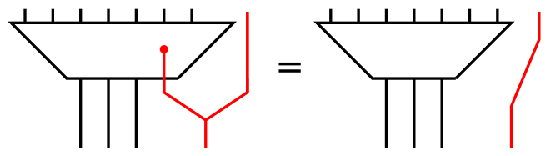} \end{equation}
(See \cite[(3.4)]{EKho}.) \end{proof}

\subsection{Coset subquotients and parabolic subgroups}\label{SS_coset_quot}

Let $\SB(W)$ denote the diagrammatic category of Soergel bimodules attached to a Coxeter system $(W,S)$ with a realization $\hg$. For $I \subset S$, not necessarily finitary, we write
$\SB(W_I)$ for the corresponding category attached to $(W_I,I)$ with realization $\hg$. In particular, both categories are $R$-bilinear for the same polynomial ring $R$. There is a fully
faithful inclusion $\SB(W_I) \subset \SB(W)$. This gives rise to an action of $\SB(W_I)$ on $\SB(W)$ induced from the monoidal structure on $\SB(W)$.

Let $C$ be a right coset in $W / W_I$. Then $C$ has an element $x$ of minimal length. We write $\le C$ for the ideal in
$W$ consisting of all elements $y$ for which there exists $z \in C$ with $y \le z$. Then $\le C$ is itself a union of cosets in $W / W_I$. We write $< C$ for $(\le C) \setminus C$.

\begin{defn} Let $\SB_C(W)$ denote the cellular subquotient category $\SB_{\le C} / \SB_{< C}$, as in Section \ref{subsec-cellsubquot}. \end{defn}

Since $x$ is a minimal element modulo $< C$, it was shown in \cite[first paragraph of Section 6.5]{EWGr4sb} that there is a canonical indecomposable isomorphism class in $\SB / \SB_{<
C}$ attached to $x$. Namely, for any two reduced expressions $\un{x}$ for $x$, the corresponding objects $\un{x} \in \SB / \SB_{< C}$ are canonically isomorphic. We simply refer to this
object abusively as $x$.

\begin{prop}\label{Prop:quot_isom} Consider the functor $\phi_x \co \SB(W_I) \to \SB_C(W)$ induced by acting (on the right) on the object $x \in \SB_C(W)$. It induces an equivalence of
Karoubi envelopes, which is obviously compatible with the $\SB(W_I)$ action on both sides. \end{prop}

\begin{proof} Fix a reduced expression $\un{x}$ of $x$, and an arbitrary expression $\un{w}$ in $W_I$. For any $y \in W_I$ we claim that there is a bijection between subexpressions $\eb
\subset \un{w}$ which express $y$, and subexpressions $\fb$ of the concatenation $\un{x}\un{w}$ which expresss $xy$. Moreover, this bijection sends $\eb$ to the subexpression $\fb$ which
contains all of $\un{x}$, and agrees with $\eb$ on $\un{w}$. After all, if $\fb$ contained a proper subset of $\un{x}$, expressing an element $z < x$, then $\fb$ must express an element
in the coset $zW_I$, which is disjoint from the subset $xW_I$.  In the remainder of this proof, $\eb$ will be a subexpression of $\un{w}$, and $\fb$ will be the corresponding subexpression of $\un{x}\un{w}$ under this bijection.

Next we observe that the decorations (U0, U1, D0, D1) on the simple reflections in $\un{w}$ coming from the subexpression $\eb$ agree with the decorations on the
corresponding simple reflections in the $\un{w}$ portion of $\un{x}\un{w}$ coming from $\fb$. This amounts to the fact that the Bruhat order on the coset $xW_I$ agrees with a Bruhat
order on $W_I$ itself. Moreover, the decorations on the $\un{x}$ portion of $\un{x}\un{w}$ are all U1, because $\fb$ contains all of $\un{x}$ and $\un{x}$ is a reduced expression.

It now follows from the algorithm for constructing light leaves (see Section \ref{subsec-lightleaves}) that $LL_{\fb} = \id_{\un{x}} \ot LL_{\eb}$. This is because the map $\phi_{\rm{last}}$ associated to U1 is just the
identity map, see \eqref{eq:LLconst}.

On objects, $\phi_x$ sends an expression $\un{w}$ in $W_I$ to the concatenation $\un{x}\un{w}$. Let $\un{w}$ and $\un{w}'$ be two expressions in $I$. Now $\Hom_{\SB(W_I)} (\un{w},
\un{w'})$ has its cellular basis (as a right $R$-module) consisting of double leaves $\LLL_{\eb, \eb', y}$, where $y \in W_I$, and $\eb \subset \un{w}$ and $\eb' \subset \un{w}'$ are
subexpressions which express $y$. Meanwhile, $\Hom_{\SB}(\un{x}\un{w}, \un{x}\un{w}')$ has a basis as a right $R$-module given by double leaves $\LLL_{\fb, \fb', z}$ for $z \in W$, and
$\fb\subset \un{x}\un{w}$, $\fb' \subset \un{x}\un{w}'$ subexpressions for $z$. Any such double leaf lies within $\SB_{\le C}$, and the double leaves which are killed in the passage to
$\SB_C$ are precisely those with $z \notin C$. Thus $\Hom_{\SB_C}(\un{x}\un{w}, \un{x}\un{w}')$ has a basis given by $\LLL_{\fb, \fb', xy}$ for $y \in W_I$. By definition, the functor
$\phi_x$ sends $\LLL_{\eb, \eb',y}$ to $\id_{\un{x}} \otimes \LLL_{\eb, \eb', y}$. Therefore, by the bijection between subsequences $\eb$ for $y$ and subsequences $\fb$ for $xy$, and the
related observations above, the functor $\phi_x$ is fully faithful, sending a basis of morphisms to a basis of morphisms.

We claim that $\phi_x$ is essentially surjective after passage to the Karoubi envelope. This follows from the classification of indecomposables in both categories. The indecomposables in
the Karoubi envelope of $\SB_C$ are given by the images of the indecomposables in $\SB$ parametrized by elements of $C$; they are described as the images of the ``top idempotents" inside
a reduced expression $\un{x}\un{w}$, where the top idempotent is the unique indecomposable idempotent which does not factor through lower cells (see \cite[Theorem 6.25]{EWGr4sb}).
Meanwhile, the indecomposables in the Karoubi envelope of $\SB(W_I)$ are the images of the top idempotents inside $\un{w}$. Because $\phi_x$ is fully faithful, these top idempotents
agree, and the indecomposables are matched up by $\phi_x$. \end{proof}

\begin{remark} The proof of \cite[Theorem 6.25]{EWGr4sb} uses cellular methods to classify indecomposables in $\SB$, though it does not directly discuss indecomposables in cellular
subquotients. Analogous arguments were used in \cite[Proposition 2.24]{ELauda} to classify the indecomposable objects in any object-adapted cellular category, and both $\SB_C$ and
$\SB(W_I)$ are object-adapted cellular categories. \end{remark}

\begin{prop}\label{Prop:quot_iso_left} Let $C$ be a left coset in $W_I \setminus W$ with minimal element $x$. Define $\SB_C(W)$ analogously, and consider the functor $\SB(W_I) \to \SB_C(W)$ induced by acting on the left on the object $x \in \SB_C(W)$. Then this induces an equivalence of Karoubi envelopes, compatible with the $\SB(W_I)$ action. \end{prop}

\begin{proof} There is a covariant, monoidally-contravariant autoequivalence of $\SB$ which flips a diagram horizontally. Applying this to the set of light leaves morphisms, one obtains
an alternate object-adapted cellular structure on $\SB$. The light leaves morphisms and double leaves basis were defined by reading an expression and a subexpression from left to right,
in the usual fashion. This alternate structure defines light leaves by a similar algorithm that reads from right to left instead.

Now there is a bijection between subexpressions of $\un{w}$ for $y \in W_I$, and subexpressions of $\un{w} \un{x}$ for $yx$. One can repeat the entire proof above in this mirrored
context. \end{proof}

\subsection{Localization}\label{subsec-localization}

We mention one more result from \cite{EWGr4sb} that we will need. Let $Q$ denote the fraction field of the polynomial ring $R$.

\begin{thm}\label{Thm:Soergel_base_change} After base change to $Q$, the category $\SB \ot_R Q$ is semisimple. This is also true after base change to any field where each polynomial in the $W$-orbit of the simple
roots is invertible. \end{thm}

The reader should think that this is a consequence of Proposition \ref{prop:whencqh}. The cellular pairings on $\SB$ are not non-degenerate over $R$, but they become non-degenerate over
$Q$, which is easy to prove.

The proof of this fact in \cite{EWGr4sb} operates in a different fashion. In \cite[Chapter 4]{EWGr4sb} a category $\DC^{\std}$ is constructed diagrammatically, and is proven to be
isomorphic to the $R$-linear 2-groupoid of $W$. Actually, this is the main result of \cite{EWFenn}, where further discussion of 2-groupoids can be found. The base change of $\DC^{\std}$
from $R$ to $Q$ is therefore semisimple. After all, the additive closure (the category obtained by allowing for finite direct sums) of any 2-groupoid, linear over a field, is semisimple
(because every nonzero map between indecomposable objects is an isomorphism).

In \cite[Section 5.4]{EWGr4sb}, and specifically in \cite[Theorem 5.16 and Proposition 5.22]{EWGr4sb}, it is proven diagrammatically that the category $\SB \ot_R Q$ has Karoubi envelope
equivalent to the additive closure of $\DC^{\std} \ot_R Q$. Hence it is semisimple. In fact, the proof in \cite{EWGr4sb} that the Karoubi envelope of $\SB \ot_R Q$ is the additive
closure of a 2-groupoid does not require passage to the full fraction field of $R$. One need only invert homogeneous polynomials, or even just the $W$-orbit of the simple roots.

%% file: SoergelO.tex
\section{Categories $\mathcal{O}$ from Soergel bimodules} \label{sec-hwnonsense}
\subsection{Highest weight categories}\label{SS_HW}
Here we recall mostly standard facts regarding highest weight categories (over complete
regular local rings).
\subsubsection{Finite posets}\label{SSS_HW_fin}
Let $\Lambda$ be a finite poset and $\Ring$ be a regular complete local ring. Following
Rouquier, \cite[Section 4.1]{rouqqsch} we introduce the notion of a highest weight
(a.k.a. quasi-hereditary) category  with poset $\Lambda$ over $\Ring$. This a Noetherian abelian
$\Ring$-linear category $\Cat_{\Ring}$ with  enough projectives
$P_\Ring(\lambda)$  indexed by $\lambda\in\Lambda$ and a collection of
{\it standard objects} $\Delta_\Ring(\lambda)$ subject to the following conditions:
\begin{itemize}
\item[(HW0)] $\operatorname{Hom}_{\Cat_\Ring}(P_\Ring(\lambda),P_\Ring(\lambda'))$ is a free
finitely generated $\Ring$-module.
\item[(HW1)] $\operatorname{Hom}_{\Cat_\Ring}(\Delta_\Ring(\lambda),\Delta_\Ring(\lambda'))$
is zero unless $\lambda\leqslant\lambda'$. Moreover, if $\lambda=\lambda'$, then
this $\Ring$-module is free of rank $1$.
\item[(HW2)] The object $P_\Ring(\lambda)$ admits an epimorphism onto $\Delta_\Ring(\lambda)$
whose kernel is filtered by $\Delta_\Ring(\lambda')$ with $\lambda'>\lambda$.
\end{itemize}
Note that (HW0) implies that $\Cat_\Ring$ is isomorphic to the category of (finitely generated) modules
over an $\Ring$-algebra $A_\Ring$ that is a free $\Ring$-module of finite rank. Also note that
the $A_\Ring$-modules corresponding to $\Delta_\Ring(\lambda)$ are free over $\Ring$. For any regular complete
local ring $\Ring'$ that is an $\Ring$-algebra, the base change $\Cat_{\Ring'}:=\Ring'\otimes_\Ring A_{\Ring}\operatorname{-mod}$ is again a highest weight category. We will write $P_{\Ring'}(\lambda)$ for $\Ring'\otimes_\Ring P_\Ring(\lambda),\Delta_{\Ring'}(\lambda)$ for $\Ring'\otimes_\Ring\Delta_\Ring(\lambda)$.

Below we will need the following easy and standard lemma.

\begin{Lem}\label{Lem:dir_sum_filtration}
Suppose that $\Lambda$ is a linearly ordered poset and let $\Delta_\Ring(\lambda),\lambda\in \Lambda,$ be a collection of
objects in $\Cat_R$ satisfying (HW1). Further, let $M_1,M_2$ be two objects in $\Cat_R$ such that
$M_1\oplus M_2$ has a $\Lambda$-filtration $F_{\leqslant \lambda}\subset M_1\oplus M_2, (\lambda\in \Lambda),$
such that $F_{\leqslant \lambda}/F_{<\lambda}\cong \Delta_\Ring(\lambda)^{\oplus ?}$.
Then each $M_1,M_2$ has such a filtration.
\end{Lem}

The following easy lemma reduces checking that the category is highest weight to the case
when the category is defined over a field.

\begin{Lem}\label{Lem:hw_field}
Let $\Ring$ be a regular complete local ring, $\F$ its residue field, and $A_\Ring$ be an
$\Ring$-algebra that is free of finite rank as an $\Ring$-module. Suppose that $A_\F\operatorname{-mod}$
is a highest weight category with poset $\Lambda$ and standard objects $\Delta_{\F}(\lambda)$.
Then $A_\Ring\operatorname{-mod}$ is a highest weight category over $\Ring$ with poset $\Lambda$.
The standard objects $\Delta_\Ring(\lambda)$ are unique flat deformations of the $\Delta_{\F}(\lambda)$'s.
\end{Lem}

Let $L_{\F}(\lambda)$ denote the head of $P_{\F}(\lambda)$ and $I_\F(\lambda)$
be the injective hull of $L_{\F}(\lambda)$. Let $I_\Ring(\lambda)$ denote a unique
$\Ring$-flat deformation of $I_{\F}(\lambda)$. We remark that being $\Ring$-flat has
a purely categorical meaning: an object $M$ is $\Ring$-flat if and only if $\Hom(P_\Ring, M)$
is flat over $\Ring$ for all projective objects $P_{\Ring}$.   Then $I_{\Ring}(\lambda)$ is an
indecomposable relatively $\Ring$-injective object in $\Cat_{\Ring}$ (meaning that
$\Ext^i(M,I_\Ring(\lambda))=0$ for any object $M\in\Cat_{\Ring}$ that is flat over
$\Ring$ and any $i>0$).


In $\Cat_\Ring$ we have costandard objects $\nabla_\Ring(\lambda),\lambda\in \Lambda$
(=standard objects for the highest weight  category $\Cat_\Ring^{opp}$).
Recall that, by definition, $\nabla_\Ring(\lambda)$ is the intersection of the kernels
of all homomorphisms $I_\Ring(\lambda)\rightarrow I_\Ring(\lambda')$, $\lambda'>\lambda$. The object
$\nabla_\Ring(\lambda)$ is flat over $\Ring$.

Recall that an object   $M\in \Cat_{\F}$ is costandardly filtered if and only
if $\operatorname{Ext}^i(\Delta_\F(\lambda),M)=0$ for all $\lambda\in \Lambda$
and $i>0$. Similarly, an object $M\in \Cat_{\Ring}$ is costandardly filtered if
and only if it is $\Ring$-flat and $\operatorname{Ext}^i(\Delta_\Ring(\lambda),M)=0$
for all $\lambda\in \Lambda$ and $i>0$. Further,   $\operatorname{Hom}_{\Cat_\Ring}(\Delta_\Ring(\lambda),\nabla_\Ring(\lambda'))$
is a free $\Ring$-module of rank $\delta_{\lambda\lambda'}$.

We also have indecomposable tilting objects $T_\Ring(\lambda), \lambda\in \Lambda$. Recall that an
object $M\in \Cat_\Ring$ is called tilting if it is both standardly and costandardly filtered.
The object $T_\Ring(\lambda)$ admits a monomorphism $\Delta_\Ring(\lambda)\rightarrow
T_\Ring(\lambda)$ whose cokernel is filtered by $\Delta_\Ring(\lambda')$ with
$\lambda'<\lambda$.

We can define the Ringel dual category $\Cat_\Ring^\vee$ as
$\operatorname{End}_{\Cat_\Ring}(T_\Ring)\operatorname{-mod}$,
where $T_\Ring:=\bigoplus_\lambda T_\Ring(\lambda)$. This category is highest weight over $\Ring$
with poset $\Lambda^{opp}$ and standard objects $\Delta^\vee_\Ring(\lambda):=\operatorname{Hom}_{\Cat_\Ring}(\Delta_{\Ring}(\lambda),T_\Ring)$.
Note that $\operatorname{Hom}_{\Cat_\Ring}(P_\Ring(\lambda), T_\Ring)$ is the tilting
object $T^\vee_\Ring(\lambda)\in \Cat_\Ring^\vee$. 

If $\Lambda_0\subset \Lambda$ is a poset ideal, then the Serre span
$\Cat_{\Ring,\Lambda_0}$ of $\Delta_{\Ring}(\lambda), \lambda\in \Lambda_0$, has
a natural structure of a highest weight category with poset $\Lambda_0$
and standard objects $\Delta_\Ring(\lambda),\lambda\in \Lambda_0$. The costandard
objects are $\nabla_\Ring(\lambda)$ and the indecomposable tilting objects are
$T_\Ring(\lambda)$.

If $\Lambda^0\subset \Lambda$ is a poset coideal, then the quotient category
$\Cat^{\Lambda^0}_\Ring:=\Cat_{\Ring}/\Cat_{\Ring,\Lambda\setminus \Lambda^0}$ is a highest
weight category with poset $\Lambda^0$. Its standard objects are the images
of $\Delta_\Ring(\lambda), \lambda\in \Lambda^0$, under the quotient morphism $\pi:\Cat_\Ring\twoheadrightarrow
\Cat_\Ring^{\Lambda^0}$.
The left adjoint $\pi^!$ gives rise to an identification of the category
$(\Cat^{\Lambda^0}_\Ring)^{\Delta}$ of standardly filtered objects to
the full subcategory of $\Cat_\Ring^\Delta$ consisting of objects filtered
by $\Delta_\Ring(\lambda), \lambda\in \Lambda^0$. Note that the indecomposable
tiltings in $\Cat^{\Lambda_0}$ are the images under $\pi$ of the indecomposable tiltings
$T_\Ring(\lambda)\in \Cat_\Ring, \lambda\in \Lambda^0$.

Note that we have natural identifications $\Cat_{\Ring,\Lambda_0}^\vee\cong (\Cat_\Ring^\vee)^{\Lambda_0}$
and $(\Cat_\Ring^{\Lambda^0})^\vee\cong (\Cat^\vee_\Ring)_{\Lambda^0}$.

By a {\it naive duality} of a highest weight category $\Cat_{\F}$  we mean
an  equivalence $\bullet^*: \Cat_{\F}\xrightarrow{\sim} \Cat_{\F}^{opp}$
with $\Delta_{\F}(\lambda)^*=\nabla_\F(\lambda)$. It is not automatic that such
a duality exists but it does in many natural examples. Note that $L_\F(\lambda)^*\cong L_\F(\lambda)$  so $\bullet^*$
is the identity on $K_0$.
 Note that $P_{\F}(\lambda)^*\cong I_\F(\lambda),
I_\F(\lambda)^*\cong P_\F(\lambda),T_\F(\lambda)^*\cong T_\F(\lambda)$.
Note also that we do not require that the naive duality equivalence is an involution.

\begin{Lem}\label{Lem:Ring_naive_dual}
If $\Cat_{\F}$ admits a naive duality then so does $\Cat_\F^\vee$.
\end{Lem}
\begin{proof}
The duality $\bullet^*$ induces an anti-automorphism, say $\sigma$, of the algebra $\operatorname{End}(T_\F)$
that preserves the indecomposable central idempotents. This gives rise to a contravariant equivalence
of $\operatorname{End}(T_{\F})\operatorname{-mod}, M\mapsto (\,^\sigma M)^*$.
This equivalence switches $P^\vee_{\F}(\lambda)$ and $I^\vee_{\F}(\lambda)$ for all $\lambda\in \Lambda$.
Hence it switches $\Delta^\vee_{\F}(\lambda)$ and $\nabla^\vee_{\F}(\lambda)$.
\end{proof}

\subsubsection{Coideal finite posets}\label{SS_hw_coid_fin}
Now let $\Lambda$ be an infinite poset. We say that $\Lambda$ is {\it coideal finite} if for every
$\lambda\in \Lambda$, the set $\{\lambda'\in \Lambda| \lambda'\geqslant \lambda\}$ is finite.
We will usually denote coideal finite posets by $\Lambda^+$.
Similarly, one introduces the notion of an {\it ideal finite} poset. Note that $\Lambda^+$
is coideal finite if and only if $(\Lambda^+)^{opp}$ is ideal finite.

By a highest weight category over $\Ring$ with poset $\Lambda^+$ we mean an abelian $\Ring$-linear category $\Cat^+_\Ring$ such that
\begin{itemize}
\item The indecomposable projectives $P_\Ring(\lambda)$ in $\Cat^+_\Ring$ are indexed by $\Lambda$ and satisfy (HW0).
\item Every object in $\Cat^+_\Ring$ is a direct limit of its subobjects that are quotients of finite sums of projectives.
\item For every $M\in \Cat^+_\Ring$ and every $\lambda\in \Lambda^+$, the $\Ring$-module $\Hom_{\Cat^+_\Ring}(P_\Ring(\lambda),M)$ is finitely generated.
\item  There are standard objects
$\Delta^+_\Ring(\lambda)$ (that are a part of the structure) subject to axioms (HW1)-(HW2)
from \S\ref{SSS_HW_fin}.
\end{itemize}
Consider the algebra $A_\Ring:=\bigoplus_{\lambda,\lambda'}\Hom_{\Cat_\Ring}(P_\Ring(\lambda'),P_\Ring(\lambda))$
that is an associative (generally, non-unital) $\Ring$-algebra with a distinguished collection of idempotents $e_\lambda,
\lambda\in \Lambda^+$.  The functor $M\mapsto \bigoplus_{\lambda}\Hom(P_\Ring(\lambda),M)$
identifies $\Cat^+_\Ring$ with the full subcategory of $A_\Ring$-modules $M$ satisfying the following two
conditions
\begin{enumerate}
\item locally unital (meaning that for every $m\in M$ there is a finite subset $\Lambda(m)\subset \Lambda$
with $\sum_{\lambda\in \Lambda(m)} e_\lambda m=m$)
\item and such that $e_\lambda M$ is a finitely generated module over $\Ring$ for every $\lambda\in \Lambda$.
\end{enumerate}
We denote this category by $A_\Ring\operatorname{-mod}$
so that $A_R\operatorname{-mod}\cong \Cat^+_{\Ring}$. We can do a base change for $\Cat^+_\Ring$ similarly to \S\ref{SSS_HW_fin} getting highest weight categories $\Cat^+_{\Ring'}$.

The definition is arranged in such a way that, for any finite coideal $\Lambda^0\subset \Lambda$, the quotient
$\Cat_\Ring^{+\Lambda^0}$ defined by the projectives $P_\Ring(\lambda), \lambda\in \Lambda^0,$
is a highest weight category over $\Ring$ with poset $\Lambda^0$ in the sense of \S\ref{SSS_HW_fin}. This is the category of all finitely generated
modules over the algebra $e_{\Lambda^0}A_\Ring e_{\Lambda^0}$, where $e_{\Lambda^0}$ is the idempotent corresponding
to $\Lambda^0$.

Thanks to the equivalence $\Cat^+_\Ring\cong A_\Ring\operatorname{-mod}$, the object $$I_\Ring(\lambda):=
\bigoplus_{\lambda'} \Hom_\Ring(e_\lambda A_\Ring e_{\lambda'},\Ring)$$ is an indecomposable relatively $\Ring$-injective
module (in the same sense as in \S\ref{SSS_HW_fin}).
We can define the costandard object $\nabla_\Ring(\lambda)\subset I_\Ring(\lambda)$ similarly
to the case of finite posets.  By the construction and axioms (HW1),(HW2), we see that
$\operatorname{Ext}^i(\Delta_\Ring(\lambda),\nabla_\Ring(\lambda'))$ is a free $\Ring$-module
of rank $\delta_{i0}\delta_{\lambda\lambda'}$. We also see that an object in
$\Cat_\Ring^+$  is costandardly filtered if and only if it is flat over $\Ring$,  embeds into a finite sum
of $I_{\Ring}(\lambda)$'s and $\operatorname{Ext}^i_{\Cat^+_\Ring}(\Delta_{\Ring}(\lambda),M)=0$
for all $i>0$ and all $\lambda$.

\subsubsection{Ideal finite posets}\label{SSS_ideal_hw}
Now let $\Lambda^-$ be an ideal finite poset.

By a highest weight category over $\Ring$ with poset $\Lambda^-$ we mean a Noetherian $\Ring$-linear abelian
category $\Cat^-_\Ring$ with standard objects $\Delta^-_\Ring(\lambda)$ such that
\begin{itemize}
\item Any object in $\Cat^-_\Ring$ admits a nonzero homomorphism from some $\Delta^-_\Ring(\lambda)$.
\item For every finite poset ideal $\Lambda_0\subset \Lambda^-$, the Serre subcategory
$\Cat_{\Ring,\Lambda_0}\subset \Cat^-_\Ring$ spanned by $\Delta^-_\Ring(\lambda)$
is highest weight over $\Ring$ with poset $\Lambda_0$ and standard objects
$\Delta_\Ring^-(\lambda),\lambda\in \Lambda_0$.
\end{itemize}
 In other words, $\Cat^-_\Ring$ is the union
of its highest weight subcategories. Note that $\Cat^-_\Ring$ has costandard
and indecomposable tilting objects indexed  by $\Lambda^-$ but does not have
enough projectives and injectives, in general. We can also define the notion
of an $\Ring$-flat object.

We have a version of the Ringel duality between highest weight categories for
$\Lambda^+$ and $\Lambda^-=(\Lambda^+)^{opp}$, we  realize $\Cat_\Ring^-$
as the category whose Ringel dual is $\Cat_\Ring^+$. Namely, take a finite coideal
$\Lambda^0\subset \Lambda^+$ and form the algebra
$$B_\Ring^{\Lambda_0}:=\operatorname{End}_{\Cat_\Ring^{\Lambda_0}}(\bigoplus_{\lambda\in \Lambda^0}T_\Ring^{\Lambda^0}(\lambda))$$
Note that, for a coideal inclusion $\Lambda^0\subset \Lambda^1$, we have a natural epimorphism
$B_\Ring^{\Lambda^1}\twoheadrightarrow B_\Ring^{\Lambda^0}$. So the algebras $B_\Ring^{\Lambda^0}, \Lambda^0\subset \Lambda^+$ form a projective system. For $\Cat_\Ring^-$ we take the category of topological modules over
$\varprojlim_{\Lambda^0} B_\Ring^{\Lambda^0}$ with discrete topology. It is straighforward to check that
$\Cat_\Ring^-$ is a highest weight category with poset $\Lambda^-$. Moreover, for every poset
ideal $\Lambda^0\subset \Lambda^-$, we have $\Cat_\Ring^{+\Lambda^0}=(\Cat^-_{\Ring,\Lambda^0})^\vee$
and we have equivalences $(\Cat_\Ring^-)^{\Delta}\cong ((\Cat_\Ring^+)^{\Delta})^{opp}$ and
$\Cat_\Ring^-\operatorname{-tilt}\cong (\Cat_\Ring^+\operatorname{-proj})^{opp}$.

\subsection{Standardly stratified categories}
\subsubsection{Standardly stratified structures}
Here we are going to recall the notion of a standardly stratified structure introduced
in \cite[Section 2]{LW} (for categories of modules for finite dimensional algebras over a field)
in a slightly more general setting (posets of the form $\Lambda^\pm$ and categories
over $\Ring$).

First, let us give the definition for categories of the form $A_\Ring\operatorname{-mod}$,
where $A_\Ring$ is an $\Ring$-algebra that is a free finite rank $\Ring$-module. Let $\Lambda$
be the indexing set for the indecomposable projectives. The additional structure involved
in the definition of a standardly stratified category is a partial preorder $\preceq$ on
$\Lambda$, equivalently an epimorphism $\varrho$ from $\Lambda$ to a poset $\Xi$.

There are three axioms. Here is the first one.
\begin{itemize}
\item[(SS0)] For a poset coideal $\Xi^0\subset \Xi$, let $e_{\Xi^0}$ denote an idempotent
in $A_\Ring$ corresponding to the projectives labelled by $\varrho^{-1}(\Xi^0)$. Then
$A_\Ring/A_\Ring e_{\Xi^0}A_\Ring$ is free over $\Ring$.
\end{itemize}
Note that this condition is vacuous if $\Ring$ is a field.

For a poset ideal $\Xi_0\subset \Xi$, let $\Cat_{\Ring,\Xi_0}$ denote the full subcategory
in $A_\Ring\operatorname{-mod}$ consisting of all modules annihilated by $e_{\Xi\setminus \Xi_0}$.
We write $\Cat_{\Ring,\preceq \xi}$ for $\Cat_{\Ring,\Xi_0}$ with $\Xi_0:=\{\xi'\in\Xi| \xi'\preceq\xi\}$.
The notation $\Cat_{\Ring,\prec\xi}$ has the similar meaning. So  we can form the quotient
category $\Cat_{\Ring,\xi}:=\Cat_{\Ring,\preceq \xi}/\Cat_{\Ring,\prec \xi}$, note that it is still
equivalent to the category of modules over a free, finite rank $\Ring$-algebra. Let
$\pi_\xi: \Cat_{\Ring,\leqslant \xi}\twoheadrightarrow \Cat_{\Ring,\xi}$ denote the
quotient functor. For general reasons, $\pi_{\xi}$ has both a left adjoint functor, to be denoted by
$\Delta_{\Ring,\xi}$ or $\pi_{\xi}^!$, and a right adjoint functor, denoted by $\nabla_{\Ring,\xi}$ or $\pi_{\xi}^*$.

Here is our second axiom.
\begin{itemize}
\item[(SS1)] The functor $\pi_{\xi}^!$ is exact.
\end{itemize}

Now pick $\lambda\in \Lambda$ and let $\xi:=\varrho(\lambda)$. Let $\iota_{\leqslant \xi}$
denote the inclusion $\Cat_{\Ring,\leqslant \xi}\hookrightarrow \Cat_\Ring$ and let $\iota_{\leqslant \xi}^!$
be the left adjoint functor. We write $P_\Ring(\lambda)$ for the indecomposable projective in $\Cat_{\Ring}$
indexed by $\lambda$ and $P_{\Ring,\xi}(\lambda)$ for $\pi_{\xi}\circ \iota_{\leqslant \xi}^!(P_\Ring(\lambda))$,
this is an indecomposable projective in $\Cat_{\Ring,\xi}$. Define the standard object $\Delta_\Ring(\lambda)$
(for our standardly stratified structure) as $\pi_{\xi}^!(P_{\Ring,\xi}(\lambda))$.

Here is our third axiom.
\begin{itemize}
\item[(SS2)] The object $P_{\Ring}(\lambda)$ admits an epimorphism onto $\Delta_\Ring(\lambda)$
and the kernel has a filtration by $\Delta_\Ring(\lambda')$ with $\lambda\prec \lambda'$.
\end{itemize}

Note that, for a regular complete local $\Ring$-algebra $\Ring'$, the category $\Cat_{\Ring'}:=\Ring'\otimes_\Ring A_\Ring\operatorname{-mod}$
is standardly stratified with the same pre-order.  We also have the following analog of Lemma \ref{Lem:hw_field}.

\begin{Lem}\label{Lem:SS_field}
Let $A_\Ring$ be an $\Ring$-algebra that is free of finite rank over $\Ring$ and let $\F$ denote the
residue field of $\Ring$. Suppose that the category $A_\F\operatorname{-mod}$ is standardly
stratified with respect to a pre-order $\preceq$ on $\Lambda$. Then $A_\Ring\operatorname{-mod}$
is standardly stratified with respect to the same pre-order.
\end{Lem}
\begin{proof}
Let us prove that (SS0) holds for $A_\Ring\operatorname{-mod}$, everything else is easy.
We may assume that each indecomposable projective occurs in $A_\Ring$ exactly once, i.e., $A_\Ring=\operatorname{End}_{\Cat_\Ring}(\bigoplus_\lambda P_\Ring(\lambda))^{opp}$.
Set $\Xi_0:=\Xi\setminus \Xi^0$ and let $\Cat_{\Xi_0}$ be the Serre subcategory
in $\Cat:=A_{\F}\operatorname{-mod}$ that is the kernel of the quotient functor defined by $P(\lambda)$
with $\varrho(\lambda)\in \Xi^0$. We write $\varrho_{\Xi_0}$ for the inclusion
$\Cat_{\Xi_0}\hookrightarrow \Cat$. The objects $\iota_{\Xi_0}^! P(\lambda'),
\lambda'\in \varrho^{-1}(\Xi_0)$, are projective in $\Cat_{\Xi_0}$.  By writing  projective resolutions of
 these objects in $\Cat$, we see  that they  have no higher self-extensions in $\Cat$.
So each $\iota_{\Xi_0}^! P(\lambda')$ uniquely deforms to an $\Ring$-flat object in $A_\Ring\operatorname{-mod}$, let $P_\Ring^{\Xi_0}(\lambda')$ denote this deformation. Then we have a natural epimorphism
$A_\Ring/A_\Ring e_{\Xi^0}A_\Ring\twoheadrightarrow \operatorname{End}_{\Cat_\Ring}(\bigoplus_{\lambda'}P_\Ring^{\Xi_0}(\lambda'))$
that is an isomorphism after specializing to $\F$. Since the target is flat over $\Ring$, we see that
this epimorphism is an isomorphism and the source is flat as well.
\end{proof}

By \cite[Lemma 2.4]{LW}, $\Cat_\F^{opp}$ is standardly stratified with the standard objects $\nabla_\F(\lambda)$.
This and Lemma \ref{Lem:SS_field} imply that $\Cat_\F^{opp}$ is standardly stratified with the standard objects $\nabla_\Ring(\lambda)$.


Let us introduce some more notation and terminology. We write $\nabla_\Ring(\lambda)$
for $\pi_{\xi}^*(I_{\Ring,\xi}(\lambda))$, where $I_{\Ring,\xi}(\lambda)$ stands for the indecomposable
relatively $\Ring$-injective object in $\Cat_{\Ring,\xi}$ indexed by $\lambda$, where $\xi=\varrho(\lambda)$.
The object $\nabla_\Ring(\lambda)$ is called costandard.  Further, let $L_{\F,\xi}(\lambda)$
denote the simple object in $\Cat_{\F,\xi}$ indexed by $\lambda$. We get
the proper standard object $\bar{\Delta}_\F(\lambda):=\pi_{\xi}^!(L_{\F,\xi}(\lambda))$
and the proper costandard object $\bar{\nabla}_{\F}(\lambda):=\pi_\xi^*(L_{\F,\xi}(\lambda))$.
Note that by \cite[Lemma 2.4]{LW} we have
\begin{equation}\label{eq:Ext_vanish}\dim \operatorname{Ext}^i(\Delta_\F(\lambda),\bar{\nabla}_\F(\lambda'))=\dim \operatorname{Ext}^i(\bar{\Delta}_\F(\lambda),\nabla_\F(\lambda'))=\delta_{i0}\delta_{\lambda\lambda'}.
\end{equation}
Consider the subcategories $\Cat_{\F}^{\Delta}\subset \Cat_{\F}^{\bar{\Delta}}\subset \Cat_{\F}$
of standardly and proper standardly filtered objects and, similarly, the subcategories
$\Cat_{\F}^{\nabla}\subset \Cat_{\F}^{\bar{\nabla}}\subset \Cat_{\F}$ of costandardly
and proper costandardly filtered objects. Note that $M\in \Cat_{\F}^{\bar{\nabla}}$
if and only if $\operatorname{Ext}^i(\Delta_\F(\lambda),M)=0$ for all $\lambda$,
\cite[Lemma 2.4]{LW}.

We also have the notions of standardly stratified categories for coideal finite and ideal
finite posets (we require that all fibers of $\varrho$ are finite). The definitions
mirror those for highest weight categories given in Section \ref{SS_HW}. We are not going to
provide them.

\subsubsection{Compatible standardly stratified structures on highest weight categories}
Let $\Ring$ be as before, $A_\Ring$ be an $\Ring$-algebra that is a free finite rank $\Ring$-module.
We write $\Cat_\Ring$ for $A_\Ring\operatorname{-mod}$ and let $\Lambda$ denote the indexing
set for the indecomposable projectives. Let $\leqslant,\preceq$ be two pre-orders
that give standardly stratified structures on $\Cat_\Ring$. We say that $\leqslant$
refines $\preceq$ if
\begin{itemize}
\item $\lambda\leqslant \lambda'$ implies $\lambda\preceq \lambda'$.
\item Every equivalence class for $\preceq$ is an interval for $\leqslant$.
\end{itemize}

Let $\xi$ be an equivalence class for $\preceq$. Note that each $\Cat_{\Ring,\xi}$ becomes
a standardly stratified category whose poset consists of the equivalence
classes for $\leqslant$ in $\xi$. Pick such an equivalence class $\vartheta$.
Let $\pi_{\vartheta}^\xi$ denote the quotient functor $\Cat_{\Ring,\xi,\leqslant \vartheta}\rightarrow
\Cat_{\Ring,\vartheta}$. Clearly, the functors $\pi_\vartheta$ and $\pi_{\vartheta}^\xi\circ \pi_\xi$
from $\Cat_{\Ring,\leqslant \vartheta}\twoheadrightarrow \Cat_{\Ring,\vartheta}$ are naturally
isomorphic. By adjointness, we get the following isomorphisms of functors.
\begin{equation}\label{eq:stand_compos}
\begin{split}
& \Delta_\vartheta\cong \Delta_\xi\circ \Delta^\xi_{\vartheta},\\
& \nabla_\vartheta\cong \nabla_\xi\circ \nabla_\vartheta^\xi.
\end{split}
\end{equation}

When $\leqslant$ defines a highest weight structure, we will say that the standardly
stratified structure defined by $\preceq$ is compatible with the highest weight structure.
So every subquotient category $\Cat_{\Ring,\xi}$ becomes highest weight with respect
to the order $\leqslant$ on $\varrho^{-1}(\xi)$ so that $\Delta_\Ring(\lambda)=\Delta_{\Ring,\xi}(\Delta_\Ring^\xi(\lambda))$,
where $\Delta_\Ring^\xi(\lambda)$ denote the standard object in $\Cat_{\Ring,\xi}$
indexed by $\lambda$.  In particular, $\Cat_\Ring^{\Delta,\preceq}
\subset \Cat_\Ring^{\Delta,\leqslant}$ and $\Cat_{\F}^{\Delta,\leqslant}\subset
\Cat_{\F}^{\bar{\Delta}}$.

The following lemma gives a criterium for a pre-order $\preceq$ refined by $\leqslant$
to give a compatible standardly stratified structure. The proof can be found in
\cite[Section 3.5]{CWR} in the case of the categories over a field, the general case follows
then from Lemma \ref{Lem:SS_field}.

\begin{Lem}\label{Lem:compat_SS_equiv}
The following are equivalent.
\begin{itemize}
\item[(a)] $\preceq$ gives a compatible standardly stratified structure.
\item[(b)] For any $\xi$ in $\Xi$ and any $\lambda\in \Lambda$, the object
$\pi_\xi\circ\iota^!_{\leqslant \xi}P_\Ring(\lambda)$ is projective in $\Cat_{\Ring,\xi}$
and the object $\pi_\xi\circ \iota^*_{\leqslant \xi}  I_\Ring(\lambda)$
is relatively $\Ring$-injective in $\Cat_{\Ring,\xi}$.
\item[(c)] For any $\xi$, the functors $\pi_\xi^!,\pi_{\xi}^*$ are exact.
\end{itemize}
\end{Lem}

Note that (b) still makes sense when $\Lambda=\Lambda^+$ is coideal finite and (a) and (b)
are still equivalent. In fact, a pre-order $\preceq$ gives a compatible standardly
stratified structure on $\Lambda^+$, when so does the restriction of $\preceq$
to any finite coideal that is the union of equivalence classes for $\preceq$.

\subsection{Categories $\,_I\mathcal{O}^+_{\tilde{\Ring}}(W)$}
In this section, $\F=\F_p$ is the finite field of $p$ elements (where $p$ is prime).
Let $p:=\operatorname{char}\F$. We set $\tilde{\Ring}:=\Z_p[[\h]]$.

Let $(W,S)$ be the affine type $A$ Coxeter system and $\h$ be such as
in Section \ref{SS_affine_prelim}. Our goal is to define analogs of
blocks of the BGG category $\mathcal{O}$ over $\Z_p$ and their deformations
over $\tilde{\Ring}$. Our approach traces back to Bernstein and Gelfand,
\cite{BG}, who related Harish-Chandra bimodules to the category $\mathcal{O}$
and then to Soergel, \cite{Soergel}, who used Soergel bimodules instead of
Harish-Chandra bimodules. Soergel worked over $\C$ with algebraic Soergel
bimodules and already defined category $\mathcal{O}$. In our situation,
we have a more mysterious category $\DG$ of diagrammatic Soergel bimodules
and use them to define the category $\mathcal{O}$ more or less following
Soergel's recipe (here and below we will abuse the notation and write $\DG$
for the Karoubi envelope of the category of singular Bott-Samelson bimodules).
We then use results from Section \ref{sec-cell}
to establish some basic properties of our category $\mathcal{O}$ that mirror
standard properties of the BGG category.

We will work in the deformed setting (working over $\tilde{\Ring}$ instead of $\Z_p$).
For the BGG categories $\OCat$ this will mean the following. Instead of dealing
with modules over a semisimple Lie algebra $\mathfrak{g}$ we will look at modules
over $\mathfrak{g}\otimes \C[x_i]_{i\in I}$ (here $I$ stands for the simple root system of
$\mathfrak{g}$) such that the simple coroot $\alpha_i^\vee$ acts diagonalizably with
eigenvalues in $\Z+x_i$.

We write $\DG=\bigoplus \DG(I,J)$ for the Karoubi envelope of
the monoidal 2-category defined in Section \ref{subsec-diagaffinetype}
(the summation is taken over all finitary subdiagrams $I,J$), where now we take $R=\Z_p[\h]$.
Let us recall some properties of $\DG$ that will be used below.
\begin{itemize}
\item[(S1)] Any object in $\DG(I,J)$ is a finitely generated graded $R^I$-$R^J$-bimodule, where we write $R=\Z_p[\h]$
and $R^I$ stand for the subalgebra of $W_I$-invariants in $R$.
\item[(S2)] The 1-morphisms are generated by $\mathfrak{P}^{I}_{J}\in \DG(J,I)$ (shifted restriction)
and $\mathfrak{P}_{I}^{J}\in \DG(I,J)$ (induction), where  $I\subsetneq J$ are finitary subsets of $S$.
\item[(S2$'$)] The 1-morphisms in $\DG(\varnothing,\varnothing)$ are generated
by $\mathfrak{P}^i_\varnothing\mathfrak{P}_i^\varnothing$ for $i\in S$.
\item[(S3)] We have a natural isomorphism $\mathfrak{P}^I_K\cong \mathfrak{P}^J_K\mathfrak{P}^I_J$
for $I\subsetneq J\subsetneq K$.
\item[(S4)] The 1-morphisms $\mathfrak{P}_I^J$ and $\mathfrak{P}^I_J$ are biadjoint.
\item[(S5)] $\DG$ categorifies the Hecke algebroid. In particular, the
ungraded (complexified) $K_0$ of $\DG(I,J)$ coincides with $e_I \C W e_J$, where $e_I,e_J$
denote the averaging idempotents in $W_I,W_J$, respectively. We have
$[\mathfrak{P}^I_J]=e_J$  and $[\mathfrak{P}_I^J]=\sum_{w\in W_J}we_I$,
for $I\subsetneq J$.
\end{itemize}
In fact, as long as we have these properties everything below in this section
related to categories $\mathcal{O}$ is going to work.

We will write $\mathfrak{P}_I$ for $\mathfrak{P}_I^\varnothing$ and $\mathfrak{P}_I^*$
for $\mathfrak{P}^I_{\varnothing}$.

%

\subsubsection{Bruhat posets}\label{SSS_Bruhat}
Recall the (inverse) Bruhat order on $W$: $w\leqslant w'$ if there are reflections $s_1,\ldots s_k$
such that $w'=s_ks_{k-1}\ldots s_1 w$ and $l(w)>l(s_1w)>l(s_2s_1w)>\ldots>l(w')$.
We remark that $1$ is the maximal element, while the minimal element exists if and only if
$W$ is finite. Also note that all coideals $\{w'\in W| w'\geqslant w\}$ are finite.
Set $\Lambda^+:=(W,\leqslant)$.

We need some quotient posets of $\Lambda^+$. First, let $I\subset S$ be a finitary subset.
Consider the left coset space $W_{I}\setminus W$.
The following lemma shows that the order $\leqslant$ descends to $W_{I}\setminus W$.

\begin{Lem}\label{Lem:Bruhat_preposet_left}
The $W_{I}$-cosets are intervals for $\leqslant$. Moreover, if $w_1^{-1} w_1', w_2^{-1}w_2'\in W_{I}$,
and $w_1\leqslant w_2, w_2'\leqslant w_1'$, then $w_1\sim_{I} w_2$.
\end{Lem}
\begin{proof}
Let $\mathfrak{t}$ denote the Cartan subalgebra of the Kac-Moody algebra with Weyl group $W$.
Pick an element $\rho^\vee\in \mathfrak{t}$ whose pairing with all simple roots equals $1$. The map
$W\mapsto W\rho$ is an embedding. The inverse Bruhat order on $W\rho$ can be described in terms of
$W\rho$ as follows: we have $\lambda>\lambda'$ if there are reflections $s_1,\ldots,s_k$
such that $\lambda'=s_k\ldots s_1 \lambda$ and $s_{i-1}\ldots s_1\lambda-s_is_{i-1}\ldots s_1\lambda$
is a positive multiple of a positive real root. Note that if $w'\in W_I$ and $w'\lambda<\lambda$,
then $\lambda-w'\lambda$ is a nonnegative linear combination of simple roots in $I$. So if
$w'w\rho<u\rho<w\rho$ with $w'\in W_I$, then $u\in W_I w$. This shows that $W_I w$ is an
interval.

Now let us show the second claim of this lemma. We may assume that $w_1<w_2$ and $w_2'<w_1'$.
The differences $w_2\rho-w_1\rho, w_1'\rho-w_2'\rho$ are positive linear combinations of
positive roots. The differences $w_1\rho-w_1'\rho, w_2\rho-w_2'\rho$ are linear combinations
of simple roots in $I$. This implies that at least one of the differences $w_2\rho-w_1\rho,
w_1'\rho-w_2'\rho$ is a linear combination of simple roots from $I$. This shows
that all four elements $w_1,w_2,w_1',w_2'$ are equivalent with respect to $\sim_I$.
\end{proof}

So we get the poset $\,_{I}\Lambda^+:=(W_I\setminus W, \leqslant)$.

Now let $J\subset S$ be another finitary subset.
We define a pre-order $\preceq=\preceq_{J}$ on $\,_{I}\Lambda^+$.
First, we define an equivalence relation $\sim_{J}$ on $W$ by setting $W_Iw\sim_{J}W_Iw'$
if $w'\in W_IwW_J$. The following lemma shows that
$\leqslant$ refines $\preceq_J$ and that  $\leqslant$ descends to
a partial order on the set of equivalence classes for $\sim_{J}$.

\begin{Lem}\label{Lem:Bruhat_preposet_right}
The equivalence classes in $\,_{I}\Lambda^+$
are intervals for $\leqslant$. Moreover, if $\lambda_1\sim_{J} \lambda_1', \lambda_2\sim_{J}\lambda_2'$,
and $\lambda_1\leqslant \lambda_2, \lambda_2'\leqslant \lambda_1'$, then $\lambda_1\sim_{J} \lambda_2$.
\end{Lem}
\begin{proof}
This is proved similarly to Lemma \ref{Lem:Bruhat_preposet_left} by embedding $W/W_J$
to $\mathfrak{t}$ by $\lambda\mapsto \lambda \rho_J$, where $\rho_J\in \mathfrak{t}$ is an element equal
$1$ on the simple roots from $S\setminus J$ and vanishing of the simple roots from $J$.
\end{proof}

So we get a pre-poset $(\,_{I}\Lambda^+, \preceq_{J})$ and the induced
poset $(\,_{I}\Lambda^+_J,\leqslant)$, where $\,_{I}\Lambda^+_J:=
W_I\setminus W/W_J$.

\subsubsection{Categories $\,_I\mathcal{O}_{\tilde{\Ring}}^+(W)$: construction}\label{SSS_Soergel_O_constr}
For $\lambda\in\,_{I}\Lambda^+$, let $\mathfrak{B}_\lambda$ be the corresponding indecomposable
object in $\DG(I,\varnothing)$.

Consider the idempotented algebra $A_R:=\left(\bigoplus_{\lambda,\mu\in \,_{I}\Lambda^+}\Hom_{\DG}(\mathfrak{B}_\lambda,\mathfrak{B}_\mu)\right)^{opp}$. Note that
$A_R$ is an $R^I\otimes R$-algebra. So we get the category $A_R\operatorname{-mod}$
as in \S\ref{SS_hw_coid_fin}.

\begin{Lem}\label{Lem:Alg_A_freeness}
We have the following.
\begin{enumerate}
\item
The algebra $A_R$ is a projective module over $R$.
\item The category of the projectives in  $A_R\operatorname{-mod}$ is
naturally identified with $\DG(I,\varnothing)$.
\end{enumerate}
\end{Lem}
\begin{proof}
In the case when $I=\varnothing$, the $R$-module $\Hom_{\DG}(\mathfrak{B}_\lambda,\mathfrak{B}_\mu)$
has a basis produced in
Section \ref{subsec-lightleaves}. For an arbitrary finitary $I$ we can argue as
follows. We need to prove that $\Hom_{\DG}(B_1,B_2)$ is free over $R$ for $B_1,B_2\in \DG(I,\varnothing)$.
Any object in $\DG(I,\varnothing)$ is a direct summand of an object of the form $\mathfrak{P}_I B$
for $B\in \DG(\varnothing,\varnothing)$. But $\Hom_{\DG}(\mathfrak{P}_I B, B_2)=\Hom_{\DG}(B, \mathfrak{P}_I^* B_2)$
as right $R$-modules.

(2) is straightforward.
\end{proof}

Abusing the notation we will write $\mathfrak{B}_\lambda$ for the indecomposable projective
$A_R$-module corresponding to $\mathfrak{B}_\lambda$.

Consider the algebra $A_{\tilde{\Ring}}:=A_R\otimes_R\tilde{\Ring}$.

\begin{Lem}\label{Lem:A_indec_proj}
The indecomposable projective $A_{\tilde{\Ring}}$-modules are precisely $\mathfrak{B}_\lambda\otimes_R\tilde{\Ring}$.
\end{Lem}
\begin{proof}
Note that $B_R:=\operatorname{End}_{\DG}(\mathfrak{B}_\lambda)$ is a finitely generated $R$-algebra
and $\operatorname{End}_{A_{\tilde{\Ring}}}(\mathfrak{B}_\lambda\otimes_R\tilde{\Ring})=
B_R\otimes_R\tilde{\Ring}$.
This algebra is indecomposable if and only $B_R/(\h,p)B_R$  is. The algebra $B_R/(\h,p)B_R$ is graded.
So the radical and also the center of the quotient by the radical are graded.
The center is the sum of some field extensions of $\F_q$ and so contains no nilpotent elements.
It follows that it is concentrated in degree $0$. A nontrivial idempotent can be lifted to
a degree $0$ idempotent in  $B_R\otimes_R\tilde{\Ring}$. But the homogeneous elements in $B_R\otimes_R\tilde{\Ring}$
lie in $B_R$. So we get a non-trivial idempotent in $B_R$, a contradiction.
\end{proof}

Set $\,_I\mathcal{O}_{\tilde{\Ring}}^+(W):=A_{\tilde{\Ring}}\operatorname{-mod}$.


%

\subsubsection{Main properties}
Here we will list basic properties of the category $\bigoplus_I\, _{I}\OCat^+_{\tilde{\Ring}}(W)$, where the summation
is taken over all finitary  subsets $I\subset S$.

We will need three groups of  properties from $\bigoplus_{I} \,_I\mathcal{O}^+_{\tilde{\Ring}}(W)$:
highest weight structure,   categorical actions and  filtrations.

{\it Highest weight structure}.  We require the following two properties.
\begin{itemize}
\item[(O1.1)] $\,_\varnothing\mathcal{O}^+_{\tilde{\Ring}}(W)$
is a highest weight category over $\tilde{\Ring}$ with poset $\,_{\varnothing}\Lambda^+=W$.
\item[(O1.2)] $\,_I\mathcal{O}^+_{\tilde{\Ring}}(W)$ has indecomposable projectives indexed by
$W_I\setminus W$ and satisfies the first three conditions from \S\ref{SS_hw_coid_fin}.
\end{itemize}

{\it Categorical action}. We require that
\begin{itemize}
\item[(O2.1)]$\bigoplus_{I} \,_I\mathcal{O}^+_{\tilde{\Ring}}(W)$
is an $\tilde{\Ring}$-linear 1-representation of $\DG$.
\item[(O2.2)] The  complexified $K_0$ of $\bigoplus_{I} \,_I\mathcal{O}^+_{\tilde{\Ring}}(W)\operatorname{-proj}$
equals to  $\bigoplus_{J}\C (W_J\setminus W)$ as a module over $K_0^{\C}(\DG):=\bigoplus_{I,J}\C(W_I\setminus W/W_J)$.
\item[(O2.3)] The standard objects in $\,_\varnothing\mathcal{O}^+_{\tilde{\Ring}}(W)$ correspond to the standard basis
in $\C W$. The action of $\DG(\varnothing,\varnothing)$ preserves
the full subcategory  of standardly filtered objects $\,_\varnothing\mathcal{O}^+_{\tilde{\Ring}}(W)^\Delta
\subset \,_\varnothing\mathcal{O}^+_{\tilde{\Ring}}(W)$.
\end{itemize}

Note that (O2.2) implies that the map $K_0(\,_\varnothing\mathcal{O}^+_{\tilde{\Ring}}(W))\operatorname{-proj})
\rightarrow K_0(\,_I\mathcal{O}^+_{\tilde{\Ring}}(W)\operatorname{-proj})$ given by ``averaging'' indecomposable object in $\mathfrak{P}_I\in \DG(I,\varnothing)$ coincides with the natural projection $\C W\twoheadrightarrow \C (W_I\setminus W)$.

Let $C\in W/W_J$. Then we can consider the highest weight quotient $\,_\varnothing\mathcal{O}^+_{\tilde{\Ring}}(W)_{\geqslant C}$ of
$\,_\varnothing\mathcal{O}^+_{\tilde{\Ring}}(W)$ corresponding to the poset coideal
$\{w\in W| wW_J\succeq_J C\}$. We have the highest weight subcategory
$\,_\varnothing\mathcal{O}^+_{\tilde{\Ring}}(W)_{C}\subset \,_\varnothing\mathcal{O}^+_{\tilde{\Ring}}(W)_{\geqslant C}$ corresponding
to the poset ideal $C$. We have the quotient functor $\iota_C^!:
\,_\varnothing\mathcal{O}^+_{\tilde{\Ring}}(W)^\Delta_{\geqslant C}
\rightarrow \,_\varnothing\mathcal{O}^+_{\tilde{\Ring}}(W)^\Delta_{C}$. Similarly,
for $C'\in W_J\setminus W$, we can form the subquotient category
$\,_\varnothing\mathcal{O}^+_{\tilde{\Ring}}(W)_{C'}$.


{\it Filtration}. Pick a finitary subset $J\subset S$. Then we have the following:
\begin{itemize}
\item[(O3.1)] For each   $C\in W/W_J$,  the subquotient
$\,_\varnothing\mathcal{O}^+_{\tilde{\Ring}}(W)_C$ is equivalent to $\OCat^+_{\tilde{\Ring}}(W_J)$,
the base change of $\OCat^+_{\tilde{\Ring}_J}(W_J)$ from $\tilde{\Ring}_J$ to $\tilde{\Ring}$.
Here we write $\tilde{\Ring}_J$ for the analog of $\tilde{\Ring}$ for $W_J$.
\item[(O3.2)] Let $w$ be the shortest element in $C$ and let $\mathfrak{B}_w\in \DG(\varnothing,\varnothing)$
be the indecomposable Soergel bimodule corresponding to $w$. Then the functor
$\iota_{C}^!\circ \mathfrak{B}_w$ is an equivalence  $\OCat^+_{\tilde{\Ring}}(W)^\Delta_{W_J}
\xrightarrow{\sim} \OCat^+_{\tilde{\Ring}}(W)^\Delta_{C}$.
\item[(O3.3)] Let $C'\in W_J\setminus W$. Then there is an equivalence
$\,_\varnothing\OCat^+_{\F}(W)_{C'}\xrightarrow{\sim}\,_\varnothing\OCat^+_{\F}(W_J)$ that
sends $\Delta_{\F}(uw)$ to $\Delta_{\F}(u)$, where $u\in W_J$ and $w$
is the shortest element in $C'$.
\end{itemize}


\begin{Prop}\label{Prop:Cat_O_prop_check}
The category $\bigoplus_I\, _{I}\OCat^+_{\tilde{\Ring}}(W)$ has the properties
stated above.
\end{Prop}
\begin{proof}
Recall cellular structure on $\DG(\varnothing,\varnothing)$ from Section \ref{subsec-lightleaves}.
From Sections \ref{SS_cell_hw} and \ref{subsec-cellsubquot} it follows that $\,_\varnothing\OCat_{\F}^+(W)$
is a highest weight category with coideal finite poset. By a coideal finite analog of Lemma
\ref{Lem:hw_field}, we get (O1.1).

(O1.2) is a direct  consequence of the construction of $\bigoplus_I \,_I\OCat_{\tilde{\Ring}}^+(W)$.

The category $\DG$ acts on $\bigoplus_I\, _{I}\OCat^+_{\tilde{\Ring}}(W)\operatorname{-proj}$
by the construction of the latter. The action is by bi-adjoint functors (by (S2) combined with (S4))
and hence  extends to an action on $\bigoplus_I\, _{I}\OCat^+_{\tilde{\Ring}}(W)$. This shows (O2.1).

(O2.2) follows from (S5).

The standard objects in $\,_\varnothing \OCat^+_{\tilde{\Ring}}(W)$ come from the standard
bimodules, see Section \ref{subsec-lightleaves}.
Now Proposition \ref{Prop:stand_action_preserv} implies (O2.3).

Let us proceed to (O3.1)-(O3.3). The first two properties follow from
Proposition \ref{Prop:quot_isom}. 
Property (O3.3) follows from Proposition \ref{Prop:quot_iso_left}.
\end{proof}

\begin{Lem}\label{Lem:act_costrand_preserv}
$\DG(\varnothing,\varnothing)$ preserves
the full subcategory  of costandardly filtered objects $\,_\varnothing\mathcal{O}^+_{\tilde{\Ring}}(W)^\nabla
\subset \,_\varnothing\mathcal{O}^+_{\tilde{\Ring}}(W)$.
\end{Lem}
\begin{proof}
Note that since $\DG(\varnothing,\varnothing)$ acts by biadjoint functors and preserves the subcategory
of projective objects, it also preserves the subcategory of flat objects. Also we know that it preserves
the subcategory of standardly filtered objects. Again, since it acts by biadjoint functors, it preserves
the subcategory of costandardly filtered objects.
\end{proof}


Finally, let us record a standard corollary of (O2.2).

\begin{Lem}\label{Lem:SB_action}
Every indecomposable projective in $\bigoplus_I\,_I \OCat_{\tilde{\Ring}}(W)$
is a direct summand in $B.\Delta_{\tilde{\Ring},\varnothing}(1)$, where $B\in \DG$.
\end{Lem}
\begin{proof}
The action of $\DG$ on $\bigoplus_I\,_I \OCat_{\tilde{\Ring}}(W)$ maps projectives to projectives.
The object $\Delta_{\tilde{\Ring},\varnothing}(1)$ is projective. Since the objects of the form
$B.\Delta_{\tilde{\Ring},\varnothing}(1)$ generated $K_0$ of the exact category
$\bigoplus_{I}\,_I\OCat^+_{\tilde{\Ring}}(W)\operatorname{-proj}$, our claim follows.
\end{proof}

\subsubsection{Special case: type $A_1$}\label{SSS_A1_Soergel}
Before moving further, let us describe what kind of categories we get in type $A_1$,
here one can easily check that the categories are specified uniquely by our axioms.

Here we work over the base ring $\tilde{\Ring}:=\Z_p[[x]]$. We will have two categories,
the regular category $\,_{\varnothing}\OCat_{\tilde{\Ring}}(W)$ and the singular category $\,_{S}\OCat_{\tilde{\Ring}}(W)$.
The latter is just $\tilde{\Ring}\operatorname{-mod}$. The category $\,_{\varnothing}\OCat_{\tilde{\Ring}}(W)$
is the category of modules over the quotient of the path algebra $\tilde{\Ring}Q/(ab=x)$,
where $Q$ is the quiver with two vertices $1,2$ and two arrows $a,b$ with $t(a)=h(b)=1, t(b)=h(a)=2$.
It is a classical fact that we get a highest weight category here (whose specialization to $\C$
is the principal block of the category $\mathcal{O}$ for $\mathfrak{sl}_2$).

\subsection{Highest weight structures for $\, _{I}\OCat^+_{\tilde{\Ring}}(W)$}
We have defined all blocks of our category $\mathcal{O}$ and have equipped the principal
block with a highest weight structure. Now we are going to equip singular blocks with
highest weight structures. Our approach here is again inspired by the BGG categories
$\mathcal{O}$. There the standard objects in the singular blocks are obtained from the
standards in a regular block by applying translation functors to walls.
We have such functors in our setting, they are given by bimodules $\mathfrak{P}_I\in \DG(I,\varnothing)$.

So we are going to prove the following result.

\begin{Prop}\label{Prop:sing_hw}
There is a highest weight structure on $\,_I\mathcal{O}^+_{\tilde{\Ring}}(W)$
over $\tilde{\Ring}$ with poset $\,_I\Lambda^+=W_I\setminus W$ and standard objects $\mathfrak{P}_I\Delta_{\tilde{\Ring}}(w)$,
where $w$ is the largest element in its $W_I$-coset.  The action of $\DG$
preserves the subcategory of standardly filtered objects. The costandard objects
are $\mathfrak{P}_I\nabla_{\tilde{\Ring}}(w)$.
\end{Prop}

The proof will be given after a  lemma.

\begin{Lem}\label{Lem:Verma_coinc}
For any finitary $I\subset S$, and any $w,w'$ with $w'\in W_Iw$,
we have $\mathfrak{P}_I\Delta_{\tilde{\Ring}}(w)\cong \mathfrak{P}_I\Delta_{\tilde{\Ring}}(w')$.
\end{Lem}
\begin{proof}
The proof is in several steps.

{\it Step 1}.
Let us reduce the proof to showing that $\mathfrak{P}_i\Delta_{\tilde{\Ring}}(s_i w)=\mathfrak{P}_i\Delta_{\tilde{\Ring}}(w)$. Suppose that we know that. Recall that
by (S2), for any $i\in I$, $\mathfrak{P}_I$ factorizes as $\mathfrak{P}_I^i \mathfrak{P}_i$.
So $\mathfrak{P}_i\Delta_{\tilde{\Ring}}(s_i w)=\mathfrak{P}_i\Delta_{\tilde{\Ring}}(w)$
implies $\mathfrak{P}_I\Delta_{\tilde{\Ring}}(s_i w)=\mathfrak{P}_I\Delta_{\tilde{\Ring}}(w)$.
From here we deduce that $\mathfrak{P}_I\Delta_{\tilde{\Ring}}(w')=\mathfrak{P}_I\Delta_{\tilde{\Ring}}(w)$
as long as $W_Iw'=W_Iw$.

{\it Step 2}. Pick $w\in W$ with $l(s_i w)>l(w)$ so that $w\geqslant s_i w$.
Form the highest weight subquotient of $\,_\varnothing\mathcal{O}^+_{\F}(W)$
corresponding to the interval $C':=\{s_i w,w\}$. Namely, consider the highest weight quotients
$\,_\varnothing\mathcal{O}^+_{\F}(W)_{\geqslant C'}\twoheadrightarrow \,_\varnothing\mathcal{O}^+_\F(W)_{>C'}$
and the corresponding kernel $\,_\varnothing\mathcal{O}^+_\F(W)_{C'}$. By (O3.3),  $\,_\varnothing\mathcal{O}^+_\F(W)_{C'}\cong \,_\varnothing\mathcal{O}^+_{\F}(W_i)$.
The right hand side was described in \S\ref{SSS_A1_Soergel}. In particular, that description
implies the existence of a nonzero homomorphism $\Delta_{\F}(s_i)\rightarrow \Delta_\F(1)$.

{\it Step 3}. The subcategory of objects in $\,_\varnothing\mathcal{O}^+_{\F}(W)$
filtered by $\Delta_{\F}(s_iw),\Delta_{\F}(w)$ is equivalent to $\,_\varnothing\mathcal{O}^+_{\F}(W_i)^{\Delta}$,
by (O3.3). So we get a nonzero homomorphism $\varphi:\Delta_{\F}(s_i w)\rightarrow
\Delta_{\F}(w)$. Suppose that we know  that the induced homomorphism
$\mathfrak{P}_i.\Delta_{\F}(s_i w)\rightarrow
\mathfrak{P}_i.\Delta_{\F}(w)$ is an isomorphism. We claim that $\mathfrak{P}_i\Delta_{\tilde{\Ring}}(w)\cong \mathfrak{P}_i.\Delta_{\tilde{\Ring}}(w')$.
Indeed, $$\operatorname{Ext}^1(\mathfrak{P}_i.\Delta_{\F}(w), \mathfrak{P}_i.\Delta_{\F}(w))=
\operatorname{Ext}^1(\mathfrak{P}_i^*\mathfrak{P}_i.\Delta_{\F}(w),\Delta_{\F}(w))=0$$
(the last equality holds because $\mathfrak{P}_i^*\mathfrak{P}_i.\Delta_{\F}(w)$ is a projective
object in $\,_\varnothing\mathcal{O}^+_{\F}(W)_{C'}^\Delta$)
so $\mathfrak{P}_i.\Delta_{\F}(w)$ admits  at most one flat deformation over $\tilde{\Ring}$.
Since both $\mathfrak{P}_i.\Delta_{\tilde{\Ring}}(s_i w),\mathfrak{P}_i.\Delta_{\tilde{\Ring}}(w)$
are flat over $\tilde{\Ring}$, we will be done.
So it remains to prove $\mathfrak{P}_i.\Delta_{\F}(s_i w)\xrightarrow{\sim}
\mathfrak{P}_i.\Delta_{\F}(w)$.

{\it Step 4}. We claim $\mathfrak{P}_i^*\mathfrak{P}_i\Delta_{\F}(s_iw)\cong
\mathfrak{P}_i^*\mathfrak{P}_i\Delta_{\F}(w)$. It is easy to see that both objects lie
in $\,_\varnothing\mathcal{O}^+_{\F}(W)_{C'}^\Delta$. Note that
$\mathfrak{P}_i^*\mathfrak{P}_i\Delta_{\F}(s_i)\cong
\mathfrak{P}_i^*\mathfrak{P}_i\Delta_{\F}(1)$ in $\,_\varnothing\mathcal{O}^+_{\F}(W_i)$.
It follows from (O3.3) that
$\mathfrak{P}_i^*\mathfrak{P}_i\Delta_{\F}(s_iw)\cong
\mathfrak{P}_i^*\mathfrak{P}_i\Delta_{\F}(w)$.

{\it Step 5}. We have $\Delta_\F(w)\hookrightarrow \mathfrak{P}_i^*\mathfrak{P}_i\Delta_{\F}(w)$.
Step 4 gives an inclusion $\psi:\Delta_\F(w)\hookrightarrow \mathfrak{P}_i^*\mathfrak{P}_i\Delta_{\F}(s_iw)$
and hence, by adjunction, a homomorphism  $\psi':\mathfrak{P}_i\Delta_\F(w)\rightarrow \mathfrak{P}_i\Delta_{\F}(s_iw)$.
The  homomorphism $\psi'\circ(\mathfrak{P}_i\varphi):\mathfrak{P}_i\Delta_{\F}(s_iw)\rightarrow \mathfrak{P}_i\Delta_{\F}(s_iw)$
corresponds via the adjunction to the homomorphism $\psi\circ \varphi:\Delta_{\F}(s_iw)\rightarrow \Delta_{\F}(w)
\hookrightarrow \mathfrak{P}_i^*\mathfrak{P}_i\Delta_{\F}(s_iw)$. So $\psi'\circ(\mathfrak{P}_i\varphi)\neq 0$.
But $$\operatorname{End}(\mathfrak{P}_i\Delta_{\F}(s_iw))=\operatorname{Hom}(\mathfrak{P}_i^*\mathfrak{P}_i\Delta_{\F}(s_iw),
\Delta_{\F}(s_iw)).$$
The latter space is one-dimensional because of the short exact sequence
$0\rightarrow \Delta_{\F}(w)\rightarrow \mathfrak{P}_i^*\mathfrak{P}_i\Delta_{\F}(s_iw)
\rightarrow \Delta_{\F}(s_i w)\rightarrow 0$. So
$\psi'\circ(\mathfrak{P}_i\varphi):\mathfrak{P}_i\Delta_{\F}(s_iw)\rightarrow \mathfrak{P}_i\Delta_{\F}(s_iw)$
is an isomorphism. It follows that $\mathfrak{P}_i\Delta_{\F}(s_i w)$ is a direct summand in
$\mathfrak{P}_i\Delta_{\F}(w)$. But the endomorphism algebra of the latter is also one-dimensional
similarly to the argument above in this step. We conclude that $\mathfrak{P}_i\Delta_{\F}(s_i w)=\mathfrak{P}_i\Delta_{\F}(w)$, which finishes the proof.
\end{proof}

\begin{proof}[Proof of Proposition \ref{Prop:sing_hw}]
Let us show that the objects $\mathfrak{P}_I\Delta_{\tilde{\Ring}}(w)$ satisfy (HW1).
Indeed, we have $\operatorname{Hom}(\mathfrak{P}_I\Delta_{\tilde{\Ring}}(w),
\mathfrak{P}_I\Delta_{\tilde{\Ring}}(w'))=
\operatorname{Hom}(\Delta_{\tilde{\Ring}}(w),
\mathfrak{P}_I^*\mathfrak{P}_I\Delta_{\tilde{\Ring}}(w'))$. If the latter
is nonzero, then $W_I w\leqslant W_I w'$,  it is a free rank one module over $\tilde{\Ring}$ else
(because $\mathfrak{P}_I^*\mathfrak{P}_I\Delta_{\tilde{\Ring}}(w)$
is filtered with $\Delta_{\tilde{\Ring}}(w)$ as a subobject and all
other filtration components having labels less than $w$).

Now let us check (HW2). Note that the projective
$\mathfrak{P}_I P_{\tilde{\Ring}}(w)$ has a required filtration,
thanks to Lemma \ref{Lem:Verma_coinc}. By (O2.2),
these projectives generate the $K_0$ group of
$\,_I\OCat^+_{\tilde{\Ring}}(W)$.  So every indecomposable projective in
$\,_I\mathcal{O}^+_{\tilde{\Ring}}(W)$ is a direct summand
in some $\mathfrak{P}_I P_{\tilde{\Ring}}(w)$. Applying Lemma \ref{Lem:dir_sum_filtration}
to $M_1\oplus M_2=\mathfrak{P}_I P_{\tilde{\Ring}}(w)$, where $M_1$ is the indecomposable projective
labelled by $W_Iw$,
we see that (HW2) holds.

Note that every indecomposable object in $\DG(I,J)$ is a summand in an object of the form
$\mathfrak{P}_I B\mathfrak{P}_J^*$. Indeed, the objects of the latter form generate
$K_0(\DG(I,J))$ (at least over the rationals).
Now the claim that the $\DG$-action on $\bigoplus_I\,_I\mathcal{O}^+_{\tilde{\Ring}}(W)$
preserves the subcategory of standardly filtered objects follows
from (O2.3).

The objects $\mathfrak{P}_I\nabla_{\tilde{\Ring}}(w)$ are flat over $\Ring$. It is easy to
see that $\operatorname{Ext}^i(\mathfrak{P}_I\Delta_{\tilde{\Ring}}(w), \mathfrak{P}_I\nabla_{\tilde{\Ring}}(w'))=0$
for $i>0$ and that $\operatorname{Hom}(\mathfrak{P}_I\Delta_{\tilde{\Ring}}(w), \mathfrak{P}_I\nabla_{\tilde{\Ring}}(w'))=\delta_{W_Iw,W_I w'}$. So these objects
are indeed costandard.
\end{proof}

We will also need a semisimplicity result for localizations of categories
$\,_{I}\mathcal{O}_{\Ring}^+(W)$.

\begin{Lem}\label{Lem:O_ss}
Let $Q$ be a field with a homomorphism from $\Ring$ such that the images of all real roots
are invertible. Then the base change $\,_{I}\mathcal{O}_{Q}^+(W)$ is semisimple.
\end{Lem}
\begin{proof}
For $I=\varnothing$, this follows from Theorem \ref{Thm:Soergel_base_change}.
For general $I$, let us notice that the standard and the costandard objects in
$\,_{I}\mathcal{O}_{Q}^+(W)$ coincide.
\end{proof}

\subsection{Standardly stratified structures on $\, _{\varnothing}\OCat^+_{\tilde{\Ring}}(W)$}
We start to prepare to construct parabolic subcategories in our categories $\mathcal{O}$.
For this we will need to produce standardly stratified structures on the principal block
$\, _{\varnothing}\OCat^+_{\tilde{\Ring}}(W)$.

Let us explain how such structures appear for the classical BGG categories.
Let $\mathfrak{g}$ be a semisimple Lie algebra over $\C$ and let $\OCat(\mathfrak{g})$ be the integral
part of the whole BGG category (so that the simples in $\OCat(\mathfrak{g})$ are labelled
by the weight lattice of $\mathfrak{g}$). Now pick a standard Levi subalgebra $\mathfrak{l}\subset \mathfrak{g}$
corresponding to a subset $J$ in the simple root system of $\mathfrak{g}$. Let $\OCat(\mathfrak{l})$
stand for the category $\mathcal{O}$ for $\mathfrak{l}$, the simples are again labelled
by the weight lattice  of $\mathfrak{g}$. We have two functors $\OCat(\mathfrak{l})\rightarrow
\OCat(\mathfrak{g})$. First, there is the induction functor $\Delta_{\mathfrak{l}}:\OCat(\mathfrak{l})\rightarrow
\OCat(\mathfrak{g})$ that sends $M$ to $U(\mathfrak{g})\otimes_{U(\mathfrak{p})}M$. Here, as usual, we write
$\mathfrak{p}$ for the standard parabolic subalgebra with Levi subalgebra $\mathfrak{l}$,
the algebra $\mathfrak{p}$ acts on $M$ via the epimorphism $\mathfrak{p}
\twoheadrightarrow \mathfrak{l}$. Next, there is the coinduction functor $\nabla_{\mathfrak{l}}:
\OCat(\mathfrak{l})\rightarrow \OCat(\mathfrak{g}), M\mapsto \Hom_{U(\mathfrak{p})}(U(\mathfrak{g}),M)$.
Both functors are exact. It is easy to see that there is a standardly stratified structure
on $\mathcal{O}(\mathfrak{g})$ with $\gr \mathcal{O}(\mathfrak{g})=\mathcal{O}(\mathfrak{l})$ where $\Delta_{\mathfrak{l}}$
and $\nabla_{\mathfrak{l}}$ are the standardization and costandardization functors.
The corresponding pre-order is obtained from the standard order on the weight lattice by declaring that two weights
whose difference lies in the root lattice of $\mathfrak{l}$ are equivalent.

 We are going to emulate the constructions in the previous paragraph for the
 category $\, _{\varnothing}\OCat^+_{\tilde{\Ring}}(W)$ (with suitable modifications).
Pick a finitary subset $J\subset S$.
We want to prove that the Bruhat order on $W/W_J$ defines a standardly
stratified structure on $\,_\varnothing\OCat_{\tilde{\Ring}}(W)$ that is compatible
with the highest weight structure.

Let $C$ be a coset in $W/W_J$
and let $w_0\in C$ be the shortest element. For $w\in C$, consider the projective object
$P_{C,\tilde{\Ring}}(w)$ in $\OCat_{\tilde{\Ring}}(W)_C$. Under the identification
of  $\OCat_{\tilde{\Ring}}(W)_C$ with $\OCat_{\tilde{\Ring}}(W_J)$ from (O3.2),
the object $P_{C,\tilde{\Ring}}(w)$ corresponds to
$\tilde{\Ring}\otimes_{\tilde{\Ring}_J}P_{\tilde{\Ring}_J}(ww_0^{-1})$. We set
$\Delta_{\tilde{\Ring},J}(w):=\pi_C^!(P_{C,\tilde{\Ring}}(w))$. We can define
the object $\nabla_{\tilde{\Ring},J}(w)$ in a similar fashion, as $\pi_C^*(I_{C,\tilde{\Ring}}(w))$.

\begin{Prop}\label{Prop:filtr_preserv}
The action of $\DG(\varnothing,\varnothing)$ on $\OCat_{\tilde{\Ring}}(W)$
preserves the full subcategory of modules that admit a filtration
by $\Delta_{\tilde{\Ring},J}(w)$'s.
\end{Prop}
\begin{proof}
Thanks to (S2$'$), it is enough to prove that $\mathfrak{P}_i^*\mathfrak{P}_i$ maps
$\Delta_{\tilde{\Ring},J}(w)$ to an object admitting a filtration
by $\Delta_{\tilde{\Ring},J}(?)$'s. We have three cases to consider:
$s_iC\succ_JC, s_iC=C, s_iC\prec_J C$.

{\it Case 1: $s_i C\prec_J C$}. We claim that $$0\rightarrow \Delta_{\tilde{\Ring},J}(w)\rightarrow\mathfrak{P}_i^*
\mathfrak{P}_i\Delta_{\tilde{\Ring},J}(w)\rightarrow \Delta_{\tilde{\Ring},J}(s_iw)\rightarrow 0.$$

Let $w=w'w''$ be the decomposition, where $w''\in W_J$ and $w'$ is shortest in $wW_J$.
Since $s_i C\prec_J C$, we see that $s_i w'$ is shortest in $s_iwW_J$, so $s_i w=(s_iw')w''$
is the similar decomposition for $s_i w$. Condition (O3.2) applied to  $C= w' W_J$ implies
that   $\mathfrak{B}_{w'}P_{\tilde{\Ring}}(w'')\twoheadrightarrow \Delta_{\tilde{\Ring},J}(w)$
with kernel filtered by standards whose labels $w_1$ satisfy $w_1 W_J\succ_J w'W_J$.
Similarly, (O3.2) applied to $s_i C=(s_iw')W_J$ implies that $\mathfrak{P}_i^*\mathfrak{P}_i\mathfrak{B}_{w'}
P_{\tilde{\Ring}}(w'')\twoheadrightarrow \Delta_{\tilde{\Ring},J}(s_iw)$ with kernel
filtered by standards whose labels $w_2$ satisfy $w_2 W_J\succ_J s_iw'W_J$. Note that $w_1 W_J, s_i w_1 W_J\succ_J
s_i w' W_J$ because $s_i C\prec_J C$. Thanks to this we get an epimorphism
$\mathfrak{P}_i^*\mathfrak{P}_i\Delta_{\tilde{\Ring},J}(w)\twoheadrightarrow \Delta_{\tilde{\Ring},J}(s_iw)$.
By (O2.2) its kernel
has the same standard subquotients (with multiplicities) as $\Delta_{\tilde{\Ring},J}(w)$.
It remains to show that the kernel is actually isomorphic to $\Delta_{\tilde{\Ring},J}(w)$.

We have an adjunction co-unit morphism
$\Delta_{\tilde{\Ring},J}(w)\rightarrow \mathfrak{P}_i^*
\mathfrak{P}_i\Delta_{\tilde{\Ring},J}(w)$. Note that a filtration on
$\Delta_{\tilde{\Ring},J}(w)$ with standard quotients $\Delta_{\tilde{\Ring}}(u)$ induces a filtration on
$\mathfrak{P}_i^*\mathfrak{P}_i\Delta_{\tilde{\Ring},J}(w)$ with quotients
$\mathfrak{P}_i^*\mathfrak{P}_i\Delta_{\tilde{\Ring}}(u)$. The induced homomorphism
$\Delta_{\tilde{\Ring}}(u)\rightarrow\mathfrak{P}_i^*\mathfrak{P}_i\Delta_{\tilde{\Ring}}(u)$
has image in $\Delta_{\tilde{\Ring}}(u)$ and hence is injective. We deduce that
$\Delta_{\tilde{\Ring},J}(w)\rightarrow \mathfrak{P}_i^*
\mathfrak{P}_i\Delta_{\tilde{\Ring},J}(w)$ is injective. The same is true after specializing
to the closed point. The image of $\Delta_{\F,J}(w)$ lies in the kernel of
$\mathfrak{P}_i^*\mathfrak{P}_i\Delta_{\F,J}(w)
\twoheadrightarrow \Delta_{\F,J}(s_iw)$, for example, by the ordering considerations. Since $[\mathfrak{P}_i^*\mathfrak{P}_i\Delta_{\F,J}(w)]=
[\Delta_{\F,J}(w)]+[\Delta_{\F,J}(s_iw)]$, we get
$$0\rightarrow \Delta_{\F,J}(w)\rightarrow\mathfrak{P}_i^*
\mathfrak{P}_i\Delta_{\F,J}(w)\rightarrow \Delta_{\F,J}(s_iw)\rightarrow 0.$$
This yields an exact sequence in the beginning of the step.

{\it Case 2: $s_i C\succ_J C$}. Note that, by the previous case,
$$0\rightarrow \Delta_{\tilde{\Ring},J}(s_iw)\rightarrow\mathfrak{P}_i^*
\mathfrak{P}_i\Delta_{\tilde{\Ring},J}(s_iw)\rightarrow \Delta_{\tilde{\Ring},J}(w)\rightarrow 0.$$
By applying an adjunction to $\Delta_{\tilde{\Ring},J}(s_iw)\rightarrow\mathfrak{P}_i^*
\mathfrak{P}_i\Delta_{\tilde{\Ring},J}(s_iw)$ we get a morphism
 $\mathfrak{P}_i\Delta_{\tilde{\Ring},J}(s_iw)\rightarrow \mathfrak{P}_i\Delta_{\tilde{\Ring},J}(w)$.
Now consider a filtration of $P_{\tilde{\Ring}_J}(w'')$ (where $w''$ has the same meaning as in Step 1)
with standard quotients. These gives filtrations $F_k\subset \Delta_{\tilde{\Ring},J}(w)$
and $F'_k\subset \Delta_{\tilde{\Ring},J}(s_i w)$ with the property that if $F_j/F_{j-1}=\Delta_{\tilde{\Ring}}(w_j)$,
then $F'_j/F'_{j-1}=\Delta_{\tilde{\Ring}}(s_iw_j)$. It follows from Lemma \ref{Lem:Verma_coinc}
that $\mathfrak{P}_i\Delta_{\tilde{\Ring}}(s_i w_j)\xrightarrow{\sim}
\mathfrak{P}_i\Delta_{\tilde{\Ring}}(w_j)$. Therefore our morphism
$\mathfrak{P}_i\Delta_{\tilde{\Ring},J}(s_iw)\rightarrow \mathfrak{P}_i\Delta_{\tilde{\Ring},J}(w)$
is an isomorphism after passing to the associated graded objects. Hence it is an isomorphism.
We conclude that
$$\mathfrak{P}_i^*\mathfrak{P}_i\Delta_{\tilde{\Ring},J}(s_iw)\xrightarrow{\sim}
\mathfrak{P}_i^*\mathfrak{P}_i\Delta_{\tilde{\Ring},J}(w).$$
Now we are done by Step 1.

{\it Case 3: $s_i C=C$}. Note that the multiplication by $s_i$ preserves
the poset coideal of all elements $\succeq_J C$ in $W/W_J$. So the functor
$\mathfrak{P}_i^*\mathfrak{P}$ preserves the categories
$\mathcal{O}_{\tilde{\Ring}}(W)_{C}^{\Delta}\subset \mathcal{O}_{\tilde{\Ring}}(W)_{\succeq C}^{\Delta}$.
The object $\Delta_{\tilde{\Ring},J}(w)$
is a projective object in the subcategory $\mathcal{O}_{\tilde{\Ring}}(W)_{C}^{\Delta}$
and therefore so is $\mathfrak{P}_i^*\mathfrak{P}_i
\Delta_{\tilde{\Ring},J}(w)$. The  class of the latter in $K_0$ is a linear combination
of $\Delta_{\tilde{\Ring}}(u)$ with $u\in C$. But any indecomposable
projective in $\mathcal{O}_{\tilde{\Ring}}(W)^{\Delta}_{C}$  is of the
form $\Delta_{\tilde{\Ring},J}(w')$ with $w'\in C$ (this follows from the definition of
the objects $\Delta_{\tilde{\Ring}}(w')$).

This finishes the proof.
\end{proof}

\begin{Cor}\label{Cor:stand_str_O}
The Bruhat order on $W/W_J$ gives rise to a standardly stratified structure
on $\OCat_{\tilde{\Ring}}(W)$  with standard objects $\Delta_{\tilde{\Ring},J}(w)$
and costandard objects $\nabla_{\tilde{\Ring},J}(w)$. The action
$\DG(\varnothing,\varnothing)$ on $\OCat_{\tilde{\Ring}}(W)$
preserves standardly filtered and costandardly filtered objects.
\end{Cor}
\begin{proof}
Using (1) of Lemma
\ref{Lem:SB_action} combined with Proposition \ref{Prop:filtr_preserv} we see that
every projective is included as a direct summand into a module filtered by
$\Delta_{\tilde{\Ring},J}(w)$ (in order).

To prove (SS2) for the indecomposable projectives in $\,_\varnothing\OCat^+_{\tilde{\Ring}}(W)$,
we need an analog of Lemma
\ref{Lem:dir_sum_filtration} in this situation. So suppose that $M_1\oplus M_2$
has a filtration with quotients of the form $\Delta_{\tilde{\Ring},J}(w)$ (in the decreasing order).
By Lemma \ref{Lem:dir_sum_filtration},
both $M_1,M_2$ have filtrations by $\Delta_{\tilde{\Ring}}(w)$'s. Let $C$ be a maximal
$W_J$-coset such that $\Delta_{\tilde{\Ring}}(w)$ occurs as a subquotient in a filtration of
$M_1\oplus M_2$ for some $w\in C$. Then we have maximal subobjects $M_i'\subset M_i, i=1,2,$
such that $M_i'$ is filtered by $\Delta_{\tilde{\Ring}}(w),w\in C$, and $M_i/M_i'$
is filtered by $\Delta_{\tilde{\Ring}}(u)$ with $uW_J\not\succeq_J C$. The projection of
$M_1'\oplus M_2'$ to the subquotient $\,_\varnothing \OCat^+_{\tilde{\Ring}}(W)_C$
is projective. In other words, $M_1'\oplus M_2'$ is a direct sum of the objects
$\Delta_{\tilde{\Ring},J}(w), w\in C$. It follows that each $M_1',M_2'$ is such
a direct sum. Now we can apply induction (on $C$). We conclude that (SS2) holds for
$\,_\varnothing\OCat^+_{\tilde{\Ring}}(W)$.

That every relatively $\Ring$-injective is costandardly filtered is proved analogously.
Now we use an equivalence of (a) and (b) in Lemma \ref{Lem:compat_SS_equiv} to deduce that $\OCat_{\tilde{\Ring}}(W)$
is standardly stratified. The claim that the $\DG(\varnothing,\varnothing)$-action
preserves the standardly filtered objects is a direct consequence of Proposition
\ref{Prop:filtr_preserv}. The claim about costandard objects is analogous.
\end{proof}

Note that the pre-order $\preceq_J$ refines $\preceq_{J'}$ provided $J\subset J'$.
From (\ref{eq:stand_compos}) we deduce that
\begin{equation}\label{eq:Delta_stages}
\Delta_{\tilde{\Ring},J}=\Delta_{\tilde{\Ring},J'}\circ \underline{\Delta}_{\tilde{\Ring},J},
\end{equation}
where we write $\underline{\Delta}_{\tilde{\Ring},J}$ for the standardization
functor for the standardly stratified structure $\preceq_{J}$ on
$\gr_{\preceq_{J'}}\OCat^+_{\tilde{\Ring},J}$.
A similar property holds for the costandard objects.

We will also need a property of certain projective objects in $\,_\varnothing\OCat_{\F}(W)$.

\begin{Lem}\label{Lem:antidom_proj}
Pick $i\in S$. Let $w\in W$ be such that  $l(w s_i)\geqslant l(w)$.
Then the projective object $P_{\tilde{\Ring}}(w)$ is filtered by $\Delta_{\tilde{\Ring},i}(w')$
with $l(w' s_i)\geqslant l(w')$. A similar claim holds for relatively $\tilde{\Ring}$-injective objects.
\end{Lem}
\begin{proof}
Note that the classes of all objects $\Delta_{\tilde{\Ring},i}(w')$ generate the submodule
$\C W e_i$ in $K^{\C}_0(\,_\varnothing\OCat_{\tilde{\Ring}}(W)^{\Delta})$. Indeed, by the definition
of the objects $\Delta_{\tilde{\Ring},i}(w')$, it is enough to check this for $W=W_i$, where the claim
is clear.

Furthermore, if a linear combination of the classes  $[\Delta_{\tilde{\Ring},i}(w')]$
belongs to $\C W e_i$, then all elements $w'$ appearing with nonzero coefficients
satisfy $l(w' s_i)\geqslant l(w')$. Again, this reduces to the case $W=W_i$, where this is clear.

Classes of objects of the form $B.\Delta_{\tilde{\Ring},i}(s_i)$, where $B\in \DG(\varnothing,\varnothing)$,
generate $\C W e_i$. Note that the  objects $B.\Delta_{\tilde{\Ring},i}(s_i)$
are projective. By the previous paragraph, they are filtered by objects of the
form $\Delta_{\tilde{\Ring},i}(w')$ with $l(w's_i)\geqslant l(w')$.
Take the indecomposable object $B_w$ in $\DG(\varnothing,\varnothing)$ corresponding to $w$. So
$P_{\tilde{\Ring}}(w)$ is a direct summand in $B_w.\Delta_{\tilde{\Ring},i}(s_i)$
and hence is also filtered by the objects $\Delta_{\tilde{\Ring},i}(w')$.
\end{proof}

\subsection{Categories $\,_{I}\mathcal{O}_{\Ring,J}^+(W)$}
We are proceeding to defining parabolic subcategories $\mathcal{O}$ in
$\,_{I}\mathcal{O}_{\Ring,J}^+(W)$, where $\Ring$ is a suitable quotient
of $\tilde{\Ring}$. Again, let us start by recalling how parabolic
subcategories in the usual BGG categories are constructed. Let $\mathfrak{p}
\subset \mathfrak{g}$ be a standard parabolic subalgebra with Levi subalgebra
$\mathfrak{l}$. Then we can consider the  full  subcategory $\OCat^{\mathfrak{p}}(\mathfrak{g})\subset \OCat(\mathfrak{g})$
consisting of all modules where $\mathfrak{p}$ acts locally finitely.
Equivalently, $\OCat^{\mathfrak{p}}(\mathfrak{g})$ is spanned by  the simples whose highest weights are non-negative on all
positive coroots in $\mathfrak{l}$.  The category $\mathcal{O}^{\mathfrak{p}}(\mathfrak{g})$
is highest weight, the standard (resp., costandard) objects are obtained by applying the induction
functor $\Delta_{\mathfrak{l}}$ (resp., coinduction functors
$\nabla_{\mathfrak{l}}$) to the finite dimensional modules in $\OCat(\mathfrak{l})$.

We can also consider a deformed version of $\OCat^{\mathfrak{p}}(\mathfrak{g})$. The simple coroots
in $\mathfrak{l}$ still need to act by integers so the deformed version will be
defined over the ring $\C[x_i]_{i\in I}/(x_j)_{j\in J}$.

Below we fix a finitary subset $J\subset S$  and set $\Ring:=\tilde{\Ring}/(x_j)_{j\in J}$.

\subsubsection{Definition of $\,_I\mathcal{O}_{\F,J}^+(W)$}
Consider the subset $\,_{I}\Lambda^+(J)\subset \,_{I}\Lambda^+$ consisting of the maximal (=shortest)
elements in all right $W_{J}$ cosets consisting of $|W_{J}|$ elements.
Note that the restriction of $\varrho: \,_{I}\Lambda^+\twoheadrightarrow \,_{I}\Lambda^+_J$
to $\,_{I}\Lambda^+(J)$ is an injection.

The following lemma is straightforward.

\begin{Lem}\label{Lem:parab_defn}
For an object $M\in \,_I\OCat_{\F}(W)$ the following is equivalent:
\begin{enumerate}
\item $\Hom_{\OCat_\F}(P_\F(\lambda), M)=0$ for all $\lambda\not\in \,_{I}\Lambda^+(J)$.
\item $\Hom_{\OCat_\F}(M,I_\F(\lambda))=0$ for all $\lambda\not\in \,_{I}\Lambda^+(J)$.
\item For any finite coideal $\Lambda^0\subset \,_{I}\Lambda^+(J)$ that is the union of equivalence classes for $\preceq_J$,
the simple $L_{\F}(\lambda)$ is not a composition factor of $\pi_{\Lambda^0}(M)$
when $\lambda\not\in \,_{I}\Lambda^+(J)$.
\end{enumerate}
\end{Lem}

For $\,_{I}\mathcal{O}_{\F,J}^+(W)$ we take the full subcategory of
$\,_{I}\mathcal{O}_{\F}^+(W)$ consisting of the objects satisfying the
conditions of the previous lemma. Let us give examples of objects in
$\,_{\varnothing}\mathcal{O}_{\F,J}^+(W)$.

\begin{Lem}\label{Lem:parab_st_cost}
For the shortest element $w$ in its right $W_J$-coset, we have $\bar{\Delta}_{\F,J}(w),
\bar{\nabla}_{\F,J}(w)\in \,_\varnothing\mathcal{O}_{\F,J}^+(W)$.
\end{Lem}
\begin{proof}
We will check that $\Hom_{\OCat_\F}(P_\F(w'), \bar{\nabla}_{\F,J}(w))=0$
when $w'$ is not the shortest element in its $W_J$-coset.
The analog of (\ref{eq:Delta_stages}) implies that  $\bar{\nabla}_{\F,J}(w)$ embeds into
$\bar{\nabla}_{\F,j}(w)$ for any $j\in J$. So it is enough to show
that  $\Hom_{\OCat_\F}(P_\F(w'), \bar{\nabla}_{\F,j}(w))=0$
when $l(w' s_j)\geqslant l(w')$ and $l(w s_j)<l(w)$. Lemma \ref{Lem:antidom_proj}
reduces this to checking that  $\Hom_{\OCat_\F}(\Delta_{\F,j}(w'), \bar{\nabla}_{\F,j}(w))=0$.
But this follows from the fact that $\Delta_{\F,j}(w'),\bar{\nabla}_{\F,j}(w)$
are standard and proper costandard objects for a standardly stratified structure
and $w'\neq w$.

Similarly, we can prove that $\Hom_{\OCat_\F}(\bar{\Delta}_{\F,J}(w),I_\F(w'))=0$.
\end{proof}

\subsubsection{Classes of $\bar{\Delta}_{\F,J}(\lambda), \bar{\nabla}_{\F,J}(\lambda)$
in $K_0$}
\begin{Lem}\label{Lem:parab_K0}
The class of $\bar{\Delta}_{\F,J}(w)$ in $K_0(\,_\varnothing\OCat^+_\F(W)^{\Delta})$
equals $\sum_{u\in W_J}(-1)^{l(u)}wu$. The similar claim is true for
the class of $\bar{\nabla}_{\F,J}(w)$ in $K_0(\,_\varnothing\OCat^+_\F(W)^{\nabla})$.
\end{Lem}
\begin{proof}
Both equalities $[\bar{\Delta}_{\F,J}(w)]=\sum_{u\in W_J}(-1)^{l(u)}wu$
and  $[\bar{\nabla}_{\F,J}(w)]=\sum_{u\in W_J}(-1)^{l(u)}wu$
boil down to $[L_\F(1)]=\sum_{w\in W}(-1)^{l(w)}w$ for $W=W_J$.
Clearly, the coefficient of $1$ in $[L_{\F}(1)]$ is $1$. On the other hand,
the proof of Lemma \ref{Lem:antidom_proj} shows that the classes of
$P_{\F}(w')$ with $w' s_j\geqslant w'$ span $\C W e_j$. So $[L_\F(1)]$
is orthogonal to all these classes in the pairing, where the elements
of $W$ form an orthonormal basis. Our claim follows.
\end{proof}

\subsubsection{$\DG(\varnothing,\varnothing)$-module structure}
We consider the full subcategory  $\,_{\varnothing}\OCat^+_{\F,J}(W)^{\Delta}
\subset \,_{\varnothing}\OCat^+_{\F,J}(W)$ consisting of all objects that
admit a filtration by $\bar{\Delta}_{\F,J}(w)$'s (with $w$ shortest in its
right $W_J$-coset). Similarly, we have
the full subcategory  $\,_{\varnothing}\OCat^+_{\F,J}(W)^{\nabla}
\subset \,_{\varnothing}\OCat^+_{\F,J}(W)$.

\begin{Lem}\label{Lem:parab_action_princ}
The action of $\DG(\varnothing,\varnothing)$ on $\,_{\varnothing}\OCat^+_{\F}(W)$ preserves
$\,_{\varnothing}\OCat^+_{\F,J}(W)$ and
$\,_{\varnothing}\OCat^+_{\F,J}(W)^{\Delta}$.
\end{Lem}
\begin{proof}
$\DG(\varnothing,\varnothing)$ preserves $\,_\varnothing\OCat^+_{\F,J}(W)$ because
it preserves the subcategory in $\,_\varnothing\OCat^+_{\F}(W)\operatorname{-proj}$
spanned by $P_{\F}(w)$, where $w$ is not the shortest element in its right $W_J$-coset.
Indeed, there is $j\in J$ such that
$w s_j> w$. By the proof of Lemma \ref{Lem:antidom_proj},
$B.P_\F(w)$ admits a filtration by $\Delta_{\F,j}(w')$ with $w's_j< w'$
that implies that $\DG(\varnothing,\varnothing)$ preserves the subcategory of the projectives of interest.

It follows from Corollary \ref{Cor:stand_str_O} and (\ref{eq:Ext_vanish}) that
$$\Ext^i(\bar{\Delta}_{\F,J}(w), \nabla_{\F,J}(w'))=\delta_{ww'}\delta_{i0}.$$
Now $\DG(\varnothing,\varnothing)$ preserves $\,_\varnothing\OCat^+_{\F,J}(W)^{\Delta}$  because
it preserves the category of  costandardly filtered objects for $\preceq_J$.
\end{proof}

\subsubsection{Singular blocks}
Now we proceed to the category $\bigoplus_I \,_I\OCat_{\F,J}(W)$.
Our goal is to prove the following proposition.

\begin{Prop}\label{Lem:parab_action_sing}
The following is true:
\begin{enumerate}
\item $\mathfrak{P}_I\bar{\Delta}_{\F,J}(w)$ and
$\mathfrak{P}_I\bar{\nabla}_{\F,J}(w)$ lie in $\,_I\OCat^+_{\F,J}(W)$
and are zero if the stabilizer of $wW_J$ in $W_I$ is nontrival.
\item The action of $\DG$ on $\bigoplus_I \,_I\OCat^+_{\F}(W)$
preserves  $\bigoplus_I \,_I\OCat^+_{\F,J}(W)$.
\item The objects $\mathfrak{P}_I\bar{\Delta}_{\F,J}(w)$ and
$\mathfrak{P}_I\bar{\nabla}_{\F,J}(w)$ depend only on $W_I w$.
\item
The indecomposable projectives in  $\,_I\OCat^+_{\F,J}(W)$
are indexed by $\,_{I}\Lambda^+(J)$ and are direct summands
of the objects of the form $\mathfrak{P}_IB.\Delta_{\F,J}(1)$ with
$B\in \DG(\varnothing,\varnothing)$.
\end{enumerate}
\end{Prop}

If the stabilizer of $W_Iw$ in $W_J$ is trivial, we write
$\bar{\Delta}_{\F,J}(\lambda),\bar{\nabla}_{\F,J}(\lambda)$ for $\mathfrak{P}_I\bar{\nabla}_{\F,J}(w)$,
where $\lambda=W_I w$.

In the proof we will need a lemma that is a singular analog of Lemma
\ref{Lem:antidom_proj}.

\begin{Lem}\label{Cor:proj_classes}
Pick $i\in S$. For $\mu\in W_I\setminus W$, the following is equivalent:
\begin{enumerate}
\item $\mu s_i\leqslant \mu$,
\item  the class of $P_{\F}(\mu)$ lies in $e_I\C W e_i$.
\end{enumerate}
\end{Lem}
\begin{proof}
Let us show that (2) implies (1). If $\mu s_i>\mu$, then standard composition
subquotients of $P_{\F}(\mu)$ are $\Delta_{\F}(\mu)$ and objects of the form
$\Delta_{\F}(\mu')$ with $\mu'W_i<\mu W_i$. This description implies that
$[P_{\F}(\mu)]$ cannot lie in $e_I\C W e_i$.

Suppose (1) holds. Then $P_{\F}(\mu)$ is a direct summand in
$\mathfrak{P}_I P_{\F}(w)$, where $\mu=W_I w$ and $l(ws_i)> l(w)$.
The object $P_{\F}(w)$ admits a filtration by $\Delta_{\F,i}(w')$,
where $w's_i>w'$. So $\mathfrak{P}_I P_{\F}(w)$ has a filtration by
$\mathfrak{P}_I\Delta_{\F,i}(w')$. The latter has a filtration
with  sub $\Delta_{\F}(W_I w')$ and quotient $\Delta_{\F}(W_Iw's_i)$.
It is enough to prove that the filtration does not split provided $W_I w's_i\neq W_Iw'$.
Assume the contrary. Then $\mathfrak{P}_I^*\mathfrak{P}_I\Delta_{\F,i}(w')=
\mathfrak{P}_I^*\Delta_{\F}(W_I w')\oplus \mathfrak{P}_I^*\Delta_{\F,i}(W_Iw's_i)$.
The first summand has a filtration by all $\Delta_{\F}(uw')$ with $u\in W_I$
(each occurring with multiplicity $1$). In particular, $\Delta_{\F}(w's_i)$ is
not among the composition factors, so $\mathfrak{P}_I^*\Delta_{\F}(W_I w')$ does
not lie in $\,_\varnothing\OCat^+_\F(W)^{\Delta,i}$.
This contradicts Proposition \ref{Prop:filtr_preserv}.
\end{proof}

\begin{proof}[Proof of Proposition \ref{Lem:parab_action_sing}]
Let us prove (1) for $\bar{\nabla}_{\F,J}(w)$ (the case of $\bar{\Delta}$ is similar).
The claim that $\mathfrak{P}_I\bar{\nabla}_{\F,J}(w)=0$ when some nontrivial
element of $W_J$ stabilizes $W_Iw$ follows from a computation in $K_0$.
In order to prove that $\mathfrak{P}_I\bar{\Delta}_{\F,J}(w)\in \,_I\OCat^+_{\F,J}(W)$
we need to check that $\Hom(P_{\F}(\mu),\mathfrak{P}_I\bar{\nabla}_{\F,J}(w))=0$
when $\mu s_j\leqslant \mu$ for some $j\in J$.
This will follow if we check that $[\mathfrak{P}_I^*P_{\F}(\mu)]\in \C W e_j$.
The latter is a consequence of  Lemma \ref{Cor:proj_classes}.

A similar argument proves that $\bigoplus_I \,_I\OCat^+_{\F,J}(W)$
is stable under the action of $\DG$ establishing (2).

Let us prove (3). Let $w$ be of minimal length in $wW_J$ and
the stabilizer of $wW_J$ in $W_I$ is trivial. Set $\lambda=W_I w$.
Let us first check that
\begin{itemize}
\item[(i)] $\mathfrak{P}_I \bar{\Delta}_{\F,J}(w)$ is the maximal quotient
of $\Delta_{\F}(\lambda)$ lying in $\,_I\OCat^+_{\F,J}(W)$ for $\lambda\in \,_I\Lambda^+(J)$
and
\item[(ii)] $\Delta_{\F}(\mu)$ has no quotients in $\,_I\OCat^+_{\F,J}(W)$ when $\mu\not\in
\,_I\Lambda^+(J)$.
\end{itemize}
(ii) is easy. Let us check (i). The object $\Delta_\F(\lambda)=\mathfrak{P}_I\Delta_{\F}(w)$
has a filtration by $\mathfrak{P}_I\bar{\Delta}_{\F,J}(wu)$ for $u\in W_J$, where
$\mathfrak{P}_I\bar{\Delta}_{\F,J}(w)$ is the top quotient.
Note that all other subquotients admit a surjection from $\Delta_{\F}(W_Iwu)$
and $W_Iwu\not\in \,_I\Lambda^+(J)$ and (i) follows.

Note that (i) implies that $\mathfrak{P}_I\bar{\Delta}_{\F,J}(w)$ depends only on $W_Iw$. The proof
for $\bar{\nabla}$ is completely analogous.

Finally, let us proceed to (4). For general reasons,  every projective $P_{\F}(\lambda)$ in $\,_I\OCat^+_{\F}(W)$ has
the  maximal quotient lying in $\,_I\OCat^+_{\F,J}(W)$ and this quotient is an indecomposable
projective there (to be denoted by $P_{\F,J}(\lambda)$) when $\lambda\in \,_I\Lambda^+(J)$.
Since $\,_I\OCat^+_{\F,J}(W)$ is a Serre subcategory in $\,_I\OCat^+_{\F}(W)$,
there are no other indecomposable projectives. Note that
$P_{\F}(\lambda)\twoheadrightarrow \bar{\Delta}_{\F,J}(\lambda)$.

Now consider the indecomposable Soergel bimodule $\mathfrak{B}_w\in \DG(I,\varnothing)$  corresponding
to $w$ that is shortest in $W_IwW_J$. The object $\mathfrak{P}_I\mathfrak{B}_w \bar{\Delta}_{\F,J}(1)$ is projective and
has a filtration by $\bar{\Delta}_{\F,J}(\mu)$'s with quotient $\bar{\Delta}_{\F,J}(\lambda)$.
So $P_{\F,J}(\lambda)$ is a summand of $\mathfrak{P}_IB \bar{\Delta}_{\F,J}(1)$.  This finishes
the proof of (4).
\end{proof}

\subsubsection{Highest weight structure}
Let us prove some corollaries of Proposition \ref{Lem:parab_action_sing}.

\begin{Cor}\label{Cor:parab_hw}
The category $\,_I\OCat^+_{\F,J}(W)$  is highest weight with poset $\varrho(\,_I\Lambda^+(J))$
and standard objects $\Delta_{\F,J}(\lambda)$. The costandard objects are $\nabla_{\F,J}(\lambda)$.
\end{Cor}
\begin{proof}
(HW1) is checked easily. Lemma \ref{Lem:dir_sum_filtration} shows that
$\,_I\OCat^+_{\F,J}(W)^{\Delta}$ is closed under taking direct summands.

(HW2) now follows from (4) of Proposition \ref{Lem:parab_action_sing}
combined with (2) of Lemma \ref{Lem:parab_action_princ}.

The claim about costandard objects
will follow if we  check $\operatorname{Ext}^i_{\OCat^+_{\F,J}}(\bar{\Delta}_{\F,J}(\lambda), \bar{\nabla}_{\F,J}(\lambda'))$ vanishes when $i=1$ and when $i=0$ and $\lambda\neq \lambda'$.
Note that these Ext's are the same as in $\OCat^+_{\F}$. If $\lambda=\lambda'$ and $i=0$,
then, since $$\Delta_{\F}(\lambda)\twoheadrightarrow \bar{\Delta}_{\F,J}(\lambda)\twoheadrightarrow
L_{\F}(\lambda)\hookrightarrow \bar{\nabla}_{\F,J}(\lambda)\hookrightarrow \nabla_{\F}(\lambda)$$
we have $\dim\Hom_{\OCat_{\F}^+}(\bar{\Delta}_{\F,J}(\lambda),
\bar{\nabla}_{\F,J}(\lambda))=1$. Now let $i=1$ or $\lambda\neq \lambda'$,
then we have, for $\lambda=W_I w, \lambda'=W_I w'$,
$$\operatorname{Ext}^i_{\OCat^+_{\F}}(\bar{\Delta}_{\F,J}(\lambda), \bar{\nabla}_{\F,J}(\lambda'))=
\operatorname{Ext}^i_{\OCat^+_{\F}}(\bar{\Delta}_{\F,J}(w), \mathfrak{P}_I^*\bar{\nabla}_{\F,J}(\lambda')).$$
Using the standardly stratified structure on $\OCat^+_\F$ it is easy to see that the
right hand side is zero. So we have proved that the costandard objects are precisely
the objects $\nabla_{\F,J}(\lambda)$.
\end{proof}

\begin{Cor}\label{Cor:parab_action_K0}
The following is true.
\begin{enumerate}
\item
The action of $\DG$ on $\bigoplus_I \,_I\OCat^+_{\F,J}(W)$ preserves
$\bigoplus_I \,_I\OCat^+_{\F,J}(W)^{\Delta}$ (the category of all objects
that admit a filtration by $\bar{\Delta}_J(\lambda)$'s) and
$\bigoplus_I \,_I\OCat^+_{\F,J}(W)^{\nabla}$.
\item The $K_0$ of the exact category $\bigoplus_I \,_I\OCat^+_{\F,J}(W)^{\Delta}$
is identified with the $K_0(\DG)$-module $\bigoplus_I e_I \C W e^-_J$,
where $e^-_J$ is the sign idempotent in $\C W_J$. The same is true for
$\bigoplus_I \,_I\OCat^+_{\F,J}(W)^{\nabla}$.
\end{enumerate}
\end{Cor}
\begin{proof}
Let us prove (1). Recall, (S2), that  the category $\DG$ is generated by
the projection bimodules $\mathfrak{P}_{I,j}:=\mathfrak{P}^I_{I\cup \{j\}}\in
\DG(I\cup\{j\},I)$ for $j\not\in I$ and their adjoints
$\mathfrak{P}^*_{I,j}\in \DG(I,I\cup\{j\})$. Clearly,
$\mathfrak{P}_{I,j}$ maps $\Delta_{\F,J}(W_Iw)$
to $\Delta_{\F,J}(W_{I\cup\{j\}}w)$ or 0. It also maps
$\nabla_{\F,J}(W_Iw)$ to $\nabla_{\F,J}(W_{I\cup\{j\}}w)$ or
0. The last sentence and Corollary \ref{Cor:parab_hw} imply
that $\mathfrak{P}^*_{I,j}\bar{\Delta}_{\F,J}(W_{I\cup j}w)$
is standardly filtered. This completes the proof of (1).

(2) follows from the formula for $[\bar{\Delta}_{\F,J}(W_Iw)]=e_I
[\bar{\Delta}_{\F,J}(w)]$ and the formula for
$[\bar{\Delta}_{\F,J}(w)]$ from Lemma \ref{Lem:parab_K0}. The $\nabla$
case is similar.
\end{proof}

\subsubsection{Definition and properties of $\,_{I}\mathcal{O}_{\Ring,J}^+(W)$}\label{SSS_parab_posit_deform}
Here we define the category $\,_{I}\mathcal{O}_{\Ring,J}^+(W)$ that will be shown
to be highest weight over $\Ring$ with poset $\varrho(\,_I\Lambda^+(J))$
(where, recall, $\varrho$ is the projection $\,_{I}\Lambda^+\twoheadrightarrow \,_{I}\Lambda^+_J$)
and to specialize to $\,_{I}\mathcal{O}_{\F,J}^+(W)$ after the base change from
$\Ring$ to $\F$.

We start with an analog of Lemma \ref{Lem:parab_defn} whose proof is left to the reader.

\begin{Lem}\label{Lem:parab_defn_deform}
For an object $M\in \,_I\OCat_{\Ring}(W)$ the following is equivalent:
\begin{enumerate}
\item $\Hom_{\OCat_\Ring}(P_\Ring(\lambda), M)=0$ for all $\lambda\not\in \,_{I}\Lambda^+(J)$.
\item For any finite coideal $\Lambda^0$ that is the union of equivalence classes for $\preceq_J$,
the simple $L_{\F}(\lambda)$ is not a composition factor of $\F\otimes_\Ring\pi_{\Lambda^0}(M)$
when $\lambda\not\in \,_{I}\Lambda^+(J)$.
\end{enumerate}
When $M$ is flat over $\Ring$ the conditions (1),(2) are equivalent to
$\Hom_{\OCat_\Ring}(M,I_\Ring(\lambda))=0$ for all $\lambda\not\in \,_{I}\Lambda^+(J)$.
\end{Lem}

For $\,_{I}\mathcal{O}_{\Ring,J}^+(W)$ we take the full subcategory of $\,_I\OCat^+_{\Ring}(W)$
of all objects satisfying the equivalent conditions of Lemma \ref{Lem:parab_defn_deform}.

\begin{Lem}\label{Lem:barDeltaR}
Let $\lambda\in \,_I\Lambda^+(J)$. There is a largest quotient, to be denoted
by $\bar{\Delta}_{\Ring,J}(\lambda)$, of $\Delta_\Ring(\lambda)$ lying in $\,_I\OCat^+_\Ring(W)$.
This object is flat over $\Ring$ and its specialization to $\F$ coincides
with $\Delta_{\F,J}(\lambda)$.
\end{Lem}
\begin{proof}
The functor $\Delta_{\tilde{\Ring},J}$ preserves $\Ring$-flat modules.
Set $\bar{\Delta}_{\Ring,J}(w)=\Delta_{\tilde{\Ring},J}(\Ring\otimes_\F L_{\F,J}(1))$
(here $w$ is shortest in its right $W_J$-coset) and $\bar{\Delta}_{\Ring,J}(W_Iw)=
\mathfrak{P}_I \bar{\Delta}_{\Ring,J}(w)$. So the object $\bar{\Delta}_{\Ring,J}(W_Iw)$
is flat over $\Ring$ and it is also a quotient of  $\Delta_\Ring(W_I w)$.
The proof that this quotient is maximal is similar to the proof of (3) of
Proposition \ref{Lem:parab_action_sing}.
\end{proof}

\begin{Prop}\label{Prop:parab_hw_deform}
The following is true:
\begin{enumerate}
\item The indecomposable projectives in $\,_{I}\mathcal{O}_{\Ring,J}^+(W)$ are indexed by $\,_{I}\Lambda^+(J)$
and are direct summands of the objects of the form $\mathfrak{P}_IB\bar{\Delta}_{\Ring,J}(1)$
for $B\in \DG(\varnothing,\varnothing)$.
\item The category $\,_{I}\mathcal{O}_{\Ring,J}^+(W)$ is highest weight over $\Ring$ with poset
$\,_{I}\Lambda^+(J)$ and standard objects $\bar{\Delta}_{\Ring,J}(\lambda)$. It specializes
to $\,_{I}\mathcal{O}_{\F,J}^+(W)$.
\end{enumerate}
\end{Prop}
\begin{proof}
The proof of (1) repeats that of (3) in Lemma \ref{Lem:parab_action_sing}.

Let us check axioms of a highest weight category over $\Ring$ listed in \S\ref{SS_hw_coid_fin}.
We need to check that $\operatorname{Hom}_{\OCat_\Ring}(P_{\Ring,J}(\lambda), P_{\Ring,J}(\lambda'))$
is free over $\Ring$. But the last Hom coincides with
$\operatorname{Hom}_{\OCat_\Ring}(P_{\Ring}(\lambda), P_{\Ring,J}(\lambda'))$ that is a direct summand
in $\operatorname{Hom}_{\OCat_\Ring}(P_{\Ring}(\lambda), B.\Delta_{\Ring,J}(1))=
\operatorname{Hom}_{\OCat_\Ring}(B^*.P_{\Ring}(\lambda), \Delta_{\Ring,J}(1))$, where $B^*$ denote the
biadjoint of $B$. The module $$\operatorname{Hom}_{\OCat_\Ring}(B^*.P_{\Ring}(\lambda), \Delta_{\Ring,J}(1))$$ is free
of finite rank over $\Ring$ because $B^*.P_{\Ring}(\lambda)$ is projective in
$\OCat_{\Ring}^+(W)$ and $\Delta_{\Ring,J}(1)$ is free of finite rank over $\Ring$. The second and the third axioms
from \S\ref{SS_hw_coid_fin} hold for $\,_{I}\mathcal{O}_{\Ring,J}^+(W)$ because they hold for
the ambient category $\,_{I}\mathcal{O}_{\Ring}^+(W)$. Finally, the axioms concerning
the standard objects are verified as in Corollary \ref{Cor:parab_hw}. The claim about specialization
is straightforward.
%
\end{proof}

\subsubsection{Definition and properties of $\,_{I}\mathcal{O}_{\Ring,J}^-(W)$}\label{SSS_neg_level_parab}
Above we have defined a ``positive level'' version of the parabolic category $\mathcal{O}$.
We need a ``negative level'' version. To define it, we get an inspiration from a classical result from the
representation theory of Kac-Moody algebras: the classical positive and negative level categories $\mathcal{O}$
are Ringel dual.

Here we are going to define the category $\,_{I}\mathcal{O}_{\Ring,J}^-(W)$ and list
its basic properties. By definition, this is a highest weight category over $\Ring$ with poset $\,_I\Lambda^-(J)$
that is opposite to $\,_I\Lambda^+(J)$. It is defined as the  category whose Ringel
dual is $\,_{I}\mathcal{O}_{\Ring,J}^+(W)$ so that $\,_{I}\mathcal{O}_{\Ring,J}^-(W)^{\Delta}\cong
(\,_{I}\mathcal{O}_{\Ring,J}^+(W)^{\Delta})^{opp}$ and $\,_{I}\mathcal{O}_{\Ring,J}^-(W)\operatorname{-tilt}
\cong (\,_{I}\mathcal{O}_{\Ring,J}^+(W)\operatorname{-proj})^{opp}$. Let $\Delta^-_{\Ring,J}(\lambda)$ denote the standard  objects in $\,_{I}\mathcal{O}_{\Ring,J}^-(W)$.

The following lemma lists basic properties of the category $\bigoplus_I \,_{I}\mathcal{O}_{\Ring,J}^-(W)$.

\begin{Lem}\label{Lem:O_neg_prop}
The following is true:
\begin{enumerate}
\item The category $\bigoplus_I \,_{I}\mathcal{O}_{\Ring,J}^-(W)$ carries an $\Ring$-linear action of $\DG$.
\item The action  preserves the subcategory
$\bigoplus_I \,_{I}\mathcal{O}_{\Ring,J}^-(W)^{\Delta}$.
\item The action induced on $K_0(\bigoplus_I \,_{I}\mathcal{O}_{\Ring,J}^-(W)^{\Delta})$
is the $\bigoplus_{I,I'} e_I \C W e_{I'}$-action on $\bigoplus_{I'}e_{I'}\C W e_J$
and the classes  $[\Delta^-_{\Ring,J}(\lambda)]$ correspond to the standard basis
vectors.
\item Any indecomposable tilting in $\bigoplus_I \,_{I}\mathcal{O}_{\Ring,J}^-(W)$
is a direct summand in an object of the form $B.\Delta^-_{\Ring,J}(1)$.
\end{enumerate}
\end{Lem}
\begin{proof}
Note that $\bigoplus_I\,_{I}\mathcal{O}_{\Ring,J}^-(W)^{\Delta}$ carries an action of $\DG$ that comes from the
identification with $(\bigoplus_I\,_{I}\mathcal{O}_{\Ring,J}^+(W)^{\Delta})^{opp}$. Recall that the collection of functors
arising from the action is closed with respect to taking biadjoints. So the action is by exact functors.
Since any object in $\bigoplus_I\,_{I}\mathcal{O}_{\Ring,J}^-(W)$
 admits a resolution by standardly filtered objects,
 the action extends to the whole category $\bigoplus_I\,_{I}\mathcal{O}_{\Ring,J}^-(W)$. This finishes the proof of (1).

The proof of (2) is similar to that of part (1) of Corollary \ref{Cor:parab_action_K0}.
To prove (3) we use the identification  $\bigoplus_I\,_{I}\mathcal{O}_{\Ring,J}^-(W)^{\Delta}
\cong(\bigoplus_I\,_{I}\mathcal{O}_{\Ring,J}^+(W)^{\Delta})^{opp}$ and (2) of  Corollary \ref{Cor:parab_action_K0}.

(4) follows from the identification
$\,_{I}\mathcal{O}_{\Ring,J}^-(W)\operatorname{-tilt}
\cong (\,_{I}\mathcal{O}_{\Ring,J}^+(W)\operatorname{-proj})^{opp}$
and (1) of Proposition \ref{Prop:parab_hw_deform}. 
\end{proof}

%% file: FockUnique4.tex
\section{Uniqueness of categorified Fock spaces}\label{S_Fock_unique}
\subsection{Highest weight categorical $\tgl_e$-actions}
Let $\FM$ be a base field of characteristic $p$. Pick  integers $e>1, \ell\geqslant 1$.
Form the formal power series ring $\Ringg:=\FM[[z_0,\ldots,z_{\ell-1}]]$.

\subsubsection{Highest weight categorical $\slf_2$-actions}
Let us start by recalling the notion of a highest weight $\slf_2$-categorification over a field
that appeared in \cite{str}.  We start with an $\FM$-linear finite length abelian category $\Cat$
equipped with a categorical  $\slf_2$-action (in the sense of Chuang and Rouquier, \cite{ChuRou06})
with functors $E',F'$. We further assume that $\Cat$ is a highest weight category with ideal finite
poset $\Lambda$. A highest weight $\slf_2$-categorification is an additional structure (so called {\it hierarchy structure}) on $\Lambda$ that satisfies some axioms that ensure a compatibility between the action and the highest weight structure.

We are going to explain the first part of a hierarchy structure, a family structure, in
detail.  This is a decomposition of $\Lambda$ into the disjoint union
of subsets $\Lambda=\bigsqcup_a \Lambda_a$ together with an identification
$\sigma_a:\{+,-\}^{n_a}\xrightarrow{\sim}\Lambda_a$. We equip the set $\{+,-\}^{n_a}$
with the dominance order: $t=(t_1,\ldots,t_n)\leqslant t'=(t_1',\ldots,t_n')$ if the total numbers
of $+$'s in $t,t'$ are the same, while for each $k\leqslant n$, the number of $+$'s in
$(t'_1,\ldots,t'_k)$ does not exceed that in $(t_1,\ldots,t_k)$. We require $\sigma_a^{-1}$
to be increasing.

The second part of a hierarchy structure is a so called {\it splitting structure} from
\cite[Section 3.1]{str} that to each $a$ assigns a decomposition
$\Lambda=\Lambda^a_{<}\sqcup \underline{\Lambda}^-_a\sqcup \underline{\Lambda}^+_a
\sqcup \Lambda^a_{>}$, where $\Lambda^a_{<}, \Lambda^a_{<}\sqcup \underline{\Lambda}^-_a,
\Lambda^a_{<}\sqcup \underline{\Lambda}^-_a\sqcup \underline{\Lambda}^+_a$ are poset ideals.
We are not going to give a precise definition: it is
very technical and we do not need it explicitly, in fact, all splitting structures
we need have appeared in \cite{VV_proof}.

Let us explain the axioms of a highest weight categorical $\slf_2$-action.
Pick $t\in \{+,-\}^{n_a}$. Let $t^1,\ldots,t^r$ be all elements obtained from $t$
by changing one $+$ to a $-$ and let $\bar{t}^1,\ldots, \bar{t}^s$ be all elements
obtained from $t$ by switching one $-$ to a $+$. The axioms we need are the following two
conditions:
\begin{itemize}
\item $F' \Delta_\FM(\sigma_a(t))$ is filtered by $\Delta_\FM(\sigma_a(t^i)), i=1,\ldots,r$,
each  with multiplicity $1$.
\item $E' \Delta_\FM(\sigma_a(t))$ is filtered by $\Delta_\FM(\sigma_a(\bar{t}^i)), i=1,\ldots,s$,
each  with multiplicity $1$.
\end{itemize}

Let us explain some consequences of the axioms from \cite{str}. Fix a family $a$.
Then $\Lambda_a$ is a poset interval and the corresponding highest weight subqotient
$\Cat_{\Lambda_a}$ inherits the categorical $\slf_2$-action and categorifies $(\C^2)^{\otimes n_a}$.
So, as an $\slf_2$-module, $K_0^{\C}(\Cat)$ is the direct sum of several (in general, infinitely many) copies
of tensor products of $\C^2$'s.
Further, for the interval $\underline{\Lambda}^a:=\underline{\Lambda}^a_-\sqcup \underline{\Lambda}^a$,
the highest weight subquotient $\Cat_{\underline{\Lambda}}$ inherits a categorical $\slf_2$-action as well.
The $\slf_2$-module $K_0^{\C}(\Cat_{\underline{\Lambda}^a})$ is of the form $M\otimes \C^2$, where
$M=K^\C_0(\Cat_{\underline{\Lambda}^a_-})\cong K^\C_0(\Cat_{\underline{\Lambda}^a_+})$. Moreover,
the subspace  $K^\C_0(\Cat_{\underline{\Lambda}^a_-})\subset K^\C_0(\Cat_{\underline{\Lambda}^a})$
is identified with $M\otimes v_-$, while the quotient space $K^\C_0(\Cat_{\underline{\Lambda}^a_+})$
is identified with $M\otimes v_+$, where $v_+,v_-$ denote the lowest and the highest vectors in $\C^2$.

In \cite[Section 4.2]{str} we had another axiom (iii) that has placed some restrictions on
the base field $\FM$. We are not going to impose this axiom here.

\subsubsection{Definition}
Now let $\Cat_\Ringg$ be a highest weight category over $\Ringg$ with ideal finite  poset $\Lambda$.
Let us define an {\it $\Ringg$-deformed  highest weight categorical action} of $\tgl_e$ on $\Cat_\Ringg$.
The data of this action is a pair of endo-functors $E,F:\Cat_\Ringg\rightarrow \Cat_\Ringg$ with fixed
adjunction $\mathsf{1}\rightarrow EF, FE\rightarrow \mathsf{1}$ and endomorphisms
$X\in \End(F), T\in \End(F^2)$ subject to the following two conditions:
\begin{itemize}
\item[(HWA1)] The data above defines an $\Ringg$-linear action of the 2-Kac-Moody category $\mathcal{U}(\tgl_e)$
on $\Cat_\Ringg$ mapping standard objects to standardly filtered ones.
\item[(HWA2)] For every $i\in \Z/e\Z$, the functors $E_i,F_i$ give a highest weight
categorical action of $\slf_2$ on $\Cat_\FM:=\FM\otimes_{\Ringg}\Cat_{\Ringg}$.
\end{itemize}
It follows from (HWA1) that $E_i,F_i$ map costandard objects in $\Cat_\Ringg$
to costandardly filtered objects, this is completely analogous to Lemma
\ref{Lem:act_costrand_preserv}.

Let us proceed to two important examples. One is the category $\bigoplus_{n=0}^\infty
\mathcal{S}^{e,\underline{r}}_{\Ringg}(n)\operatorname{-mod}$, the other is the category
$\OCat^{-e}_\Ringg(J)$ defined below.

%

\subsection{Example: cyclotomic  Schur algebras}\label{SS_cycl_Schur}
\subsubsection{Fock space}
Pick residues $r_0,\ldots,r_{\ell-1}$ modulo $e$ and let us write
$\underline{r}$ for $(r_0,\ldots,r_{\ell-1})$. Let us define the
Fock space $\mathcal{F}^{e,\underline{r}}$ as the vector space
with basis formed by the $\ell$-multipartitions. For $n\geqslant 0$, we write $\mathcal{F}^{e,\underline{r}}(n)$
for the span of the basis elements indexed by the multipartitions of $n$.

Let $\mathcal{P}_\ell$
denote the set of $\ell$-multipartitions. We can represent a multipartition
$\lambda$ as a collection $(\lambda^{(0)},\ldots,\lambda^{(\ell-1)})$
of Young diagrams. In particular, we can talk about addable and removable
boxes in $\lambda$. We define a (shifted) residue of a box $b$ with coordinates
$(x,y)$ in the $i$th diagram as $x-y+r_i\in \Z/e\Z$. By a $j$-box we mean a
box with shifted residue $j$.

Now we can define a $\tgl_e$-action on $\mathcal{F}^{e,\underline{r}}$.
Let us write $|\lambda\rangle$ for the basis element labelled by $\lambda$.
We set $e_j|\lambda\rangle=\sum_\mu |\mu\rangle$, where the summation is
over all $\ell$-multipartitions $\mu$ obtained from $\lambda$ by removing
a $j$-box. Similarly,    $f_j|\lambda\rangle=\sum_\nu |\nu\rangle$, where
the summation is over all $\ell$-multipartitions $\nu$ obtained from $\lambda$ by adding
a $j$-box. The Cartan subalgebra of $\tgl_e$ acts as follows: $\alpha_i^\vee |\lambda\rangle=(a_i-r_i)\lambda$,
where $a_i$ (resp., $r_i$) stands for the number of addable (resp., removable) $i$-boxes,
and $d|\lambda\rangle=|\lambda||\lambda\rangle$ (and the usual generator of the  center of $\hat{\mathfrak{gl}}_e$
acts by $\ell$). The action of the unit matrix from $\mathfrak{gl}_e$ can be by an arbitrary scalar.
So we get a level $\ell$ integrable highest weight representation of $\tilde{\gl}_e$ with highest weight
(for $\tilde{\sl}_e$) equal to $\sum_{i=1}^{\ell}\omega_{r_i}$. This representation is reducible.

Let us proceed to higher level Schur algebras  $\mathcal{S}^{e,\underline{r}}_{\Ringg}(n)$
from \cite{DJM}.

\subsubsection{Cyclotomic Hecke algebras}
We assume that either $e$ is coprime to $p$ or $e=p$.
Let $\mathcal{H}^{e,\underline{r}}_\Ringg(n)$ stand for the cyclotomic Hecke algebra (for primitive
$e$th root of unity $q$)  for $e$ coprime to $p$.
Recall that this is the quotient of the (degenerate) affine Hecke algebra  $\mathcal{H}^{e,aff}_{\Ringg}(n)$
(with standard generators $X_1,\ldots,X_n,T_1,\ldots,T_{n-1}$)
by the cyclotomic relation $\prod_{i=0}^{e-1} (X_1-(1+z_i)q^{r_i})$ (for $e$ coprime to $p$)
or $\prod_{i=0}^{e-1} (X_1-r_i-z_i)$ (for $e=p$).  We also consider the level
$1$ cyclotomic KLR algebra $\mathcal{H}^{e}_\Ringg(n)$ (where the cyclotomic
relation corresponds to $r_1=0$).

Let $\Cat^{\underline{r}}_\Ringg(n)$ denote the category of finitely generated
$\mathcal{H}^{e,\underline{r}}_\Ringg(n)$-modules. Thanks to results of Rouquier, \cite[5.3.7, 5.3.8]{Rouq-2KM}
and Brundan and Kleshchev, \cite{BK_KLR}, the category $\Cat^{\underline{r}}_{\Ringg}:=
\bigoplus_{n=0}^\infty \Cat^{\underline{r}}_{\Ringg}(n)$ carries an action of the 2-Kac-Moody category
$\UC(\tgl_e)$ so that the specialization $\Cat^{\underline{r}}_{\FM}$ is a minimal $\tgl_e$-categorification
with highest weight $\sum_{i=1}^\ell \omega_{r_i}$.  We remark that the algebras
$\mathcal{H}^{e,\underline{r}}_{\Ringg}(n)$ are symmetric, in particular,
the categories $\Cat^{\underline{r}}_{\Ringg}$ carry a contravariant duality $\bullet^*$
that  fixes the projectives.

Let us discuss induction functors for Hecke algebras. Pick a partition $\mu=(\mu_1,\ldots,\mu_k)$
of $m\in \Z_{>0}$. Then we can consider the induction functor
$$\Ind^{\mu}_n: \mathcal{H}^{e,\underline{r}}_{\Ringg}(n)\operatorname{-mod}
\boxtimes \bigotimes_{i=1}^k \mathcal{H}^e_{\Ringg}(\mu_i)\operatorname{-mod}\rightarrow
\mathcal{H}^{e,\underline{r}}_{\Ringg}(n+m)\operatorname{-mod}.$$
It is intertwined by the dualities.



\subsubsection{Young modules}
To define cyclotomic Schur algebras we need Young modules
$Y^{\underline{r}}_\Ringg(\lambda)\in \Cat^{\underline{r}}_\Ringg(n)$. Set
$\underline{r}':=(r_0,\ldots,r_{\ell-2})$. We have
a natural epimorphism $\pi:\mathcal{H}^{e,\underline{r}}_{\Ringg}(m)\twoheadrightarrow
\mathcal{H}^{e,\underline{r}'}_{\Ringg}(m)$ for any $m$. The Young modules are defined inductively
(with respect to $\ell$). If $\lambda^{(\ell-1)}=\varnothing$, then we set
$Y^{\underline{r}}_\Ringg(\lambda):=\pi^* Y^{\underline{r}'}_\Ringg(\lambda)$. In the general case,
consider a composition $\mu=(\mu_0,\mu_1,\ldots,\mu_k)$, where
$\lambda':=(\lambda^{(0)},\ldots,\lambda^{(\ell-2)},\varnothing), \mu_0:=|\lambda'|,
\mu_i:=\lambda^{(\ell-1)}_i$.  Consider $\Ind_{n}^\mu(Y^{\underline{r}}_\Ringg(\lambda')\boxtimes \bigotimes_{i=1}^{k}\mathsf{1}_{\mu_i})$,
where $\mathsf{1}_{\mu_i}$ is the trivial module over the level $1$ Hecke algebra $\mathcal{H}^e_\Ringg(\mu_i)$
(i.e., the one-dimensional module where all the generators $T_i$ act by $q$).
Then there is a unique direct summand $Y_\Ringg^{\underline{r}}(\lambda)$ in this induced module
which does not appear in the modules induced from $Y^{\underline{r}}_\Ringg(\lambda'')\boxtimes \bigotimes_{i=1}^{k}\mathsf{1}_{\mu''_i}$ with $|\mu''|<|\mu|$ or $|\mu''|=|\mu|$
and $\mu''>\mu$ with respect to the dominance ordering. This is the Young module we need.
We note that $Y^{\underline{r}}_\Ringg(\lambda)\cong Y^{\underline{r}}_\Ringg(\lambda)^*$
because the induction is intertwined by the dualities.

\subsubsection{Cyclotomic Schur algebras}\label{SSS_Quiv_Schur}
We define the algebra
$$\mathcal{S}^{e,\underline{r}}_\Ringg(n)=\bigoplus_{\lambda\in \mathcal{P}_\ell(n)}\End_{\mathcal{H}^{e,\underline{r}}_{\Ringg}(n)}(Y^{\underline{r}}_\Ringg(\lambda))^{opp},$$
where we write $\mathcal{P}_\ell(n)$ for the set of $\ell$-multipartitions of $n$.
Further, consider the category $\OCat_{\Ringg,\underline{r}}^S(n):=\mathcal{S}^{e,\underline{r}}_\Ringg(n)\operatorname{-mod}$.
It is known, \cite{DJM}, that this is a highest weight category over $\Ringg$ with poset
$\mathcal{P}_\ell(n)$, where the order is introduced in the following way: we say that $\lambda\leqslant \mu$
if for all $m,k\geqslant 0$, we have $\sum_{j=0}^{m-1}|\lambda^{(j)}|+\sum_{i=1}^k|\lambda^{(m)}_i|
\leqslant \sum_{j=0}^{m-1}|\mu^{(i)}|+\sum_{i=1}^k |\mu^{(m)}_i|$. A family structure is
introduced as follows. All $i$-families (for $i=0,\ldots,e-1$,) consist of the $\ell$-multipartitions that differ
only in $i$-boxes. To get the map $\sigma^{-1}:\Lambda_a\rightarrow \{+,-\}^{m_a}$
we read all addable and removable boxes in the decreasing order (i.e., we start with the
first diagram and list boxes right to left, then we do the same with the second diagram
and so on). We write a $+$ for an addable box and a $-$ for a removable one.

\begin{Lem}\label{Lem:schur_hw_cat} The category $\OCat_{\Ringg,\underline{r}}^S:=\bigoplus_{n=0}^\infty \OCat_{\Ringg,\underline{r}}^S(n)$ carries
a highest weight $\Ringg$-deformed categorical $\tgl_e$-action categorifying $\mathcal{F}^{e,\underline{r}}$.
\end{Lem}
\begin{proof}
A hierarchy structure on $\mathcal{P}_\ell(n)$ appeared in \cite[Section 3.2]{str}
(one needs to consider $r_0,\ldots,r_{\ell-1}$ there that are very far away from each
other). A categorical $\tgl_e$-action was defined by Wada in \cite{Wada}. This implies (HWA1).
In  \cite[Section 3]{Wada} he also checked that the functors $E_i,F_i$ preserve
standardly filtered objects.
(HWA2) follows. Wada also checked that the module being categorified is
$\mathcal{F}^{e,\underline{r}}$.
\end{proof}

Let us now examine a compatibility of the action of the endomorphism $X$ on $F_i\Delta^S(\lambda)$
with a filtration by standards (where the standards appear in the decreasing order, as usual).

\begin{Lem}\label{Lem:X_filtrn1}
The action of $X$ on $F_i\Delta^S_{\Ringg}(\lambda)$ preserves the filtration.
Moreover, $X$ acts on the filtration quotient $\Delta^S(\mu)$ by $z_k$
if $\mu\setminus \lambda$ lies in $k$th diagram.
\end{Lem}
\begin{proof}
Any endomorphism of $F_i\Delta^S(\lambda)$ preserves the filtration so
we only need to prove the claim about eigenvalues. Let us change the
base to $\operatorname{Frac}(\Ringg)$.  The category $\Cat^{\underline{r}}_{\operatorname{Frac}(\Ringg)}$
is equivalent to $\bigotimes_{k=0}^{\ell-1} \Cat^{r_k}_{\operatorname{Frac}(\Ringg)}$ so that
each functor $F_i$ splits as $\bigoplus_{k=0}^{\ell-1}F_{i,k}$, where $F_{i,k}$
stands for the categorification functor on $\Cat^{r_k}_{\operatorname{Frac}(\Ringg)}$.
Note that the dot acts on $F_{i,k}$ with a single eigenvalue $z_k$. 
Now the equivalence $\Cat^{\underline{r}}_{\operatorname{Frac}(\Ringg)}\cong \bigotimes_{k=0}^{\ell-1} \Cat^{r_k}_{\Ringg}$ induces a labelling preserving highest weight
equivalence $\OCat^S_{\operatorname{Frac}(\Ringg),\underline{r}}\cong
\bigotimes_{k=0}^{\ell-1} \OCat^S_{\operatorname{Frac}(\Ringg),r_k}$.
Still, $F_i=\bigoplus_k F_{i,k}$ and the dot acts on $F_{i,k}$ with eigenvalue $z_k$.
Our claim follows.
\end{proof}

The specialized (from $\Ringg$ to $\F$; we will usually omit the subscript when we consider specializations to $\F$)
algebra $S^{e,\underline{r}}(n)$ coincides with
$$\bigoplus_{\lambda\in \mathcal{P}_\ell(n)}\End_{\mathcal{H}^{e,\underline{r}}(n)}(Y^{\underline{r}}(\lambda))^{opp}$$
It follows that the quotient functor $\pi^S:\OCat^S\twoheadrightarrow \Cat$ is fully faithful on
the projectives. We note that the projective defining $\pi^S$ (it equals $\bigoplus_{n}F^n \mathsf{1}$)
is also injective (and hence tilting), where we write $\mathsf{1}$ for the indecomposable
object in $\OCat^S(0)$.

Another property of $\OCat^S$ that we are going to use is that this category admits a naive duality functor,
see \S\ref{SSS_HW_fin}.
Namely, each $Y(\lambda)$ is self-dual with respect to the anti-involution
$\bullet^*$. So the algebra $S^{e,\underline{r}}(n)$ comes with an
anti-involution that gives rise to an equivalence $\OCat^S\cong (\OCat^S)^{opp}$.
Since $\pi^S(I^S(\lambda))=Y(\lambda)$, we see that
the anti-involution produces a labeling preserving functor.




\subsection{Example: category $\OCat^{-e}_\Ringg(J)$}\label{SSS_Soegel_O_HW_cat}
Take the Dynkin diagram $S$ of type $\tilde{A}_m$. We fix
$J\subset S^{fin}$ (where $S^{fin}$ is the system of Dynkin roots in $S$),
the system of simple Dynkin roots.
The set of such $J$ is in a natural bijection with compositions
$m=m_1+\ldots+m_\ell$ (the sum of positive integers). Recall the ring $\tilde{\Ring}=\Z_{p}[[\mathfrak{h}]]$,
where $\mathfrak{h}$ has basis $x_1,\ldots,x_m,y$ and its quotient $$\Ring=\tilde{\Ring}/(x_1-x_2,\ldots,x_{m_1-1}-x_{m_1},
x_{m_1+1}-x_{m_1+2},\ldots, x_{m-1}-x_m)$$
corresponding to $J$. The map $\Ring\rightarrow \Ringg$ given by $\Z_p\twoheadrightarrow \FM_p\hookrightarrow \FM,
x_i\mapsto z_j$ for $1+\sum_{k=0}^{j-1}m_k\leqslant i\leqslant \sum_{k=0}^{j}m_k$ and $y\mapsto 0$
gives a ring homomorphism.
So we can consider the category $\,_I\OCat^-_{\Ringg,J}(W):=\Ringg\otimes_{\Ring}\,_I\OCat^-_{\Ringg,J}(W)$.

We will produce the category $\OCat^{-e}_\Ringg(J)$ that will be the direct sum
$\bigoplus_{I\in \mathbb{O}} \,_I\OCat^-_{\Ringg,J}(W)$, where $I$ are taken from a certain
multiset $\mathbb{O}$ of finitary subsets to be produced below in this section.

\subsubsection{Combinatorics}\label{SSS_combinatorics_XJ}
Take the lattice $X:=\Z^m$ and recall the right $W$-action on $X$ from
Section \ref{SS_affine_prelim}, which is as follows: the reflection
$s_\alpha$ for $\alpha=\epsilon_i-\epsilon_j+n\delta$ sends $x=(x_1,\ldots,x_n)$
to $x'=(x_1',\ldots,x_n')$, where $x_k':=x_k$ for $k\neq i,j$, $x_i':=x_j+ne,
x_j':=x_i-ne$. For $\mathbb{O}$ we take the set of $W$-orbits on $X$
for this action (or more precisely, the set of standard parabolic stabilizers
for such orbits). Note that the set $\tilde{\Lambda}^{n,e}$ is the fundamental region for
the action of $W$ on $X$ so $\mathbb{O}$ is identified with $\tilde{\Lambda}^{n,e}$.

We want to describe the poset of $\OCat^{-e}_\Ringg(J)$ in a more elementary language.
We start with relating $\bigsqcup_{I\in \mathbb{O}}W_I\setminus W$ to the lattice
$X$ and introduce a partial order on $X$. We say that $x\leqslant x'$ if there is
a sequence $x^0=x, x^1,\ldots, x^k=x'$ and positive real roots $\alpha^0,\ldots, \alpha^{k-1}$
with $x^{i+1}=s_{\alpha^i}x^i$ for all $i$ and $x^i_k+ne-x^i_\ell>0$ if $\alpha^i=\epsilon_\ell-\epsilon_k+n\delta$.

For $I\in \mathbb{O}$, we choose a unique element $x=x^I$ in the intersection of the  corresponding
$W$-orbit with the set $\tL(n,e)$ from Section \ref{SS_affine_prelim}, i.e.,
$x^I=(x_1,\ldots,x_m)$ satisfies $x_1\leqslant x_2\leqslant\ldots\leqslant x_m\leqslant x_1+e$. The stabilizer of $x^I$
is $W_I$. The following lemma is straightforward.

\begin{Lem}\label{Lem:combin1}
A unique $W$-equivariant map $\bigsqcup_{I\in \mathbb{O}}W_I\setminus W\rightarrow X$
sending each $W_I1$ to $x^I$ is an isomorphism of posets (where the poset structure on
$W_I\setminus W$ was introduced in \S\ref{SSS_Bruhat}).
\end{Lem}

Set $X(J)$ to be the set of all elements  $(x_1,\ldots,x_m)\in X$ such that
$$x_1>x_2>\ldots>x_{m_0}, x_{m_0+1}>\ldots>x_{m_0+m_1},\ldots,x_{m_0+\ldots+m_{\ell-2}+1}>\ldots>x_m.$$
Let us introduce a partial order on  $X(J)$. For an element $x\in X$ define $x^+$ to be
a formal symbol $\emptyset$ if $x\not\in X(J)W_J$ and the element in $X(J)\cap xW_J$
else. Now define the partial order on $X(J)$ by setting $x\leqslant x'$ if there are elements
$x^0=x, x^1,\ldots,x^k=x'\in X(J)$ and positive roots    $\alpha^0,\ldots, \alpha^{k-1}$ with
$x^{i+1}=(s_{\alpha^i}x^i)^+$ for all $i$ and $x^i_k+ne-x^i_\ell>0$ if $\alpha^i=\epsilon_\ell-\epsilon_k+n\delta$.
The next lemma is also straightforward.

\begin{Lem}\label{Lem:combin2}
The isomorphism $\bigsqcup_{I\in \mathbb{O}}W_I\setminus W\xrightarrow{\sim} X$
restricts to $\bigsqcup_{I\in \mathbb{O}}\,_I\Lambda^-(J)\xrightarrow{\sim} X(J)$.
\end{Lem}

\subsubsection{Category $\OCat^{-e}_\Ringg(J)$}\label{SSS_full_affine_cat}
We set $\OCat^{-e}_\Ringg(J):=\bigoplus_{I\in \mathbb{O}} \,_I\OCat^-_{\Ringg,J}(W)$.
We remark that the complexified $K_0$ of $\OCat^{-e}_\Ringg(J)$ is naturally identified
with $\Lambda^{m_0}\C^{\Z}\otimes \Lambda^{m_1}\C^{\Z}\otimes\ldots\otimes \Lambda^{m_{\ell-1}}\C^{\Z}$.

The poset $X(J)$ carries a hierarchy structure as defined in \cite[Section 3.2]{str}.
Let us explain how the family structure looks like in this case. We present the elements
of $X(J)$ as  $\ell$-tuples of {\it virtual Young diagrams} with the number of rows equal
to $m_0,\ldots,m_{\ell-1}$ (the difference between a virtual and the usual diagram
is that the former is infinite to the left). We can present a virtual diagram with $m_i$
rows as a weakly decreasing collection of integers (that may be positive or negative).
To an element $x\in X(J)$ we assign the $\ell$-tuple corresponding to the following
collections of integers $(x_1-m_0,x_2+1-m_0,\ldots, x_{m_0}-1), (x_{m_0+1}-m_1,\ldots,
x_{m_0+m_1}-1),\ldots$.

Then we can speak about shifted contents of boxes in these diagrams: to a box $b$ in $x$th
row and $y$th column of the $i$th diagram we assign the integer $c(b):=x-y+m_i$.
So we can speak about $j$-boxes in a virtual Young diagram for $j\in \Z/e\Z$.

Again, two $\ell$-tuples of virtual Young diagrams lie in the same $i$-family
if they are different only in $i$-boxes. We impose the following order on $j$-boxes:
$b<b'$ if $c(b)<c(b')$ or $c(b)=c(b')$ and $i<i'$, where $i,i'$ are the numbers of
diagrams, where $b,b'$ appear. Then we produce maps $\Lambda_a\rightarrow \{+,-\}^{n_a}$
as in \S\ref{SSS_Quiv_Schur}.

\begin{Prop}\label{Lem:hw_cat_action}
The category $\OCat^{-e}_\Ringg(J)$ carries a highest weight categorical
$\tgl_e$-action making it into a  categorifiication
of $\Lambda^{m_0}\C^{\Z}\otimes \Lambda^{m_1}\C^{\Z}\otimes
\ldots\otimes\Lambda^{m_{\ell-1}}\C^\Z$.
\end{Prop}
\begin{proof}
Let $\DC_0(I,K)$ denote the specialization of $\DC(I,K)$ to $y=0$.
Consider the 2-category $\DC_{\mathbb{O}}:=\bigoplus_{I,K\in \mathbb{O}}\DC_0(I,K)$. By the main result of
Section \ref{sec-quantumaction-affine}, we have a 2-functor $\UC(\tilde{\gl}_e)\rightarrow
\DC_{\mathbb{O}}$.
The category $\DC_{\mathbb{O}}$ naturally
acts on  $\OCat^{-e}_\Ringg(J)$ and so $\OCat^{-e}_{\Ringg}(J)$ becomes a $\UC(\tilde{\gl}_e)$-module.
Recall, Proposition \ref{Lem:parab_action_sing}, that $\DC_{\mathbb{O}}$ preserves the standardly
filtered objects in  $\OCat^{-e}_\Ringg(J)$. Hence so does $\UC(\tilde{\gl}_e)$. This establishes (HWA1).

Recall also that on the level of $K_0$, the action of $\DC(I,K)$ on $\,_K\OCat^-_{\Ringg,J}(W)$
corresponds to the multiplication map $e_I \C W e_K \otimes e_K \C W e^-_J\rightarrow
e_I \C W e_J^-$. Using the description of the images of the generating 1-morphisms
of $\UC(\tilde{\gl}_e)$ in $\DC_{\mathbb{O}}$ (see Section \ref{subsec:e_bigger_2}
for $e>2$ and in Section  \ref{subsec-e2} for $e=2$) we see that
$K_0^{\C}(\OCat^{-e}_\Ringg(J))$ is $\tgl_e$-equivariantly identified
with $\Lambda^{m_0}\C^{\Z}\otimes \Lambda^{m_1}\C^{\Z}\otimes
\ldots\otimes\Lambda^{m_{\ell-1}}\C^\Z$. Combining this with the fact that
the $\UC(\tilde{\gl}_e)$-action preserves the standardly filtered objects, we
see that (HWA2) holds.
%
\end{proof}

We will need an analog of Lemma \ref{Lem:X_filtrn1} for the category $\OCat_{\Ringg}^{-e}(J)$.

\begin{Lem}\label{Lem:X_filtrn2}
The endomorphism $X$ of $F_i\Delta_{\Ringg}(\lambda)$
preserves a decreasing filtration by standards on $\OCat_{\Ringg}^{-e}(J)$.
Moreover, $X$ acts on the filtration quotient $\Delta^S(\mu)$ by $z_k$
if $\mu\setminus \lambda$ lies in the $k$th diagram.
\end{Lem}
\begin{proof}
Consider the Ringel dual $\OCat_{\Ringg}^{+e}(J)$ of $\OCat_{\Ringg}^{-e}(J)$.
The 2-category $\mathcal{U}(\tgl_e)$ acts on $\OCat_{\Ringg}^{+e}(J)$ and it
is enough to prove the similar claim for the action on $\OCat_{\Ringg}^{+e}(J)$
because $\OCat_{\Ringg}^{-e}(J)^{\Delta}\cong\OCat_{\Ring}^{+e}(J)^{\Delta,opp}$
and the categorical action respects this identification.

First of all, let us reduce the claim  to $J=\varnothing$.

Form the ring $\tilde{\Ringg}=\F[[x_1,\ldots,x_m]]$ that surjects onto $\Ringg$
similarly to the beginning of the present section. Then we have the category $\OCat_{\tilde{\Ringg}}^{+e}(\varnothing)$
and its base change $\OCat_{\Ringg}^{+e}(\varnothing)$. Note that $\OCat_{\Ringg}^{+e}(J)$
is a full subcategory of $\OCat_{\Ringg}^{+e}(\varnothing)$. For $x\in X$, let
$\Delta_{\tilde{\Ringg},\varnothing}(x)$ denote the standard object in
$\OCat^{+e}_{\tilde{\Ringg}}(\varnothing)$ corresponding to $x$ and $\Delta_{\Ringg,\varnothing}(x):=\Ringg\otimes_{\tilde{\Ringg}}\Delta_{\tilde{\Ringg},\varnothing}(x)$.
Then recall that, for $x\in X(J)$, $\Delta_{\Ringg}(x)$ is the maximal
quotient of $\Delta_{\Ringg,\varnothing}(x)$ lying in $\OCat_{\Ringg}^{+e}(J)$,
see \S\ref{SSS_parab_posit_deform}.
Moreover, $F_i \Delta_{\Ringg}(x)$ is the maximal quotient of
$F_i\Delta_{\Ringg,\varnothing}(x)$ lying in $\OCat_{\Ringg}^{+e}(J)$.
This follows from the observation that the decreasing filtration by
standards on  $F_i\Delta_{\Ringg}(x)$ is obtained from  that on
$F_i\Delta_{\Ringg,\varnothing}(x)$, the latter is a consequence of  Lemma \ref{Lem:hw_cat_action}.
More precisely,  $F_i\Delta_{\Ringg}(x)$ is filtered precisely by the maximal nonzero
quotients of the standards appearing in the filtration of
$F_i\Delta_{\Ringg,\varnothing}(x)$ that lie in $\OCat_{\Ringg}^{+e}(J)$.
The endomorphism $X$ of $F_i\Delta_{\Ringg}(x)$ is restricted from the similar endomorphism
of $F_i\Delta_{\Ringg,\varnothing}(x)$.  This reduces the claim to the case
$J=\varnothing$.

Now consider the base change $\OCat_{\operatorname{Frac}(\tilde{\Ringg})}(\varnothing)$.
It is a semisimple category by Lemma \ref{Lem:O_ss}. We claim that $X$ acts on
the summand $\Delta_{\operatorname{Frac}(\tilde{\Ringg}),\varnothing}(x+\epsilon_k)$
by $x_k$. This follows from the description of the image of the dot, see Section \ref{subsec-polysredux}.
This finishes the proof.
\end{proof}

\subsection{Restricted highest weight categorifications of Fock spaces}\label{SS_Fock_restr_categ}
Now we fix integers $s_0,\ldots,s_{\ell-1}$.  We introduce
a partial order on $\Part_\ell(n)$ that will depend on $s_0,\ldots,s_{\ell-1}$.
First, we introduce an order on the $i$-boxes completely analogously to Section \ref{SSS_Soegel_O_HW_cat}
and declare that boxes with different contents mod $e$ are not comparable.
Now we set $\lambda\leqslant \lambda'$ if one can order
boxes $b_1,\ldots,b_n$ of $\lambda$ and $b_1',\ldots,b_n'$ of $\lambda'$
in such a way that $b_i\leqslant b_i'$.

\subsubsection{Definition}\label{SSS_restr_categ_defn}
Let us proceed to specifying what kind of a categorical action we want. Let us consider
a category $\Cat_\Ringg$  decomposed as $\Cat_\Ringg=\bigoplus_{i=0}^N \Cat_\Ringg(i)$
such that every $\Cat_\Ringg(i)$ is equivalent to the category of
finitely generated modules over an $\Ringg$-algebra that is finitely generated
as an $\Ringg$-module.
By a restricted $\Ringg$-deformed $\tgl_e$-action on
$\Cat_\Ringg$ we mean  functors $F:\Cat_\Ringg(j-1)\rightarrow \Cat_\Ringg(j)$
and $E:\Cat_\Ringg(j)\rightarrow \Cat_\Ringg(j-1), j=0,\ldots,N$ together with endomorphisms
$X\in \End(E), T\in \End(E^2)$. The conditions are that $E,F$ are biadjoint and
$X,T$ give homomorphisms of the KLR algebra $\mathcal{H}^{e,aff}_\Ringg(d)\rightarrow \End(E^d)$
for any $d$.

Now we are ready to define a restricted $\Ringg$-deformed categorification of $\mathcal{F}^{e,\underline{s}}$
(the Fock space representation itself depends only on $\underline{s}$ modulo  $e$ but the category will depend on
$\underline{s}$ itself).
This is a highest weight $\Ringg$-linear category
$\OCat_\Ringg(\leqslant N)=\bigoplus_{i=0}^N \OCat_\Ringg(i)$ equipped
with a restricted $\Ringg$-deformed categorical $\tgl_e$-action
such that
\begin{itemize}
\item The order is that introduced above (for $s_0,\ldots,s_{\ell-1}$).
\item The functors $E_i,F_i, i\in \Z/e\Z,$ map $\Delta_\Ringg(\lambda)$ to
standardly filtered objects for all $\lambda$.
\item Under the identification of $K_0(\OCat_{\FM})$ with $\mathcal{F}^{e,\underline{s}}(\leqslant N):=
\bigoplus_{i=0}^N \mathcal{F}^{e,\underline{s}}(i)$
given by $[\Delta_\FM(\lambda)]\mapsto |\lambda\rangle$ the operators
$[E_i],[F_i]$ become the $\tgl_e$-generators $e_i,f_i$.
\item The action of $X$ on the filtration quotient $\Delta_{\Ringg}(\mu)$
in $F_j\Delta_{\Ringg}(\lambda)$ is by $z_k$, where the box $\mu\setminus \lambda$
appears in the $k$th partition $\mu^{(k)}$ of the $\ell$-multipartition $\mu$.
\end{itemize}

One example is provided by $\bigoplus_{k=0}^N\OCat^{S}_{\underline{r},\Ringg}(n)$, where $\underline{r}$
stands for $\underline{s}\mod p$. Here $s_0-s_1\gg s_1-s_2\gg \ldots \gg s_{\ell-2}-s_{\ell-1}\gg N$.

\subsubsection{Example from category $\OCat^{-e}_\Ringg(J)$}\label{SSS:cat_O_restr_catn}
We are going to produce a restricted highest weight $\Ringg$-deformed
categorical $\tgl_e$-action on a certain highest weight
subcategory of $\OCat^{-e}_\Ringg(J)$. An analogous construction appeared in
\cite[Section 5]{VV_proof}, let us briefly recall it here.

Pick $n$ such that $n<m_i$ for all $i=1,\ldots,\ell$ and set $s_i:=m_i$.
Recall, \S\ref{SSS_full_affine_cat}, that we have identified  $X(J)$, the poset
of  $\OCat^{-e}_\Ringg(J)$, with the set of virtual multipartitions. This
gives rise to the embedding $\mathcal{P}_\ell(\leqslant n)\hookrightarrow X(J)$:
we adjoin the $\ell$ components of $\lambda$ to the virtual multipartitions
corresponding to sequences of $0$'s. It is easy to see that the order on
$\mathcal{P}_{\ell}(\leqslant n)$ restricted from $X(J)$ is refined by the
order on $\mathcal{P}_\ell(n)$ introduced in \S\ref{SSS_restr_categ_defn}.


For $n'\leqslant n$, let $\mathcal{O}^{-e}_\Ringg(J,n')$ denote the highest weight subcategory in $\mathcal{O}^{-e}_\Ringg(J)$ corresponding to the ideal $\mathcal{P}_\ell(n')$. Clearly, for $n'<n$, we have
$F_i\mathcal{O}^{-e}_\Ringg(J,n')\subset \mathcal{O}^{-e}_\Ringg(J,n'+1)$. Also
$E_i \mathcal{O}^{-e}_\Ringg(J,n')\subset \mathcal{O}^{-e}_\Ringg(J,n'-1)$ for $i\neq 0$.
However,  $E_0 \mathcal{O}^{-e}_\Ringg(J,n')\not\subset \mathcal{O}^{-e}_\Ringg(J,n'-1)$.
There is a recipe to fix this called {\it categorical truncation} that was discovered in
\cite[Section 5]{str} and subsequently used in \cite{VV_proof,LW} in various contexts. We will
use it similarly to \cite[Section 5.1]{VV_proof}.

Let us recall how the categorical truncation works. Fix a {\it frozen box}, say $b$,
and consider the set $\Lambda$ of all virtual $\ell$-multipartitions $\lambda$ such that
$b$ is the smallest removable box in $\lambda$ (with  residue $0$ mod $e$, to be definite).
Then $\Lambda$ is an interval in $\mathcal{P}_{\ell}(n)$ and so $\Lambda$ gives rise to a
highest weight subquotient $\mathcal{O}^{-e}_\Ringg(J)_\Lambda$ of $\mathcal{O}^{-e}_{\Ringg}(J)$.

In \cite[Section 5.3]{str} the second named author has constructed a truncated endo-functor
$\underline{E}_0$ of $\mathcal{O}^{-e}_\Ringg(J)_\Lambda$ together with a functor
morphism $F_0^!\rightarrow \underline{E}_0\otimes_{\Ringg}\Ringg^2$ (where $F_0^!$
is a left adjoint of $F_0$) depending on an embedding
$\Ringg\hookrightarrow \Ringg^2$. Changing the base field to some finite extension of
$\F_p$, one can assume that the composition of this embedding with the projection
$\underline{E}_0\otimes_{\Ringg}\Ringg^2\twoheadrightarrow \underline{E}_0$ gives rise to $F_0^!\xrightarrow{\sim} \underline{E}_0$ (on $\mathcal{O}^{-e}_\Ringg(J,\leqslant n)$). Similarly, there is an isomorphism
$F_0^*\xrightarrow{\sim}\underline{E}_0$. So the functors $F_0, \underline{E}_0$
form a  categorical $\mathfrak{sl}_2$-action.



As in \cite[Section 5.1]{VV_proof}, we can apply the construction above repeatedly to
boxes  $(0,m_i,i)$ that are the smallest
removable $0$-boxes in any $\lambda\in \mathcal{P}_\ell(n)$ viewed as a virtual multipartition.
So the category $\OCat^{-e}_{\Ringg}(J,\leqslant n):=\bigoplus_{n'=0}^n \OCat^{-e}_{\Ringg}(J,n')$
carries a structure of a deformed highest weight categorification of the level $\ell$ Fock
space with multi-charge $(m_1,\ldots,m_\ell)$ (given by the  functors $F_i, i=0,\ldots,e-1$,
$E_j, j=1,\ldots,e-1,\underline{E}_0$).

\subsubsection{Ringel duality}\label{SSS_restr_Ringel_duality}
Now let $\OCat_{\Ringg}(\leqslant N)$ be a restricted highest weight $\Ringg$-deformed categorification of
$\mathcal{F}^{e,\underline{s}}$. Consider the Ringel dual category $\OCat_\Ringg(\leqslant N)^{\vee}$.
 By \cite[Section 7.1]{str}, the Ringel dual $\Cat^\vee$ of a category $\Cat$ with a highest weight categorical $\slf_2$-action also comes with a highest weight categorical $\slf_2$-action that makes the identification
$(\Cat^\vee)^\Delta\cong (\Cat^\Delta)^{opp}$ equivariant.
It follows that $\OCat_{\Ringg}(\leqslant N)^{\vee}$
carries a restricted $\Ringg$-deformed categorical $\tgl_e$-action. Using the identification
$(\OCat_{\Ringg}(\leqslant N)^{\vee})^\Delta\cong (\OCat_{\Ringg}(\leqslant N)^{\Delta})^{opp}$, we see
that our new restricted $\Ringg$-deformed categorical $\tgl_e$-action still categorifies
$\mathcal{F}^{e,\underline{s}}$. 
Note, however, that $\OCat_{\Ringg}(\leqslant N)^{\vee}$ is not a restricted highest
categorification of $\mathcal{F}^{e,\underline{s}}$, because the highest weight order for
$\OCat_{\Ringg}(\leqslant N)^{\vee}$ is opposite that of a highest weight categorification.
To fix this, we relabel the functors $E_i,F_i$ using the involution $i\mapsto -i$ on $\Z/e\Z$
and also change the labelling of standard objects in $\OCat_{\Ringg}(\leqslant N)^{\vee}$,
we now declare that the standard $\Delta^\vee_{\Ringg}(\lambda)$ corresponds to
$\Delta_{\Ringg}(\lambda^*)$ under the identification
$(\OCat_{\Ringg}(\leqslant N)^{\vee})^\Delta\cong (\OCat_{\Ringg}(\leqslant N)^{\Delta})^{opp}$,
where $\lambda^*:=(\lambda^{(\ell-1)t},\ldots, \lambda^{(0)t})$. The category
$\OCat_{\Ringg}(\leqslant N)^{\vee}$ becomes a restricted highest weight categorification of
$\mathcal{F}^{e,\underline{s}^*}(\leqslant N)$, where $\underline{s}^*:=(-s_{\ell-1},\ldots,-s_0)$.



\subsection{Equivalence theorems}\label{SS_equiv_thms}
\subsubsection{Main result}
The following is the main result of this section. Let $n,e,\ell,s_0,\ldots,s_{\ell-1}$ have the same meaning as before
and let $\OCat_\Ringg(\leqslant N)$ be a restricted $\Ringg$-deformed highest weight $\tgl_e$-categorification of
$\mathcal{F}^{e,\underline{s}}$.

\begin{Thm}\label{Thm:asymp_Fock_uniqueness} Assume that $e$ is coprime to $p$ or $e=p$.
Further, assume $N\gg s_0-s_1\gg s_1-s_2\gg\ldots\gg s_{\ell-2}-s_{\ell-1}\gg n$. Then we have an $\Ringg$-linear equivalence
$\OCat_\Ringg(n)\xrightarrow{\sim}\OCat^S_{\underline{s},\Ringg}(n)\operatorname{-mod}$
with $\Delta_\Ringg(\lambda)\mapsto \Delta^S_\Ringg(\lambda)$.
\end{Thm}

\begin{Rem}
We could remove the restriction on $e$ by considering quiver Schur algebras, \cite{SW}, instead
of classical ones. However, this would have complicated our story even further without giving
applications to classical representation theory.
\end{Rem}

We will apply this theorem to the category $\OCat_\Ringg(n)$ produced in \S\ref{SSS:cat_O_restr_catn}
in order to relate the multiplicities in the representation theoretic type A categories
to $p$-Kazhdan-Lusztig polynomials.

The proof of Theorem \ref{Thm:asymp_Fock_uniqueness} is going to be based on a much more
technical result. To state it we need to describe some crystal structures on $\mathcal{P}_\ell$.

\subsubsection{Crystal structures}
Recall, \cite[Section 4.4]{VV_proof}, that we have two $\tgl_e$-crystal structures on $\Part_\ell$ (collections
of maps $\tilde{f}_i,\tilde{e}_i:\Part_\ell\rightarrow \Part_\ell\sqcup\{0\}$ and
$\tilde{f}^*_i,\tilde{e}^*_i:\Part_\ell\rightarrow \Part_\ell\sqcup\{0\}$) produced from
the multi-charge $(s_0,\ldots,s_{\ell-1})$. In order to define the crystal operators
$\tilde{e}_i,\tilde{f}_i$, where $i\in \Z/e\Z$, we use the following combinatorial
recipe. We list all addable and removable boxes in the decreasing order with respect to $\leqslant$.
We write a ``$+$'' for an addable box and a ``$-$'' for a removable one. The resulting sequence
is called the {\it $i$-signature} of $\lambda$. Then we apply the following
cancellation procedure: when we have two consecutive elements $-+$ we remove them.
In the end we get a sequence, where all the $+$'s are to the left and all $-$'s are
to the right. By definition, $\tilde{e}_i\lambda$ is the multi-partition
obtained from $\lambda$ by deleting the box corresponding to the leftmost $-$
(or zero if there is no such a box). Similarly, $\tilde{f}_i\lambda$ is obtained
from $\lambda$ by adding the box corresponding to the rightmost $+$ (or is zero if
there are no $+$'s).

This crystal has a categorical meaning: it is the crystal of a highest weight categorification
of $\mathcal{F}^{e,\underline{s}}$, by \cite{cryst}.

The ``dual'' crystal operators $\tilde{f}_i^*, \tilde{e}_i^*$ are defined in the similar
fashion but instead of removing the consecutive $-+$ we remove consecutive $+-$.
We say that $\lambda$ is {\it singular}, resp., {\it cosingular}, if $\tilde{e}_i\lambda=0$
for all $i$, resp., $\tilde{e}_i^* \lambda=0$ for all $i$.

Let us give an example of computation. Let $e=2, s_0=11, s_1=0$.

\begin{picture}(120,40)
\put(2,1){\line(0,1){14}}
\put(9,1){\line(0,1){14}}
\put(16,1){\line(0,1){14}}
\put(2,1){\line(1,0){14}}
\put(2,8){\line(1,0){14}}
\put(2,15){\line(1,0){14}}
\put(42,1){\line(0,1){28}}
\put(49,1){\line(0,1){28}}
\put(56,1){\line(0,1){7}}
\put(63,1){\line(0,1){7}}
\put(42,1){\line(1,0){21}}
\put(42,8){\line(1,0){21}}
\put(42,15){\line(1,0){7}}
\put(42,22){\line(1,0){7}}
\put(42,29){\line(1,0){7}}
\put(65,3){$b_4$}
\put(17,3){$b_1$}
\put(11,10){$b_2$}
\put(4,17){$b_3$}
\put(44,24){$b_5$}
\put(58,3){$b_1'$}
\put(51,10){$b_2'$}
\put(44,31){$b_3'$}
\end{picture}

The addable/removable boxes with residue $1$ written in the decreasing order are
$b_1,b_2,b_3,b_4,b_5$, which gives the signature $+-++-$. The reduced signature
for the usual crystal is $++-$. We will have $\tilde{f}_1\lambda=((2,2),(4,1^3))$
and $\tilde{e}_1\lambda=((2,2),(3,1^2))$. For the dual crystal, the reduced signature
is $+$ (with the second plus in the original signature surviving). So $\tilde{e}_1^*\lambda=0$
and $\tilde{f}_1^*\lambda=((2,2,1),(3,1^3))$.

The addable/removable boxes with residue $0$ written in the decreasing order are
$b_1',b_2',b_3'$. So the signature is $-++$. The reduced signature for the
usual crystal is $+$. We have $\tilde{e}_0\lambda=0, \tilde{f}_0\lambda=((2,2),(3,1^4))$.
The reduced signature for the dual crystal is $-++$. We have $\tilde{e}_0^*\lambda=((2,2),(2,1^3))$
and $\tilde{f}_0^*\lambda=((2,2),(3,2,1^2))$.

\subsubsection{Technical results}
Let $\OCat^i_\Ringg:=\bigoplus_{j=0}^N \OCat^i_\Ringg(j), i=1,2,$ be two restricted $\Ringg$-deformed
highest weight categorifications of $\mathcal{F}^{e,\underline{s}}$.

Consider the projectives $P^i_\Ringg:=F^n \Delta_\Ringg^i(\varnothing)
\in \OCat_\Ringg^i(n)$, where $n\leqslant N$.
The proof of the next lemma repeats that of \cite[Proposition 5.1(2)]{VV_proof}.

\begin{Lem}\label{Lem:End_Hecke_iso}
We have $\operatorname{End}_{\OCat^i_\Ringg(n)}(P^i_\Ringg)^{opp}\cong \mathcal{H}^{e,\underline{s}}_{\Ringg}(n)$.
\end{Lem}

Following \cite[Section 7.1]{VV_proof}, we also consider the bigger projectives $\bar{P}^i_\Ringg$ defined
as follows. If $e>2$, then we set $\bar{P}^i_\Ringg:=P^i_\Ringg\oplus \bigoplus_\lambda
F^{n-1}P^i_\Ringg(\lambda)$, where the sum is taken over all singular (with respect
to the $\tgl_e$-crystal on the set of partitions) multipartitions $\lambda$ of $1$.
If $e=2$, then we set $\bar{P}^i_\Ringg:=P_\Ringg\oplus \bigoplus_{\lambda} F^{n-1} P^i_\Ringg(\lambda)
\oplus F^{n-2}P_\Ringg(\nu)$, where the range of $\lambda$ is the same as before and
$\nu=(\varnothing,\ldots,\varnothing, (2))$.

We set $\Cat_\Ringg(n):=\operatorname{End}_{\OCat^i_\Ringg(n)}(P^i_\Ringg)^{opp}\operatorname{-mod}$.
Let $\bar{\Cat}_\Ringg^i(n)$
denote the category of right modules over $\operatorname{End}_{\OCat^i_\Ringg(n)}(\bar{P}^i_\Ringg)^{opp}$.
We have natural quotient functors $\bar{\pi}^i_\Ringg:\OCat^i_\Ringg\twoheadrightarrow \bar{\Cat}^i_\Ringg\operatorname{-mod},
\underline{\pi}^i_\Ringg:\bar{\Cat}^i_\Ringg\twoheadrightarrow \Cat_\Ringg$ and $\pi^i_\Ringg:\OCat^i_\Ringg\twoheadrightarrow
\Cat_\Ringg$. Note that $\pi^i_\Ringg=\underline{\pi}^i_\Ringg\circ \bar{\pi}^i_\Ringg$.

Let $\K$ denote the fraction field of $\Ringg$. We are going to consider the base changes of
the categories of interest to $\K$ and to $\FM$ (we will replace the subscript $\Ringg$
with $\K,\FM$ in these cases).

We introduce the following assumptions.

\begin{itemize}
\item[(I)]  We have an equivalence $\bar{\Cat}^1_\Ringg(n)\xrightarrow{\sim}\bar{\Cat}^2_\Ringg(n)$
intertwining the functors $\underline{\pi}^i_\Ringg$ and preserving the labels of the
indecomposable projective objects.
\item[(F$_{-1}$)] The functor $\pi^i_\FM:\OCat^i_{\FM}(n)\twoheadrightarrow \Cat_{\FM}(n)$ is {\it $(-1)$-faithful}
(=faithful on standardly filtered objects) and also faithful on costandard objects.
\item[(F$_0$)] The functor $\bar{\pi}^i_{\K}:\OCat^i_{\K}(n)
\twoheadrightarrow \bar{\Cat}^i_{\K}(n)$ is {\it $0$-faithful}(=fully faithful on standardly
filtered objects).
\item[(S)] Further, under an equivalence from (I), we have
$\bar{\pi}^1_{\K}(\Delta^1_{\K}(\lambda))\cong \bar{\pi}^2_{\K}(\Delta^2_{\K}(\lambda))$
for all $\lambda\in \mathcal{P}_\ell(n)$.
\end{itemize}

Then we have the following results.

\begin{Thm}\label{Thm:equi_techn}
Suppose that (F$_{-1}$),(F$_{0}$),(I) and (S) hold. Identify $\bar{\Cat}^1_\Ringg,\bar{\Cat}^2_\Ringg$
using the equivalence in (I). Then there is a highest weight equivalence
$\OCat^1_\Ringg(n)\xrightarrow{\sim}\OCat^2_\Ringg(n)$ that intertwines the functors
$\bar{\pi}^1_\Ringg,\bar{\pi}^2_\Ringg$. On the categories of standardly filtered objects,
the equivalence is given by $(\bar{\pi}^2_\Ringg)^*\circ \bar{\pi}^1_\Ringg$.
\end{Thm}

\begin{Prop}\label{Prop:unique_check_conditions} The assumptions
(F$_{-1}$),(F$_0$),(I),(S) hold in the following two cases:
\begin{enumerate}
\item when $\OCat_\Ringg(\leqslant N)=\OCat^S_{\underline{s},\Ringg}(\leqslant N)^\vee$,
\item when $N\gg s_1-s_2\gg s_2-s_3\gg\ldots \gg s_{\ell-1}-s_\ell\gg n$.
\end{enumerate}
\end{Prop}

We will prove Theorem \ref{Thm:equi_techn} in the next section and Proposition
\ref{Prop:unique_check_conditions} in the remaining five sections. Now let us
finish the proof of Theorem \ref{Thm:asymp_Fock_uniqueness} modulo these two results.

\begin{proof}[Proof of Theorem \ref{Thm:asymp_Fock_uniqueness}]
Recall, \S\ref{SSS_restr_Ringel_duality}, that $\OCat_\Ringg(\leqslant N)^\vee$ is  a restricted $\Ringg$-deformed
highest weight categorification of $\mathcal{F}^{e,\underline{s}^*}$.
By Theorem \ref{Thm:equi_techn} and Proposition \ref{Prop:unique_check_conditions},
we see that $\OCat_{\Ringg}(\leqslant n)^\vee\cong \OCat^{S}_\Ringg(\leqslant n)^\vee$. From
here we deduce the required equivalence $\OCat_\Ringg(\leqslant n)\cong \OCat^S_\Ringg(\leqslant n)$.
\end{proof}

\subsection{Proof of Theorem \ref{Thm:equi_techn}}
Before proceeding to the proof of Theorem \ref{Thm:equi_techn} (that is inspired by that
of \cite[Theorem 3.4]{VV_proof}), we will need the following three technical lemmas.

\begin{Lem}\label{Lem:proj_incl}
Let $\tilde{P}^i_\Ringg$ be a projective object in $\OCat^i_\Ringg(n)$. Then there is
$d>0$ and an inclusion $\tilde{P}^i_\Ringg\hookrightarrow (P^i_\Ringg)^{\oplus d}$
with standardly filtered cokernel.
\end{Lem}
The proof repeats that of \cite[Proposition 3.7]{VV_proof}.

Let us order the labels $\lambda_1,\ldots,\lambda_m$ for $\OCat^i_\Ringg(n)$ in such a
way that $\lambda_i>\lambda_j\Rightarrow i<j$. For a standardly filtered
object $M\in \OCat^i_\Ringg(n)$ and a positive integer $k$, define $M_{\leqslant k}$
as the unique standardly filtered subobject in $M$ filtered by $\Delta^i_\Ringg(\lambda_j), j\leqslant k,$
and such that $M/M_{\leqslant k}$ is filtered by $\Delta^i_\Ringg(\lambda_j)$ with $j>k$.

\begin{Lem}\label{Lem:coinc_image}
The following is true.
\begin{enumerate}
\item Let $M^i\in \OCat^i_\Ringg(n)$ be a standardly filtered object.
If $\bar{\pi}^1_\Ringg(M^1)\cong \bar{\pi}^2_\Ringg(M^2)$, then
$\bar{\pi}^1_\Ringg(M^1_{\leqslant k})\cong \bar{\pi}^2_\Ringg(M^2_{\leqslant k})$.
\item We have $\bar{\pi}^1_\Ringg(\Delta^1_\Ringg(\lambda))\cong
\bar{\pi}^2_\Ringg(\Delta^2_\Ringg(\lambda))$ for any $\lambda\in \mathcal{P}_\ell(n)$.
\end{enumerate}
\end{Lem}
\begin{proof}
First of all, (S) means that (2) is true if we replace $\Ringg$ with $\K$.
Let us deduce that (1) is true with $\Ringg$ replaced with $\K$.
The image of the natural homomorphism $$\bar{\pi}^i_\K(\Delta_\K^i(\lambda_1))\otimes_\K\Hom_\K(\Delta_\K^i(\lambda_1),M^i_\K)\rightarrow
\bar{\pi}^i_\K(M^i_\K)$$ coincides with $\bar{\pi}^i_\K((M^i_\K)_{\leqslant 1})$ because $\Hom(\Delta^i_\K(\lambda_1),\Delta^i_\K(\lambda_i))=0$ and because $\bar{\pi}^i_\K$
is 0-faithful. The cokernel coincides with $\bar{\pi}^i_\K(M^i_\K/(M^i_\K)_{\leqslant 1})$.
Then we get (1) for $\K$ by induction on $k$.

Now let us deduce (1) for $\Ringg$. Let us first check that $M^i_{\leqslant k}=M^i\cap (M^i_\K)_{\leqslant k}$.
The inclusion $M^i_{\leqslant k}\subset M^i\cap (M^i_{\K})_{\leqslant k}$ is straightforward.
Clearly, $M^i/(M^i \cap (M^i_{\K})_{\leqslant k})$ is torsion-free over $\Ringg$. On the other hand,
we have $(M^i/M^i_{\leqslant k})_{\K}=M^{i}_{\K}/(M^i_{\K}\cap (M^i_{\K})_{\leqslant k})$.
This implies $M^i_{\leqslant k}=M^i\cap (M^i_\K)_{\leqslant k}$.

Since $\bar{\pi}^i_\Ringg$ is an exact functor (the multiplication by
an idempotent), we get
$$\bar{\pi}^i_\Ringg(M^i_{\leqslant k})=\bar{\pi}^i_\Ringg(M^i)\cap \bar{\pi}^i_\Ringg((M^i_\K)_{\leqslant k})=
\bar{\pi}^i_\Ringg(M^i)\cap \bar{\pi}^i_\K((M^i_\K)_{\leqslant k}).$$
From here we get (1) for $\Ringg$ as in \cite[Lemma 4.48]{rouqqsch}. (2) follows as
in that lemma as well.
\end{proof}

We consider the categories
$\OCat^{i\Delta}_\Ringg(n)$ of all standardly filtered objects in $\OCat^i_\Ringg(n)$. We also
consider a category $\Cat^\Delta_\Ringg(n)$ defined as the full subcategory in $\Cat_\Ringg(n)$
consisting of all modules $M$ that admit a filtration $M=M_k\supset M_{k-1}\supset M_{k-2}\supset\ldots
\supset M_0=\{0\}$ such that $M_i/M_{i-1}$ is the direct sum of several copies of
$\bar{\pi}^1_\Ringg(\Delta^1_\Ringg(\lambda_i))\cong \bar{\pi}^2_\Ringg(\Delta^2_\Ringg(\lambda_i))$.

\begin{Lem}\label{Lem:stand_filt_emb}
We have the following.
\begin{enumerate}
\item
The functor $\bar{\pi}^i_\Ringg:\OCat^{i\Delta}_\Ringg(n)\rightarrow \bar{\Cat}^\Delta_\Ringg(n)$
is a fully faithful embedding. Moreover, $(\bar{\pi}^i_R)^*$ is the left inverse of
$\bar{\pi}^i_\Ringg|_{\OCat^{i\Delta}_\Ringg(n)}$.
\item If the images of $\OCat^{1}_\Ringg(n)\operatorname{-proj},
\OCat^2_\Ringg(n)\operatorname{-proj}$ in $\bar{\Cat}^{\Delta}_\Ringg(n)$ coincide, then
Theorem \ref{Thm:equi_techn} holds.
\end{enumerate}
\end{Lem}
\begin{proof}
The claim that $\bar{\pi}^i_{\Ringg}$ is $0$-faithful is a special case of \cite[Proposition 2.18]{RSVV}.
The claim that $(\bar{\pi}^i_\Ringg)^*\circ \bar{\pi}^i_\Ringg$ is the identity on the standardly filtered
objects is equivalent to  $\bar{\pi}^i_\Ringg$ being $0$-faithful, see \cite[Proposition 4.40]{rouqqsch}.

Let us prove (2). The functors $\bar{\pi}^i_{\Ringg}, i=1,2,$ are fully faithful on the projective
objects by (1). So there is an equivalence $\OCat^1_{\Ringg}(n)\operatorname{-proj}\xrightarrow{\sim}
\OCat^2_{\Ringg}(n)\operatorname{-proj}$ intertwining the functors $\bar{\pi}^i_{\Ringg}$. This equivalence extends to
an equivalence $\OCat^1_{\Ringg}(n)\xrightarrow{\sim}\OCat^2_{\Ringg}(n)$ as needed in Theorem \ref{Thm:equi_techn}.
\end{proof}

\begin{proof}[Proof of Theorem \ref{Thm:equi_techn}]
To simplify the notation, we write $\OCat^{i\Delta}_\Ringg$ for $\OCat^{i\Delta}_\Ringg(n)$,
and $\Cat^{i\Delta}_\Ringg$ for $\Cat^{i\Delta}_\Ringg(n)$.

{\it Step 1}. Let $M\in \OCat^{i\Delta}_\Ringg$ be such that
$M$ is included into $(P^i_\Ringg)^{\oplus m}$ with standardly filtered cokernel.  Then $R^1(\bar{\pi}_\Ringg^i)^*\circ \bar{\pi}_\Ringg^i M=0$. This is proved as
\cite[Lemma 3.6]{VV_proof}.

{\it Step 2}. Let $M$ be such as in the previous step and $N\in \OCat^{i\Delta}_\Ringg$.
We claim that if $Q\in \Cat_\Ringg$ is included into an exact sequence
$0\rightarrow \bar{\pi}^i_\Ringg M\rightarrow Q\rightarrow \bar{\pi}^i_\Ringg N\rightarrow 0$, then $Q\in \bar{\pi}^i_\Ringg(\OCat^{i\Delta}_\Ringg)$ and hence $(\bar{\pi}^i_\Ringg)^*Q\in \OCat^{i\Delta}_\Ringg$.
Indeed, we have an exact sequence $0\rightarrow M\rightarrow
(\bar{\pi}^i_\Ringg)^* Q\rightarrow N\rightarrow R^1(\bar{\pi}^i_\Ringg)^*\circ \bar{\pi}^i_\Ringg(M)$ in $\OCat^i_\Ringg$. From this sequence and the previous step, we see that $(\bar{\pi}^i_\Ringg)^*
Q\in \OCat^{i\Delta}_\Ringg$. It follows that $Q\in \bar{\pi}^i_\Ringg(\OCat^{i\Delta}_\Ringg)$.

{\it Step 3}.  Suppose that there is an equivalence $\iota_1:(\OCat^1_\Ringg)^{\Delta,\leqslant r}\xrightarrow{\sim}
(\OCat^2_\Ringg)^{\Delta,\leqslant r}$ intertwining the functors $\bar{\pi}^1_\Ringg,\bar{\pi}^2_\Ringg$. Since
$\bar{\pi}^1_\Ringg,\bar{\pi}^2_\Ringg$ are $0$-faithful, we see that $\iota_1:=(\bar{\pi}^2_\Ringg)^*\circ\bar{\pi}^1_\Ringg$. Let $\iota_2$ denote the
inverse of $\iota_1$ so that $\iota_2=(\bar{\pi}^1_\Ringg)^*\circ\bar{\pi}^2_\Ringg$.
In the subsequent steps we will see that there is an equivalence
$(\OCat^1_\Ringg)^{\Delta,\leqslant r+1} \xrightarrow{\sim} (\OCat^2_\Ringg)^{\Delta,\leqslant r+1}$
intertwining $\bar{\pi}^1_\Ringg,\bar{\pi}^2_\Ringg$. The claim of the theorem will follow by induction
with trivial base $r=-1$.

{\it Step 4}. We note that if $\psi: M\rightarrow M'$ be a monomorphism in $\OCat^{i\Delta}_\Ringg(n)$ with
standardly filtered cokernel, then, for any $j$, $\psi|_{M_{\leqslant j}}: M_{\leqslant j}\rightarrow
M'_{\leqslant j}$ also has standardly filtered cokernel (that coincides with $(M'/M)_{\leqslant j}$).
For an object $N\in \Cat_\Ringg^{\Delta}(n)$, let $N_{\leqslant j}$ denote the maximal
subobject filtered by $\bar{\pi}^i_\Ringg(\Delta^i_\Ringg(\lambda_s))$ with $s\leqslant j$
such that the quotient is filtered by $\bar{\pi}^i_\Ringg(\Delta^i_\Ringg(\lambda_s))$ with $s> j$.
Since $\Hom(\bar{\pi}^i_\Ringg(\Delta^i_\Ringg(\lambda_s)),\bar{\pi}^i_\Ringg(\Delta^i_\Ringg(\lambda_t)))=0$
for all $s<t$, the object $N_{\leqslant j}$ is well defined. As in the proof
of Lemma \ref{Lem:coinc_image}, $\bar{\pi}^i_\Ringg(M_{\leqslant j})=\bar{\pi}^i_\Ringg(M)_{\leqslant j}$.

{\it Step 5}. Let $a,b$ be two different elements of $\{1,2\}$.
Let $M^a\in \OCat^{a\Delta}_{\Ringg}$ be an object that can
be included into $(P^a_\Ringg)^{\oplus m}$ with standardly filtered cokernel.
We claim that $\iota_a(M^a_{\leqslant r})\in \OCat^{b}_\Ringg$ is included
into $(P^b_\Ringg)^{\oplus m}$ with standardly filtered cokernel. Indeed, from the inclusion
$M^a\subset (P^a_\Ringg)^{\oplus m}$ with standardly filtered cokernel,
we deduce the inclusion  $M^a_{\leqslant r}\subset
(P^a_{\Ringg,\leqslant r})^{\oplus m}$ with standardly
filtered cokernel. By Step 4,   $\bar{\pi}^a_\Ringg(P^a_{\Ringg,\leqslant r})=
\bar{\pi}_\Ringg^a(P^a_\Ringg)_{\leqslant r}=
\bar{\pi}_\Ringg^b(P^b_\Ringg)_{\leqslant r}=\bar{\pi}^b_\Ringg(P^b_{\Ringg,\leqslant r})$.
By Step 3, $\iota_a(M^2_{\leqslant r}) \hookrightarrow (P^b_\Ringg)^{\oplus m}_{\leqslant r}$
with standardly filtered cokernel. Our claim follows.

{\it Step 6}. In particular, let $\tilde{P}^a_\Ringg\in \OCat^a_\Ringg\operatorname{-proj}$.
By Lemma \ref{Lem:stand_filt_emb}, $\tilde{P}^a_\Ringg$ is included
to $(P^a_\Ringg)^{\oplus m}$ with standardly filtered cokernel. From Step 5, we see
that $\iota_a(\tilde{P}^a_{\Ringg,\leqslant r})\in \OCat^{b\Delta}_\Ringg$ is included
into $(P^b_\Ringg)^{\oplus m}$ with standardly filtered cokernel.

{\it Step 7}. Let us prove now that, for any projective $\tilde{P}^b_\Ringg\in
\OCat^b_\Ringg$, we have  $(\pi^a_\Ringg)^*\circ \pi^b_{\Ringg}((\tilde{P}^b_\Ringg)_{\leqslant r+1})
\in \OCat^{a\Delta,\leqslant r+1}_\Ringg$. We have an exact sequence
$0\rightarrow \bar{\pi}^b_\Ringg(\tilde{P}^b_{\Ringg,\leqslant r})\rightarrow \bar{\pi}^b_\Ringg(\tilde{P}^b_{\Ringg,\leqslant r+1})\rightarrow
\bar{\pi}^b_\Ringg(\Delta^b_\Ringg(\lambda_{r+1})^{\oplus m})\rightarrow 0$ for some $m\geqslant 0$.
By Step 6, the object $M:=\iota_b(\tilde{P}^b_{\Ringg,\leqslant r})$ satisfies the condition
of Step 1. Our claim follows now from Step 2.

{\it Step 8}. We claim that $(\bar{\pi}^a_\Ringg)^*\circ \bar{\pi}^b_\Ringg(P^b_\Ringg(\lambda_{r+1}))$
(standardly filtered by Step 7) is projective in $\OCat^a_\Ringg$. First of all, note that
$$P^a_\Ringg(\lambda_{r+1})\oplus \tilde{P}^a_\Ringg\twoheadrightarrow (\bar{\pi}^a_\Ringg)^*\circ \bar{\pi}^b_\Ringg(P^b_\Ringg(\lambda_{r+1})),$$
where $\tilde{P}^a_\Ringg$ is a projective cover  of $\ker[(\bar{\pi}^a_\Ringg)^*\circ \bar{\pi}^b_\Ringg(P^b_\Ringg(\lambda_{r+1}))\twoheadrightarrow \Delta^a_\Ringg(\lambda_{r+1})]$.
The projective $\tilde{P}^a_\Ringg$ lies in $(\OCat^a_\Ringg)^{\Delta,\leqslant r}$.
We have an exact sequence
\begin{equation}\label{eq:ex_seq_1} 0\rightarrow M^a\rightarrow P^a_\Ringg(\lambda_{r+1})\oplus \tilde{P}^a_\Ringg\rightarrow
(\bar{\pi}^a_\Ringg)^*\circ \bar{\pi}^b_\Ringg(P^b_\Ringg(\lambda_{r+1}))\rightarrow 0,\end{equation}
where $M^a\in (\OCat^a_\Ringg)^{\Delta\leqslant r}$.
Applying $\bar{\pi}^a_\Ringg$ to (\ref{eq:ex_seq_1}), we get an exact sequence
\begin{equation}\label{eq:ex_seq_2} 0\rightarrow \bar{\pi}^a_\Ringg(M^a)\rightarrow \bar{\pi}^a_\Ringg(P^a_\Ringg(\lambda_{r+1})\oplus \tilde{P}^a_\Ringg)
\rightarrow \bar{\pi}^b_\Ringg(P^b_\Ringg(\lambda_{r+1}))\rightarrow 0.\end{equation}
The exact sequence (\ref{eq:ex_seq_1}) combined with Step 7 imply that the object $M^a$ is
included into $(P^a_\Ringg)^{\oplus m}$ with standardly filtered cokernel.
By Step 5, the object $M^b:=(\bar{\pi}^b_\Ringg)^*\circ\bar{\pi}^a_\Ringg(M^a)$ is included
into $(P^b_\Ringg)^{\oplus m}$  with standardly filtered cokernel. By
Step 1, $R^1\bar{\pi}_\Ringg^{b*}\circ\bar{\pi}^a_\Ringg(M^a)=0$. By Step 2,
$(\bar{\pi}^b_\Ringg)^*\circ\bar{\pi}^a_\Ringg(P^a_\Ringg(\lambda_{r+1})\oplus \tilde{P}^a_\Ringg)$
is in $(\OCat^b_\Ringg)^{\Delta,\leqslant r+1}$ and we have the following exact sequence
\begin{equation}\label{eq:ex_seq_3}
0\rightarrow M^b\rightarrow (\bar{\pi}^b_\Ringg)^*\circ\bar{\pi}^a_\Ringg(P^a_\Ringg(\lambda_{r+1})\oplus \tilde{P}^a_\Ringg)
\rightarrow P^b_\Ringg(\lambda_{r+1})\rightarrow 0.
\end{equation}
This sequence splits, and hence so does (\ref{eq:ex_seq_2}).   Since the functor $\bar{\pi}^a_{\Ringg}$ is
0-faithful, the sequence (\ref{eq:ex_seq_1}) splits as well. The claim in the beginning of this step follows.

{\it Step 9}. Step 8 implies that the image of $\OCat^a_\Ringg\operatorname{-proj}\cap \OCat^{a,\Delta,\leqslant r+1}_\Ringg$
under $\bar{\pi}^a_\Ringg$ is contained in the image of $\OCat^b_\Ringg\operatorname{-proj}\cap \OCat^{b\Delta,\leqslant r+1}_\Ringg$.
So these images coincide. We get an equivalence between these two categories that then extends
to $\OCat^{a\Delta,\leqslant r+1}_\Ringg\xrightarrow{\sim} \OCat^{b\Delta,\leqslant r+1}_\Ringg$.
\end{proof}

In the following five sections we will prove Proposition \ref{Prop:unique_check_conditions}.

\subsection{Checking $(-1)$-faithfulness}\label{SS_-1_faith}
In this section we will prove (F$_{-1}$) for the categories in Proposition \ref{Prop:unique_check_conditions}.

\subsubsection{Case $N\gg s_0-s_1\gg s_1-s_2\gg\ldots\gg s_{\ell-2}-s_{\ell-1}\gg n$}
In \cite[Section 6]{VV_proof} we have found a sufficient condition for
the functor $\pi:\OCat(n)\twoheadrightarrow \Cat(n)$ to be $(-1)$-faithful
(on standardly filtered and on costandardly filtered objects). This criterium is
as follows.

In \cite[Section 6.1]{VV_proof}, for $w\in \hat{S}_e$ and a singular
$\lambda$ (resp., cosingular $\mu$) we have defined
$w\lambda\in \mathcal{P}_\ell$ (resp., $w^*\mu\in \mathcal{P}_\ell$). The rule is that if
$\lambda'$ satisfies $\tilde{e}_i \lambda'=0$, then $\sigma_i\lambda':=\tilde{f}_i^{d}\lambda'$,
where we write $\sigma_i$ for the simple reflection corresponding to $i\in \{0,\ldots,e-1\}$
and $d$ is maximal with the property that $\tilde{f}_i^d\lambda'\neq 0$.
As we have seen in {\it loc. cit.}, this gives rise to a well-defined element $w\lambda$
when $\lambda$ is singular. The element $w^*\mu$ is defined similarly
using the dual crystal structure.

Now, following \cite[Section 6.1]{VV_proof},  for a singular
$\lambda$ and cosingular $\mu$, consider the following condition:
\begin{itemize}
\item[$(\mathfrak{C}_{\lambda\mu})$] There is $w\in W$ such that $w\lambda\not\leqslant w^*\mu$.
\end{itemize}

\begin{Prop}\label{Prop:comb_cond1}
Suppose $N\gg s_0-s_{1}\gg s_1-s_2\gg\ldots \gg s_{\ell-2}-s_{\ell-1}\gg n$.
Then $(\mathfrak{C}_{\lambda\mu})$ is satisfied for any singular $\lambda$ and
any cosingular $\mu$ with $|\lambda|=|\mu|\leqslant n$.
\end{Prop}
\begin{proof}
As in \cite[Section 6.3]{VV_proof}, we will use  elements $w$ of the form
$C_{a,m}:=\sigma_{a+1-m}\ldots \sigma_{a-1}\sigma_a$, where $\sigma_i$ denotes the simple
reflection in $\hat{S}_e$ corresponding to $i$.  Note that if $\lambda$ is
singular, then so is $\lambda^{(\ell-1)}$ (in the crystal on $\Part_1$). It follows
that $\lambda^{(\ell-1)}$ is divisible by $e$ (partwise).
Consider the following cases:

{\it Case 1}: $\mu^{(\ell-1)}\neq \varnothing$. Let $k$ denote the number of rows in $\mu^{(\ell-1)}$.
We set $a=s_{\ell-1}-k$ (this is the content of the smallest addable box in
$\mu$). We will apply an element $C_{a,m}$ to $\lambda$ and $\mu$.
The smallest box in $C_{a,m}\mu$ has content $a+1-m$ (each time we apply
a reflection, we add a box in the first column of $\mu^{(\ell-1)}$). On the other hand,
the $(C_{a,m}\lambda)^{(\ell-1)}$ cannot have more that $m$ rows (we add boxes
there starting from the first row). So when $m$ is bigger than  the number of rows
in $\lambda^{(\ell-1)}$, we get  $C_{a,m}\lambda\not\leqslant C_{a,m}^*\mu$.

{\it Case 2}: $\mu^{(\ell-1)}=\varnothing$ but the largest addable box in
$\lambda^{(\ell-1)}$ does not survive in the reduced $s_{\ell-1}$-signature of
$\lambda$. Pick $a:=s_{\ell-1}$. When we apply $s_{a}$ to $\lambda$ we do not affect
$\lambda^{(\ell-1)}$.  Now take $m$ bigger than the number of rows in
$\lambda^{(\ell-1)}$. Then the smallest  box in $C_{a,m}^*\mu$ is $(1,m,\ell-1)$,
while the number of rows in $(C_{a,m}\lambda)^{(\ell-1)}$ is less than
$m$. We get $C_{a,m}\lambda\not\leqslant C_{a,m}^*\mu$.

{\it  Case 3}: Now suppose that the largest addable box in $\lambda^{(\ell-1)}$
survives in the reduced $s_{\ell-1}$-signature. This is equivalent to say that
$\lambda'=(\lambda^{(0)},\ldots,\lambda^{(\ell-2)})$ is singular. In particular,
$\lambda^{(\ell-2)}$ is divisible by $e$. Set $\nu_m:=C_{0,m}\varnothing$
(in the crystal of $\mathcal{P}_1$). Take $a=s_{\ell-2}$. Then $(C_{a,m}\lambda)^{(\ell-1)}=\nu_{m-k}+\lambda^{(\ell-1)}$
(component-wise addition), and $(C_{a,m}^*\mu)^{(\ell-1)}=\nu_{m-k}$, where $k$
is the remainder of $s_{\ell-1}-s_{\ell-2}$ under the division by $e$, see
\cite[Section 6.3]{VV_proof}. If $\mu^{(\ell-2)}\neq \varnothing$ or the largest addable box of $\lambda^{(\ell-2)}$
does not survive in the reduced signature, then the smallest box in $(C_{a,m}^*\mu)^{(\ell-2)}$
will be  smaller than any box in $(C_{a,m}\lambda)^{(i)}, i<\ell-1$. The content of the smallest
box in $(C_{a,m}^*\mu)^{(\ell-2)}$ is smaller than or equal to $s_{\ell-2}-m$. On the other hand,
the boxes in $(C_{a,m}\lambda)^{(\ell-1)}\setminus (C_{a,m}^*\mu)^{(\ell-1)}=(\nu_m+\lambda^{(\ell-2)})\setminus \nu_m$
all have shifted content that is larger than $s_{\ell-1}+\lfloor m/(k'e)\rfloor-k'$, where $k'$
is the number of rows in $\lambda^{(\ell-1)}$.
We conclude that $C_{a,m}\lambda\not\leqslant C_{a,m}^*\mu$ for $m>e(s_{\ell-2}-s_{\ell-1}+k')$.

{\it Case 4}:
So we only need to consider the case
when $(\lambda^{(0)},\ldots,\lambda^{(\ell-3)})$ is singular. We continue as before. At some point
we get $\mu^{(i)}\neq 0$ or $(\lambda^{(0)},\ldots,\lambda^{(i-1)})$ is not singular. Then we
are done.
\end{proof}

Using results of \cite[Section 6.2]{VV_proof}, we deduce the following corollary.

\begin{Cor}\label{Cor:-1_faithf}
Let $\OCat(\leqslant N)=\bigoplus_{j=0}^N \OCat(j)$ be a restricted categorification of
$\mathcal{F}^{e,\bs}$. If $n,\bs,N$ satisfy the conditions of
Proposition \ref{Prop:comb_cond1}, then the functor $\pi:\OCat(n)\twoheadrightarrow
\Cat(n)$ is $(-1)$-faithful on the standardly filtered and on the costandardly
filtered objects.
\end{Cor}

\subsubsection{Ringel dual Schur category}
\begin{Prop}\label{Prop:-1_faith_Schur}
The functor $\pi^{S\vee}$ is faithful on standardly filtered objects and
on costandardly filtered objects.
\end{Prop}
\begin{proof}
Consider the quotient functor $\pi^{S\vee}:\OCat^{S}_\FM(n)^\vee\twoheadrightarrow \Cat_\FM(n)$.
It is fully faithful on the tilting objects.

Recall that the category $\OCat^S_{\FM}$
admits a naive duality functor, see \S\ref{SSS_Quiv_Schur}. By Lemma \ref{Lem:Ring_naive_dual},
so does the category $(\OCat^S_{\FM})^\vee$. The naive duality functor preserves all
simple objects.

Let $L$ be a simple killed by $\pi^{S\vee}$ that lies in the socle of some
standard object $\Delta^\vee(\lambda)$. Thanks to the naive duality, $L$ also lies  in the head of the costandard
$\nabla^\vee(\lambda)$. It follows that there is an endomorphism of $T^\vee(\lambda)$
that is killed by $\pi^{S\vee}$. We get a contradiction.
\end{proof}

\subsection{Equivalence of extended quotients}\label{SS_Ext_quot_equiv}
In this section we establish property (I) for restricted $\Ringg$-deformed highest weight
$\tgl_e$-categorifications $\OCat_\Ringg(\leqslant N)$, where $N,s_0-s_1,\ldots,s_{\ell-2}-s_{\ell-1}\gg 0$
provided property (F$_{-1}$) holds.
The proof is inspired by that of a similar result in \cite[Section 7]{VV_proof}.
Much of that proof carries over to our situation, in particular, the case of $e>2$.
So we assume that $e=2$. The most essential and difficult step in the proof for $e=2$ is
to describe the projective object $P^i_\Ringg(\nu)$.



\subsubsection{Inductive description}
Pick $n\ll N, s_0-s_1,\ldots,s_{\ell-2}-s_{\ell-1}$ and write $\OCat_\Ringg$ for $\OCat_\Ringg^i, i=1,2$.
Let $\mathcal{P}^+_{\ell-1}(n)\subset \mathcal{P}_\ell(\leqslant n)$ denote the subsets
of all $\ell$-multipartitions $(\lambda^{(0)},\ldots,\lambda^{(\ell-2)},\varnothing)$.
This is a poset co-ideal. So we can consider the highest weight quotient
$\OCat^{+}_\Ringg(\leqslant n)$ of $\OCat_\Ringg(\leqslant n)$. Similarly, let $\mathcal{P}^-_{\ell-1}(n)$
denote  the subset of $\mathcal{P}_{\ell}(n)$ consisting of all $\ell$-partitions
$(\varnothing, \lambda^{(1)},\ldots,\lambda^{(\ell-1)})$. This is a poset ideal.
So we can consider the highest weight subcategory $\OCat^{-}_\Ringg(\leqslant n)$.
Set $\underline{s}^+=(s_0,\ldots,s_{\ell-2}), \underline{s}^-=(s_1,\ldots,s_{\ell-1})$.

\begin{Prop}\label{Prop:trunc_Schur}
The category $\OCat^{\pm}_\Ringg(\leqslant n)$ is a restricted $\Ringg$-deformed highest weight
$\tgl_e$-categorification of $\mathcal{F}^{e,\underline{s}^{\pm}}$.
\end{Prop}
\begin{proof}
To prove the claim about $\OCat^-_\Ringg(\leqslant n)$ we use an argument
similar to \S\ref{SSS:cat_O_restr_catn}, in particular, a suitable version of
the  categorical truncation explained there. Namely, we freeze the $(1,1)$-box
in $\lambda^{(0)}$ and truncate the functor $F_{s_0}$. To deal with
$\OCat^+_\Ringg(\leqslant n)$ we use the dual argument freezing the $(1,1)$-box
in $\lambda^{(\ell-1)}$ and truncate the functor $F_{s_{\ell-1}}$.
\end{proof}

Now we are going to describe the image of $\OCat^{-}_\Ringg(n')\subset \OCat_\Ringg(n')$
(where $n'\leqslant n$) under the quotient functor $\pi_\Ringg:\OCat_\Ringg(n')\twoheadrightarrow
\Cat_\Ringg(n')$. Note that $\mathcal{H}^{q,\underline{s}^-}(n')$ is naturally the quotient
of $\mathcal{H}^{q,\underline{s}}(n')$ so we get the inclusion
$\Cat^-_\Ringg(n')\hookrightarrow \Cat_\Ringg(n')$.

\begin{Lem}\label{Lem:induct_Schur_minus}
The inclusions $\OCat^{-}_\Ringg(n')\subset \OCat_\Ringg(n'), \Cat^-_\Ringg(n')\hookrightarrow \Cat_\Ringg(n')$
intertwine the quotient functors $\pi^-_\Ringg,\pi_\Ringg$.
\end{Lem}
\begin{proof}
The proof repeats that of \cite[Lemma 4.9]{LW}.
\end{proof}

We want an analog of this result for $\OCat^+_\Ringg(n')$. Here we have the
inclusion $\OCat^+_\Ringg(n')^{\Delta}\subset \OCat_\Ringg(n)^{\Delta}$. Let us write
$\Cat^+_\Ringg(n')$ for the quotient of $\OCat^+_\Ringg(n')$  defined analogously
to  the  quotient $\Cat^-_\Ringg(n')$. The proof
of the next lemma again repeats the proof of \cite[Lemma 4.9]{LW}.

\begin{Lem}\label{Lem:induct_Schur_plus}
The inclusions $\OCat^{+}_\Ringg(n')^{\Delta}\subset \OCat_\Ringg(n')^{\Delta}, \Cat^+_\Ringg(n')\hookrightarrow \Cat_\Ringg(n')$
intertwine the  functors $\pi^+_\Ringg,\pi_\Ringg$.
\end{Lem}

\subsubsection{Standardly filtered (co)kernel}
Let $\OCat_\Ringg(\leqslant N)$ be as before. We assume that it is one of the categories
in Proposition \ref{Prop:unique_check_conditions}, in particular, (F$_{-1}$) holds, see Section \ref{SS_-1_faith}.  We write $\OCat(\leqslant N)$
for the specialization to $\FM$.

\begin{Prop}\label{Prop:ker_condns}
The kernel (resp., cokernel) of $T+1$ in $(F_0F_1\oplus F_1F_0)\Delta(\varnothing)$ is standardly filtered
so that $\Delta(\lambda)$ appears in the filtration if and only if $\lambda$ is not a column
(resp., not a row) and with multiplicity $1$.
\end{Prop}
\begin{proof}
The kernel of an epimorphism of standardly filtered objects is standardly filtered, so it is
enough to prove the claim about the cokernel.  The proof is in several steps.

{\it Step 1}. We start with some preparatory remarks.

{\it Step 1.1}. First, let us observe that $(F_0F_1\oplus F_1F_0)\Delta(\varnothing)$
is standardly filtered with standard quotients of the form $\Delta(\mu)$, where $\mu$
is a column or a row in a single diagram (these factors appear with multiplicity $1$)
or $\mu$ consists of two boxes in different diagrams (these appear with multiplicity $2$).
We have a filtration $M_{i,j}\subset (F_0F_1\oplus F_1F_0)\Delta(\varnothing)$ such that
$M_{i,j}\subsetneq M_{i',j'}$ if $j<j'$ or $j=j'$ and $i<i'$ and the quotient of $M_{i,j}$
by the smaller filtration terms is filtered by standard objects that have one box in the $i$th
diagram and the other box in the $j$th diagram.
Any endomorphism of
$(F_0F_1\oplus F_1F_0)\Delta(\varnothing)$ preserves this filtration.

{\it Step 1.2}. Now let $M$ be  a filtration subquotient of $(F_0F_1\oplus F_1F_0)\Delta(\varnothing)$
for the filtration from Step 1.1. We remark that $\operatorname{im}(T+1)|_M=M/\ker(T+1)|_M$
admits a filtration whose successive quotients are quotients of the objects $\Delta(\mu)$
appearing in $M$. We will prove that in the case when $M=(F_0F_1\oplus F_1F_0)\Delta(\varnothing)$,
this filtration on $M/\ker (T+1)|_M$ has quotients $\Delta(\mu)$, where $\mu$ is either a
column or two boxes in different diagrams, all occur with multiplicity $1$.

{\it Step 1.3}. The kernel of $(T+1)$ in $(F_0F_1\oplus F_1F_0)\Delta_{\K}(\varnothing)$
is the sum $\bigoplus_\mu \Delta_{\K}(\mu)$, where the sum is taken over
all $\mu$ that are either a single row or two boxes in different partitions (recall
that $\K=\operatorname{Frac}(\Ringg)$). Indeed, the proof reduces
to the case of $\ell=1$, where it is easy. It follows
that the class of $\ker(T+1)$ in $K_0(\OCat)$ has the form
$\sum_\mu [\Delta(\mu)]+C$, where $C$ is some effective class, and
$\mu$ is as above. Therefore the class of $\operatorname{im}(T+1)\cong(F_0F_1\oplus F_1F_0)\Delta(\varnothing)/
\ker(T+1)$ has the form $\sum_{\mu'} [\Delta_{\FM_q}(\mu')]-C$, where $\mu'$ runs over the
set of all diagrams that are either columns or consist of two boxes in different diagrams.

{\it Step 1.4}. In $\pi((F_0F_1\oplus F_1F_0)\Delta(\varnothing))$ (that is a
direct summand in the Hecke algebra $\mathcal{H}^{2,\underline{r}}(2)$) we have
$\operatorname{im}(T+1)=\ker(T+1)$. Since $\pi$ is an exact functor, it annihilates
$\ker(T+1)/\operatorname{im}(T+1)$. Equivalently, the number of simple composition
factors not annihilated by $\pi$ in $\operatorname{im}(T+1)$ and $\ker(T+1)$
(taken in $(F_0F_1\oplus F_1F_0)\Delta(\varnothing)$)
coincide.

{\it Step 1.5}. The algebra $\mathcal{H}(2):=\mathcal{H}^{2,\underline{s}}_{\FM}(2)$
has two simple modules in the block with
residues $0,1$: $\tilde{f}_1\tilde{f}_0\varnothing$ and $\tilde{f}_0\tilde{f}_1\varnothing$
(we ignore the case where all residues are the same, in which case there is a single
simple module; this case is easier than what we consider). The labels $\tilde{f}_1\tilde{f}_0\varnothing,
\tilde{f}_0\tilde{f}_1\varnothing$ are as follows. One of them is the column in the $\ell-1$th diagram.
The other consists of two boxes, in diagrams $x,\ell-1$, where $x$ is maximal with $s_x-s_{\ell-1}$ odd.

The module $\pi(\Delta(\mu))$ is simple if and only if
$\mu$ is a single row or a single column, it is isomorphic to $\tilde{f}_{i+1}\tilde{f}_i\varnothing$,
where $i$ is the residue of the box with coordinates $(0,0)$ that is present in the diagram of $\mu$.
If $\mu$ is two boxes located in the different diagrams, then $\pi(\Delta(\mu))$
has two different composition factors, each with multiplicity $1$. Both cases are proved by
applying the  functors $E_0,E_1$ to $\pi(\Delta(\mu))$ because $E_0 \pi(\Delta(\mu)),
E_1 \pi(\Delta(\mu))$ are simple or zero.

{\it Step 1.6}. Recall that $\pi$ is faithful on both standard and costandard objects.

{\it Step 2}. In this step we will analyze the action of $T+1$ on various subquotients
$M$ as in Step 1.

{\it Step 2.1}. Let $\mu$ be the diagram with boxes (of two different residues) in two
different partitions. We claim that $\pi(\Delta(\mu))$ is indecomposable. Indeed,
using Lemmas \ref{Lem:induct_Schur_minus},\ref{Lem:induct_Schur_plus}, we
reduce the proof to the case when the partitions are $0,\ell-1$. In this case,
$\Delta(\mu)$ is the sub in the standard filtration of $F_{s_0}\Delta(\mu')$,
where $\mu'$ is the single box in the $\ell-1$th diagram. The object $\pi(\Delta(\mu'))$
is simple. It follows that $\pi(F_{s_0}\Delta(\mu'))=F_{s_0}\pi\Delta(\mu')$
has simple socle. Therefore $\pi(\Delta(\mu))$ has simple socle, which finishes the proof.
Similarly, $\pi(\nabla(\mu))$ has simple socle, hence it is indecomposable as well.

{\it Step 2.2}. Consider the subquotient $M=\Delta(\mu)^{\oplus 2}$, where $\mu$
is a multipartition with $2$ boxes in different diagrams. We claim that $M$
is free over $\FM[T]/((T+1)^2)$, equivalently, $\ker(T+1)|_{M}=\operatorname{im}(T+1)|_M
\cong \Delta(\mu)$. In order to show this, consider
$$\Hom(M,\nabla(\mu))=\Hom(F^2\Delta(\varnothing),\nabla(\mu))=
E^2 \nabla(\mu)=\pi(\nabla(\mu)).$$
The operators $X_1,X_2$  on the 2-dimensional module $\pi(\nabla(\mu))$
are diagonalizable.  The module $\pi(\nabla(\mu))$ is indecomposable by Step 2.1.
So the action of $T+1$ on $\pi(\nabla(\mu))$ is nonzero and our claim follows.

{\it Step 2.3}. Consider now $M$ filtered with $\Delta(\lambda_2),\Delta(\lambda_{11})$,
where $\lambda_2$ is a row and $\lambda_{11}$ is a column located in the same diagram.
Note that $T+1$ annihilates $\Delta(\lambda_2)$. For $\ell=1$,  $M=(F_0F_1\oplus F_1F_0)\Delta(\varnothing)$.
 So in this case, the kernel of $(T+1)$ in $\pi(M)$  coincides
with its image. Using Lemmas \ref{Lem:induct_Schur_minus},\ref{Lem:induct_Schur_plus}
and Step 1.4, we  see that the same is the case for the general $\ell$ as well.

The socle of $\Delta(\lambda_{11})$ is simple. Indeed, $\pi(\Delta(\lambda_{11}))$ is simple
and the socle cannot contain simples killed by $\pi$ because $\pi$ is faithful on standard
objects. This socle is not contained in $\ker(T+1)/\Delta(\lambda_2)$
(Step 1.4 combined with Step 1.6). We conclude that $\ker(T+1)=\Delta(\lambda_2)$.

{\it Step 3}.
Now we are going to prove the claim in Step 1.2 that will finish the proof
of the lemma. Recall that the filtration  we consider on $(F_0F_1\oplus F_1F_0)\Delta(\varnothing)$
induces a filtration on $(F_0F_1\oplus F_1F_0)\Delta(\varnothing)/\ker (T+1)$ whose successive
quotients are quotients of $\Delta(\mu)$ (when $\mu$ is supported at a single
diagram) or of $\Delta(\mu)^{\oplus 2}$.
More precisely, the quotients are $Q_j:=\operatorname{F_j}/(\ker(T+1)_{\operatorname{F}_j}+ \operatorname{F}_{j-1})$,
where $\operatorname{F}_j$ stand for the filtration terms from Step 1.1.
If $F_j/F_{j-1}=\Delta(\mu)^{\oplus 2}$, then we see that $Q_j\twoheadrightarrow \Delta(\mu)$.

Now suppose that $\operatorname{F}_j/\operatorname{F}_{j-1}$
is the extension  of $\Delta(\lambda_{11})$ by $\Delta(\lambda_2)$.
So $Q_j\twoheadrightarrow \Delta(\lambda_{11})$. Using Step 1.3,
we complete the proof.
\end{proof}

\subsubsection{Proof of equivalence}\label{SSS_ext_equiv_proof}
Note that the functors $\pi_\Ringg^i:\OCat^i_\Ringg\rightarrow \Cat_\Ringg$ are fully faithful on the
projective objects. Indeed,  by our assumptions the functor $\pi^i_{\FM}:\OCat^i_\FM\twoheadrightarrow
\Cat_\FM$ is $(-1)$-faithful on both standardly and costandardly filtered objects.
Lemma \ref{Lem:proj_incl} together with a standard argument (see, e.g., the proof
of \cite[Theorem 5.1]{LW}) shows that $\pi_\Ringg^i$ is fully faithful on the projective objects.
Now arguing as in \cite[Section 7.2]{VV_proof} we see that to establish (I) it is enough
to check that $\pi_\Ringg^1(P^1_\Ringg(\nu))=\pi_\Ringg^2(P^2_\Ringg(\nu))$.

Let $Q^i_\Ringg(\nu)$ be the kernel of $(T+1)$ in $(F_0F_1\oplus F_1F_0)\Delta^i_\Ringg(\varnothing)$.
It follows from Proposition \ref{Prop:ker_condns} (and the argument of Step 1.3 of
its proof) that $Q^i_\Ringg(\nu)$ is standardly filtered object (the filtration subquotients
are $\Delta^i_\Ringg(\mu)$, where $\mu$ is not a column, with multiplicity $1$).

\begin{Cor}\label{Cor:QP}
We have $Q^i_\Ringg(\nu)=P^i_\Ringg(\nu)$.
\end{Cor}
\begin{proof}
Again, in the proof we suppress the superscript. Also we note that it is enough
to prove the isomorphism after specializing to $\FM$ (because $P^i_\FM(\nu)$ is projective and hence
has ho higher self-extensions). We omit the subscript as well.

Let $R(\nu)$ be the image of $(T+1)$ in $(F_0F_1\oplus F_1F_0)\Delta(\varnothing)$.
It is standardly filtered and is isomorphic to the coimage.
Because of this, the object $R(\nu)$ is still standardly filtered
in the opposite category (that is still a restricted highest weight $\tgl_e$-categorification
of the same Fock space). It follows that
$\Ringg(\nu)$ is costandardly filtered and hence tilting in the initial category.

Now consider the Ringel dual category $\OCat^\vee$.
As we have seen in \S\ref{SSS_restr_Ringel_duality}, it is also a highest weight
restricted $\tgl_e$-categorification of a Fock space. There $Q(\nu)$
plays the role of $\Ringg(\nu)$ and so is tilting in the Ringel dual category.
It follows that it is projective in
the original category.

The proof that $Q(\nu)$ is indecomposable is now done as in the proof of
\cite[Proposition 7.5]{VV_proof}, see Steps 4 and 5 there, in particular.
\end{proof}

\subsection{Companion elements and condition $\widetilde{\mathfrak{C}}$}\label{SS_companion}
In order to test $0$-faithfulness we will need a stronger condition than $(\mathfrak{C}_{\lambda\mu})$.

\subsubsection{Condition}
For this we have to consider certain {\it companion} multi-partitions
together with $w\lambda$. Here $\lambda$ is singular.

In the definition of companions, we need to use a  construction with crystals.
Pick $\lambda\in \mathcal{P}_\ell$ and a residue $j$. Then we can form
the $j$-signature of $\lambda$ (equal to $\sigma_a^{-1}(\lambda)$, where $\lambda\in \Lambda_a$).
By a {\it marked pair} $p$, we mean  a pair of the form $-,+$ that are erased together, the set of
such pairs is well-defined. We write $\lambda[p]$ for the multipartition, where we move the box
corresponding  to the $-$ in $p$ to the position corresponding to the $+$ in $p$.
So on the level of the $j$-signature, the $-$ and the $+$ in the pair $p$ get swapped.

The companions of $w\lambda$ (relative to a reduced expression of $w$) are defined
inductively. As the base, let us define the companion elements of $\lambda$. Those are either  
\begin{itemize}
\item singular
partitions that are less than $\lambda$ and do not lie in the $i$-family of $\lambda$ for any $i$
\item or the 
partitions of the form $\lambda[p]$, where $p$ is a marked pair for some residue $i$.
\end{itemize}
Let, as before, $w=\sigma_{i_k}\ldots \sigma_{i_1}$.
Now suppose that $\bar{\xi}^1,\ldots,\bar{\xi}^m$ are the companions of $\xi:=\sigma_{i_q}\ldots \sigma_{i_1}\lambda$
and let $i:=i_{q+1}$. The companions of $\sigma_i\xi$ are determined as follows. If $\tilde{e}_i\bar{\xi}^j=0$
and $\bar{\xi}^j$ is not in the $i$-family of $\xi$,
then, by definition, $\sigma_i\bar{\xi}^j$ is a companion of $\sigma_i\xi$. Also the multipartitions obtained as
$(\sigma_i\xi)[p]$, where $p$ is a marked pair in the $i$-signature of $\sigma_i\xi$, 
are  companions of $\sigma_i\xi$. By definition, all companions of $\sigma_i\xi$ are obtained in one of these two ways.
It is unclear to us whether the companions of $w\lambda$ are independent of the choice of a reduced
expression for $w$ -- their construction obviously depends on the choice.
Roughly speaking, companions of $w\lambda$ are labels that can appear in the head of the radical of $\Delta(w\lambda)$,
a precise statement is Proposition \ref{Prop:head_ker}. Let us point out that all companions
of $\lambda$ are not $\geqslant \lambda$.


Now let $\lambda$ be singular and $\mu$ be cosingular. We consider the following condition.
\begin{itemize}
\item[$(\widetilde{\mathfrak{C}}_{\lambda,\mu})$] There is $w\in \hat{S}_e$ and a reduced expression $w=\sigma_{i_k}\sigma_{i_{k-1}}\ldots \sigma_{i_1}$ with the following property: 
Let $\xi$ be  a companion of $w\lambda$ constructed from the reduced expression above. Then $\xi\not\preceq w^*\mu$.
\end{itemize}

\subsubsection{Checking for $\ell=1$}
\begin{Prop}\label{Prop:0_faithf_level1}
Let $\ell=1$. Then $(\widetilde{\mathfrak{C}}_{\lambda,\mu})$ holds for all singular $\lambda$ and cosingular $\mu$
with $|\lambda|=|\mu|>r$, where
\begin{itemize}\item $r=0$ if $e>2$,
\item and $r=2$ if $e=2$.
\end{itemize}
\end{Prop}
To prove the proposition, we need to analyze the structure of the companion elements.
Recall the element $C_{0,n}=\sigma_{1-n}\ldots \sigma_1\sigma_0\in \hat{S}_e$.
The proof of the following lemma is straightforward.

\begin{Lem}\label{Lem:pairing}
Let $i$ be a residue mod $e$ and $n>0$. Then the marked pairs of $i$-boxes in $C_{0,n}\lambda$ are precisely
the following: the removable $i$-box originally contained in $\lambda$ lying in the row with number $> n$
and the addable box lying in the next row.
\end{Lem}

\begin{proof}[Proof of Proposition \ref{Prop:0_faithf_level1}]
Lemma \ref{Lem:pairing} shows that for a large enough $n$ there are no marked pairs in $\lambda$. Actually
all boxes with given residue will be either addable (if the residue is $n$) or removable (if the residue
is $n-1$). For large $n$, the number of boxes in the first row of $C_{0,n}^*\mu$ equals $\lceil n/(e-1)\rceil$.
We also see that, on each step, we add a box to the first row of any companion of $\lambda$.
The number of boxes in the first row of a companion of $\lambda$ is always $\geqslant \lceil n/(e-1)\rceil$.
The equality can only occur if $\lambda_1=e$ and the companion is $C_{0,n}(e-1,1)$. If this is the case and 
$e>2$, then we can take $w=1$. And for $e=2$, the partition $\lambda=(e)$ does not appear. 
\end{proof}

\begin{Rem}\label{Rem:prod_extension}
All results  in this section can be trivially extended to the $\tgl_e^{\ell}$-crystals $\mathcal{P}_1^\ell$.
\end{Rem}

\subsection{Checking 0-faithfulness}\label{SS_0_faith}
In this section we use results of the previous one to establish (F$_0$). Here $\OCat_\Ringg(\leqslant N)$
is a restricted $\Ringg$-deformed highest weight categorification of $\mathcal{F}^{e,\underline{s}}$
for some $s$. We pick $n\ll N$.

\subsubsection{Representation theoretic meaning of companion elements}
Let $\lambda\in \Part_\ell(n)$ be a singular element.
Pick a Weyl group element $w$ with reduced decomposition $ \sigma_{i_k}\ldots \sigma_{i_1}$.
Fix a residue $i$ and an $i$-family $\Lambda_a\subset \mathcal{P}_{\ell}(\leqslant N)$
containing $\tau=w\lambda$ with $l(\sigma_i w)>l(w)$.
Then  $\upsilon:=\sigma_iw\lambda=\tilde{f}_i^{h_{i,+}(\tau)}\tau$
makes sense. Let $\tau^1,\ldots,\tau^k$ be all elements in $\mathcal{P}_{\ell}\setminus \Lambda_a$ such that
$\tilde{e}_i \tau^j=0, L(\tau^j)$ lies in the head of $\ker[\Delta(\tau)\rightarrow L(\tau)]$.
We assume that the $i$-families containing the elements $\tau^1,\ldots,\tau^k$ are in $\mathcal{P}_{\ell}(\leqslant N)$
(this can be achieved by increasing $N$).
Set $\upsilon^j:=\upsilon[p_j]$, where $p_1,\ldots,p_l$ are all marked pairs in
the $i$-signature of $\upsilon$.

\begin{Prop}\label{Prop:head_ker}
If $L(\xi)$ lies in the head of $\ker[\Delta(\upsilon)\rightarrow L(\upsilon)]$, then $\xi$ is one
of the elements $\sigma_i\tau^1,\ldots,\sigma_i\tau^k, \upsilon^1,\ldots,\upsilon^l$.
\end{Prop}

In the proof we will need two auxiliary lemmas.

\begin{Lem}\label{Lem:big_head}
Let $\xi\in \mathcal{P}_{\ell}(N)$ be such that
\begin{itemize} \item $\tilde{f}_i\xi\neq 0$ or $\xi$ lies in the same $i$-family as $\sigma_i\tau$,
\item and $L(\xi)$ lies in the head of $\ker[\Delta(\sigma_i\tau)\rightarrow L(\sigma_i\tau)]$.
\end{itemize}
Then $\xi$ belongs to
the same $i$-family as $\sigma_i\tau$ and is obtained from $\sigma_i\tau$  by switching  $-$ and  $+$
in a marked pair.
\end{Lem}
\begin{proof}
Set $\hat{\tau}=\sigma_i\tau$.

Suppose that $\tilde{f}_i\xi\neq 0$.
It follows that $L(\tilde{f}_i\xi)$ lies in the head of $F_i\Delta(\hat{\tau})$. But all labels
of the simples in that head belong to the same family as $\hat{\tau}$.  So we can assume that $\xi$
lies in the same $i$-family as $\hat{\tau}$.

Consider the highest weight sub $\OCat_{\leqslant \Lambda_a}$ corresponding to all labels that
are less than or equal to some element in $\Lambda_a$ and its highest weight quotient $\OCat_{\Lambda_a}$.
The latter is a basic highest weight $\slf_2$-categorification (=single family) in the sense of \cite{str}. It is enough
to prove the following claim: if $L(\xi)$ lies in the head of $\ker[\Delta(\hat{\tau})\rightarrow L(\hat{\tau})]$
taken in $\OCat_{\Lambda_a}$, then $\xi$ is obtained from $\hat{\tau}$ by switching boxes in a marked pair.
Modulo $\xi<\hat{\tau}$, the condition that $L(\xi)$ appears in the kernel is equivalent to
$\Ext^1(L(\hat{\tau}),L(\xi))\neq 0$. Indeed, $\Delta(\hat{\tau})$ is projective in the highest weight
subcategory with labels $\leqslant \xi$, so it surjects onto any nontrivial extension of $L(\hat{\tau})$
by $L(\xi)$.

Now let $\Cat$ be a basic highest weight $\slf_2$-categorification. We claim that if $\lambda,\mu$
are two labels with $\lambda>\mu$ and $\Ext^1(L(\mu),L(\lambda))\neq 0$, then $\mu$ is obtained from
$\lambda$ by switching the box in a marked pair. The Ext vanishing in the previous paragraph will
follow from here by considering $\OCat_{\Lambda_a}^{opp}$ that is also a basic highest weight categorification.

Let $K$ denote the kernel of $\Delta(\mu)\rightarrow L(\mu)$. Since $\lambda>\mu$, we have
$\Hom(K,L(\lambda))=0$. So $\Ext^1(L(\mu),L(\lambda))\hookrightarrow \Ext^1(\Delta(\mu),L(\lambda))$.
The latter is nonzero if and only if $P(\lambda)$ occurs in the degree $1$ part of the minimal
projective resolution of $\Delta(\mu)$. Now we are done by \cite[Theorem 6.1]{str}.
\end{proof}

The next result follows from the work of Chuang and Rouquier, \cite[Section 6]{ChuRou06}.

\begin{Lem}\label{Lem:sing_iso}
Let $\Cat=\bigoplus_{d\in \Z}\Cat_d$ be an $\slf_2$-categorification with functors $E_i,F_i$.
Let $\Cat_d^{sing}$ denote the Serre subcategory of all objects in $\Cat_d$ annihilated by $E_i$
(for $d\geqslant 0$) and annihilated by $F_i$ (for $d\leqslant 0$). The functor
$F_i^{(d)}:\Cat_{d}\rightarrow \Cat_{-d}, d\geqslant 0,$ restricts to an
abelian equivalence $\Cat_{d}^{sing}\xrightarrow{\sim}\Cat_{-d}^{sing}$
with quasi-inverse $E_i^{(d)}$.
\end{Lem}

We denote the equivalences of this lemma by $\Theta_i$.

\begin{proof}[Proof of Proposition \ref{Prop:head_ker}]
The case when $\xi\in \Lambda_a$ follows from Lemma \ref{Lem:big_head}.
So we consider the case when  $\xi\not\in \Lambda_a$:  we need to check that $\xi=\sigma_i\tau^j$ for some
$j$.  By Lemma \ref{Lem:big_head}, $\tilde{f}_i\xi=0$ and so $\sigma_i\xi$ makes sense.
Next, $\sigma_i\xi\not\geqslant\tau$.  Indeed, since the $i$-families of $\xi, \tau$ are different,
no element of the family of $\xi$ can be bigger than or equal to  $\tau$ (to see this one compares
the boxes with residue different from $i$ in $\xi,\tau$).

Fix an epimorphism $\ker[\Delta(\sigma_i\tau)\rightarrow L(\sigma_i\tau)]\twoheadrightarrow L(\xi)$.
Let $M_{\sigma_i\tau}$ denote the quotient of $\Delta(\sigma_i\tau)$
by the kernel of that epimorphism. So we have a non-split exact sequence
$0\rightarrow L(\xi)\rightarrow M_{\sigma_i\tau}\rightarrow
L(\sigma_i\tau)\rightarrow 0$. In particular, $M_{\sigma_i\tau}$
is annihilated by $F_i$. Let $M_\tau:= \Theta_i M_{\sigma_i\tau}$. Thanks to Lemma \ref{Lem:sing_iso},
we have a non-split exact sequence $0\rightarrow L(\sigma_i\xi)\rightarrow  M_{\tau}\rightarrow L(\tau)\rightarrow 0$.
Consider the Serre subcategory generated by $L(\zeta)$ with $\zeta\leqslant \tau$. Then $\Delta(\tau)$ is a projective
in this subcategory, and $M_\tau$ is an object there. So the projection $\Delta(\tau)\rightarrow L(\tau)$
factors through $\Delta(\tau)\twoheadrightarrow M_\tau$.  We see that $L(\sigma_i\xi)$ lies in the head of
$\ker[\Delta(\tau)\twoheadrightarrow L(\tau)]$ and conclude that $\sigma_i\xi=\tau^j$ for some $j$
and hence $\xi=\sigma_i\tau^j$.
%
\end{proof}

Now we can explain a representation theoretic meaning of companion elements.

\begin{Cor}\label{Cor:companion_RT}
Let $\lambda$ be a singular element in $\mathcal{P}_{\ell}(n)$ and $w\in \hat{S}_e$.
Fix a reduced expression  $w=\sigma_{i_k}\ldots \sigma_{i_1}$. Suppose that $\mathcal{P}_{\ell}(\leqslant N)$
contains the whole $i_j$-family of $\sigma_{i_{j-1}}\ldots \sigma_{i_1}\lambda'$, for any singular $\lambda'$
with $|\lambda'|=|\lambda|$ and any $j\leqslant k$. Further, let $\xi\in \mathcal{P}_\ell(|w\lambda|)$
be such that $L(\xi)$ lies in the head of $\ker[\Delta(w\lambda)\rightarrow L(w\lambda)]$. Then
$\xi$ is a companion element of $w\lambda$ (associated to the fixed reduced expression
of $w$).
\end{Cor}
\begin{proof}
The proof is by induction on the length $\ell(w)$. The induction step is provided by Proposition \ref{Prop:head_ker}.
The base of induction is $w=1$. The straightforward analog of Lemma \ref{Lem:big_head} (with $E_i$ instead of $F_i$)
applied to all $i$ shows that $\xi$  is a companion of $\lambda$.
\end{proof}

\subsubsection{Checking (F$_0$)}
\begin{Prop}\label{Prop:Ext1_vanish}
Suppose that $(\widetilde{\mathfrak{C}}_{\lambda,\mu})$ holds. Then $\Ext^1(L(\lambda),T(\mu))=0$.
\end{Prop}
\begin{proof}
Let $w\in \hat{S}_\ell$ be such that $(\widetilde{\mathfrak{C}}_{\lambda,\mu})$ works for $w$.
Set $\hat{\lambda}=w\lambda, \hat{\mu}=w^*\mu$. As we have seen in the proof of \cite[Proposition 6.3]{VV_proof}, $\Ext^1(L(\lambda),T(\mu))=\Ext^1(L(\hat{\lambda}),T(\hat{\mu}))$.

Let $K:=\ker[\Delta(\hat{\lambda})\rightarrow L(\hat{\lambda})]$.
By Corollary \ref{Cor:companion_RT}, if $L(\xi)\subset \head(K)$ (the head of $K$), then $\xi$ is a companion of
$\hat{\lambda}$. In particular, $\xi\not\preceq \hat{\mu}$. So $L(\xi)$ is not a constituent of
$T(\hat{\mu})$. It follows that $\Hom(K,T(\hat{\mu}))=0$. Also $\Ext^1(\Delta(\hat{\lambda}), T(\hat{\mu}))=0$
because $T(\hat{\mu})$ is tilting. From the exact sequence
$$\Hom(K,T(\hat{\mu}))\rightarrow \Ext^1(L(\hat{\lambda}), T(\hat{\mu}))\rightarrow \Ext^1(\Delta(\hat{\lambda}), T(\hat{\mu}))$$
we deduce that $\Ext^1(L(\hat{\lambda}), T(\hat{\mu}))=0$.
\end{proof}

\begin{Rem}\label{Rem:prod_extension1}
An analog of Proposition \ref{Prop:Ext1_vanish} holds for a highest weight $\tgl_e^{\oplus \ell}$-categorification of $\mathcal{F}^{e,\underline{s}}$.
\end{Rem}

Now let us check that (F$_0$) holds for any restricted highest weight categorification
$\OCat_\Ringg(\leqslant N)$ of $\mathcal{F}$, where the multi-charge satisfies $s_i-s_{i+1}\gg 0$,
and the functor $\bar{\pi}_\Ringg:\OCat_\Ringg\twoheadrightarrow \bar{\Cat}_\Ringg$.
We assume that $n\ll N$.

\begin{Prop}\label{Prop:0_faith}
The functor $\bar{\pi}_{\K}:\OCat_{\K}(\leqslant n)
\twoheadrightarrow \Cat_{\K}(\leqslant n)$ is $0$-faithful.
\end{Prop}
\begin{proof}
The category $\OCat_{\K}(\leqslant N)$ is a restricted
$\tgl_e^{\ell}$-categorification of $\mathcal{F}^{e,\underline{s}}$
(after we pass from $\Ringg$ to $\K$ each of the functors $E_i,F_i$
decomposes into the sum of $\ell$ summands).
The object $Q_{\K}(\nu)$ introduced in \S\ref{SSS_ext_equiv_proof} has a filtration
with successive quotients $\Delta_{\K}(\mu)$, where $\mu$
is an $\ell$-partition of $2$ that is not a column. By block considerations, we
see that $\Delta_{\K}(\lambda)$ is a direct summand of
$Q_{\K}(\nu)$ for any row $\lambda$.  Remark \ref{Rem:prod_extension}
shows that $(\tilde{\mathfrak{C}}_{\lambda\mu})$ holds
for  any $\tgl_e^{\ell}$-singular $\lambda$ with more than two boxes.
Applying Remark \ref{Rem:prod_extension1}, we see that
$\Ext^1(L_{\K}(\lambda),T_{\K}(\mu))=0$
for any $\tgl_e^\ell$-singular $\lambda$ such that
$L_{\K}(\lambda)$ is killed by $\bar{\pi}_\K$
and any $\tgl_e^\ell$-cosingular $\mu$. Arguing as in the proof of
\cite[Proposition 6.3]{VV_proof}, we see that $\bar{\pi}_{\K}$ is $0$-faithful.
\end{proof}

\subsection{Equivalence for $\ell=1$}\label{SS_equiv_ell1}
In this section we will prove the following theorem.

\begin{Thm}\label{Thm:level1_unique}
Let $N\gg n$ and let $\OCat^1(\leqslant N),\OCat^2(\leqslant N)$ be two restricted highest weight $\tgl_e$-categorifications
of the level 1 Fock space $\mathcal{F}^e$. Then there is a labeling preserving
equivalence $\OCat^1(n)\xrightarrow{\sim} \OCat^2(n)$ of highest weight categories
intertwining the quotient functors to $\Cat(n)$.
\end{Thm}
\begin{proof}
We prove this by induction on $n$. The cases $n=0,n=1$ are vacuous. In the case $n=2$, the two categories
are either semisimple (for $e>2$) or are
equivalent to the length two block of the basic highest weight $\mathfrak{sl}_2$-categorification
in the sense of \cite{str} (for $e=2$; the functors acting are $E_1,F_1$).

Now assume that $\OCat^1(<n)\cong \OCat^2(<n)$ (a labeling preserving highest weight equivalence
intertwining the quotient functors to $\Cat(<n)$).
Since the quotient functors are fully faithful on the projectives,  the equivalence
$\OCat^1(<n)\cong \OCat^2(<n)$ is strongly $\tgl_e$-equivariant. Let $\bar{\Cat}^i(n)$ denote the quotient
of $\OCat^i(n)$ by the Serre subcategory spanned by the singular simples. We can argue as in \cite[Section 7.2]{VV_proof}
and show that $\OCat^1(<n)\cong \OCat^2(<n)$ implies $\bar{\Cat}^1(\leqslant n)\cong \bar{\Cat}^2(\leqslant n)$
(a strongly equivariant equivalence, in particular, it preserves the labels of the projectives).
We note that the functor $\OCat^i(n)\twoheadrightarrow \Cat^i(n)$ is 0-faithful, this follows from
Proposition \ref{Prop:0_faith}.
From Theorem \ref{Thm:equi_techn} we see that an equivalence $\OCat^1(n)\xrightarrow{\sim}
\OCat^2(n)$ will follow if we check that $\bar{\pi}^1(\Delta^1(\lambda))\cong \bar{\pi}^2(\Delta^2(\lambda))$
for all $\lambda\in \mathcal{P}(n)$.

In the proof we will assume that $\bar{\pi}^1(\Delta^1(\lambda'))\cong \bar{\pi}^2(\Delta^2(\lambda'))$
for all $\lambda'> \lambda$ (in the dominance ordering). Let $\underline{\lambda}$ be the partition
of $n$ obtained from $\lambda$ by removing the minimal (=leftmost removable) box.
Let $i$ be the residue of this box. Consider the object $F_i\Delta^s(\underline{\lambda}), s=1,2$.

{\it Case 1}. Suppose that $\lambda$ is  the minimal label of a standard occurring in     $F_i\Delta^s(\underline{\lambda})$,
then we can argue as follows. Let us order the labels of standard subquotients in the decreasing way:
$\lambda_1>\ldots>\lambda_k=\lambda$. Let $Q^1_j$ denote the quotient of $F_i\Delta^s(\underline{\lambda})$
filtered with $\Delta^s(\lambda_j),\ldots,\Delta^s(\lambda_k)$. We use the increasing induction
on $j$ to prove that $\pi^1(Q^1_j)\cong \pi^2(Q^2_j)$. The case $j=1$ is an isomorphism
$$\pi^1(F_i\Delta^1(\underline{\lambda}))\cong F_i\pi^1(\Delta^1(\underline{\lambda}))
\cong F_i\pi^2(\Delta^2(\underline{\lambda}))\cong \pi^2(F_i\Delta^2(\underline{\lambda})).$$
Now suppose that we have the required isomorphism for $j$ and let us prove it for $j+1$.
Note that the space $\Hom(\pi^s(\Delta^s(\lambda_j)), \pi^s(Q_j))=\Hom(\Delta^s(\lambda_j), Q_j)$
is one-dimensional and $\pi^s(Q_{j+1})$ is the quotient of $\pi^s(Q_j)$ by the image
of $\pi^s(\Delta^s(\lambda_j))$. By our assumption on $\lambda$, we know that
$\pi^1(\Delta^1(\lambda_j))=\pi^2(\Delta^2(\lambda_j))$. It follows that
$\pi^1(Q_{j+1})=\pi^2(Q_{j+1})$. In particular, we get $\pi^1(\Delta^1(\lambda))=\pi^2(\Delta^2(\lambda))$.

{\it Case 2}. Let us consider the case when $\lambda$ is not the minimal partition obtained from $\underline{\lambda}$
by adding an $i$-box. This is equivalent to the condition that the addable box in the first column
of $\lambda$ is an $i$-box. Let $\lambda^\circ$ be the partition obtained from $\underline{\lambda}$ by adding
the box in the first column. Similarly to the above, we get an isomorphism $\pi^1(Q^1)\cong \pi^2(Q^2)$,
where $Q^s$ is the quotient of $F_i\Delta^s(\underline{\lambda})$ with
$$0\rightarrow \Delta^s(\lambda)\rightarrow Q^s\rightarrow \Delta^s(\lambda^\circ)\rightarrow 0.$$
Here we need to consider two cases: a) $\lambda$ is not singular, and b) $\lambda$ is singular.

{\it Case 2.a}. We see that $\pi^s(L^s(\lambda))\neq 0$. We claim that $\pi^s(L^s(\lambda))=\mathsf{head}(\pi^s(\Delta^s(\lambda)))$. Assume the converse, we have some other simple $\pi^s(L^s(\lambda'))$ in the head.
Consider the surjection $\pi^s(\Delta^s(\lambda))\twoheadrightarrow \pi^s(L^s(\lambda))\oplus
\pi^s(L^s(\lambda'))$. Applying the right adjoint functor $\pi^{s*}$   to this surjection,
we get a morphism
\begin{equation}\label{eq:morphism}\Delta^s(\lambda)=\pi^{s*}\circ\pi^s(\Delta^s(\lambda))\rightarrow
\pi^{s*}\circ\pi^s(L^s(\lambda))\oplus \pi^{s*}\circ\pi^s(L^s(\lambda')).\end{equation}
The first equality is a consequence of the claim that $\pi^s$ is $0$-faithful.
The cokernel of (\ref{eq:morphism}) is killed by $\pi^s$. So the image of
$\Delta^s(\lambda)$ coincides with $L^s(\lambda)\oplus L^s(\lambda')$, which is
impossible.

All simples that occur in $\Delta^s(\lambda^\circ)$ have labels smaller than $\lambda$.
It follows that $\pi^s(\Delta^s(\lambda))$ is the minimal (i.e., contained in every other such object) subobject of $\pi^s(Q)$ that has $\pi^s(L^s(\lambda))$ as a composition factor. From $\pi^1(Q^1)\cong \pi^2(Q^2)$, it follows that
$\pi^1(\Delta^1(\lambda))=\pi^2(\Delta^2(\lambda))$.

{\it Case 2.b}. Here $\pi^s(L^s(\lambda))=0$, equivalently, $\lambda$ is divisible by $e$.
Let $\tilde{L}$ be the maximal singular quotient of $\Delta(\lambda)$. 
Arguing as in Proposition \ref{Prop:head_ker}, the only simples that can appear in the head
of the kernel of $\Delta^s(\lambda)\twoheadrightarrow \tilde{L}$  are $L^s(\nu)$, where
$\nu=\lambda[p]$ for some marked pair $p$ in the signature for some residue $i\in \Z/e\Z$.
Similarly to Case a), we see that any simple in the head of $\pi^s(\Delta^s(\lambda))$
are of the form $\pi^s(L^s(\nu))$. By the structure of basic highest weight $\slf_2$-categorifications,
all $L^s(\nu)$ indeed occur in the composition series of $\Delta^s(\lambda)$
(see, e.g., \cite[Section 7.5]{str}). On the other hand, all simples $\pi^s(L(\mu))$ in $\pi^s(\operatorname{rad}\Delta^s(\lambda^\circ))$ satisfy $\mu<\nu$ for any $\nu$ as above.
So we get the following characterization of $\pi^s(\Delta^s(\lambda))$: it is the minimal
submodule in the kernel of $\pi^s(Q^s)\twoheadrightarrow \pi^s(L^s(\lambda^\circ))$
that contains all $\pi^s(L^s(\nu))$ as composition factors. It follows that
$\pi^1(\Delta^1(\lambda))=\pi^2(\Delta^2(\lambda))$ and completes the proof of
the theorem.
\end{proof}

Again, this theorem trivially generalizes to the case
of $\tgl_e^{\oplus \ell}$-categorifications of $\mathcal{F}_1^{\otimes\ell}$.
This has the following corollary.

\begin{Cor}
Suppose that $s_i-s_{i+1}\gg 0$ for all $i$ and $N\gg n$.
Let $\OCat^1_\Ringg(\leqslant N),\OCat^2_\Ringg(\leqslant N)$
be two restricted $\Ringg$-deformed $\tgl_e$-categorifications of $\mathcal{F}^{e,\underline{s}}$.
Then property (I) holds.
\end{Cor}
\begin{proof}
Apply the analog of Theorem \ref{Thm:level1_unique} to see that we have a labelling
preserving equivalence $\OCat^1_{\K}\xrightarrow{\sim} \OCat^2_{\K}$.
The   functors $\bar{\pi}^i_{\K}$ are defined by the projectives
whose labels are precisely those contained in the $\tgl_e$-crystal components having a multi-partition
with two boxes. So the equivalence $\OCat^1_{\K}\xrightarrow{\sim}
\OCat^2_{\K}$ intertwines the functors $\bar{\pi}^i_{\K}$.
This completes the proof.
\end{proof}

Proposition \ref{Prop:unique_check_conditions} is now fully proved.

%% file: decomp_pKL.tex
\section{Application to decomposition numbers} \label{sec-application}
In this relatively short section we will give applications of Theorem \ref{Thm:asymp_Fock_uniqueness}
to computing the decomposition numbers in the representation theoretic categories
mentioned in the introduction. When we deal with highest weight categories we talk about
multiplicities (of simples in standards) rather than about decomposition numbers.

\subsection{p-Kazhdan-Lusztig polynomials and multiplicities in Soergel categories $\mathcal{O}$}
We start by recalling p-Kazhdan-Lusztig polynomials and their parabolic analogs.
Recall that the indecomposable objects $\mathfrak{B}_w$ in $\DC(\varnothing,\varnothing)$ (here we consider this
category over $\mathbb{F}_p$)
are labelled by the elements  $w\in W$ and that $K_0(\DC(\varnothing,\varnothing))$
is identified with the Hecke algebra (over $\Z[v^{\pm 1}]$). By the p-KL
polynomial $P^p_{y,w}(v)$ we mean the coefficient of the standard basis vector
$T_w$ in the class $[\mathfrak{B}_y]$.

Now let us relate the p-KL polynomials to the multiplicities in the categories
$\,_I\mathcal{O}_{\F,J}^+(W)$. This relation is basically the same as in
the usual Kac-Moody categories $\mathcal{O}$. Pick finitary subsets $I,J\subset S$
and $\alpha\in \,_I\Lambda^+(J)$. Pick $w\in W$ that is minimal (=longest) in the coset $\alpha$.
Then $w$ is shortest in its right $W_J$-coset.

\begin{Prop}\label{Prop:positive_O_mult}
Let $\alpha,\beta\in\,_I\Lambda^+(J)$, and $w,y\in W$ be minimal elements in $\alpha,\beta$.
Then the multiplicity of $L^+_{\F}(\beta)$ in $\Delta^+_{\F,J}(\beta)$ equals
$P_{\beta,\alpha}^{p,J}(1):=\sum_{u\in W_J}(-1)^{\ell(u)}P^p_{y,wu}(1)$.
\end{Prop}
\begin{proof}
Recall, \S\ref{SSS_Soergel_O_constr},  that the indecomposable projective $P^+_{\F}(y)$ in $\,_\varnothing\OCat^+_\F(W)$
is obtained as $\mathfrak{B}_y\Delta^+_{\F}(1)$. So $[P^+_\F(y):\Delta_\F^+(w)]=P^p_{y,w}(1)$.

Let us check that  $[\Delta_\F^+(w):L_{\F}^+(u)]=[\nabla_\F^+(w):L_{\F}^+(u)]$ for all $u$.
Note that  the category $\,_\varnothing\OCat^+_{\operatorname{Frac}(\tilde{\Ring})}(W)$
is semisimple, see Lemma \ref{Lem:O_ss}. So we have
$[\Delta_\F(w)]=[\Delta_{\operatorname{Frac}(\tilde{\Ring})}(w)]=[\nabla_{\operatorname{Frac}(\tilde{\Ring})}(w)]=
[\nabla_\F(w)]$ in any finite highest quotient of $\OCat^+_{\tilde{\Ring}}(W)$.
This implies the claim in the beginning of the paragraph.

It follows that
$$[\Delta_\F^+(w):L^+_{\F}(y)]=[\nabla_\F^+(w):L^+_{\F}(y)]=[P^+_\F(y):\Delta_\F^+(w)]=P^p_{y,w}(1).$$
Now   the formula for  $[\Delta^+_{\F,J}(\alpha):L^+_{\F}(\beta)]$ is derived in
the same way as for the usual singular parabolic categories $\mathcal{O}$.
\end{proof}

Now let us deal with the negative level categories $\mathcal{O}$. The observation that the Ringel duality functor
maps $T^-_{\F,J}(\beta)$ to $P^+_{\F,J}(\beta)$  and $\Delta^-_{\F,J}(\alpha)$ to $\Delta^+_{\F,J}(\alpha)$
gives the following statement.

\begin{Cor}\label{Cor:negative_O_mult}
We have $[T^-_{\F,J}(\beta):\Delta^-_{\F,J}(\alpha)]=P^{p,J}_{\beta,\alpha}(1)$.
\end{Cor}


\subsection{Multiplicities for cyclotomic Schur and Hecke  algebras}
Here explain how to compute the multiplicities in the categories  $\mathcal{S}^{e,\underline{r}}_\FM(n)\operatorname{-mod}$ and the decomposition numbers in  $\mathcal{H}^{e,\underline{r}}_\FM(n)\operatorname{-mod}$
via the numbers $P_{\mu,\lambda}^{p,J}(1)$.

Recall the categories $\OCat^{-e}_\FM(J,n)$ from \S\ref{SSS:cat_O_restr_catn}. Here we consider $W=\hat{S}_m$ and $J\subset S$
recovered from the $\ell$-tuple $(m_1,\ldots,m_\ell)$ (as explained in
\S\ref{SSS_combinatorics_XJ}) with $m_1+\ldots+m_\ell=m$
and $m_i>n$. The labelling set for $\OCat^{-e}_{\FM}(J,n)$ is the set of
$\ell$-multipartitions of $n$. For this category we can compute the multiplicities of standards
in tiltings using Corollary \ref{Cor:negative_O_mult}. But we cannot compute the multiplicities
of simples in standards directly.

Now define the category $\OCat^{+e}_\FM(J,n)$ as the highest weight category that is the Ringel dual
of  $\OCat^{-e}_{\FM}(J,n)$. So the category $\OCat^{+e}_\FM(J,n)$ is a highest
weight  quotient of $\bigoplus_{I\in \mathbb{O}}\,_I\OCat^+_{\F,J}(W)$. Let us elaborate
on how to compute the multiplicities (of simples in standards) for this category.

We can view an $\ell$-partition $\lambda$ of $n$ as a virtual multipartition
with rows of heights $m_1,\ldots,m_\ell$. In its turn, a virtual multipartition
encodes an element of $X(J)$, see \S\ref{SSS_combinatorics_XJ} and hence
a $W$-orbit $O_\lambda\in \mathbb{O}$, a finitary subset $I(\lambda)\subset S$
with  $O_\lambda=W/W_{I(\lambda)}$ and an element $\alpha(\lambda)\in W/W_{I(\lambda)}$.

The following claim is a direct consequence of Proposition \ref{Prop:positive_O_mult}.

\begin{Cor}\label{Cor:mult_posit_level}
We have $[\Delta^{+e}(\lambda):L^{+e}(\mu)]=P^{p,J}_{\alpha(\mu),\alpha(\lambda)}(1)$
if $O_\lambda=O_\mu$. If $O_\lambda\neq O_\mu$, then
$[\Delta^{+e}(\lambda):L^{+e}(\mu)]=0$.
\end{Cor}

Now let us explain how to express the multiplicities in $\mathcal{S}^{e,\underline{r}}_\FM(n)\operatorname{-mod}$,
where $\underline{r}=(r_0,\ldots,r_{\ell-1})$ with $r_i\in \Z/e\Z$.
Pick the numbers $m_1,\ldots,m_\ell\in \Z$ with
\begin{align}\label{eq:cond_m_congr}&m_i\equiv -r_{\ell-i}\mod e,\\\label{eq:cond_m_order}
&m_1-m_2\gg m_2-m_3\gg\ldots \gg m_\ell\gg n.
\end{align}

Set $m=m_1+\ldots+m_\ell$ and consider the group $W$ of type
$\hat{S}_m$.

Combining \S\ref{SSS_restr_Ringel_duality} with Theorem \ref{Thm:asymp_Fock_uniqueness}
we get an equivalence $\OCat^{+e}_\FM(J,n)\rightarrow \mathcal{S}^{e,\underline{r}}_\FM(n)\operatorname{-mod}$ that maps $\Delta^{+e}(\lambda^*)$  to $\Delta^S(\lambda)$,
where $\bullet^*$ is the involution of the set of $\ell$-multipartitions introduced
in \S\ref{SSS_restr_Ringel_duality}. Corollary \ref{Cor:mult_posit_level} implies
the following claim.

\begin{Cor}\label{Cor:mult_Schur}
Suppose that $(m_1,\ldots,m_\ell)$ satisfy (\ref{eq:cond_m_congr}) and (\ref{eq:cond_m_order}).
Then we have $[\Delta^{S}(\lambda):L^{S}(\mu)]=P^{p,J}_{\alpha(\mu^*),\alpha(\lambda^*)}(1)$
provided $O_{\lambda^*}=O_{\mu^*}$. If $O_{\lambda^*}\neq O_{\mu^*}$, then
$[\Delta^S(\lambda):L^S(\mu)]=0$.
\end{Cor}

Since we have the quotient functor $\pi:\mathcal{S}^{e,\underline{r}}_\FM(n)\operatorname{-mod}
\twoheadrightarrow \mathcal{H}^{e,\underline{r}}_\FM(n)\operatorname{-mod}$, we can deduce
a formula for the decomposition numbers of $\mathcal{H}^{e,\underline{r}}_\FM(n)\operatorname{-mod}$
from Corollary \ref{Cor:mult_Schur}. However, we can do better and weaken (\ref{eq:cond_m_order}).

\begin{Cor}\label{Cor:mult_Hecke}
Suppose that $(m_1,\ldots,m_\ell)$ satisfy (\ref{eq:cond_m_congr}) and $m_i>n$ for all $i$.
If $\pi(L^S(\mu))\neq 0$, then we have $[\pi(\Delta^{S}(\lambda)):\pi(L^{S}(\mu))]=P^{p,J}_{\alpha(\mu^*),\alpha(\lambda^*)}(1)$
if $O_{\lambda^*}=O_{\mu^*}$. If $O_{\lambda^*}\neq O_{\mu^*}$, then
$[\pi(\Delta^S(\lambda)):\pi(L^S(\mu))]=0$.
\end{Cor}
Note that $\pi(\Delta^{S}(\lambda))$ is the dual Specht module labelled by $\lambda$.
\begin{proof}
Under our assumptions, we have a restricted categorical $\tgl_e$-action on
$\OCat^{+e}_\FM(J,\leqslant n)$ categorifying $\mathcal{F}^{e,\underline{r}}(\leqslant n)$.
Rouquier's uniqueness result, \cite[5.1.2]{Rouq-2KM}, yields a quotient functor
$\OCat^{+e}_\FM(J,\leqslant n)\twoheadrightarrow \bigoplus_{i=0}^n\mathcal{H}^{e,\underline{r}}_\FM(i)\operatorname{-mod}$.
The class  $[\Delta^{+e}(\lambda^*)]$ in $K_0$ coincides with $[\Delta^-(\lambda)]$,
which then coincides with $[\Delta^S(\lambda)]$ (both are equal to $|\lambda\rangle$).
It follows that  $[\pi(\Delta^{S}(\lambda))]=[\pi(\Delta^{+e}(\lambda^*))]$.
Now we are done by Proposition \ref{Prop:positive_O_mult}.
\end{proof}

\subsection{Multiplicities for rational representations and quantum group modules}
In the remainder of the section we will be dealing with the categories of finite dimensional
$U_{q,\FM}(\gl_m)$-modules (where $q$ is a root of $1$) and also $\GL_m(\FM)$-modules.
We will concentrate on the case of $\GL_m(\FM)$, the case of quantum groups is analogous.
Let $\operatorname{Pol}_{m,d}$ denote the category of polynomial representations
of $\GL_m(\FM)$ of degree $d$, this is a highest weight subcategory of $\operatorname{Rep}(\GL_m(\FM))$.

\begin{Lem}\label{Lem:polyn_reln}
For $m<m'$, the category $\operatorname{Pol}_{m,d}$ is a highest weight quotient of
$\operatorname{Pol}_{m',d}$ corresponding to labels with at most $m$ rows.
\end{Lem}
\begin{proof}
Consider the functor $\tau$ that takes the invariants for the torus $T_{m,m'}$ of all diagonal matrices with
first $m$ entries equal to $1$. It has a right inverse that is an induction functor.
So $\tau$ is a quotient functor.
Moreover, $\tau$ maps $\Delta_{m'}(\lambda)$ to $\Delta_m(\lambda)$ if $\lambda$ has at most
$m$ rows and to $0$ else. This implies the claim.
\end{proof}

Now note that every rational representation becomes polynomial after twisting by a sufficiently
large power of $\det$. In other words, $\operatorname{Rep}(\GL_m)$ is the direct limit of
its highest weight subcategories equivalent to $\operatorname{Pol}_{m,d}$ for various $d$.
Together with Lemma \ref{Lem:polyn_reln}, this reduces the computation of the multiplicities in $\operatorname{Rep}(\GL_m)$ to those in the categories $\operatorname{Pol}_{m',d}$ with $m'>d$, i.e., to the multiplicities
in the category $\mathcal{O}^S_{p,0}$.